\documentclass[reqno]{amsart}
\usepackage[utf8]{inputenc}
\usepackage{mathrsfs}
\usepackage{mathabx}
\usepackage{bm,bbm}
\usepackage{xcolor}

\usepackage{amsfonts}
\DeclareMathAlphabet{\mathpgoth}{OT1}{pgoth}{m}{n}
\DeclareMathAlphabet{\mathesstixfrak}{U}{esstixfrak}{m}{n}
\DeclareMathAlphabet{\mathboondoxfrak}{U}{BOONDOX-frak}{m}{n}
\usepackage{amsmath,amssymb}
\usepackage[bbgreekl]{mathbbol}
\usepackage{yfonts,mathtools}

\usepackage{tikz-cd}
\usepackage[all,cmtip]{xy}

\numberwithin{equation}{section}
\usepackage{soul}
\usepackage{comment}

\usepackage{hyperref}
\hypersetup{colorlinks}
\definecolor{darkred}{rgb}{0.5,0,0}
\definecolor{darkgreen}{rgb}{0,0.5,0}
\definecolor{darkblue}{rgb}{0,0,0.5}
\hypersetup{colorlinks, linkcolor=darkblue, filecolor=darkgreen, urlcolor=darkred, citecolor=darkblue}
\makeatletter 
\@addtoreset{equation}{section}
\makeatother  

\numberwithin{equation}{section}

\newtheorem{thm}{Theorem}[section]
\newtheorem{cor}[thm]{Corollary}

\newtheorem{prop}[thm]{Proposition}
\newtheorem{lemma}[thm]{Lemma}
\theoremstyle{definition}
\newtheorem{defn}[thm]{Definition}
\theoremstyle{remark}
\newtheorem{rem}[thm]{Remark}

\newtheorem{addendum}[thm]{Addendum}

\newtheorem{example}[thm]{Example}

\newcommand{\beq}{\begin{equation}}
\newcommand{\eeq}{\end{equation}}
\newcommand{\beqn}{\begin{equation*}}
\newcommand{\eeqn}{\end{equation*}}
\newcommand{\ov}{\overline}
\newcommand{\mb}{\mathbb}
\newcommand{\mc}{\mathcal}
\newcommand{\mf}{\mathfrak}

\newcommand{\coloru}{\color{olive}}

\newcommand{\wt}{\widetilde}
\newcommand{\wh}{\widehat}

\newcommand{\uds}[1]{\underline{\smash{#1}}}

\renewcommand{\subset}{\subseteq}

\newcommand{\ev}{{\rm ev}}

\newcommand\rightthreearrow{%
        \mathrel{\vcenter{\mathsurround0pt
                \ialign{##\crcr
                        \noalign{\nointerlineskip}$\rightarrow$\crcr
                        \noalign{\nointerlineskip}$\rightarrow$\crcr
                        \noalign{\nointerlineskip}$\rightarrow$\crcr
                }%
        }}%
}

\makeatletter
\newcommand{\colim@}[2]{%
  \vtop{\m@th\ialign{##\cr
    \hfil$#1\operator@font colim$\hfil\cr
    \noalign{\nointerlineskip\kern1.5\ex@}#2\cr
    \noalign{\nointerlineskip\kern-\ex@}\cr}}%
}
\newcommand{\colim}{%
  \mathop{\mathpalette\colim@{\rightarrowfill@\textstyle}}\nmlimits@
}
\makeatother

\newcommand{\Gammait}{{\mathit{\Gamma}}}

\title[Transversality on Orbifolds and ${\mb Z}$-valued GW Invariants]{A New Transversality Condition on Orbifolds and Integer-Valued Gromov--Witten Type Invariants}

\author{Shaoyun Bai}
\address{Department of Mathematics, MIT, 77 Massachusetts Avenue, Cambridge, MA 02139-4307, USA}
\email{shaoyunb@mit.edu}

\author{Guangbo Xu}
\address{Department of Mathematics, Rutgers University, Hill Center 
 -  Busch Campus,
110 Frelinghuysen Road, Piscataway, NJ 08854-8019, USA}

\email{gx49@math.rutgers.edu}

\thanks{S. B. is supported by NSF DMS-2404843. G. X. is supported by NSF DMS-2345030.}

\date{\today}

\begin{document}

\maketitle

\begin{abstract}
Following the proposal of Fukaya--Ono \cite{Fukaya_Ono_integer} and the exploration by B. Parker \cite{BParker_integer}, we introduce a new transversality condition, the {\it FOP transversality condition}, for sections of orbifold vector bundles ${\mc E} \to {\mc U}$ when both ${\mc E}$ and ${\mc U}$ have ``normal complex structures." This notion allows one to define various integral virtual cycles on moduli spaces of pseudoholomorphic curves. Two immediate applications in symplectic topology are the definition of integer-valued Gromov--Witten type invariants  in all genera  for general compact symplectic manifolds using the global Kuranishi chart constructed by Abouzaid--McLean--Smith \cite{AMS}\cite{AMS2} and Hirschi--Swaminathan \cite{Hirschi_Swaminathan_2022} , and an alternative proof of the cohomological splitting theorem for Hamiltonian fibrations over $S^2$ with integer coefficients from \cite{AMS}.
\end{abstract}






\setcounter{tocdepth}{1}
\tableofcontents

\section{Introduction}

Gromov--Witten invariants are generally rational numbers but not integers. A conceptual reason for the non-integrality is that we need to consider curves with symmetries. In the algebraic geometric situation, this means the moduli spaces are stacks rather than schemes. In the technical construction of symplectic Gromov--Witten invariants (see  \cite{Li_Tian},\cite{Fukaya_Ono},\cite{Ruan_1999},\cite{Cieliebak_Mohnke},\cite{Pardon_virtual},\cite{HWZ-GW},\cite{Siebert_virtual}, and \cite{AMS}), it is also the consequence of using multi-valued perturbations. For the same reason, Hamiltonian Floer homology of a general symplectic manifold is also defined over ${\mb Q}$ but not over ${\mb Z}$ (see \cite{Fukaya_Ono},\cite{Liu_Tian_2}, and \cite{Pardon_virtual}). This defect limits our applications. For example, regarding the Arnold conjecture, the Hamiltonian Floer chain complex can only detect rational homology classes of the symplectic manifold\footnote{The recent advancement of the Arnold conjecture by Abouzaid--Blumberg \cite{abouzaid2021arnold} strengthens the lower bound to the total Betti number over finite fields using more involved stable homotopy theory; afterwards we use the method developed in this paper to obtain a stronger lower bound including contributions from all characteristics.}. One cannot either define quantum Steenrod operations for general symplectic manifolds as one needs ${\mb Z}/p$ counts of genus zero curves (see \cite{Wilkins_2020} and \cite{Seidel_Wilkins}).

The main motivation of this paper is to define integer-valued Gromov--Witten invariants and prove integral Arnold conjecture, following a 
proposal of Fukaya--Ono \cite{Fukaya_Ono_integer}. Roughly speaking, because moduli spaces of pseudoholomorphic curves carry a ``normal complex structure,'' one can use a specified type of {\it single-valued} perturbations to separate curves with different automorphism groups; as a result, one could count only curves with trivial automorphism group and obtain well-defined integral invariants.

The main contribution of this paper is to rigorously realize Fukaya--Ono's proposal in an abstract, finite-dimensional setting which can be applied to symplectic geometric problems. In particular, we define a new notion of transversality for sections on orbifolds. The most important topological consequence  is the existence of integer-valued Euler cycles of orbifold vector bundles. We include two geometric applications. The first one is the definition of integer-valued Gromov--Witten type invariants for general compact symplectic manifolds. As far as we know, this is the first general construction of integer-valued curve-counting invariants for \emph{all} symplectic manifolds, which in particular include all smooth projective varieties over complex numbers. The second application is a proof of principle: we offer an alternative proof of Abouzaid--McLean--Smith's cohomological splitting result for Hamiltonian fibrations \cite{AMS} not relying ideas inspired by Floer homotopy theory.

In the subsequent work \cite{Bai_Xu_Arnold}, we prove the Arnold conjecture over ${\mb Z}$ by further developing the ideas in this paper. We expect this paper to serve as a starting point of constructing refinements of curve-counting invariants for general symplectic manifolds. 
Moreover, the methods here should also be useful for tackling longstanding conjectures concerning quantitative lower bound for fixed/intersection points in various settings, beyond the integral version of homological Arnold conjecture.


\subsection{Statement of the main technical result}

In this paper we define a new notion of transversality for sections of orbifold vector bundles. This notion, called the {\it FOP transversality condition}, originates from the proposal of Fukaya--Ono \cite{Fukaya_Ono_integer} and which was further elaborated by B. Parker \cite{BParker_integer}. One of the upshots of this notion is to define integral virtual fundamental classes when the usual notion of transversality fails in the presence of symmetry. 

Before we get into details we would like to discuss the general philosophy. It is not hard to realize that transversality and symmetry are two conflicting features: the former is generic while the latter is specific. However we need to reconcile them when discussing transversality for sections of orbifold vector bundles, the model for regularizing moduli spaces in symplectic topology. The usual transversality theorem can provide a regular zero locus only away from orbifold points. Indeed, a generic section (which preserves the symmetry) results in a singular zero locus stratified by the isotropy types of points in the orbifold. 

One way to reconcile the two conflicting aspects is to drop the symmetry by using ``multivalued'' sections while sacrificing the integrality of the fundamental cycles. This is the method used in \cite{Li_Tian, Fukaya_Ono} which essentially produces cycles dual to the (rational-valued) orbifold Euler class. The method proposed by Fukaya--Ono \cite{Fukaya_Ono_integer} and worked out in this paper is, in another direction, to enhance the notion of transversality without dropping the symmetry or losing the integrality. One can see from Theorem \ref{thm11} below that this new notion is almost as flexible as traditional transversality. Moreover, even though the zero locus is still singular, it has a stratification whose strata lead to integral cycles which can be used to construct enumerative invariants. An essential requirement for discussing this new notion, which is the ``normal complex structure,'' is satisfied by moduli spaces of pseudoholomorphic curves, thanks to the nature of the Cauchy--Riemann operator.

\subsubsection{The FOP transversality condition}

To state the transversality theorem in more precise terms, we need certain notational preparations. Consider an effective\footnote{All orbifolds in this paper are effective unless otherwise declared.} orbifold ${\mc U}$ and an orbifold vector bundle ${\mc E} \to {\mc U}$. The pair $({\mc U}, {\mc E})$ has local charts of the form $(G, U, E)$ where $G$ is a finite group, $U$ is a smooth effective $G$-manifold, and $E \to U$ is a $G$-equivariant vector bundle. 
Each point $x \in U$ can be labelled by the isotropy subgroup $G_x \subset G$ and two representations of $G_x$: the tangent space $V_x = T_x U$ and the fiber $W_x = E_x$. Such a triple $(G_x, V_x, W_x)$ is called the {\it isotropy type} of $x$, whose isomorphism class only depends on the corresponding point in the orbifold ${\mc U}$. In this way the orbifold ${\mc U}$ is stratified by isotropy types $\gamma = (G, V, W)$. Let ${\mc U}_\gamma \subset {\mc U}$ be the set of points with isotropy type $\gamma$. In this way ${\mc U}$ becomes a (Thom--Mather) stratified space (Definition \ref{defn_Thom_Mather}) whose strata ${\mc U}_\gamma$ consist of points with fixed isotropy types. Among them the top stratum of ${\mc U}$ is a dense open subset ${\mc U}_{\rm free}$ of isotropy free points corresponding to the isotropy type $(
\{1\}, {\mb R}^m, {\mb R}^n)$.

The usual transversality for sections can be achieved within each stratum, but may fail in the normal directions. If an isotropy type $\gamma$ is represented by $(G, V, W)$, then one can decompose
\begin{align}\label{basic_1}
&\ V = V_G \oplus \check V_G,\ &\ W = W_G \oplus \check W_G
\end{align}
where $V_G, W_G$ resp. $\check V_G, \check W_G$ are the direct sums of irreducible trivial resp. nontrivial subrepresentations of $V, W$. If we restrict the tangent bundle $T{\mc U}$ and the bundle ${\mc E}$ to ${\mc U}_\gamma$, then one has corresponding decompositions
\begin{align*}
    &\ T{\mc U}|_{{\mc U}_\gamma} = T{\mc U}_\gamma \oplus {\mc N}_\gamma,\ &\ {\mc E}|_{{\mc U}_\gamma} =  {\mc E}_\gamma \oplus \check {\mc E}_\gamma.
\end{align*}
The ``normal'' components ${\mc N}_\gamma$ and $\check{\mc E}_\gamma$ are roughly bundles with fibers being $\check V_G$ and $\check W_G$ specified above. Notice that for any section ${\mc S}: {\mc U} \to {\mc E}$, when restricted to ${\mc U}_\gamma$, the $\check{\mc E}_\gamma$-component must vanish. If we decompose
\beqn
{\mc S}|_{{\mc U}_\gamma} = {\mc S}_\gamma \oplus \check{\mc S}_\gamma,
\eeqn
then usually the best we can hope is the transversality of $  {\mc S}_\gamma: {\mc U}_\gamma \to {\mc E}_\gamma$. The derivative of ${\mc S}$ in the direction ${\mc N}_\gamma$ may not be enough to cover the missing direction $\check {\mc E}_\gamma$, causing the failure of transversality.

A {\it normal complex structure} on ${\mc U}$ resp. ${\mc E}$ consists of $G$-invariant complex structures on fibers of ${\mc N}_\gamma$ resp. $\check{\mc E}_\gamma$ for all isotropy types $\gamma$ which are ``compatible'' between adjacent strata (see Section \ref{sec:orbifolds} for the rigorous definition). In particular, if ${\mc U}$ and ${\mc E}$ are normally complex, then the isotropy types appearing are those $(G, V, W)$ with $\check V_G$, $\check W_G$ being complex representations; such isotropy types are called {\it normally complex isotropy types}. 

For the purpose of dealing with infinite-dimensional situations, we need to introduce the {\it stabilization} of isotropy types. A normally complex isotropy type $(G, V, W)$ can be stabilized to $(G, V \oplus R, W \oplus R)$ where $R$ is a normally complex representation (cf. Definition \ref{defn:rep-normal-complex}). Let $\uds {\mb \Gamma}^{\rm NC}$ be the set of all complex isotropy types modulo stabilization, which has induced partial order induced from inclusions of groups and representations. Let $\uds {\mb\Gamma}_k^{\rm NC} \subset \uds {\mb \Gamma}^{\rm NC}$ be the subset of elements whose virtual dimension ${\rm dim}_{\mb R} (V)- {\rm dim}_{\mb R} (W)$ is equal to $k\in {\mb Z}$. There is a special maximal element $\gamma_k^* \in \uds {\mb \Gamma}_k^{\rm NC}$ represented by $(\{1\}, {\mb R}^m, {\mb R}^n)$ with $m-n = k$.

Our main technical result is stated as the following blackbox-type theorem. 

\begin{thm}\label{thm11}
Suppose ${\mc U}$ is a normally complex orbifold without boundary and ${\mc E} \to {\mc U}$ is a normally complex vector bundle. Let $\Gamma({\mc U}, {\mc E})$ be the space of smooth sections. Then there is a $C^0$-dense\footnote{In general it is not $C^1$-dense. See Remark \ref{rem_density}.} subset $\Gamma^{\rm FOP}({\mc U}, {\mc E}) \subset \Gamma({\mc U}, {\mc E})$ whose elements are called \emph{FOP transverse sections} satisfying the following properties.
\begin{enumerate}

\item {\bf (Classical Transversality)} If ${\mc U}$ is a manifold, FOP transversality is equivalent to classical transversality.

\item {\bf (Locality)} The restrictions of FOP transverse sections to open subsets are still FOP transverse. 

\item {\bf (Extension Property I)} For any pair of closed subsets ${\mc Y} \subset {\mc Y}'$ and open neighborhoods ${\mc V} \subset {\mc U}$ of ${\mc Y}$ and ${\mc V}' \subset {\mc U}$ of ${\mc Y}'$, if ${\mc S} \in \Gamma^{\rm FOP}({\mc V}, {\mc E}|_{{\mc V}})$, then there exists ${\mc S}' \in \Gamma^{\rm FOP}({\mc V}', {\mc E}|_{{\mc V}'})$ which agrees with ${\mc S}$ near ${\mc Y}$.

\item {\bf (Extension Property II)} Let ${\mc X} \subset {\mc U}$ be a suborbifold with ordinary normal bundle (Definition \ref{defn_ordinary_normal}), which implies that ${\mc X}$ and ${\mc E}|_{\mc X}$ are normally complex. Let ${\mc S}_{\mc X}: {\mc X} \to {\mc E}|_{{\mc X}}$ be an FOP transverse section. Then there exists an FOP transverse extension of ${\mc S}_{\mc X}$ to ${\mc U}$.

\item {\bf (Product Property)} Let ${\mc U}'$ be another normally complex  orbifold and ${\mc E}' \to {\mc U}'$ be a normally complex  vector bundle. Then the product map 
\beqn
\Gamma({\mc U}, {\mc E}) \times \Gamma({\mc U}', {\mc E}') \to \Gamma({\mc U} \times {\mc U}', {\mc E}\boxplus {\mc E}')
\eeqn
sends products of FOP transverse sections to FOP transverse sections.

\item {\bf (Stabilization Property)} If $\pi_{\mc F}: {\mc F} \to {\mc U}$ is a normally complex orbifold vector bundle,\footnote{The total space ${\mc F}$ has a canonical concordance class of normal complex structures.} then the stabilization map 
\beqn
\Gamma({\mc U}, {\mc E}) \to \Gamma({\mc F},  \pi_{\mc F}^* {\mc E} \oplus \pi_{\mc F}^* {\mc F}),\quad \quad {\mc S} \mapsto \pi_{\mc F}^* {\mc S} \oplus \tau_{\mc F},
\eeqn
where $\tau_{\mc F}: {\mc F} \to \pi_{\mc F}^* {\mc F}$ is the tautological section, sends FOP transverse sections to FOP transverse sections.

\item {\bf (Stratified Regularity)} For each ${\mc S}\in \Gamma^{\rm FOP}({\mc U}, {\mc E})$ and each normally complex isotropy type $\gamma$,  if we write ${\mc S}|_{{\mc U}_\gamma}$ as the direct sum $ {\mc S}_\gamma \oplus \check{\mc S}_\gamma$, then $ {\mc S}_\gamma: {\mc U}_\gamma \to   {\mc E}_\gamma$ is transverse in the usual sense.
\end{enumerate}
In addition, the structure of the zero loci of FOP transverse sections can be described as follows. For each integer $k$, there exists a countable partially ordered set ${\mf Z}_k^{\rm univ}$ (called the set of {\it universal strata}), an order-preserving map $\rho_k: {\mf Z}_k^{\rm univ} \to \uds{\mb \Gamma}_k^{\rm NC}$, and a strictly order-preserving\footnote{Namely, $\theta_1 < \theta_2$ implies $n_k (\theta_1) < n_k(\theta_2)$.} map $n_k : {\mf Z}_k^{\rm univ} \to 2{\mb Z}$ (cf. Definition \ref{defn:n-map}) (all independent of ${\mc U}$ and ${\mc E}$) satisfying the following conditions.
\begin{enumerate}

\item For any maximal element $\theta \in {\mf Z}_k^{\rm univ}$, if $\rho_k (\theta) \in \uds{\mb \Gamma}_k^{\rm NC}$ is represented by a triple $(G, V, W)$, then $n_k (\theta) = - {\rm dim}_{\mb R} (\check{V}_G) + {\rm dim}_{\mb R} (\check{W}_G)$.

\item There is a special maximal $\theta_k \in {\mf Z}_k^{\rm univ}$ with $\rho_k (\theta_k) = \gamma_k^* \in \uds{\mb \Gamma}_k^{\rm NC}$.

\item Suppose ${\rm dim} ({\mc U}) - {\rm rank} ({\mc E}) = k$. For ${\mc S}\in \Gamma^{\rm FOP}({\mc U}, {\mc E})$, ${\mc S}^{-1}(0)$ admits a partition
\beq\label{zero_set_partition}
{\mc S}^{-1}(0) = \bigsqcup_{\theta \in {\mf Z}_k^{\rm univ}} {\mc S}^{-1}(0)_\theta
\eeq
satisfying the following conditions.
\begin{enumerate}

\item Each ${\mc S}^{-1}(0)_\theta$ is a $k + n_k (\theta)$-dimensional smooth submanifold of ${\mc U}_{\rho_k (\theta)}$. 

\item After removing empty pieces, the partition \eqref{zero_set_partition} makes ${\mc S}^{-1}(0)$ a Thom--Mather stratified space (Definition \ref{defn_Thom_Mather}) whose natural partial order among strata coincides with the partial order of ${\mf Z}_k^{\rm univ}$. 
\end{enumerate}
\end{enumerate}



\end{thm}

\begin{rem}
The FOP transversality condition arises from the original proposal of Fukaya--Ono \cite{Fukaya_Ono_integer} and the crucial contribution of B. Parker \cite{BParker_integer}. Many difficult points in proving the desired properties of FOP transverse perturbations have their roots in \cite{BParker_integer}. Besides these works, the notion of normally polynomial perturbations introduced in \cite{Fukaya_Ono_integer} was used in \cite{FOOO_Chapter8}. 
\end{rem}

\begin{rem}\label{rem_boundary}
Theorem \ref{thm11} (as well as a few other technical results proved in this paper) is stated for orbifolds without boundary. In applications, one needs to consider orbifolds with boundary. In different scenarios, the corresponding extensions are addressed specifically regarding the context. However the general principle is very simple: on an orbifold with boundary, it is easy to construct collar neighborhoods of the boundary. The structures we need to construct (such as an FOP transverse section), can be first chosen on the boundary, then extended to the collar neighborhood by a trivial product. Properties such as the {\bf (Extension Property I)} of Theorem \ref{thm11} allow us to further extend to whole interior of the orbifold.   
\end{rem}

We briefly explain the spirit of the above conditions. Conditions (1)---(5) show that the FOP transversality is almost as flexible as the classical transversality on manifolds. In particular, the condition {\bf (Extension Property II)} is useful when considering incidence conditions for evaluation maps into manifolds; for example, moduli spaces of pseudoholomorphic curves in a symplectic manifold satisfying constraints at marked points are suborbifolds with ordinary normal bundles inside the moduli spaces of curves without such constraints. Regarding Condition (6), any finite-dimensional reduction of an infinite-dimensional problem is only well-defined up to stabilization; this condition is stating that the transversality condition is invariant under stabilization. Condition (7) says that the FOP transversality condition implies the stratified transversality.

\subsubsection{FOP pseudocycles}

To define integral cycles, it is convenient to introduce the notion of derived orbifolds. They are models for regularizing various kinds of moduli spaces of pseudoholomorphic curves used in symplectic geometry.

\begin{defn}\label{def:dorbchart}
A \emph{derived orbifold} 
is a triple ${\mc D} = ({\mc U}, {\mc E}, {\mc S})$ where ${\mc U}$ is an orbifold, ${\mc E} \rightarrow {\mc U}$ is a  vector bundle and ${\mc S}: {\mc U} \rightarrow {\mc E}$ is a continuous 
section. ${\mc D}$ is said to be \emph{compact} if ${\mc S}^{-1}(0)$ is compact; is \emph{normally complex} if both ${\mc U}$ and ${\mc E}$ are normally complex; is \emph{oriented} if the virtual vector bundle $T{\mc U} - {\mc E}$ is oriented. The {\it virtual dimension} of ${\mc D}$ is defined as ${\rm dim} ({\mc U}) - {\rm rank} ({\mc E})$. 
\end{defn}

Given a compact, oriented, and normally complex derived orbifold ${\mc D}= ({\mc U}, {\mc E}, {\mc S}_0)$ of virtual dimension $k$, consider an FOP transverse section ${\mc S} \in \Gamma^{\rm FOP}({\mc U}, {\mc E})$ which is $C^0$-close to ${\mc S}_0$. Assuming compactness and orientation, by Theorem \ref{thm11}, ${\mc S}^{-1}(0)$ is a compact Thom--Mather stratified space with oriented strata. Moreover, as the dimensions of any pair of adjacent strata differ by an even number, the closure $\ov{{\mc S}^{-1}(0)_\theta}$ is an oriented pseudomanifold (Definition \ref{defn_TM_pseudomanifold}), hence carries a fundamental class
\beqn
[{\mc S}^{-1}(0)_\theta] \in H_{k + n_k(\theta)} ( \ov{{\mc S}^{-1}(0)_\theta}; {\mb Z}) \to H_{k + n_k(\theta)} ({\mc U}; {\mb Z}).
\eeqn
We call this the {\it $\theta$-th FOP Euler cycle} of the derived orbifold ${\mc D}$. As a theorem we prove that it is independent of the choice of the FOP transverse perturbations.

\begin{thm}
Given a compact, oriented, and normally complex derived orbifold ${\mc D} = ({\mc U}, {\mc E}, {\mc S}_0)$ of virtual dimension $k$, for any $C^0$-small FOP transverse perturbation ${\mc S}$ of ${\mc S}_0$, for any universal stratum $\theta \in {\mf Z}_k^{\rm univ}$, the integral class $[{\mc S}^{-1}(0)_\theta]$ is independent of the choice of ${\mc S}$. 
\end{thm}

\begin{rem}
The reader would be curious about the relation between the rational-valued Euler cycles and the various FOP Euler cycles as it would imply a decomposition of the Gromov--Witten invariants into enumerative invariants in different isotropy types. It seems plausible to connect them by further perturbing an FOP transverse section via multivalued transverse sections. This will be discussed in forthcoming work.
\end{rem}

\subsection{Intuitions about the FOP transversality}

\begin{figure}[h]
    \centering
\includegraphics[width=\linewidth]{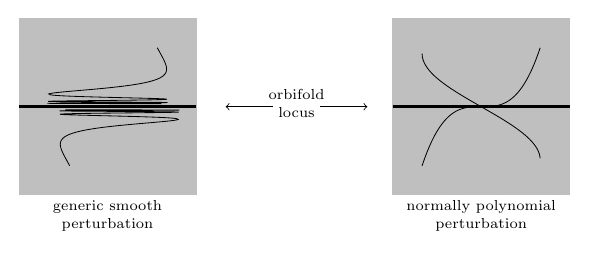}
    \caption{Smooth perturbations and FOP transverse perturbations.}
    \label{fig_FOP}
\end{figure}

We would like to provide an intuitive understanding of the FOP transversality condition. Let ${\mc E} \to {\mc U}$ be an orbifold vector bundle. As we have mentioned, the best transversality condition one can hope for a generic smooth section ${\mc S}: {\mc U} \to {\mc E}$ is the stratified transversality. As a result, the zero locus is the union of smooth manifolds contained in each stratum of the isotropy stratification. However, in the normal direction of each stratum, the transversality is typically lost. As a result, we do not know how the frontier of the zero locus in a higher stratum behaves when approaching to a lower stratum, which could be {\it a priori} pathological (see the left picture of Figure \ref{fig_FOP}). 

One idea to regularize the frontier of a stratum of the zero locus is to consider perturbations in a category more restricted than the smooth category. Fukaya--Ono's suggestion in \cite{Fukaya_Ono_integer} essentially says that one can require the section to be algebraic in normal directions of each stratum. Then the frontier of the zero locus of each stratum looks like a smooth family of algebraic varieties (see the right picture of Figure \ref{fig_FOP}). Although the ordinary transversality still fails, this type of singular space can be characterized using languages of Whitney stratified space or Thom--Mather stratified space. 

Using the normally polynomial type section, the zero locus then becomes a Thom--Mather stratified space. In order to support a fundamental class rather than a chain, one needs to guarantee that adjacent strata differ in dimensions at least by 2. The assumption of the normal complex structure makes this work, as the frontier of each stratum looks like a family of complex algebraic varieties.

\subsection{Applications in symplectic topology}

\subsubsection{Integral Gromov--Witten type invariants}
Let $(X, \omega)$ be a compact symplectic manifold. Let $A \in H_2(M; {\mb Z})$ be a homology class. Let $g, n \geq 0$ be nonnegative integers. Upon choosing a compatible almost complex structure $J$ on $X$, one obtains the moduli space $\ov{\mc M}{}_{g,n}(X, J, A)$ of $J$-holomorphic stable maps of genus $g$, degree $A$, with $n$ marked points. Let $k$ be its virtual dimension. The following theorem is proved in Section \ref{section_GW}.

\begin{thm}\label{thm16}
For each element $\theta \in {\mf Z}_k^{\rm univ}$, there is a symplectic deformation invariant 
\beq\label{eqn13}
[\ov{\mc M}_{g, n}(X, J, A)]_{\theta}^{vir} \in H_{k + n_k (\theta)} ( \ov{\mc M}_{g, n} \times X^n; {\mb Z}).
\eeq
When $\theta$ is the maximal element $\theta_k$ provided in Theorem \ref{thm11}, this class ``morally'' counts $J$-holomorphic curves in class $A$ whose automorphism group is trivial. More precisely, when all elements of $\ov{\mc M}{}_{g,n}(X, J, A)$ have trivial automorphism group (hence the ordinary virtual fundamental class is integral),  the above virtual class for $\theta_k$ coincides with the ordinary virtual class. 
\end{thm}



\begin{rem}
When $(X, \omega)$ is semi-positive, $g = 0$, and $n \geq 3$, the numerical invariants (correlation function) defined via the class \eqref{eqn13} corresponding to the maximal element $\theta_k$ should coincide with the usual genus zero Gromov--Witten invariant of $(X, \omega)$ defined by Ruan \cite{Ruan_96} and Ruan--Tian \cite{Ruan_Tian, Ruan_Tian_97} (see also \cite{McDuff_Salamon_2004}). A proof of this claim will be provided in forthcoming work. 
\end{rem}

\begin{rem}
One may expect the integral virtual class $[\ov{\mc M}{}_{g,n}(X, J, A)]_{\theta_k}^{vir}$ to define a cohomological field theory satisfying the Kontsevich--Manin axioms. However, there are a few notable exceptions. First, as the FOP perturbation scheme is sensative to the symmetry group, the $S_n$-invariance axiom (with respect to the permutation of marked points) may not be true. Second, the axioms related to forgetting marked points (the fundamental class axiom and the divisor axiom) may not hold; this is because when forgetting a marked point, the isotropy types of curves may change and it is difficult to find FOP transverse perturbations which respect forgetting marked points. Third, the splitting axiom must be modified. Indeed, the usual splitting axiom of Gromov--Witten invariants is stated using the classical Poincar\'e duality which is not quantum deformed. However, this is not a universal truth; for example, the splitting axiom of quantum K-theory (see \cite{YPLee_QK} or the theories recently constructed by Abouzaid--McLean--Smith \cite{AMS2} both use quantum deformed Poincar\'e duality. For the FOP invariants, one also needs to deform the classical Poincar\'e duality in order to have a splitting property. The discussion of Kontsevich--Manin axioms for the FOP Gromov--Witten invariants will be given in forthcoming work. 
\end{rem}

\subsubsection{Homological splitting for Hamiltonian loops}

The second main result of this paper is an alternate proof of the main theorem of \cite{AMS} which is conceptually much simpler than the original one (see Section \ref{sec:proof-AMS}).

\begin{thm}(cf. \cite[Theorem 1.1]{AMS})\label{thm15}
Let $P \to S^2$ be a smooth Hamiltonian fibration over the $2$-sphere with fiber being a closed symplectic manifold $(X, \omega)$. Then there is an isomorphism of graded abelian groups
\begin{equation}\label{eqn:intro-split}
H^*(P; \mathbb{Z}) \cong H^*(X; \mathbb{Z}) \otimes_{\mathbb{Z}} H^*(S^2; \mathbb{Z}).
\end{equation} 
\end{thm}

\begin{rem}
Although \cite{AMS} shows that the splitting \eqref{eqn:intro-split} actually holds for any complex oriented generalized cohomology theory, for the case over $\mathbb{Z}$, our short proof uses only the virtual cycle underlying the class \eqref{eqn13} indexed by the trivial group without appealing to Morava K-theories or Atiyah duality for orbifolds. 
\end{rem}

\subsection{Relations with other work}

\subsubsection{$\mathbb{Z}$-valued enumerative invariants in symplectic topology}

There are other known $\mathbb{Z}$-valued Gromov--Witten type invariants in symplectic topology. Firstly, as mentioned above, for semi-positive symplectic manifolds, Ruan--Tian's construction (\cite{Ruan_Tian}) shows that the genus-$0$ Gromov--Witten invariants with at least three insertions are indeed $\mathbb{Z}$-valued. Secondly, there are the celebrated $\mathbb{Z}$-valued Gopakumar--Vafa invariants of Calabi--Yau $3$-folds \cite{Gupakumar_Vafa_1}, \cite{Gupakumar_Vafa_2} which govern the $\mathbb{Q}$-valued Gromov--Witten invariants as demonstrated in full generality by Ionel--Parker \cite{IP-GV}. It was speculated by some experts (e.g.\ \cite[Section 6.3]{joyce2007kuranishi}) that the invariants defined using Fukaya--Ono's normally polynomial sections are related to the Gopakumar--Vafa invariants. However, recent developments \cite{doan2023counting}, \cite{bai-swaminathan}, which should be viewed as continuations of the groundbreaking work of Taubes in dimension $4$ \cite{taubes-gr}, construct $\mathbb{Z}$-valued invariants of Calabi--Yau $3$-folds by (virtually) counting \emph{embedded} pseudo-holomorphic curves and these invariants seems to be better connected with the Gopakumar--Vafa invariants. It is an intriguing question to understand the relation between these two seemingly distinct proposals for geometrically constructing the Gopakumar--Vafa invariants. Lastly, the $K$-theoretic Gromov--Witten invariants, defined by Givental and Lee \cite{YPLee_QK} for algebraic objects (which are expected to exist for general symplectic manifolds in light of \cite[Section 6.12]{AMS}), are also $\mathbb{Z}$-valued.



\subsubsection{$\mathbb{Z}$-valued invariants in algebraic geometry}

There are several types of integral enumerative
invariants constructed using algebraic geometry, most notably the Donaldson--Thomas invariants \cite{Thomas-DT} and the Pandharipande--Thomas invariants \cite{PT}. Although these invariants are closely tied with the Gromov--Witten invariants \cite{MNOP-I, MNOP-II}, their constructions depend on sheaf theory. Fukaya--Ono's proposal is differential-topological in nature, as it is still a variant of the general position argument. It would be interesting to see if the virtual classes $[\overline{\mathcal{M}}_{0,n}(X,J, A)]^{\text{vir}}_{\theta}$ admit a purely algebro-geometric interpretation. 

\subsubsection{Wasserman's theorem and stable homotopy theory}

A renowned theorem of Wasserman \cite{wasserman}, recaptured in \cite[Theorem 5.6]{pardon2023orbifold}, provides a sufficient condition for equivariant transversality to hold. Our result can be interpreted as a variant of Wasserman's theorem given the presence of normal complex structures, whose existence allows us to extract more information. 
The homotopical cobordism perspective gives far-reaching corollaries of Wasserman's theorem, see \cite[Theorem 6.2.33]{schwede}. As mentioned by \cite[Remark 5.7]{pardon2023orbifold}, it is an interesting question to understand Fukaya--Ono's proposal on the homotopical cobordism side.

\begin{addendum}
    After the first version of this paper was posted on arXiv, the results and methods herein have been used for several applications. As mentioned above, the integral counting scheme is a crucial tool for our proof of the integral Arnold conjecture \cite{Bai_Xu_Arnold}; we also obtained cohomological splitting result for Hamiltonian fibrations over general bases in \cite{BPX} using our integral Gromov--Witten type invariants. More recently, the FOP perturbation method is combined with resolution of singularities to construct complex-cobordism-valued Gromov--Witten invariants \cite{abouzaid-bai}. In forthcoming work \cite{Bai_Xu_Floer}, we will also derive applications to Hamiltonian dynamics based on constructing quantum power operations on general symplectic manifolds.
\end{addendum}

\subsection{Plan of the paper}

We start with basic notions related to effective orbifolds in Section \ref{sec:orbifolds}, especially the notion of \emph{normal complex structure} (Definition \ref{defn_normal_complex}) and \emph{normally complex sections} (Definition \ref{defn:FOP}). In Section \ref{section3}, we study the canonical Whitney stratification on the universal zero locus $Z^d(G, V, W)$ and understand its behavior under various operations. This is the technical core of this paper. The definition of the FOP transversality condition and the proof of Theorem \ref{thm11} is given in Section \ref{section4}. In Section \ref{section5}, we provide the homological consequence of the FOP transversality condition. In Section \ref{section_GW} and Section \ref{sec:proof-AMS} we provide our geometric applications and prove Theorem \ref{thm16} and Theorem \ref{thm15}. 

In Appendix \ref{appendixa} we provide the proof of technical results about straightenings on orbifolds. In Appendix \ref{appendixb} we provide technical details on Whitney stratifications supporting the proofs in Section \ref{section3}. In Appendix \ref{appendixc} we prove that zero loci of FOP transverse sections are Thom--Mather stratified spaces, which is needed for constructing the integral Euler cycles. 

\subsection{Acknowledgements}

We would like to thank Mohammed Abouzaid, Kenji Fukaya, Mark Goresky, Helmut Hofer, Eleny Ionel, John Pardon, Paul Seidel, Mohan Swaminathan, and Dingyu Yang for useful correspondences and discussions at various stages of this project. The first-named author would like to thank his Ph.D. advisor John Pardon for constant encouragement and support. The second-named author would like to thank his family for their love and support during the Covid-19 pandemic, and thank Princeton University and Chenyang Xu for hosting the second-named author's visit in Fall 2024. 

The authors would also thank the anonymous referees for carefully examining our submission and making many important suggestions.

\section{Normally Complex Orbifolds and Bundles}\label{sec:orbifolds}

In this section we recall basic notions about effective orbifolds and the prerequisite structures for defining the FOP transversality condition. In Subsection \ref{subsection21} we recall notions and facts about effective orbifolds and orbifold vector bundles. In Subsection \ref{subsection22} we specify the stratification of orbifolds (together with a vector bundle) induced by the isotropy types of points. In Subsection \ref{subsection23} we define the notion of normal complex structures (which comes from the original proposal of Fukaya--Ono \cite{Fukaya_Ono_integer}). In Subsection \ref{subsection24} we describe the auxiliary and technical structures called ``straightening'' which can be viewed as a system of normal tubular neighborhoods for the strata of the isotropy stratification. The techical construction of straightenings is deferred to Appendix \ref{appendixa}. In Section \ref{subsection25} we define the notion of normally complex sections which have more rigid behavior in normal directions to each stratum than general smooth sections.

\subsection{Orbifolds and orbifold vector bundles}\label{subsection21}

We recall the basic definition of effective orbifolds. We follow the definition of \cite[Section 1.1]{Adem_Leida_Ruan}. Working with effective orbifolds allows us to use orbifold charts exclusively without appealing to the language of groupoids. 

Let ${\mc U}$ be a Hausdorff and second countable topological space. An {\it $n$-dimensional orbifold chart} of ${\mc U}$ is a triple
\beqn
C = (G, U, \psi)
\eeqn
where $U$ is a nonempty smooth manifold, $G$ is a finite group acting effectively and smoothly on $U$, and $\psi: U \to {\mc U}$ is a $G$-invariant continuous map such that the induced map 
\beqn
\uds \psi: U/G \to {\mc U}
\eeqn
is a homeomorphism onto an open subset of ${\mc U}$. If $p \in \psi(U)$ we also say that $p$ is {\it contained} in the chart $C$. If $\psi^{-1}(p) \subset U_G$, the fixed point locus of the $G$-action, then we say that the chart $C$ is {\it centered at $p$}.

A {\it chart embedding} from another chart $C' = (G', U', \psi')$ to $C$ is a smooth open embedding $\iota: U' \hookrightarrow U$ such that
\beqn
\psi \circ \iota = \psi'.
\eeqn
It follows that (see \cite[Page 3]{Adem_Leida_Ruan}) given a chart embedding $\iota$ as above there exists a canonical group injection $G' \hookrightarrow G$ such that $\iota$ is equivariant. Therefore we often include the group injection as part of the data of a chart embedding. 

We say two charts $C_i = (G_i, U_i, \psi_i)$, $i = 1, 2$ are {\it compatible} if for each $p \in \psi_1(U_1) \cap \psi_2(U_2)$, there exists an orbifold chart $C_p = (G_p, U_p, \psi_p)$ containing $p$ and chart embeddings into both $C_1$ and $C_2$.

An {\it orbifold atlas} ${\mc A} = \{C_i\ |\ i \in I\}$ on $X$ is a collection of mutually compatible charts $C_i$ which cover ${\mc U}$. We say an atlas ${\mc A}' = \{ C_j'\ |\ j \in J\}$ is a {\it refinement} of ${\mc A}$ if for each $C_j'$ there exists a chart embedding $C_j' \hookrightarrow C_i$ for some $i \in I$. We say two orbifold atlases are {\it equivalent} if they have a common refinement. A Hausdorff and second countable topological space ${\mc U}$ together with an equivalence class of orbifold atlases is called a {\it smooth effective orbifold}. Every smooth effective orbifold has a unique maximal atlas; two atlases are equivalent if they are contained in the common maximal atlas (see \cite[Page 3]{Adem_Leida_Ruan}). It is convenient to work with the maximal atlas. From now on, an orbifold chart of a smooth effective orbifold means a chart in the maximal atlas. 

We often use $|{\mc U}|$ to denote the underlying topological space (called the {\it coarse space}) of an effective orbifold ${\mc U}$ while forgetting the orbifold structure. 

We do not define the general form of orbifold morphisms. Below are a few special cases of morphisms between orbifolds. First, a continuous function on an effective orbifold (or more generally a continuous map into a smooth manifold) is {\it smooth} if its pullback to each chart is a smooth function. On the other hand, an {\it isomorphism} of orbifolds from ${\mc U}$ to ${\mc U}'$ is a homeomorphism $f: |{\mc U}| \to |{\mc U}'|$ such that the correspondence
\beqn
(G, U, \psi) \mapsto (G, U, f \circ \psi) 
\eeqn
is a one-to-one correspondence between their maximal atlases. An {\it open embedding} from ${\mc U}$ to ${\mc U}'$ is an isomorphism from ${\mc U}$ to an open subset of ${\mc U}'$.

\begin{rem}\label{rem21}
One can see that orbifolds are all locally compact. As we also assume they are Hausdorff and second countable, they are paracompact spaces. Hence for any open cover by charts, there exists a subordinate continuous partition of unity; as one can approximate continuous functions by smooth functions on each chart, there always exist a subordinate smooth partition of unity. 
\end{rem}

\begin{rem}
We also need the notion of orbifolds with boundary. In that case, the domain of a chart $C = (G, U, \psi)$ is allowed to be a smooth manifold with boundary and the group $G$ acts by diffeomorphisms of $(U, \partial U)$. As a consequence, the fixed point locus $U_G$ is a manifold with boundary $\partial U_G = U_G \cap \partial U$.
\end{rem}

The definition of orbifold vector bundles is very similar to that of effective orbifolds. Let ${\mc U}$ be an effective orbifold, ${\mc E}$ be a topological space, and $\pi_{\mc E}: {\mc E} \to {\mc U}$ be a continuous map. A {\it bundle chart} of $\pi_{\mc E}: {\mc E} \to {\mc U}$ consists of an orbifold chart $C = (G, U, \psi)$ of ${\mc U}$, a $G$-equivariant smooth vector bundle $\pi_E: E \to U$, and a $G$-invariant continuous map $\hat \psi: E \to {\mc E}$ such that the induced map from $E/G$ to ${\mc E}$ is a homeomorphism onto $\pi_{\mc E}^{-1}(\psi(U))$ and such that the following diagram commutes:
\beqn
\vcenter{ \xymatrix{      E \ar[r]^{\hat \psi} \ar[d]_{\pi_E} & {\mc E} \ar[d]^{\pi_{\mc E}} \\
                U  \ar[r]_{\psi} &                     {\mc U}} }.
                \eeqn
In notation we will use a quadruple $ \hat C = (G, U, E, \hat \psi)$ to denote the bundle chart where the map $\psi: U \to {\mc U}$ is determined by the map $\hat\psi: E \to {\mc E}$. If $\hat C' = (G', U', E', \hat\psi')$ is another bundle chart, a {\it bundle chart embedding} from $\hat C'$ to $\hat C$ consists of an orbifold chart embedding $\iota: U' \hookrightarrow U$ (equivariant with respect to a group injection $G' \hookrightarrow G$) covered by a vector bundle embedding $\hat\iota: E' \to E$ such that 
\beqn
\hat\psi \circ \hat \iota = \hat\psi'.
\eeqn
We can similarly define the notions of compatibility between bundle charts and bundle atlases. Then an {\it orbifold vector bundle structure} over $\pi_{\mc E}: \mathcal{E} \rightarrow {\mc U}$ is defined to be an equivalence class of bundle atlases as before.  Similarly, an orbifold vector bundle has a unique maximal atlas and two atlases are equivalent if and only if they are contained in a common maximal atlas. A bundle chart  then means a chart in the maximal atlas. 

\begin{defn}
An orbifold vector bundle ${\mc E} \to {\mc U}$ is said to be {\it ordinary} if for each bundle chart $\hat C = (G, U, E, \hat\psi)$ and any $x \in U$, $G_x$ acts trivially on the fiber $E_x$.
\end{defn}

We do not attempt to deal with general suborbifolds in this paper. Instead we consider the following very restricted type of suborbifolds.

\begin{defn}\label{defn_ordinary_normal}
Let ${\mc U}$ be an effective orbifold. A {\it closed suborbifold with ordinary normal bundle} of ${\mc U}$ is a closed subset ${\mc X} \subset |{\mc U}|$ such that for each chart $C = (G, U, \psi)$ of ${\mc U}$, $\psi^{-1}({\mc X})$ is a closed submanifold of $U$ such that for each $x \in \psi^{-1}({\mc X})$, $G_x$ acts trivially on the fiber of the normal bundle at $0$. In particular, $(G, \psi^{-1}({\mc X}), \psi|_{\psi^{-1}({\mc X})})$ is an orbifold chart of ${\mc X}$ and all such charts make ${\mc X}$ itself an effective orbifold.
\end{defn}

We spell out the definition of sections of an orbifold vector bundle because of their importance in this paper.

\begin{defn}[Sections]\label{defn_section} Let ${\mc E} \to {\mc U}$ be an orbifold vector bundle.
\begin{enumerate}
\item Let $\hat C_i = (G_i, U_i, E_i, \hat\psi_i)$, $i = 1, 2$ be two bundle charts. We say that a $G_1$-equivariant section $S_1: U_1 \to E_1$ and a $G_2$-equivariant section $S_2: U_2 \to E_2$ are {\it compatible} if for any bundle chart $\hat C_0 = (G_0, U_0, E_0, \hat\psi_0)$ of $E$ and chart embeddings  $\hat\iota_i: \hat C_0 \hookrightarrow \hat C_i$, $i = 1, 2$ there holds  
\beqn
\hat \iota_1^{-1} \circ S_1 \circ \iota_1 = \hat \iota_2^{-1} \circ S_2 \circ \iota_2
\eeqn
as sections of $E_0 \to U_0$. 

\item A {\it section} of ${\mc E}$, denoted by ${\mc S}: {\mc U} \to {\mc E}$, is a collection of mutually compatible $G_i$-equivariant sections $S_i: U_i \to E_i$ for all bundle charts $\hat C_i$ belonging to the maximal atlas of ${\mc E}$.  


\end{enumerate}
\end{defn}

On each single chart $\hat C = (G, U, E, \hat\psi)$, there are a lot of $G$-equivariant sections $S: U \to E$. The existence of partitions of unity implies that any orbifold vector bundle over an effective orbifold has a lot of smooth sections. 


As a final topic of this subsection, we recall the definition of Riemannian metrics on orbifolds and connections on orbifold vector bundles.

\begin{defn}
Let ${\mc U}$ be an effective orbifold (possibly with boundary). 
\begin{enumerate}

\item A {\it Riemannian metric} on ${\mc U}$, denoted by $g^{T{\mc U}}$, is a collection $g^{TU}$ for all orbifold charts $C = (G, U, \psi)$ where $g^{TU}$ is a $G$-invariant Riemannian metric on $U$, such that every chart embedding is isometric.

\item Let ${\mc F} \to {\mc U}$ be an orbifold vector bundle. A {\it connection} resp. {\it inner product} on ${\mc F}$, denoted by $\nabla^{\mc F}$ resp. $h^{\mc F}$, is a collection $\nabla^F$ resp. $h^F$ for all bundle charts $\hat C = (G, U, F, \hat\psi)$ where $\nabla^F$ resp. $h^E$ is a $G$-invariant connection resp. inner product on $F$ such that all bundle chart embeddings preserve the connections resp. inner products. We say that $\nabla^{\mc F}$ {\it preserves} $h^{\mc F}$ if $\nabla^F$ preserves $h^F$ on each chart. 

\end{enumerate}
\end{defn}

One can use the standard way to construct Riemannian metrics, connections, and inner products using smooth partition of unity on orbifolds (see Remark \ref{rem21}).

\subsection{The isotropy stratification}\label{sec:iso-strat}\label{subsection22}

 The isomorphism class of the stabilizer of a point $x$ in an orbifold chart $C = (G, U, \psi)$ only depends on the point $\psi(x) \in {\mc U}$ (see \cite[Page 3]{Adem_Leida_Ruan}). 
One can use this isomorphism class to stratify the orbifold. Moreover, together with the data from a vector bundle, we can actually obtain a more refined decomposition.

We first spell out certain frequently used notions. 

\begin{defn}
Let $G$ be a group acting on a set $A$. A subgroup $H \subset G$ is called {\it $A$-essential}, denoted by $H \subset_A G$, if there exists $a \in A$ with isotropy group $H$. 
\end{defn}

Notice that $G \subset_A G$ if and only if the set of fixed points $A_G \neq \emptyset$; $\{e\} \subset_A G$ if and only if the action is effective. The set of $A$-essential subgroups is closed under conjugations. If $(G, U, \psi) \to (G', U', \psi')$ is a chart embedding of an orbifold, then the map $G \to G'$ sends $U$-essential subgroups of $G$ to $U'$-essential subgroups of $G'$.

\begin{defn}\label{defn_isotropy_type}(Isotropy triples and isotropy types)

\begin{enumerate}

\item An {\it isotropy triple} is a triple $(G, V, W)$ where $G$ is a finite group and $V, W$ are finite-dimensional real representations of $G$. 

\item An {\it isomorphism} of isotropy triples from $(G, V, W)$ to $(G', V', W')$ consists of a group isomorphism $G \cong G'$ and equivariant linear isomorphisms $V \cong V'$ and $W \cong W'$. An {\it isotropy type} is an isomorphism class of isotropy triples, typically denoted by $\gamma, \delta$, etc. Let ${\mb\Gamma}$ be the set of all isotropy types. 

\item A partial order on ${\mb\Gamma}$ is defined as follows. We write $\gamma \leq \delta$ if $\gamma$ can be represented by $(G, V, W)$, $\delta$ can be represented by $(H, V, W)$ such that $H$ is isomorphic to a $V$-essential subgroup of $G$ and such that $X \cong V, Y \cong W$ as representations of $H$.
\end{enumerate}
\end{defn}

\begin{defn}\label{basic_decomposition}
Suppose $W$ is a finite-dimensional real representation of a finite group $G$. Then $W$ can be decomposed as the direct sum of irreducible representations. We call the (canonical) decomposition 
\beqn
W = W_G \oplus \check{W}_G
\eeqn
where $W_G$ is the direct sum of all trivial summands  and $\check W_G$ is the direct sum of all nontrivial summands the \emph{basic decomposition} of $W$ with respect to $G$. Notice that when $W$ is a complex representation, the basic decomposition is an isomorphism of complex vector spaces.

More generally, if $U$ is a $G$-manifold and $E \to U$ is a $G$-equivariant vector bundle, then over the $G$-fixed points $U_G \subset U$, the fiberwise basic decomposition of $E|_{U_G}$ induces a decomposition of vector bundles (possibly with varying ranks on different connected components)
\beqn
E = E_G \oplus \check E_G
\eeqn
where $E_G \subset E|_{U_G}$ coincides with the set of $G$-fixed points of $E$.
\end{defn}

Now consider an effective orbifold ${\mc U}$ with an orbifold vector bundle ${\mc E} \to {\mc U}$. For each $p \in {\mc U}$, consider a bundle chart $\hat C = (U, G, E, \hat \psi)$ containing $p$ with $x \in \psi^{-1}(p)\subset U$.  Define $\gamma_p$ to be the isotropy type represented by the triple $(G_x, T_x U, E_x)$. Then $\gamma_p$ is well-defined, independent of the choice of charts and the point $p\in \psi^{-1}(x)$. 

\begin{lemma}\label{lem:iso-prestrat}
For each isotropy type $\gamma$, the subset 
\beqn
{\mc U}_\gamma:= \{ p \in {\mc U}\ |\ \gamma_p = \gamma\}
\eeqn
is a topological manifold (possibly with boundary if ${\mc U}$ has boundary) equipped with a natural smooth structure.\footnote{It is more ``correct'' to view ${\mc U}_\gamma$ as an non-effective orbifold. However we refrain from discussing non-effective cases in this paper.}
\end{lemma}

\begin{proof}
By definition, for each $p\in {\mc U}_\gamma$, there exists a bundle chart $(G, U, E, \hat\psi)$ and $x \in \psi^{-1}(p) \subset U$ such that $G = G_x$ and $\gamma$ is represented by $(G, T_x U, E_x)$. As $U_G \subset U$ is a closed submanifold, it implies that ${\mc U}_\gamma$ is a topological manifold and $\psi^{-1}: {\mc U}_\gamma \cap \psi(U) \to U_G \subset U$ is a local chart. 

To put a smooth structure on ${\mc U}_\gamma$ one needs to choose an atlas of $C^\infty$-compatible charts. Let $\hat {\mc A}_\gamma$ be the set of bundle charts $\hat C = (G, U, E, \hat \psi)$ of ${\mc E}$ such that for each $x \in U_G \subset U$, the triple $(G, T_x U, E_x)$ represents $\gamma$ (so that $\psi(U) \cap {\mc U}_\gamma \neq \emptyset$). Define the atlas
\beq
\Big\{ \psi^{-1}: {\mc U}_\gamma \cap \psi(U) \to U_G \ |\ \hat C = (G, U, E, \hat \psi) \in \hat {\mc A}_\gamma \Big\}.
\eeq
We claim that this atlas is $C^\infty$-compatible, namely, if $\hat C_i = (G_i, U_i, E_i, \hat\psi_i)$, $i = 1, 2$, are two such bundle charts of ${\mc E}$, then the map 
\beqn
\psi_2^{-1} \circ \psi_1: \psi_1^{-1} \Big( {\mc U}_\gamma \cap \psi_1(U_1) \cap \psi_2(U_2) \Big) \to U_{2, G_2}
\eeqn
is smooth. Indeed, the smoothness at each point $q$ of the domain of this transition map can be checked by using a smaller bundle chart containing $q$ which embeds into both $\hat C_1$ and $\hat C_2$. Therefore, ${\mc U}_\gamma $ is equipped with a smooth structure.
\end{proof}

Now we recall the definition of the basic concept of stratification, which is frequently used in this paper. 

\begin{defn}\label{defn_stratification}
Let $X$ be a topological space. 
\begin{enumerate}

\item A {\it stratification} of $X$ is a set ${\mf X}$ of subsets of $X$ satisfying
\begin{enumerate}
\item $X$ is equal to the disjoint union of all elements $X_\alpha \in {\mf X}$.

\item Each $X_\alpha \in {\mf X}$ (called a {\it stratum}) is locally closed and nonempty.

\item ${\mf X}$ is locally finite.

\item (Axiom of frontier) If $X_\alpha \cap \ov{X_\beta} \neq \emptyset$, then $X_\alpha \subset \ov{X_\beta}$.
\end{enumerate}

\item Notice that the last condition induces a partial order relation among strata:
\beqn
X_\alpha \leq X_\beta \Longleftrightarrow X_\alpha \subset \ov{X_\beta}.
\eeqn
A {\it top stratum} of a stratification ${\mf X}$ is a stratum $X_\alpha$ which is maximal with respect to the above partial order, or equivalently, an open stratum.

\item A homeomorphism $f: X \to X$ is said to {\it preserve} a stratification ${\mf X} = \{ X_\alpha \}$ on $X$ if ${\mf X} = \{ X_\alpha \} = \{ f^{-1}(X_\alpha )\}$ as sets.

\item Let $X$ resp. $Y$ be equipped with a stratification ${\mf X} = \{ X_\alpha \}$ resp. ${\mf Y} = \{ Y_\beta \}$. The {\it product stratification} on $X \times Y$ is the set 
\beqn
{\mf X}  \times {\mf Y} = \{ X_\alpha \times Y_\beta \ |\ X_\alpha \in {\mf X}, Y_\beta \in {\mf Y} \}.
\eeqn
\end{enumerate}
\end{defn}

The following lemma follows from inspecting the definition. 
\begin{lemma}\label{lemma_isotropy_stratification}
The set $\{ {\mc U}_\gamma\neq \emptyset\ |\ \gamma \in {\mb\Gamma} \}$ is a stratification of ${\mc U}$, called the {\it isotropy stratification} of $({\mc U}, {\mc E})$. Moreover, the partial order induced from the stratification coincides with the restriction of the partial order in ${\mb\Gamma}$. \qed
\end{lemma}


\subsection{Normal complex structures}\label{subsec:ncs}\label{subsection23}

Now we introduce the most important geometric condition which plays the central role in our construction. In applications, normal complex structures appear naturally as the Cauchy--Riemann operator has a complex linear principal symbol. 

In this paper, we frequently abbreviate the phrase ``normal(ly) complex'' by NC.

\subsubsection{Linear algebra}

We first discuss the related linear algebra.

\begin{defn}\label{defn:rep-normal-complex}
A {\it normally complex representation} of a group $G$ is a real representation $V$ together with a $G$-invariant complex structure $I^{\check V_G}$ on $\check V_G \subset V$.
\end{defn}

The following lemma is obvious. 

\begin{lemma}
Let $V$ be an NC representation of $G$.
\begin{enumerate}

\item If $H \subset G$ is a subgroup, then $V$ is naturally an NC representation of $H$.

\item Any $G$-invariant subspace $V'\subset V$ is naturally an NC representation of $G$.
\end{enumerate}
\end{lemma}

\begin{defn}\label{defn_NC_triple}(Normally complex isotropy types)
\begin{enumerate}

\item A {\it normal complex isotropy triple} (NC triple for short) consists of an isotropy triple $(G, V, W)$ where $V, W$ are NC representations of $G$.

\item An isomorphism of NC triples from $(G', V', W')$ to $(G, V, W)$ is an isomorphism of isotropy triples such that the induced isomorphisms $\check V_{G'}' \cong \check V_G$ and $\check W_{G'}' \cong \check W_G$ are complex linear. An {\it NC isotropy type} is an isomorphism class of NC triples. Let ${\mb\Gamma}^{\rm NC}$ denote the set of NC isotropy types.

\item A partial order on ${\mb\Gamma}^{\rm NC}$ is defined as follows. We write $\gamma \leq \delta$ if $\gamma$ can be represented by $(G, V, W)$, $\delta$ can be represented by $(H, V, W)$ such that $H$ is isomorphic to a $V$-essential subgroup of $G$ and such that $X \cong V, Y \cong W$ as NC representations of $H$.
\end{enumerate}
\end{defn}

\subsubsection{NC structure on charts and orbifolds}


\begin{defn}
Let $G$ be a finite group and $U$ be a $G$-manifold. When $U$ has boundary, this implies that $\partial U$ is $G$-invariant.
\begin{enumerate}

\item Let $E \to U$ be a $G$-equivariant vector bundle. A {\it normal complex structure} (NC structure for short) on $E$, denoted by ${\bm I}^E$, consists, for each $U$-essential subgroup $H \subset G$, an $H$-invariant complex structure $I^{\check E_H}$ on the  bundle $\check E_H \to U_H$ satisfying the following conditions. 
\begin{enumerate}

\item For each pair of $U$-essential subgroups $H \subsetneq K$, one has the $H$-equivariant and $I^{\check E_K}$-linear decomposition 
\beqn
\check E_K = ( E_H \cap \check E_K) \oplus \check E_H|_{U_K}.
\eeqn
We require that the restriction of $I^{\check E_K}$ on the second summand $\check E_H|_{U_K}$ coincides with $I^{\check E_H}$ restricted to $U_K$.

\item For each $g \in G$ and $H' = g H g^{-1}$, the bundle isomorphism $\check E_H \to \check E_{H'}$ induced by $g$ is complex linear.
\end{enumerate}

\item An {\it NC structure} on $U$ is an NC structure on the tangent bundle $TU \to U$.
\end{enumerate}
\end{defn}

The extension to orbifolds (with boundary) is straightforward.

\begin{defn}\label{defn_normal_complex}
\begin{enumerate}
\item An {\it NC structure} on an orbifold vector bundle ${\mc E}\to {\mc U}$, denoted by ${\mc I}^{\mc E}$, consists, for each chart $\hat C = (G, U, E, \hat\psi)$ of ${\mc E}$, an NC structure ${\bm I}^E = (I^{\check E_H})$ on $E$ satisfying the following conditions. For each chart embedding from $\hat C' = (G', U', E', \hat\phi')$ to $\hat C = (G, U, E, \hat\psi)$ given by a group injection $G' \hookrightarrow G$, an equivariant open embedding $\iota: U' \hookrightarrow U$ covered by an equivariant bundle isomorphism $\hat\iota: E' \to E$, for any $U'$-essential subgroup $H' \subseteq G'$ mapped onto $H \subseteq G$, $\hat \iota$ maps $E'|_{U_{H'}'}$ into $E|_{U_H}$, we require that the induced bundle isomorphism
\beqn
\check E_{H'}'|_{U_{H'}'} \to \check E_H |_{U_H}
\eeqn
is complex linear with respect to the complex structures $I^{\check E_{H'}'}$ and $I^{\check E_H}$. 

\item An {\it NC structure} on an orbifold ${\mc U}$ is an NC structure on  $T{\mc U}$.

\item An {\it NC vector bundle} is an orbifold vector bundle ${\mc E}$ equipped with an NC structure ${\mc I}^{\mc E}$. An {\it NC orbifold} is an  orbifold ${\mc U}$ equipped with an NC structure ${\mc I}^{T{\mc U}}$. An {\it NC pair} is a pair $({\mc U}, {\mc E})$ where ${\mc U}$ is an NC orbifold and ${\mc E}$ is an NC vector bundle over ${\mc U}$. 
\end{enumerate}
\end{defn}

For example, if a vector bundle ${\mc E}$ is equipped with a complex structure $I^{\mc E}$, then there is a naturally induced NC structure ${\mc I}^{\mc E}$ on ${\mc E}$.

Notice that when $({\mc U}, {\mc E})$ is an NC pair, the associated isotropy stratification of ${\mc U}$ can be indexed by NC isotropy types.

\subsubsection{Induced NC structures on suborbifolds with ordinary normal bundle}

Let ${\mc U}$ be an effective orbifold equipped with an NC structure ${\mc I}^{T{\mc U}}$. Let ${\mc X} \subset {\mc U}$ be a suborbifold with ordinary normal bundle (Definition \ref{defn_ordinary_normal}). We claim that there is a naturally induced NC structure ${\mc I}^{T{\mc X}}$. Indeed, given any local chart $C = (G, U, \psi)$, $X:= \psi^{-1}({\mc X}) \subset U$ is a $G$-invariant submanifold. For each $U$-essential subgroup $H \subset_U G$, the normal bundle of $X$ restricted to $X_H$ is the same as the normal bundle of $X_H$ inside $U_H$. Hence
\beqn
NX_H = NU_H|_{X_H}.
\eeqn
Therefore, $I^{NX_H}:= I^{NU_H}|_{X_H}$ defines an NC structure on ${\mc X}$. 

On the other hand, if ${\mc E} \to {\mc U}$ is a vector bundle equipped with an NC structure ${\mc I}^{\mc E}$, then the restriction ${\mc E}|_{\mc X}$ carries an induced NC structure. 

\begin{rem}
When ${\mc U}$ is an orbifold with boundary, the normal bundle to the boundary $\partial {\mc U}$ is an ordinary vector bundle. One can see that if ${\mc U}$ has an NC structure, then $\partial {\mc U}$ has an induced NC structure.
\end{rem}

\subsubsection{Nearby and concordant NC structures}

In later constructions, one often needs to modify the given NC structures. To connect nearby NC structures, one needs to define the closeness properly. First recall some basic linear algebra. Let $V$ be an $2n$-dimensional real vector space. Then a complex structure $J$ on $V$ is equivalent to a $n$-dimensional complex subspace $R \subset V \otimes {\mb C}$ such that $R \oplus \ov{R} = V \otimes {\mb C}$. If $V$ has a linear $G$-action and $J$ is $G$-invariant, then the associated spaces $R$ and $\ov{R}$ are $G$-invariant complex subspaces of $V \otimes {\mb C}$. A nearby complex structure $J'$ on $V$ corresponds to another subspace $R' \subset V \otimes {\mb C}$ which is close to $R$. Then $R' = {\rm graph}(\rho)$ for some linear map $\rho: R \to \ov{R}$ which is close to zero. If $R'$ is also $G$-invariant, then $\rho$ is $G$-equivariant.

We want to characterize the closeness without referring to any extra structure. Let ${\mc J}(V)$ be the space of complex structures on $V$.

\begin{defn}
Given $J \in {\mc J}(V)$ corresponding to a subspace $R \subset V \otimes {\mb C}$, we say that $J' \in {\mc J}(V)$ is {\it close to $J$} if the corresponding space $R' \subset V \otimes {\mb C}$ coincides with ${\rm graph}(\rho)$ for a linear map $\rho: R \to \ov{R}$ such that for all $t \in [0,1]$ and $R_t:= {\rm graph}(t \rho)$, there holds $R_t \cap \ov{R_t} = \{0\}$; in other words, $R_t$ is associated to a complex structure $J_t\in {\mc J}(V)$. We call $J_t$ the {\it $t$-interpolation} from $J$ to $J'$. 
\end{defn}

Notice that in the above definition, if both $J$ and $J'$ are invariant under a linear $G$-action on $V$, then $J_t$ is $G$-invariant. In addition, if $H \subset G$ is a $V$-essential subgroup, then one has the basic decompositions
\begin{align*}
&\ R = R_H \oplus \check R_H,\ &\ \ov{R} = \ov{R}_H \oplus \check{\ov{R}}_H.
\end{align*}
The $H$-equivariant linear map $\rho: R \to \ov{R}$ is block-diagonal with respect to this decomposition with diagonal blocks $\rho_H$ and $\check \rho_H$. Then the restriction of $J'$ to $\check V_H$ is close to $J|_{\check V_H}$ and the restriction $J_t|_{\check V_H}$ is the $t$-interpolation from $J|_{\check V_H}$ to $J'|_{\check V_H}$.

\begin{defn}\label{defn220}
Suppose ${\mc E} \to {\mc U}$ is equipped with an NC structure ${\mc I}$. We say that another NC structure ${\mc I}_0$ on ${\mc E}$ is {\it close to ${\mc I}$} if the following holds. For each bundle chart $\hat C = (G, U, E, \hat\psi)$ of ${\mc E}$ and $x \in U$, let $I^{\check E_{G_x}}$, $I_0^{\check E_{G_1}}$ be the $G_x$-invariant complex structures on $\check E_{G_x} \to U_{G_x}$. Then $I_0^{\check E_{G_x}}$ is close to $I^{\check E_{G_x}}$. 
\end{defn}

\begin{lemma}\label{lemma_NC_extension}
Let ${\mc E} \to {\mc U}$ be equipped with NC structures ${\mc I}^{\mc E}$. Let ${\mc Y} \subset {\mc U}$ be a closed subset and let ${\mc I}_0^{\mc E}$ be another NC structure on ${\mc E}$ defined near ${\mc Y}$ which is close to ${\mc I}^{\mc E}$. Then there exists an NC structure ${\mc I}_1^{\mc E}$ which agrees with ${\mc I}_0^{\mc E}$ near ${\mc Y}$ and which is close to ${\mc I}^{\mc E}$. 
\end{lemma}

\begin{proof}
Choose a smooth cut-off function $\eta: {\mc U} \to [0,1]$ supported in the domain of ${\mc I}_1^{\mc E}$ which is identically zero near ${\mc Y}$. For each bundle chart $\hat C = (G, U, E, \hat\psi)$ and each point $x \in U$, define $I_1^{\check E_{G_x}}(x)$ to be the $(1-\eta(\psi(x)))$-interpolation from $I^{\check E_{G_x}}(x)$ to $I_0^{\check E_{G_x}}(x)$. Then the collection $I_1^{\check E_{G_x}}$ defines an NC structure ${\mc I}_1^{\mc E}$ on ${\mc E}$ which is close to ${\mc I}^{\mc E}$ and which agrees with ${\mc I}_0^{\mc E}$ near ${\mc Y}$.
\end{proof}

\begin{defn}
Let ${\mc E} \to {\mc U}$ be an orbifold vector bundle and let ${\mc I}_0^{\mc E}$, ${\mc I}_1^{\mc E}$ be two NC structures on ${\mc E}$. A {\it concordance} from ${\mc I}_0^{\mc E}$ to ${\mc I}_1^{\mc E}$ is a family of NC structures ${\mc I}_t^{\mc E}$, $t \in [0, 1]$ which depends smoothly on $t$ and domain variables. 
\end{defn}

Above discussions essentially says that if ${\mc I}_1^{\mc E}$ is close to ${\mc I}_0^{\mc E}$, then there exists a canonical concordance from ${\mc I}_0^{\mc E}$ to ${\mc I}_1^{\mc E}$.

\subsubsection{Stabilization of normal complex structures}

We first describe an important model of normally complex orbifolds.

\begin{example}\label{example218}
Let $G$ be a finite group. Let $X$ be a smooth manifold and $F \to X$ be a complex vector bundle equipped with a fiberwise linear $G$-action. Then a $G$-equivariant complex structure $I^F$ on $F$ induces a canonical NC structure on the total space of $F$. Indeed, for each essential subgroup $H \subset_F G$, the fixed point set of $H$ is the total space of $ F_H \subset F$ and the normal bundle $NF_H$ is $\pi_{  F_H}^* \check F_H \to  F_H$. The restriction of $I^F$ to $\check F_H$ is pulled back to an $H$-invariant complex structure on $NF_H$, giving an NC structure on the $G$-manifold $F$. 
\end{example}

Now consider the following situation. Let $U$ be a $G$-manifold. Let $F \to U$ be a $G$-equivariant real vector bundle. Consider the total space $F$ which is a $G$-manifold. A $G$-invariant connection $\nabla^F$ provides an $G$-equivariant splitting of the exact sequence
\beqn
\xymatrix{ 0 \ar[r] & \pi_F^* F \ar[r] & TF \ar[r] &  \pi_F^* TU \ar[r] & 0 }.
\eeqn
If $H \subset G$ is an $F$-essential subgroup, then along the $H$-fixed point $F_H \subset F$, there is an induced split exact sequence
\beqn
\xymatrix{ 0 \ar[r] & \pi_{F_H}^* \check F_H \ar[r]  & NF_H \ar[r] &  \pi_{F_H}^* NU_H \ar[r] & 0}.
\eeqn
If $TU$ and $F$ are normally complex, the normal bundle $NF_H \to F_H$ carries an $H$-invariant complex structure 
\beqn
I^{NF_H} = \pi_{F_H}^* I^{\check F_H} \oplus \pi_{F_H}^* I^{NU_H}.
\eeqn
(We do not need to require that $\nabla^F$ preserves the NC structure.) Therefore, the total space $F$ carries an NC structure ${\bm I}^{TF}$ determined by ${\bm I}^{TU}$, ${\bm I}^F$, and the connection $\nabla^F$. We call ${\bm I}^{TF}$ the {\it bundle NC structure}. Notice that the situation described by Example \ref{example218} is the case when $NU_H = 0$, where a connection is not necessary. 

The above notion can be easily extends to orbifolds. Let $({\mc U}, {\mc I}^{T{\mc U}})$ be an NC orbifold and ${\mc F} \to {\mc U}$ be an NC vector bundle. Let $\nabla^{\mc F}$ be a connection on ${\mc F}$. Then there is a bundle NC structure ${\mc I}^{T{\mc F}}$. Moreover, as connections form an affine space, different choices of connections result in concordant NC structures.

One considers the analogous situation of NC structures on bundles under pullback. A typical example is as follows. Let $U$ be a $G$-manifold, $F \to U$ be a $G$-equivariant vector bundle, and $E \to U$ be a $G$-equivariant vector bundle equipped with an NC structure ${\bm I}^E$. Consider the pullback $\pi_F^* E$. For each $F$-essential subgroup $H \subset G$, one has 
\beqn
\pi_F^* E|_{F_H} = \pi_{F_H}^* E_H \oplus \pi_{F_H}^* \check E_H
\eeqn
where the second summand carries the pullback complex structure $\pi_{F_H}^* I^{\check E_H}$. This makes the pullback bundle $\pi_F^* E \to F$ an NC vector bundle. This example extends to orbifolds. Suppose ${\mc F} \to {\mc U}$ is an orbifold vector bundle and ${\mc E} \to {\mc U}$ is an NC vector bundle. Then $\pi_{\mc F}^* {\mc E} \to {\mc F}$ carries a pullback NC structure. 

\begin{defn}[Stabilization of NC pairs]
Let $({\mc U}, {\mc E})$ be an NC pair. Let ${\mc F} \to {\mc U}$ be an NC vector bundle equipped with a connection $\nabla^{\mc F}$. The {\it stabilization} of $({\mc U}, {\mc E})$ is the pair $({\mc F}, \pi_{\mc F}^* {\mc E} \oplus \pi_{\mc F}^* {\mc F})$ where ${\mc F}$ is equipped with the bundle NC structure induced from the NC structure ${\mc I}^{T{\mc U}}$ on $T{\mc U}$, the NC structure ${\mc I}^{\mc F}$, and the connection $\nabla^{\mc F}$, and where $\pi_{\mc F}^* {\mc E} \oplus \pi_{\mc F}^* {\mc F}$ is equipped with the pullback NC structure.
\end{defn}

\subsection{Straightened geometric structures}\label{subsection24}

The local models for FOP perturbations are equivariant complex polynomial maps between two complex representations. To consider this notion on an orbifold, one needs to locally linearize the orbifold and the bundle so that the notion of polynomial maps (or certain generalization) transits between different local models. This is the purpose of introducing \emph{straightened} geometric structures. A convenient way to construct such structures is to use exponential maps associated to a suitable Riemannian metric on the orbifold and parallel transport associated to a suitable connection on the vector bundle.

Since the discussion here is rather technical, we would like to begin with the case of straightening domain structures, while deferring bundle structures to the second part of this subsection.

Note: in this subsection as well as Appendix \ref{appendixa}, since the discussion is rather technical, we will restrict ourselves to the case where the manifolds or orbifolds do not have boundary. When dealing with spaces with boundary in applications, we construct collar neighborhoods of the boundary and extend structures on the boundary to interior via trivial products (cf. Remark \ref{rem_boundary}).

\subsubsection{Notions and notations}

The model case for our purpose is the total space of a vector bundle equipped with some canonical metric. 

\begin{defn}[Bundle metric]
Let $(X, g^{TX})$ be a smooth Riemannian manifold, $F \to X$ be a smooth real vector bundle equipped with an inner product $h^F$ and an $h^F$-preserving connection $\nabla^F$. The {\it bundle metric} $g^{TF}$ on $F$ is the Riemannian metric on the total space of $F$ determined as the direct sum of
\beqn
g^{TF} = \pi_F^* g^{TX} \oplus \pi_F^* h^F
\eeqn
under the isomorphism of vector bundles over $F$
\beqn
TF \cong \pi_F^* TX \oplus \pi_F^*F 
\eeqn
determined by the connection $\nabla^F$. 
\end{defn}

We fix notations related to exponential maps and tubular neighborhoods. Let $(X, g^{TX})$ be a Riemannian manifold and $S \subset X$ be a closed submanifold. Identify the normal bundle $NS$ with the orthogonal complement of $TS$ in $TM|_S$. For any smooth function $\epsilon: S \to {\mb R}_+$, denote
\beq\label{disk_bundle}
N^\epsilon S:= \Big\{ (x, v) \ |\ x \in S, v \in NS|_x,\ \| v \|^2 < \epsilon(x) \Big\}\subset NS.
\eeq
Then for $\epsilon$ small enough, the exponential map in the normal direction to $S$ defines an open embedding
\beqn
\exp^{NS}: N^\epsilon S \to X.
\eeqn
Let the image be $|N^\epsilon S|$. 

On the other hand, on the normal bundle $NS$ there is an induced inner product $h^{NS}$. Let $p^{NS}: TX|_S \to NS$ be the orthogonal projection induced from $g^{TX}$. Then there is the connection 
\beqn
\nabla^{NS}:= p^{NS} \circ \nabla^{TX}|_{NS}
\eeqn
induced from the Levi--Civita connection $\nabla^{TX}$, which preserves $h^{NS}$.
Together with the restricted Riemannian metric $g^{TS}$ on $S$ there is a canonically induced bundle metric $g^{TNS}$ on $NS$. 

\begin{defn}\label{defn226}
\begin{enumerate}
\item Let $U$ be a manifold and $S \subset U$ be a closed submanifold. A Riemannian metric $g^{TU}$ on $U$ is said to be {\it straightened along $NS$} if for $\epsilon: S \to {\mb R}_+$ sufficiently small, the exponential map $\exp^{NS}: N^\epsilon S \to |N^\epsilon S|$ is an isometry between the bundle metric on $N^\epsilon S$ and the restriction $g^{TU}|_{|N^\epsilon S|}$.

\item Let $U$ be a $G$-manifold. A $G$-invariant Riemannian metric $g^{TU}$ is said to be {\it straightened} if it is straightened along $NU_H$ for all $U$-essential subgroup $H\subset G$. $(U, g^{TU})$ is called a {\it straightened Riemannian $G$-manifold}.
\end{enumerate}
\end{defn}

In the presence of normal complex structures, we would like the NC structure to be ``fiberwise constant'' in normal directions. Let $U$ be a $G$-manifold equipped with an NC structure ${\bm I}^{TU}$. For each $U$-essential subgroup $H\subset G$, the total space of $NU_H$ carries a bundle NC structure ${\bm I}^{TNU_H}$ induced from the $H$-invariant complex structure $I^{NU_H}$ (see Example \ref{example218}).

\begin{defn}\label{defn227}
Let $(U, g^{TU})$ be a Riemannian $G$-manifold equipped with an NC structure ${\bm I}^{TU}$. We say that ${\bm I}^{TU}$ is {\it straightened with respect to $g^{TU}$} if for each $U$-essential subgroup $H \subset G$ the following conditions are satisfied.
\begin{enumerate}

\item The inner product $h^{NU_H}$ on $NU_H \to U_H$ induced from $g^{TU}$ is Hermitian with respect to $I^{NU_H}$ and the induced connection $\nabla^{NU_H}$ is $I^{NU_H}$-linear.

\item For $\epsilon: U_H \to {\mb R}_+$ sufficiently small, the exponential map $\exp^{NS}: N^\epsilon U_H \to |N^\epsilon U_H|$ sends the bundle NC structure ${\bm I}^{TNU_H}$ to the NC structure ${\bm I}^{TU}$.

\end{enumerate}
\end{defn}

Now we turn to orbifolds. The definition is straightforward.

\begin{defn}\label{straightened_metric} 
Let ${\mc U}$ be an effective orbifold.
\begin{enumerate}
    \item A Riemannian metric $g^{T{\mc U}}$ is said to be {\it straightened} if its pullback to each chart $C = (G, U, \psi)$ is a straightened Riemannian metric on $U$.

    \item An NC structure ${\mc I}^{T{\mc U}}$ is said to be {\it straightened} with respect to $g^{T{\mc U}}$ if its pullback to each chart is a straightened NC structure with respect to the chartwise metric.
\end{enumerate}
\end{defn}

Note that the above definitions are well-defined. Indeed, once ${\mc U}$ is equipped with a Riemannian metric, the chart embeddings become isometric embeddings respecting the group action. Therefore, the property of being straightened for a Riemannian metric is compatible for different charts as dictated by Definition \ref{defn226}. A similar discussion can be applied to NC structures as well.

We give the detailed proof of the following existence result for straightened metrics in Appendix \ref{appendixa}. 

\begin{prop}\label{metric_existence}
Let $({\mc U}, {\mc I}^{T{\mc U}})$ be an NC orbifold and ${\mc Y} \subset {\mc U}$ be a closed set. Let $g_0^{T{\mc U}}$ be a Riemannian metric on ${\mc U}$ and ${\mc I}_0^{T{\mc U}}$ be an NC structure on ${\mc U}$ satisfying the following conditions.
\begin{enumerate}

\item ${\mc I}_0^{T{\mc U}}$ is close to ${\mc I}^{T{\mc U}}$ (see Definition \ref{defn220}).

\item $g_0^{T{\mc U}}$ is straightened near ${\mc Y}$.

\item ${\mc I}_0^{T{\mc U}}$ is straightened near ${\mc Y}$ with respect to $g_0^{T{\mc U}}$.

\end{enumerate}
Then there exists a straightened Riemannian metric $g_1^{T{\mc U}}$, a straightened NC structure ${\mc I}_1^{T{\mc U}} $ with respect to $g_1^{T{\mc U}}$ such that
\begin{enumerate}

\item ${\mc I}_1^{T{\mc U}}$ is close to ${\mc I}^{T{\mc U}}$.

\item $(g_1^{T{\mc U}}, {\mc I}_1^{T{\mc U}})$ coincides with $(g_0^{T{\mc U}}, {\mc I}_0^{T{\mc U}})$ near ${\mc Y}$.
\end{enumerate}
\end{prop}

\begin{proof}
See Appendix \ref{subsectiona1}.
\end{proof}



\subsubsection{Straightened bundle structures}

We would also like the vector bundle to be ``straightened'' in the direction normal to fixed point loci in a way analogous to Riemannian metrics. Specifically, we consider connections and inner products on vector bundles and their behavior along normal geodesics. 

We first specify a few notations regarding behaviors of vector bundles within tubular neighborhoods of submanifolds.

\begin{defn} 
Let $(U, g^{TU})$ be a Riemannian manifold and $S \subset U$ be a closed submanifold. Let $\exp^{NS}: N^\epsilon S \to U$ be the normal exponential map. Let $E \to U$ be a real vector bundle equipped with a connection $\nabla^E$. Then the parallel transport along normal geodesics induces a bundle isomorphism
\beqn
\vcenter{ \xymatrix{     \pi_{NS}^* (E|_S)|_{N^\epsilon S} \ar[r]^-{{\rm par}^{NS}} \ar[d]_{\pi_{NS}}  & E|_{|N^\epsilon S|}\ar[d]^{\pi_E} \\
                               N^\epsilon S \ar[r]_{\exp^{NS}} &  |N^\epsilon S|} }. 
\eeqn
We say that $\nabla^E$ is {\it straightened along $NS$} with respect to $g^{TU}$ if ${\rm par}^{NS}$ sends the pullback connection $\pi_{NS}^* (\nabla^{E|_S})$ to the connection $\nabla^E$.
\end{defn}

On the other hand, the corresponding straightenedness for bundle inner products is a consequence of requiring that the connection preserves the inner product. Therefore we do not need a more specific definition. 

Now we consider the case with group actions.

\begin{defn}\label{defn231}
Let $(U, g^{TU})$ be a Riemannian $G$-manifold. Let $E \to U$ be a $G$-equivariant vector bundle equipped with a $G$-invariant connection $\nabla^E$, a $G$-invariant inner product $h^E$, and an NC structure ${\bm I}^E$. 
\begin{enumerate}

\item We say that $\nabla^E$ is {\it straightened} with respect to $g^{TU}$ if for each $U$-essential subgroup $H \subset G$, $\nabla^E$ is straightened along $NU_H$. 


\item We say that ${\bm I}^E$ is {\it straightened} (with respect to $g^{TU}$ and $\nabla^E$) if for each $U$-essential subgroup $H \subset G$, the bundle map ${\rm par}^{NS}$ sends the pullback NC structure on $\pi_{NU_H}^* (E|_{U_H})$ to the NC structure on $E|_{|N^\epsilon U_H|}$. 

\item We say that the triple $(\nabla^E, h^E, {\bm I}^E)$ is {\it normally Hermitian} if for each $H \subset_U G$, the restriction $h^E|_{\check E_H}$ is Hermitian with respect to $I^{\check E_H}$, and the induced connection $\nabla^{\check E_H}$ is a Hermitian connection.
\end{enumerate}
\end{defn}

Notice that the bundle isomorphism ${\rm par}^{NS}$ is defined by parallel transport. Then 
if for all $H \subset_U G$, $\nabla^{\check E_H}$ is complex linear with respect to $I^{\check E_H}$, then ${\bm I}^E$ is straightened. In particular, one has the following situation.
\begin{lemma}
If $E \to U$ is a $G$-equivariant complex vector space with complex structure $I^E$. Let ${\bm I}^E$ be the induced NC structure. Let $h^E$ be a Hermitian inner product and $\nabla^E$ is a Hermitian connection. Then ${\bm I}^E$ is straightened with respect to $g^{TU}$ and $\nabla^E$ and the triple $(\nabla^E, h^E, {\bm I}^E)$ is normally Hermitian.
\end{lemma}

The extension of these notions to orbifold vector bundles is straightforward. 

\begin{defn}\label{defn_bundle_straightening}
Let $({\mc U}, g^{T{\mc U}})$ be a Riemannian orbifold. Let ${\mc E} \to {\mc U}$ be an orbifold vector bundle equipped with a connection $\nabla^{\mc E}$, an inner product $h^{\mc E}$, and an NC structure ${\mc I}^{{\mc E}}$. 
\begin{enumerate}

\item We say that $\nabla^{\mc E}$ is {\it straightened} with respect to $g^{T{\mc U}}$ if its pullback to each chart $\hat C = (G, U, E, \hat\psi)$ is a straightened connection $\nabla^E$ with respect to the chartwise metric $g^{TU}$. 

\item We say that ${\mc I}^{\mc E}$ is {\it straightened} with respect to $g^{T{\mc U}}$ and $\nabla^{\mc E}$ if its pullback to each chart is straightened with respect to the chartwise pullbacks of the bundle connection.

\item We say that the triple $(\nabla^{\mc E}, h^{\mc E}, {\mc I}^{\mc E})$ is {\it normally Hermitian} if its pullback to each bundle chart $\hat C = (G, U, E, \hat\psi)$ is normally Hermitian. 
\end{enumerate}
\end{defn}

Then one has the existence result for straightened structures. Its proof is simpler than the case of straightened Riemannian metrics.

\begin{prop}\label{prop_straightened_connection}
Let $({\mc U}, g^{T{\mc U}})$ be a straightened Riemannian orbifold. Let ${\mc Y} \subset {\mc U}$ be a closed subset. Let ${\mc E} \to {\mc U}$ be a vector bundle equipped with an NC structure ${\mc I}^{\mc E}$. Let $\nabla_0^{\mc E}$ be a connection on ${\mc E}$, $h_0^{\mc E}$ be an inner product on ${\mc E}$, and ${\mc I}_0^{\mc E}$ be an NC structure on ${\mc E}$, satisfying the following conditions.
\begin{enumerate}

\item ${\mc I}_0^{\mc E}$ is close to ${\mc I}^{\mc E}$.

\item $\nabla_0^{\mc E}$ is straightened near ${\mc Y}$ with respect to $g^{T{\mc U}}$.

\item $h_0^{\mc E}$ is preserved by $\nabla_0^{\mc E}$ near ${\mc Y}$.

\item ${\mc I}_0^{\mc E}$ is straightened with respect to $g^{T{\mc U}}$ and $\nabla_0^{\mc E}$ near ${\mc Y}$.

\item $(\nabla_0^{\mc E}, h_0^{\mc E}, {\mc I}_0^{\mc E})$ is normally Hermitian near ${\mc Y}$.

\end{enumerate}
Then there exist a connection $\nabla_1^{\mc E}$, an inner product $h_1^{\mc E}$, and an NC structure ${\mc I}_1^{\mc E}$ satisfying the following conditions.
\begin{enumerate}

\item ${\mc I}_1^{\mc E}$ is close to ${\mc I}^{\mc E}$.

\item $\nabla^{\mc E}$ is straightened with respect to $g^{T{\mc U}}$.

\item $h_1^{\mc E}$ is preserved by $\nabla_1^{\mc E}$.

\item ${\mc I}_1^{\mc E}$ is straightened with respect to $g^{T{\mc U}}$ and $\nabla_1^{\mc E}$.

\item $(\nabla_1^{\mc E}, h_1^{\mc E}, {\mc I}_1^{\mc E})$ is normally Hermitian.

\item $(\nabla_1^{\mc E}, h_1^{\mc E}, {\mc I}_1^{\mc E})$ coincides with $(\nabla_0^{\mc E}, h_0^{\mc E}, {\mc I}_0^{\mc E})$ near ${\mc Y}$.
\end{enumerate}
In addition, if ${\mc E}$ is a complex vector bundle, ${\mc I}^{\mc E} $ is the NC structure induced from the complex structure, $h_0^{\mc E}$ is a Hermitian inner product, $\nabla_0^{\mc E}$ is a Hermitian connection, and ${\mc I}_0^{\mc E} = {\mc I}^{\mc E}$, then we can require that ${\mc I}_1^{\mc E} = {\mc I}_0^{\mc E} = {\mc I}^{\mc E}$, $h_1^{\mc E} = h_0^{\mc E}$, and  $\nabla_1^{\mc E}$ is Hermitian.
\end{prop}

\begin{proof}
See Appendix \ref{subsectiona2}.
\end{proof}

Choosing a system of the above ``straightened'' structures are the prerequisite for the FOP perturbation method. We introduce the following concept. 

\begin{defn}\label{defn_straightening}
Let $({\mc U}, {\mc I}^{T{\mc U}})$ be an NC orbifold and $({\mc E}, {\mc I}^{\mc E})$ be an NC vector bundle over ${\mc U}$. A {\it straightening} of $({\mc U}, {\mc E})$ is a pair $({\mc U}^\#, {\mc E}^\#)$ where 
\begin{enumerate}

\item ${\mc U}^\# = ({\mc U}, g_1^{T{\mc U}}, {\mc I}_1^{T{\mc U}})$ consists of a straightened Riemannian metric $g_1^{T{\mc U}}$ on ${\mc U}$ and an NC structure ${\mc I}_1^{T{\mc U}}$ which is straightened with respect to $g_1^{T{\mc U}}$ and which is concordant to ${\mc I}^{T{\mc U}}$, and 

\item ${\mc E}^\# = ({\mc E}, \nabla_1^{\mc E}, h_1^{\mc E}, {\mc I}_1^{\mc E})$ consists of a connection $\nabla_1^{\mc E}$ on ${\mc E}$ which is straightened with respect to $g_1^{T{\mc U}}$, $h_1^{\mc E}$ is an inner product on ${\mc E}$ which is preserved by $\nabla_1^{\mc E}$, and ${\mc I}_1^{\mc E}$ is an NC structure on ${\mc E}$ which is straightened with respect to $g_1^{T{\mc U}}$ and $\nabla_1^{\mc E}$ and which is concordant to ${\mc I}^{\mc E}$, such that $(\nabla_1^{\mc E}, h_1^{\mc E}, {\mc I}_1^{\mc E})$ is normally Hermitian.
\end{enumerate}
\end{defn}

\begin{rem}
Later we will see that the bundle inner product will not enter the definition of normally complex 
perturbation.
\end{rem}

The following corollary is a direct consequence of Proposition \ref{metric_existence} and Proposition \ref{prop_straightened_connection}. 

\begin{cor}\label{cor_straightening_extension}
Let $({\mc U}, {\mc E})$ be an NC pair and ${\mc Y} \subset {\mc U}$ be a closed subset. Suppose a straightening of $({\mc U}, {\mc E})$ is given over a neighborhood of ${\mc Y}$. Then there exists a straightening of $({\mc U}, {\mc E})$ which agrees with the given one near ${\mc Y}$.
\end{cor}

Another direct corollary to Proposition \ref{metric_existence} and Proposition \ref{prop_straightened_connection} is that two straightenings can be connected in a one-parameter family. 

\begin{cor}
Let $({\mc U}, {\mc I}^{T{\mc U}})$ be an NC orbifold without boundary. Let $({\mc E}, {\mc I}^{\mc E})$ be an NC vector bundle over ${\mc U}$. Let 
\beqn
\Big( {\mc U}_i^\# = ({\mc U}, g_i^{T{\mc U}}, {\mc I}_i^{T{\mc U}}),\ {\mc E}_i^\# = ({\mc E}, \nabla_i^{\mc E}, h_i^{\mc E}, {\mc I}_i^{\mc E})\Big),\ i = 0, 1
\eeqn
be two straightenings of $({\mc U}, {\mc E})$. Let $\tilde {\mc U} = {\mc U} \times [0, 1]$ equipped with the product NC structure ${\mc I}^{T\tilde {\mc U}}$ and $\tilde {\mc E} = {\mc E} \times [0,1]$ equipped with the product NC structure ${\mc I}^{\tilde {\mc E}}$. Then there exists a straightening $(\tilde {\mc U}^\#, \tilde {\mc E}^\#)$ of $(\tilde {\mc U}, \tilde {\mc E})$ whose restriction to ${\mc U}\times [0, \epsilon)$ resp. ${\mc U}\times (1-\epsilon, 1]$ coincides with the product of $({\mc U}_0^\#, {\mc E}_0^\#)$ resp. $({\mc U}_1^\#, {\mc E}_1^\#)$ with the interval $[0, \epsilon)$ resp. $(1-\epsilon, 1]$.
\end{cor}

\begin{proof}
We first choose a concordance of the NC structures. As ${\mc I}_i^{T{\mc U}}$ resp. ${\mc I}_i^{\mc E}$ is concordant to ${\mc I}^{T{\mc U}}$ resp. ${\mc I}^{\mc E}$, there exists an NC structure ${\mc I}_{01}^{T\tilde {\mc U}}$ on $\tilde {\mc U}$ resp. ${\mc I}_{01}^{\tilde {\mc E}}$ on $\tilde {\mc E}$ which is concordant to ${\mc I}^{T\tilde {\mc U}}$ resp. ${\mc I}^{\tilde {\mc E}}$ and which agrees with the product of ${\mc I}_i^{T{\mc U}}$ resp. ${\mc I}_i^{\mc E}$ near the corresponding boundaries. For $\epsilon$ sufficiently small, denote $\tilde {\mc Y} = {\mc U}\times ([0, \epsilon] \cup [1-\epsilon, 1])$. Then the desired straightening can be constructed using Proposition \ref{metric_existence} and Proposition \ref{prop_straightened_connection}.
\end{proof}

\subsubsection{Straightenings and stabilizations}

Recall that for an NC pair $({\mc U}, {\mc E})$ and an NC vector bundle ${\mc F} \to {\mc U}$, the stabilization is 
\beqn
({\mc F}, \pi_{\mc F}^* {\mc E} \oplus \pi_{\mc F}^* {\mc F}).
\eeqn
Upon choosing a connection $\nabla^{\mc F}$ on ${\mc F}$, the stabilization becomes an NC pair. We need to consider straightenings on the stabilization. First, similar to the manifold case, if ${\mc U}$ is equipped with a metric $g^{T{\mc U}}$ and if $h^{\mc F}$ is an inner product on ${\mc F}$ preserved by $\nabla^{\mc F}$, then on the total space ${\mc F}$ there is a canonically induced {\it bundle metric} $g^{T{\mc F}}$. There are similar pullback connections. 

\begin{prop}\label{prop_straightening_stabilization}
Let $({\mc U}^\#, {\mc E}^\#)$ be a straightened NC pair with domain metric $g^{T{\mc U}}$, the domain NC structure ${\mc I}^{T{\mc U}}$, the bundle connection $\nabla^{\mc E}$, the bundle inner product $h^{\mc E}$, and the bundle NC structure ${\mc I}^{\mc E}$. 

Let $\pi_{\mc F}: {\mc F} \to {\mc U}$ be a vector bundle equipped with an inner product $h^{\mc F}$, an $h^{\mc F}$-preserving connection $\nabla^{\mc F}$, and an NC structure ${\mc I}^{\mc F}$ such that the triple $(\nabla^{\mc F}, h^{\mc F}, {\mc I}^{\mc F})$ is normally Hermitian (Definition \ref{defn_bundle_straightening}). Suppose $\nabla^{\mc F}$ is straightened with respect to $g^{T{\mc U}}$. Consider the pair
\beqn
({\mc F}, \pi_{\mc F}^* {\mc E})
\eeqn
Then the domain metric $g^{T{\mc F}}$ as the bundle metric on ${\mc F}$, the domain NC structure as the bundle NC structure on ${\mc F}$, the bundle connection $\pi_{\mc F}^* \nabla^{\mc E}$, the bundle inner product $\pi_{\mc F}^* h^{\mc E}$, the bundle NC structure $\pi_{\mc F}^* {\mc I}^{\mc E}$ give a straightening of $({\mc F}, \pi_{\mc F}^* {\mc E})$.
\end{prop}

\begin{proof}
See Appendix \ref{subsectiona3}.
\end{proof}

\subsubsection{Suborbifolds with ordinary normal bundle}

Now we discuss extending straightenings on an suborbifold with ordinary normal bundle to the ambient orbifold.

We first reorganize the notations. Let ${\mc U}$ be an orbifold and ${\mc X} \subset {\mc U}$ be a suborbifold with ordinary normal bundle. Recall that for each chart $C = (G, U, \psi)$, $X:= \psi^{-1}({\mc X}) \subset U$ is a $G$-invariant closed submanifold. Let $\pi^{N{\mc X}}: N{\mc X} \to {\mc X}$ be the normal bundle and $N^\epsilon {\mc X}$ is the disk bundle of radius $\epsilon: {\mc X} \to {\mb R}_+$. If $g^{T{\mc U}}$ is a Riemannian metric on ${\mc U}$, then there is a normal exponential map 
\beqn
\exp^{N{\mc X}}: N^\epsilon {\mc X} \to {\mc U}
\eeqn
which is an open embedding of orbifolds. Notice that the total space $N{\mc X}$ carries a bundle metric. When ${\mc I}^{T{\mc X}}$ is an NC structure on ${\mc X}$ induced from an NC structure ${\mc I}^{T{\mc U}}$, there is a bundle NC structure on $N{\mc X}$.

In addition, if ${\mc E} \to {\mc U}$ is a vector bundle and $\nabla^{\mc E}$ is a connection on ${\mc E}$. Abbreviate ${\mc E}_{\mc X}:= {\mc E}|_{\mc X}$. Then along normal geodesics, there is a bundle map covering the normal exponential map 
\beqn
{\rm par}^{N{\mc X}}: \Big( \pi_{N{\mc X}}^* {\mc E}_{\mc X}\Big)|_{N^\epsilon {\mc X}} \to {\mc E}|_{|N^\epsilon {\mc X}|}.
\eeqn

\begin{defn} 
Let $g^{T{\mc U}}$ be a Riemannian metric on ${\mc U}$ and ${\mc I}^{T{\mc U}}$ is an NC structure on ${\mc U}$. Let $\nabla^{\mc E}$, $h^{\mc E}$, ${\mc I}^{\mc E}$ be respectively a connection, an inner product, an NC structure on ${\mc E} \to {\mc U}$.
\begin{enumerate}

\item $g^{T{\mc U}}$ is said to be {\it straightened} along $N{\mc X}$ if for $\epsilon: {\mc X} \to {\mb R}_+$ sufficiently small, the normal exponential map $\exp^{N{\mc X}}$ is an isometric open embedding.

\item ${\mc I}^{T{\mc U}}$ is said to be {\it straightened} along $N{\mc X}$ if for $\epsilon: {\mc X} \to {\mb R}_+$ sufficiently small, $\exp^{N{\mc X}}$ is an open embedding of NC orbifolds.

\item $\nabla^{\mc E}$ is said to be {\it straightened} along $N{\mc X}$ if ${\rm par}^{N{\mc X}}$ sends the pullback connection $\pi_{N{\mc X}}^* \nabla^{{\mc E}_{\mc X}}$ to $\nabla^{\mc E}$.

\item ${\mc I}^{\mc E}$ is said to be {\it straightened} along $N{\mc X}$ with respect to $g^{T{\mc U}}$ and $\nabla^{\mc E}$ if ${\rm par}^{N{\mc X}}$ sends the pullback 
$\pi_{N{\mc X}}^* {\mc I}^{{\mc E}_{\mc X}}$ to 
${\mc I}^{\mc E}$.

\end{enumerate}
\end{defn}

\begin{prop}\label{prop241}
Let $({\mc U}, {\mc E})$ be an NC pair. Let ${\mc X} \subset {\mc U}$ be a suborbifold with ordinary normal bundle and ${\mc E}_{\mc X}:= {\mc E}|_{\mc X}$. Then $({\mc X}, {\mc E}_{\mc X})$ carries an induced NC structure. Suppose $({\mc X}^\#, {\mc E}_{\mc X}^\#)$ is a straightening of $({\mc X}, {\mc E}_{\mc X})$. Then there exists a straightening $({\mc U}^\#, {\mc E}^\#)$ of $({\mc U}, {\mc E})$ whose restriction to ${\mc X}$ is $(g_0^{T{\mc X}}, {\mc I}_0^{T{\mc X}})$ such that 
\begin{enumerate}
\item The metric $g_0^{T{\mc U}}$ is straightened along $N{\mc X}$.

\item The NC structure ${\mc I}_0^{T{\mc U}}$ is straightened along $N{\mc X}$ with respect to $g_0^{T{\mc U}}$.

\item The connection $\nabla_0^{\mc E}$ is straightened along $N{\mc X}$ with respect to $g_0^{T{\mc U}}$.

\item The 
NC structure ${\mc I}_0^{\mc E}$ is straightened along $N{\mc X}$ with respect to $g_0^{T{\mc U}}$ and $\nabla_0^{\mc E}$.
\end{enumerate}
\end{prop}

\begin{proof}
We only give the proof of the extension of the straightened metric. The extensions of other structures are similar and omitted. Let $g_0^{T{\mc X}}$ be the straightened metric on ${\mc X}$. Let $g^{T{\mc U}}$ be an arbitrary Riemannian metric on ${\mc U}$ extending $g_0^{T{\mc X}}$, which gives a normal exponential map $\exp^{N{\mc X}}: N^\epsilon {\mc X} \to {\mc U}$. Moreover, choose a connection $\nabla_0^{N{\mc X}}$ on the normal bundle $N{\mc X} \to {\mc X}$ which is straightened with respect to $g_0^{T{\mc X}}$ and choose an inner product $h_0^{N{\mc X}}$ which is preserved by $\nabla_0^{N{\mc X}}$ and which is straightened with respect to $g_0^{T{\mc X}}$ and $\nabla_0^{N{\mc X}}$ (the existence of such a connection and an inner product is guaranteed by Proposition \ref{prop_straightened_connection}). The triple $(g_0^{T{\mc X}}, h_0^{N{\mc X}}, \nabla_0^{N{\mc X}})$ induces a bundle metric $g_0^{TN{\mc X}}$ on the total space $N{\mc X}$, which is a straightened metric. Then using the normal exponential map to pushforward the bundle metric to a tubular neighborhood $|N^\epsilon {\mc X}|$. This gives a straightened metric in $|N^\epsilon {\mc X}|$ which extends $g_0^{T{\mc X}}$. Then using the relative version of Proposition \ref{metric_existence} to find a straightened metric $g_0^{T{\mc U}}$ on ${\mc U}$ without altering its value near ${\mc X}$.
\end{proof}

\subsection{Normally complex sections}\label{subsection25}

Now we discuss the notion of normally complex sections originally introduced by Parker \cite{BParker_integer} which generalizes Fukaya--Ono's notion of normally polynomial sections \cite{Fukaya_Ono_integer}. Such sections will form the category of perturbations we use. 

\subsubsection{Normally complex  maps between representations}

We first discuss equivariant polynomial maps between representations. Let $G$ be a finite group and $V, W$ be complex representations of $G$. Let ${\rm Poly}_G (V, W)$ be the space of $G$-equivariant polynomial maps. For $d \geq 0$, let
\beqn
{\rm Poly}_G^d (V, W) \subset {\rm Poly}_G (V, W)
\eeqn
be the subspace of maps of degrees at most $d$. More generally, for an NC triple $(G, V, W)$, define
\beqn
\wh {\rm Poly}{}_G (V, W):= {\rm Poly}_G (\check V_G, \check W_G) \oplus W_G.
\eeqn
Its elements, typically denoted by $P$, can be viewed as $G$-equivariant maps from $V$ to $W$ which are independent of the variable $v_G\in  V_G$ and whose $ W_G$-components are constant maps. Hence there is an inclusion
\beqn
\wh{\rm Poly}{}_G (V, W) \subset C_G^\infty(V, W)
\eeqn
where the latter is the space of smooth $G$-equivariant maps. Let $\wh{\rm Poly}{}_G^d(V, W) \subset \wh{\rm Poly}{}_G(V, W)$ be the subset of maps of degree at most $d$. It is the model of maps considered in \cite{Fukaya_Ono_integer}. Parker \cite{BParker_integer} generalizes it to the following larger space of equivariant maps.

\begin{defn}
Let $(G, V, W)$ be an NC triple. A $G$-equivariant {\it normally complex map} from $V$ to $W$ is an element of the module over the ring of $G$-invariant smooth functions $C_G^\infty(V; {\mb R})$ generated by $\wh{\rm Poly}{}_G(V, W)$, denoted by 
\beqn
C_G^{\rm NC} (V, W) \subset C_G^\infty(V, W).
\eeqn
A {\it lift} of an NC map $S\in C_G^{\rm NC} (V, W)$ is a $G$-invariant smooth map 
\beqn
{\mf p} \in C_G^\infty(V, \wh{\rm Poly}{}_G^d(V, W))
\eeqn
such that $S(v) = \mf{p}(v)(v)$. We denote
\beqn
S_{\mf p}(v) = {\mf p}(v)(v).
\eeqn
\end{defn}

Notice that there is a natural diagram
\beqn
C_G^\infty(V, \wh{\rm Poly}{}_G^d(V, W)) \to \colim_d C_G^\infty( V, \wh{\rm Poly}{}_G^d(V, W)) \to C_G^{\rm NC} (V, W).
\eeqn

\begin{rem}\label{rem_density}
In general $C_G^{\rm NC} (V, W)$ is only a $C^0$-dense subspace of $C_G^\infty(V, W)$. For example, when $G = {\mb Z}_2$ and $V\cong W \cong {\mb C}$ is the nontrivial complex representation, the map $z\mapsto \ov{z}$, which is $G$-equivariant, cannot be approximated by an NC map in $C^1$-topology. 
\end{rem}

\subsubsection{Change of symmetry groups}

Let $(G, V, W)$ be an NC triple. If $H \subset G$ is a $V$-essential subgroup, then any $G$-equivariant map $f: V \to W$ is automatically $H$-equivariant. Hence there is an inclusion
\beqn
C_G^\infty(V, W) \hookrightarrow C_H^\infty(V, W).
\eeqn

\begin{lemma}\label{lemma244}
Under the above inclusion, the space $C_G^{\rm NC} (V, W)$ is mapped into $C_H^{\rm NC} (V, W)$. Moreover, there is a natural map $\mu$ (the dashed arrow below) making the following diagram commute.
\beqn
\xymatrix{   \displaystyle \colim_d C_G^\infty\left( V, \wh{\rm Poly}{}_G^d(V, W) \right) \ar[r] \ar@{.>}[d]_{\mu} & C_G^{\rm NC} (V, W) \ar[d]\\
             \displaystyle \colim_d C_H^\infty \left( V, \wh{\rm Poly}{}_H^d(V, W) \right) \ar[r] & C_H^{\rm NC} (V, W) }
\eeqn
\end{lemma}

\begin{proof}
A map ${\mf p}: V \to \wh{\rm Poly}{}_G^d(V, W)$ is determined by its graph
\beqn
{\rm graph}({\mf p}) \subset M^d(G, V, W):= V \times \wh{\rm Poly}{}_G^d(V, W).
\eeqn
We would like to define a natural map 
\beqn
\mu: M^d(G, V, W) \to M^d(H, V, W)
\eeqn
which will play an important role. Given $v \in V$ and $P \in \wh{\rm Poly}{}_G^d(V, W)$, with respect to the basic decomposition (for the $H$-action) $V = V_H \oplus \check V_H$, write $v = v_H + \check v_H$; with respect to the basic decomposition $W = W_H \oplus \check W_H$, write $P =  P_H + \check P_H$. Then define 
\beq\label{eqn22}
\mu(v, P):= (v, \check P_H(v_H, \cdot), P_H(v)) \in M^d(H, V, W).
\eeq
It is an $H$-equivariant map. Then there is a unique $H$-invariant map ${\mf q}: V \to \wh{\rm Poly}{}_H^d(V, W)$ whose graph is $\mu ({\rm graph}({\mf p}))$. By abuse of notation, define 
\beqn
\mu ({\mf p}):= {\mf q}.
\eeqn
From the definition one can see that $S_{\mf q} = P$ as maps from $V$ to $W$.
\end{proof}

\subsubsection{Family case}

We also consider the family case. 

\begin{defn}\label{defn246} Let $X$ be a manifold. 
\begin{enumerate}

    \item A {\it family NC triple} over $X$ is a triple $(G, F, E)$ consisting of a pair of real vector bundles $\pi_F: F \to X$, $\pi_E:E \to X$ with fiberwise linear $G$-actions, a $G$-invariant complex structure $I^{\check F_G}$ on $\check F_G \subset F$ and a $G$-invariant complex structure $I^{\check E_G}$ on $\check E_G \subset E$. 

    \item Given a family NC triple $(G, F, E)$, let $C_G^\infty(F, E)$ be the space of $G$-equivariant smooth bundle maps $S: F \to E$. Let 
\beqn
\wh{\rm Poly}{}_G(F, E)  \subset C_G^\infty(F, E)
\eeqn
be the subbundle whose fiber at $x \in X$ is the space $\wh{\rm Poly}{}_G(F_x, E_x)$, and $\wh{\rm Poly}{}_G^d(F, E) \subset \wh{\rm Poly}{}_G(F, E)$ be the subbundle of fiberwise polynomial maps of degree at most $d$.

\item An element $S \in C_G^\infty(F, E)$ is called a {\it normally complex} bundle map if there exists $d\geq 0$ such that it is an element of the module over $C_G^\infty(F; {\mb R})$ generated by $\Gamma(\wh{\rm Poly}{}_G^d (F, E))$, the global sections of $\wh{\rm Poly}{}_G^d(F, E)$. Denote
    \beqn
    C_G^{\rm NC} (F, E):= \colim_{d} C_G^\infty(F; {\mb R}) \cdot \Gamma( \wh{\rm Poly}{}_G^d(F, E)) \subset C_G^\infty(F, E).
    \eeqn
\end{enumerate} 
\end{defn}

When we discuss transversality, we need to choose a particular presentation of a NC bundle map. Notice that for each $d\geq 0$, there is a natural map 
\beq
\begin{aligned}
C_G^\infty \Big( F, \wh{\rm Poly}{}_G^d(F, E) \Big) &\to C_G^{\rm NC} (F, E) \\
{\mf p} &\mapsto S_{\mf p} := (F \ni v \mapsto {\mf p}(v)(v)).
\end{aligned}
\eeq

\begin{defn}\label{defn244}(cf. \cite[Definition 4.12]{BParker_integer})
A {\it lift} of an NC bundle map $S\in C_G^{\rm NC} (F, E)$ is a $G$-invariant bundle map ${\mf p}: F \to \wh{\rm Poly}{}_G^d(F, E)$ such that $S = S_{\mf p}$.
\end{defn}

It is important to consider how the NC bundle maps behave under restriction to a smaller symmetry group. Let $H \subset G$ be an $F$-essential subgroup. Then there is the natural inclusion
\beq\label{eqn23}
C_G^\infty(F, E) \to C_H^\infty(F, E).
\eeq
Lemma \ref{lemma244} implies that \eqref{eqn23} sends $C_G^{\rm NC} (F, E)$ into $C_H^{\rm NC} (F, E)$. Moreover, there is a natural map on the level of lifts.

\begin{lemma}\label{lemma247}
There is a natural arrow $\mu$ (the dotted one) making the following diagram commutes.
\beqn
\xymatrix{ \displaystyle \colim_d C_G^\infty \Big( F, \wh{\rm Poly}{}_G^d(F, E) \Big) \ar[r] \ar@{.>}[d]_\mu & C_G^{\rm NC} (F, E) \ar[d] \\
\displaystyle \colim_d C_H^\infty\Big(F, \wh{\rm Poly}{}_H^d(F, E) \Big)  \ar[r] & C_H^{\rm NC} (F, E)} 
\eeqn
\end{lemma}

\begin{proof}
One simply applies the map $\mu$ provided by Lemma \ref{lemma244} along each fiber of $F$ and the proof follows from the same argument.
\end{proof}

One can view the total space $F$ as the total space of the pullback bundle $\pi_{F_H}^* \check F_H$. Then an $H$-equivariant bundle map $S: F \to E$ can also be viewed as an $H$-equivariant bundle map from $\pi_{F_H}^* \check F_H$ to $\pi_{F_H}^* E$. Hence one has the natural identification
\beq\label{eqn24}
C_H^\infty(F, E) \cong C_H^\infty( \pi_{ F_H}^* \check F_H, \pi_{  F_H}^* E),
\eeq
which induces the identification
\beqn
\pi_{ F_H}^* \wh{\rm Poly}{}_H^d(F, E) \cong \wh{\rm Poly}{}_H^d( \pi_{ F_H}^* \check F_H, \pi_{ F_H}^* E).
\eeqn
The following lemma is readily to check.

\begin{lemma}\label{lemma248}
The map \eqref{eqn24} sends $C_H^{\rm NC} (F, E)$ to $C_H^{\rm NC} (\pi_{F_H}^* \check F_H, \pi_{F_H}^* E)$ and the following diagram commutes.
\beqn
\xymatrix{  C_H^\infty \left( F,  \wh{\rm Poly}{}_H^d(F, E) \right) 
\ar[r] \ar[d]_{\cong} &  C_H^{\rm NC} (F, E) \ar[d]  \\
C_H^\infty \left( \pi_{ F_H}^* \check F_H, \wh{\rm Poly}{}_H^d \left( \pi_{  F_H}^* \check F_H, \pi_{ F_H}^* E \right) \right) \ar[r]    & C_H^{\rm NC} \left( \pi_{ F_H}^* \check F_H, \pi_{ F_H}^* E \right)  }
\eeqn
\end{lemma}

\subsubsection{Normally complex sections}

Now we define the notion of normally complex sections over orbifolds. Let $({\mc U}, {\mc E})$ be an NC pair equipped with a straightening $({\mc U}^\#, {\mc E}^\#)$. Let $\hat C = (G, U, E, \hat \psi)$ be a bundle chart. Let $H \subset G$ be a $U$-essential subgroup. Recall that the straightening induces exponential maps
\beqn
\exp^{NU_H}: N^\epsilon U_H \to |N^\epsilon U_H| \subset U
\eeqn
and the bundle isomorphism
\beqn
{\rm par}^{NU_H}: \pi_{NU_H}^* E|_{U_H} \cong E|_{|N^\epsilon U_H|}.
\eeqn
Then any $G$-equivariant section $S: U \to E$ induces a $G$-equivariant bundle map 
\beqn
S_H^\#: N^\epsilon U_H \to E|_{U_H}.
\eeqn

\begin{defn}\label{defn:FOP}
Under the above setting, let $\hat C = (G, U, E, \hat \psi)$ be a bundle chart of ${\mc E}$ and $S: U \to E$ is a $G$-equivariant section. 

\begin{enumerate}

\item Let $H \subset G$ be a $U$-essential subgroup. We say that $S$ is {\it normally complex along $NU_H$} if for $\epsilon: U_H \to {\mb R}_+$ sufficiently small, the bundle map $S_H^\#$ is the restriction of an NC bundle map from $NU_H$ to $E|_{U_H}$ to  $N^\epsilon U_H$.

\item $S$ is said to be {\it normally complex} if it is NC along $NU_H$ for all $U$-essential subgroups $H \subset G$.

\item A smooth section ${\mc S}: {\mc U} \to {\mc E}$ is an {\it NC section} if its pullback to each bundle chart is an NC section.
\end{enumerate}
\end{defn}

The importance of straightening is stated in the following lemma, which guarantees that the compatibility of being normally complex along tubular neighborhoods of fixed point loci of different groups.

\begin{lemma}
Let $\hat C = (G, U, E, \hat\psi)$ be a bundle chart. Suppose $S: U \to E$ is NC along $NU_G$, then for any $U$-essential subgroup $H \subset G$, $S$ is NC along $NU_H$ in a neighborhood of $U_G$.
\end{lemma}

\begin{proof}
By the definition of straightening, near $U_G$, we may regard $U$ with a disk bundle of $NU_G$ equipped with the bundle metric and bundle NC structure, and identify $E$ with the pullback bundle $\pi_{NU_G}^* (E|_{U_G})$ with the pullback structures. By the definition of NC sections, $S_G^\#$ is the restriction of a section in $C_G^{\rm NC} ( NU_G, E|_{U_G})$. Let $H\subset_U G$ be an essential subgroup. In the neighborhood $|N^\epsilon U_G|$ for $\epsilon$ sufficiently small, the $H$-fixed point set is the total space of $(NU_G)_H$. By Lemma \ref{lemma247}, Lemma \ref{lemma248}, and the definition of straightenings, the restriction of $S_H^\#$ to a neighborhood of $(NU_G)_H$ is also normally complex.
\end{proof}

Next we prove some simple formal properties of the space of normally complex sections. Let $\Gamma({\mc U}, {\mc E})$ be the space of smooth sections of ${\mc E}$. Let $({\mc U}^\#, {\mc E}^\#)$ be a straightening of $({\mc U}, {\mc E})$. Let 
\beqn
\Gamma^{\rm NC} ({\mc U}^\#, {\mc E}^\#) \subset \Gamma({\mc U}, {\mc E})
\eeqn
be the subset of normally complex sections (with respect to the given straightening). 

\begin{lemma}\label{lemma252}
$\Gamma^{\rm NC} ({\mc U}^\#, {\mc E}^\#)$ is a $C^\infty({\mc U})$-submodule of $\Gamma({\mc U}, {\mc E})$.
\end{lemma}

\begin{proof}
It follows from the fact that $C_G^{\rm NC}(F, E)$ is a $C_G^\infty(F; {\mb R})$-submodule (cf. Definition \ref{defn246}).
\end{proof}

Lastly we prove the $C^0$-density of the space of normally complex sections.

\begin{prop}\label{prop_complex_dense}
$\Gamma^{\rm NC} ({\mc U}^\#, {\mc E}^\#)$ is dense in $\Gamma ({\mc E})$ with respect to the $C^0$-norm induced by any inner product on ${\mc E}$.
\end{prop}

To prove this result, we first recall a fundamental lemma proved by Fukaya--Ono \cite{Fukaya_Ono_integer}. For the convenience of the reader we provide the proof here.

\begin{lemma}[Fukaya--Ono's lemma] \cite[Lemma 5]{Fukaya_Ono_integer}\label{lemma251} 
Let $(G, V, W)$ be an NC triple. Then for any $d \geq |G|$ and any $v_0 \in V$ with isotropy subgroup $H \subseteq G$ and $w_0 \in  W_H$, there exists $P \in \wh{\rm Poly}{}_G^d (V, W)$ such that $P(v_0) = w_0$.
\end{lemma}

\begin{proof}
It is easy to see that it suffices to consider the case when $V_G = W_G = \{0\}$; in particular, $V$ and $W$ are complex representations. Then by decomposing $W$ into irreducible components, we may assume that $W$ is an irreducible (complex) representation of $G$. Define the $G$-vector space
\beqn
R:= \bigoplus_{\gamma \in G} {\mb C}\{\langle \gamma \rangle \}
\eeqn
with $G$-action defined as 
\beqn
g \left( \sum_\gamma c_\gamma \cdot \langle \gamma \rangle \right) = \sum_\gamma c_\gamma \cdot \langle \gamma g^{-1}\rangle = \sum_\gamma c_{\gamma g} \cdot \langle \gamma \rangle.
\eeqn
Since $R$ is a regular representation, there is a $G$-equivariant homomorphism $\Psi: R \to W$ and an element $r \in R$ such that
\beqn
\Psi(r) = \Psi \left( \sum_\gamma w_\gamma \cdot \langle \gamma \rangle \right) = w_0.
\eeqn
Since $w \in  W_H$, for all $h \in H$, one has 
\beqn
\Psi(h r) = h \Psi(r ) = w_0.
\eeqn
Hence by taking average over $H$, one may assume that 
\beqn
\gamma' H = \gamma'' H \Longrightarrow w_{\gamma'} = w_{\gamma''}.
\eeqn
Now we claim that one can choose a complex polynomial $f: V \to {\mb C}$ (not necessarily $G$-invariant) of degree at most $|G|$ such that
\beqn
\forall \gamma \in G,\ f(\gamma v_0) = w_\gamma.
\eeqn
Indeed, there are $n:= |G/H|$ distinct elements in the $G$-orbit of $v_0$. One can choose a linear decomposition $V = V_1 \oplus V_2$ such that $V_1$ is one-dimensional and that the projection of these $n$ distinct elements are still distinct in $V_1$. Then by Lagrange's method of interpolation, one can find a complex polynomial $f: V_1 \to {\mb C}$ of degree at most $|G|$ taking the prescribed values $w_\gamma$ at the corresponding projection image of $\gamma v_0$ in $V_1$. Extend $f$ trivially to $V$ one obtains a polynomial $f: V \to {\mb C}$ satisfying the required conditions. Now define $P: V \to W$ by 
\beqn
P(v) = \Psi \left( \sum_{\gamma\in G} f(\gamma v) \cdot \langle \gamma\rangle \right).
\eeqn
Then this is a $G$-equivariant polynomial map sending $v_0$ to $w_0$.
\end{proof}

\begin{proof}[Proof of Proposition \ref{prop_complex_dense}]
Fix an inner product on ${\mc E}$. Given ${\mc S}_0\in \Gamma({\mc U}, {\mc E})$ and $\delta > 0$, it suffices to find a section ${\mc S} \in \Gamma^{\rm NC}({\mc U}^\#, {\mc E}^\#)$ such that $\| {\mc S} - {\mc S}_0\|_{C^0}\leq \delta$. We construct ${\mc S}$ by local approximations and gluing using a partition of unity. 

For each $p \in {\mc U}$, one can find a bundle chart $\hat C = (G, U, E, \hat\psi)$ centered at $p$ satisfying the following conditions.
\begin{enumerate}
    \item The normal bundle $NU_G \to U_G$ is trivial with fiber being a representation $V$ of $G$. The tubular neighborhood $|N^\epsilon U_G|$ is the whole chart $U$, identified with a subset of $U_G \times V$.

    \item $E|_{U_G}$ is trivial with  fiber being an NC representation $W$ of $G$. Then the straightening induces a trivialization $E \cong U \times W$ and the section $S_0: U \to E$ is equivalent to a $G$-equivariant map $S_0: U \to W$.
\end{enumerate}

\vspace{0.1cm}

\noindent {\it Claim.} For $d \geq |G|$ and any
$\epsilon>0$ there is a $G$-invariant map
\beqn
{\mf p}: U_G \times V \to \wh {\rm Poly}{}_G^d(V, W)
\eeqn
providing a lift of an NC map $S_{\mf p}$ such that 
\beqn
\sup_{(y, v) \in U} \| S_0 (y, v) -S_{\mf p}(y, v) \| \leq \delta.
\eeqn

\vspace{0.1cm}

Assuming this claim is true, we prove the proposition as follows. As ${\mc U}$ is second-countable, it is Lindel\"of, hence there is a countable collection of points $p_i \in {\mc U}$ and bundle charts $\hat C_i = (G_i, U_i, E_i, \hat\psi_i)$, $i = 1, \ldots$, which satisfy the above conditions such that the collection $\psi_i(U_i)$ cover ${\mc U}$. Let $S_{0, i}: U_i \to E_i$ be the pullback of ${\mc S}_0$ to the chart $\hat C_i$. The above claim implies that there exists a $G_i$-equivariant NC section $S_i: U_i \to E_i$ such that 
\beqn
\| S_{0, i} - S_i\|_{C^0(U_i)} \leq \frac{\delta}{2^i}.
\eeqn
Then choose a smooth partition of unity $\rho_i$ subordinate to the open cover $\psi_i(U_i)$. Consider 
\beqn
{\mc S}_i:= \sum_{i=1}^\infty \rho_i S_i
\eeqn
which belongs to $\Gamma^{\rm NC}({\mc U}^\#, {\mc E}^\#)$ by Lemma \ref{lemma252}. Then one has 
\beqn
\| {\mc S}_0 - {\mc S} \|_{C^0} \leq \sum_{i=1}^\infty  \| \rho_i {\mc S}_{0, i} - \rho_i {\mc S} \|_{C^0} \leq \sum_{i=1}^\infty \| S_{0, i} - S_i \|_{C^0(U_i)} \leq  \sum_{i=1}^\infty \frac{\delta}{2^i}  = \delta.
\eeqn

Now we prove the above claim. For each $x = (y_0, v_0) \in U$, the equivariance condition on $S_0$ implies $S_0(x) \in W_{G_x} \subset W$. Then by Lemma \ref{lemma251}, for $d \geq |G|$,  there exists $P_x \in \wh {\rm Poly}{}_G^d (V, W)$ such that $P_x(v_0) = S_0 (x)$. We view $P_x$ as a map 
\beqn
P_x: U \to W,\ P_x(y, v) = P_x(v),
\eeqn
which is a $G$-equivariant NC section of $E \to U$. Then for the given $\delta > 0$, one can find an open neighborhood $O_x \subset U$ of $x$ such that
\beqn
\sup_{x'\in O_x} | S_0 (x) - S_0 (x')| + \sup_{x'\in O_x} | P_x (x) - P_x(x')| \leq \delta.
\eeqn
As $U$ is Lindel\"of and paracompact, one can find countably many points $x_j$ such that the collection $O_{x_j}$ is a locally finite cover of $U$. Choose a partition of unity $\rho_j$ on $U$ subordinate to $O_{x_j}$. Define 
\beqn
{\mf q}: U \to \wh{\rm Poly}{}_G^d (V, W),\ {\mf q}(x) = \sum_{j=1}^\infty \rho_j (x) P_{x_j}.
\eeqn
We first check that this (not necessarily $G$-invariant) map defines a section that is close to the original section $S_0$. Indeed, for each $x = (y, v) \in U$, then one has 
\beqn
|S_0 (x) - {\mf q}(x)(v)| \leq \sum_{j=1}^\infty \rho_j(x) |S_0 (x) -  P_{x_j}(x)|.
\eeqn
For $j = 1, \ldots, l$, if $x \notin O_{x_j}$, then $\rho_j (x) = 0$ and the corresponding summand above vanishes; if $x \in O_{x_j}$, then 
\beqn
|S_0 (x) - P_{x_j}(x)|\leq |S_0 (x) - S_0 (x_j)| + |P_{x_j}(x) - P_{x_j}(x_j)| \leq \delta.
\eeqn
Hence one has
\beqn
\sup_{x \in U} |S_0 (x) - {\mf q}(x)(v)| \leq \delta.
\eeqn
Now we make ${\mf q}$ invariant by setting 
\beqn
{\mf p} (x) = \frac{1}{|G|} \sum_{\tau \in G} {\mf q} ( \tau x) \in \wh{\rm Poly}{}_G^d(V, W).
\eeqn
Then ${\mf p}$ satisfies the condition stated in the above claim.
\end{proof}

\section{Whitney Stratifications on the Variety $Z$}\label{section3}

This section serves as another facet of the technical foundation of the FOP transversality condition, which depends on an intricate analysis of a special class of affine varieties, especially their Whitney stratifications. Many technical discussions we carry out here originated from Parker's work \cite{BParker_integer}. 

Before the technical discussion, we explain why such study of Whitney stratifications is necessary. As mentioned in the introduction, by doing single-valued perturbation to a section of an orbifold vector bundle, the best hope is to have a stratified zero locus whose isotropy free part has codimension two or higher frontier. The local models for such a stratified zero locus are certain ``universal'' zero locus, i.e., the varieties labelled by $Z$. The Whitney stratification on $Z$ we specify in this section will give models for stratifying the perturbed zero locus. The challenge in the discussion is that one needs to identify a kind of canonical Whitney stratifications which have sufficient functorial properties so that the corresponding transversality condition is well-defined and is an open condition on an orbifold. One is then compelled to revisit arguments of Whitney and Mather etc. about the existence and properties of canonical Whitney stratifications on algebraic varieties and adopt them to the current situation.

In Subsection \ref{subsection31} we review basics about Whitney stratifications. Some technical results stated here are proved in Appendix \ref{appendixb}. In Subsection \ref{subsection32} we specify the consideration to a special class of algebraic varieties whose Whitney stratifications are essential for defining the FOP transversality condition. The proofs of several important propositions are given at the end of this section. In Subsection \ref{subsection33}, we describe a universal set which will index the strata of the zero loci of FOP transverse sections. The rest of this section contains proofs of results stated in Subsection \ref{subsection32}.

\subsection{Basic notions and facts about Whitney stratifications}\label{subsection31}

We first recall basic notions and facts about Whitney stratifications. Proofs of some results will be given in Appendix \ref{appendixb}. 

\begin{defn}\label{defn31}
Let $M$ be a smooth manifold.

\begin{enumerate}

\item Let $X, Y \subset M$ be disjoint smooth submanifolds. We say that $(X, Y)$ satisfies {\it Whitney's condition (b)} at $y \in Y$ if the following is true. Given a sequence $x_i \in X$ converging to $y$ and $y_i \in Y$ converging to $y$, 
suppose the secant line $\ov{x_i  y_i}$ \footnote{The secant line depends on a choice of local coordinate system near $y$, but the condition does not depend on such choices.} converges to a line $L_y \subseteq T_y M$ and $T_{x_i} X$ converges to a ${\rm dim } X$-dimensional subspace $H_y \subseteq T_y M$, then $L_y \subseteq H_y$. We say that $(X, Y)$ is {\it Whitney regular} if $(X, Y)$ satisfies Whitney's condition (b) at all points $y \in Y$.

\item A 
stratification (Definition \ref{defn_stratification}) of a subset $Z \subset M$, denoted by ${\mf Z} = \{Z_\lambda \}$, is called a \emph{Whitney stratification} if all strata are smooth submanifolds and each distinct pair $(Z_\lambda, Z_\mu )$ of strata is Whitney regular.

\item Let $N$ be a smooth manifold and $f: N \to M$ be a smooth map. For a Whitney stratified subset $Z \subset M$, $f$ is said to be {\it transverse} to $Z$ with respect to a Whitney stratification ${\mf Z}$ on $Z$, if $f$ is transverse to each stratum of ${\mf Z}$. The {\it naive pullback Whitney stratification} on $f^{-1}(Z)$ is the stratification
\beqn
f^* {\mf Z}:= \Big\{ f^{-1}(Z_\alpha) \neq \emptyset\ |\ Z_\alpha \in {\mf Z} \Big\}
\eeqn
(see Lemma \ref{lemma_pullback_Whitney} proving that this partition satisfies the axiom of frontier). Notice that there is a canonical map between the sets of strata:
\beqn
f_*: f^* {\mf Z} \to {\mf Z}.
\eeqn
In particular, for any open subset $U \subset M$, the {\it restriction} of ${\mf Z}$ to $U \cap Z$ is the pullback of ${\mf Z}$ by the natural inclusion $\iota: U \hookrightarrow M$, denoted by ${\mf Z}|_{U \cap Z}$. 

\item Suppose $f: N \to M$ is transverse to ${\mf Z}$ as above. The {\it refined pullback} of ${\mf Z}$ is the Whitney  stratification (cf. Lemma \ref{lemmab3}) on $f^{-1}(Z)$
\beqn
f^\dagger {\mf Z}:= \Big\{\ {\rm connected\ components\ of\ } f^{-1}(Z_\lambda)\ |\ Z_\lambda \in {\mf Z}\Big\}.
\eeqn
In particular, denote by ${\mf Z}^\dagger$ be the refined pullback by the identity map; in other words, ${\mf Z}^\dagger$ is obtained from ${\mf Z}$ by taking connected components. Notice that there is canonically induced order-preserving map
\beqn
f_\dagger: f^\dagger {\mf Z} \to {\mf Z}.
\eeqn
\end{enumerate}
\end{defn}

Whintey stratification may not exist on an arbitrary subset. It was originally proved by Whitney \cite{Whitney_1965} in the complex analytic category that analytic sets have  Whitney stratifications. We make a few more definitions.

\begin{defn}\label{defn_strongly_analytic}
Let $M$ be a complex manifold. 
\begin{enumerate}

\item A subset $U \subset M$ is called a {\it strongly analytic subset} if both $\ov{U}$ and $\ov{U} \setminus U$ are closed analytic subsets. $U$ is called a {\it strongly analytic submanifold} if it is an analytic submanifold and a strongly analytic set.

\item A Whitney stratification ${\mf Z}$ on a subset $Z \subset M$ is called {\it strongly analytic} if its strata are all strongly analytic submanifolds.
\end{enumerate}
\end{defn}

Whitney's construction actually provides an ``optimal'' Whitney stratification. Below we define how to compare different Whitney stratifications.

\begin{defn}\label{defn_minimal_Whitney}
Let $M$ be an $m$-dimensional   manifold and $Z \subset M$.
\begin{enumerate}
    \item The {\it dimension filtration} of a Whitney stratification ${\mf Z}= \{Z_\lambda\}$ on $Z$ is
    \beqn
    Z = {\mf Z}_m \supseteq {\mf Z}_{m-1} \supseteq \cdots \supseteq {\mf Z}_0
    \eeqn
    where 
    \beqn
    {\mf Z}_k:= \bigcup_{{\rm dim}(Z_\lambda) \leq k} Z_\lambda.
    \eeqn

    \item Two Whitney stratifications ${\mf Z}$ and ${\mf Z}'$ on $Z$ are called {\it equivalent}, denoted by ${\mf Z} \equiv {\mf Z}'$, if their dimension filtrations coincide. We write ${\mf Z} < {\mf Z}'$ if there exists $k$ such that 
    \begin{align*}
        &\ {\mf Z}_l = {\mf Z}_l'\ \forall l > k,\ &\ {\mf Z}_k \subsetneq {\mf Z}_k'.
    \end{align*}
    We write ${\mf Z} \leq {\mf Z}'$ if either ${\mf Z} \equiv {\mf Z}'$ or ${\mf Z} < {\mf Z}'$.

    \item Let $\mc{WS}^\infty(Z)$ be the set of all smooth Whitney stratifications on $Z$ and let ${\mc N}(Z) \subset \mc{WS}^\infty(Z)$ be a subset. We say ${\mf Z} \in {\mc N}(Z)$ is {\it minimal} within ${\mc N}(Z)$ if for any ${\mf Z}' \in \mc{N}(Z)$, either ${\mf Z} \equiv {\mf Z}'$ or ${\mf Z} < {\mf Z}'$. 
\end{enumerate}
\end{defn}

It is easy to see that if a minimal Whitney stratification exists, then it is unique up to equivalence.

\begin{defn}\label{defn_nice_Whitney_stratification}
Let $M$ be a smooth manifold equipped with a stratification ${\mf M} = \{ M_\alpha \}$ where each stratum $M_\alpha$ is a locally closed smooth submanifold (${\mf M}$ itself may not be a Whitney stratification). We say a Whitney stratification ${\mf Z}$ on a subset $Z \subset M$ {\it respects ${\mf M}$} if each $Z \cap M_\alpha$ is a union of strata in ${\mf Z}$. Denote by 
\beqn
\mc{WS}^\infty (Z; {\mf M}) \subset \mc{WS}^\infty(Z)
\eeqn
the subset of ${\mf M}$-respecting Whitney stratifications on $Z$.
\end{defn}

The case we are interested in comes from group actions.

\begin{defn}\label{defn_action_stratification}
Let $G$ be a finite group acting on a smooth manifold $M$. The {\it action stratification} on $M$ is the stratification ${\mf M}$ indexed by $M$-essential subgroups of $G$ whose strata are 
\beqn
M_H^*:= \Big\{ x \in M\ |\ G_x = H \Big\}.
\eeqn
For each $G$-invariant subset $Z \subset M$, let 
\beqn
\mc{WS}^\infty(Z; G)\subset \mc{WS}^\infty(Z)
\eeqn
be the subset of Whitney stratifications which respects the action stratification.
\end{defn}

In the appendix we prove the following theorem, which can be viewed as a refinement of Whintey's original theorem.

\begin{prop}\label{thm_nice_Whitney}
In the situation of Definition \ref{defn_nice_Whitney_stratification}, suppose $M$ is a complex manifold, each $M_\alpha$ is a strongly analytic submanifold, and $Z \subset M$ is a closed analytic subset. Then there exists a unique minimal element ${\mf Z} \in \mc{WS}^\infty(Z; {\mf M})$ that has connected strata. In addition, ${\mf Z}$ satisfies the following properties.
\begin{enumerate}
\item For each open subset $U \subset M$, ${\mf Z}|_{U \cap Z}$ (which may not have connected strata) is minimal in $\mc{WS}^\infty(U \cap Z, {\mf M}|_U)$.

\item ${\mf Z}$ is a strongly analytic Whitney stratification (Definition \ref{defn_strongly_analytic}).

\item If $f: M \to M$ is a diffeomorphism such that $f^* {\mf M} = {\mf M}$ and $f(Z) = Z$, then $f^* {\mf Z} = {\mf Z}$.
\end{enumerate}

\end{prop}

\begin{proof}
See Appendix \ref{appendix_nice_Whitney}.
\end{proof}

\begin{rem}
The construction needed for the proof was provided in the original paper of Whitney (see \cite[Section 18--21]{Whitney_1965}). However, as Whitney's proof does not explicitly discuss the minimality (which is important for our purpose), we need to revisit his argument and reproduce the proof in Appendix \ref{appendixb}.
\end{rem}

For certain applications (for example when proving the Arnold conjecture over ${\mb Z}$ in \cite{Bai_Xu_Arnold}), we need the property that the certain Whitney stratifications respect direct product. Notice that products of Whitney stratifications are still Whitney stratifications (see \cite[(1.2)]{Topological_stability}). In the appendix we also prove the following.

\begin{prop}\label{prop_product_nice}
Let $M$ resp. $N$ be complex manifolds equipped with stratifications ${\mf M}$ resp. ${\mf N}$ by strongly analytic submanifolds. Let $S \subset M$ resp. $T \subset N$ be closed analytic sets. Let ${\mf S}$ resp. ${\mf T}$ be the minimal ${\mf M}$-respecting resp. ${\mf N}$-respecting Whitney stratification on $S$ resp. $T$ provided by Proposition \ref{thm_nice_Whitney}. Then ${\mf S}\times {\mf T}$ is a minimal ${\mf M}\times {\mf N}$-respecting Whitney stratification on $S \times T$.
\end{prop}

\begin{proof}
See Appendix \ref{appendix_product_nice}.
\end{proof}

\subsection{The canonical Whitney stratifications on the variety $Z$}\label{subsection32}

Now we use the existence results (Theorem \ref{thm_nice_Whitney}) to specify the canonical Whitney stratifications on the class of algebraic varieties $Z$. We then state the properties of such Whitney stratifications which will be proved later in this section.

\subsubsection{The variety $Z$}

Recall that for each NC triple $(G, V, W)$ (Definition \ref{defn_NC_triple}) and $d \geq 0$, there is the vector space introduced in the proof of Lemma \ref{lemma244}
\beqn
M^d:= M^d(G, V, W) = V \times \wh{\rm Poly}{}_G^d(V, W)
\eeqn
where $G$ acts on the factor $V$. There is the natural {\it evaluation map}
\beqn
\ev: M^d \to W,\ (v, P) \mapsto P(v)\in W
\eeqn
whose zero locus is a $G$-invariant subvariety 
\beqn
Z^d = Z^d(G, V, W):= \Big\{ (v, P) \in M^d(G, V, W)\ |\ \ev(v, P) = 0 \Big\}.
\eeqn
In particular, the $W_G$-component of $P$ for a point in $Z^d$ is $0$. The group $G$ acts on the factor $V$ and $Z^d$ is $G$-invariant. Hence there is an action stratification on $M^d$ (Definition \ref{defn_action_stratification}). The strata are indxed by $V$-essential subgroups $H \subset_V W$ with
\begin{align*}
&\ V_H^* = \{ v \in V_H\ |\ G_v = H\},\ &\ M^d(G, V, W)_H^* = V_H^* \times \wh{\rm Poly}{}_G^d(V, W).
\end{align*}
For each $d\geq 0$, denote
\beqn
Z^d(G, V, W)_H^*:= Z^d(G, V, W) \cap M^d(G, V, W)_H^*.
\eeqn
We first look at this induced partition on $Z^d$.

\begin{prop}\label{prop_stratified_transverse}
For all $d \geq |G|$ and $H \subset_V G$, $Z^d(G, V, W)_H^*$ is a smooth submanifold of dimension ${\rm dim}_{\mb R} (\wh{\rm Poly}{}_G^d(V, W)) + {\rm dim}_{\mb R} (V_H) - {\rm dim}_{\mb R} (W_H)$.
\end{prop}

\begin{proof}
This is a consequence of Fukaya--Ono's lemma (Lemma \ref{lemma251}). For each $H \subset_V G$ and $P \in \wh{\rm Poly}{}_G^d(V, W)$, with respect to the basic decomposition $W = \check W_H \oplus W_H$, we can write $P = \check P_H + P_H$. Then the $H$-equivariance of $P$ implies that
\beqn
\check P_H|_{ V_H} \equiv 0.
\eeqn
Then 
\beqn
Z^d(G, V, W)_H^* =  \Big\{( v_H, P) \in   V_H^* \times \wh{\rm Poly}{}_G^d(V, W)\ |\ P_H(v_H) = 0 \Big\}.
\eeqn
Then Lemma \ref{lemma251} implies that when $d \geq |G|$, $Z^d(G, V, W)_H^*$ is transversely cut out, having dimension being ${\rm dim}_{\mb R} ( \wh{\rm Poly}{}_G^d(V, W) ) + {\rm dim}_{\mb R}  (V_H) - {\rm dim}_{\mb R} (W_H)$. Indeed, the linearization of the defining equation of $Z^d(G, V, W)_H^*$ is surjective: given $(v_H, P) \in Z^d(G, V, W)_H^*$,  Lemma \ref{lemma251} ensures that for any $w_0 \in W_H$, we can find a polynomial $P$ such that $P(v_H) = w_0$, and we can linearize the equation $P_H(v_H) = 0$ along the direction speficied by $P$ along the polynomial component. So, the smoothness follows accordingly, and the dimension formula follows by a count of number of parameters.
\end{proof}

\begin{rem}
Proposition \ref{prop_stratified_transverse} implies that $Z^d(G, V, W)$ has a natural partition into smooth manifolds. However, this partition is not a stratification, as the axiom of frontier fails (see Definition \ref{defn_straightening_model}). 
\end{rem}

\subsubsection{The canonical Whitney stratification}

Recall that for the NC pair $(G, V, W)$, one has the basic decompositions
\begin{align*}
&\ V = V_G \oplus \check V_G,\ &\ W = W_G \oplus \check W_G.
\end{align*}
Moreover, recall the definition
\beqn
\wh{\rm Poly}{}_G^d(V, W) = {\rm Poly}{}_G^d(\check V_G, \check W_G) \times W_G.
\eeqn
Then 
\beqn
M^d  = V_G \times M^d(G, \check V_G, \check W_G) \times W_G
\eeqn
where $\check M^d:= M^d(G, \check V_G, \check W_G)$ is  a complex vector space and 
\beqn
Z^d = V_G \times Z^d(G, \check V_G, \check W_G) \times \{0\}
\eeqn
where $\check Z^d:= Z^d(G, \check V_G, \check W_G) \subset \check M^d$ is a complex algebraic variety with a $G$-action. In this case, by Proposition \ref{thm_nice_Whitney}, there exists a unique minimal element
\beqn
\check {\mf Z}^d \in \mc{WS}^\infty( \check Z^d; G)
\eeqn
with connected strata. Using the identification $Z^d = V_G \times \check Z^d \times \{0\}$, there is a canonically induced Whitney stratification ${\mf Z}^d$ on $Z^d$ given by 
\beqn
{\mf Z}^d = \Big\{ V_G \times \check Z_\alpha^d\times \{0\} \ |\ \check Z^d_\alpha \in \check {\mf Z}^d \Big\}
\eeqn
It is easy to check that ${\mf Z}^d$ is a minimal Whitney stratification in $\mc{WS}^\infty(Z^d; G)$. We call ${\mf Z}^d$ the {\it canonical Whitney stratification} on $Z^d$. 

\begin{example}
Consider the situation where $G = {\mb Z}_3$, $V \cong {\mb C}$ is the weight 1 representation and $W \cong {\mb C}$ is the weight 2 representation. Let $d = 5$. Then 
\beqn
{\rm Poly}{}_{{\mb Z}_3}^5 (V, W) = \{ P(v) = av^2 + bv^5\ |\ a , b \in {\mb C} \}
\eeqn
and 
\beqn
Z^5 = \{ (a, b, v) \in {\mb C}^3 \ |\ av^2 + bv^5 = 0 \}.
\eeqn
The canonical Whitney stratification has three strata. Two top strata are 
\beqn
{\mf Z}^5_{\lambda_1} = \{ (a, b,v)\ |\ v \neq 0,\ a + bv^3 = 0 \}
\eeqn
and 
\beqn
{\mf Z}^5_{\lambda_2} = \{(a, b, v)\ |\ v = 0,\ a \neq 0 \}.
\eeqn
The intersection of their closure is the lowest stratum
\beqn
{\mf Z}^5_{\lambda_3} = \{ (a, b, v)\ |\ v = 0, a = 0 \}.
\eeqn
Notice that this Whitney stratification is the same as the minimal Whitney stratification of $Z^5$ without requiring to respect the action stratification.
\end{example}

\begin{example}
Consider $G = {\mb Z}_k$ and $R \cong {\mb C}$ being the weight 1 representation. Choose $V  = R \oplus R$, $W = R$, and $d = 1$. Let the variable of $V$ be $(x, y)$. Then
\beqn
{\rm Poly}_{{\mb Z}_k}^1(V, W) = \{ P(x, y) = ax + by\ |\ a, b \in {\mb C} \}
\eeqn
and $Z^1$ is the cone $ax + by = 0$. It only singularity is the origin and its minimial Whitney stratification (without referring to the group action) has two strata: the smooth locus and the cone point. However, with respect to the group action, the canonical Whitney stratification is finer: its top stratum is 
\beqn
{\mf Z}_{\lambda_1}^1 = \{ (a, b, x, y)\ |\ (x, y) \neq (0, 0),\ ax + by = 0 \}
\eeqn
which is properly contained in the smooth locus of $Z^1$.
\end{example}

The various properties of the FOP transversality condition listed in Theorem \ref{thm11} relies on the following important properties of the canonical Whitney stratification. 

\subsubsection{Invariance}

\begin{prop}\label{prop_nice_invariance}
Given an isomorphism $(G, V, W) \cong (G', V', W')$ of NC triples which induces a linear isomorphism $f: M^d(G, V, W) \to M^d(G', V', W')$ which sends $Z^d(G, V, W)$ to $Z^d(G', V', W')$, there holds 
\beqn
f^* {\mf Z}^d(G', V', W') = {\mf Z}^d(G, V, W).
\eeqn
In particular, ${\mf Z}^d(G, V, W)$ is ${\rm Aut}(G, V, W)$-invariant, where ${\rm Aut}(G, V, W)$ denotes the group consisting of a pair of $G$-equivariant linear isomorphisms of $V$ and $W$.
\end{prop}

\begin{proof}
It is equivalent to consider the case when $(G, V, W) = (G', V', W')$ and $f$ is induced from a pair of $G$-equivariant automorphisms $f_V: V \to V$ and $f_W: W \to W$. Then $f$ induces a complex linear automorphism of the subspace $M^d(G, \check V_G, \check W_G)$ preserving the complex variety $Z^d(G, \check V_G, \check W_G) $ as well as the action stratification and $\check Z^d$. Then by Proposition \ref{thm_nice_Whitney}, $f^* {\mf Z}^d =  {\mf Z}^d$. 
\end{proof}

Proposition \ref{prop_nice_invariance} allows us to extend the canonical Whitney stratification to the family case. Let $(G, F, E)$ be a family NC triple over a manifold $X$ (Definition \ref{defn246}). The one can define the vector bundle
\beqn
M^d(G, F, E) \to X
\eeqn
and the $G$-invariant subbundle
\beqn
Z^d(G, F, E)\subset M^d(G, F, E).
\eeqn
As the structure group of $M^d(G, F, E)$ is contained in ${\rm Aut}(G, V, W)$, by Proposition \ref{prop_nice_invariance}, there is a locally trivial Whitney stratification ${\mf Z}^d(G, F, E)$ on $Z^d(G, F, E)$, which is still called the {\it canonical Whitney stratification} on $Z^d(G, F, E)$.

\subsubsection{Change of degrees}

The following proposition is necessary for showing the FOP transversality condition is well-defined and independent of the cut-off degree $d$. Consider $d \geq 0$ and $d' > d$. Then there is the natural inclusion map 
\beq\label{eqn31}
\sigma_{d'd}: M^d (G, V, W) \to M^{d'}(G, V, W),\ (v,   P) \mapsto (v,   P )
\eeq
which respects the action stratifications such that $\sigma_{d'd}^{-1} ( Z^{d'}) =  Z^d$. 

\begin{prop}\label{prop313}
For any NC triple $(G, V, W)$, there exists $d_1 = d_1(G, V, W) \geq |G|$ such that when $d > d' \geq d_1$, there holds 
\beqn
\sigma_{d'd}^* {\mf Z}^{d'} = {\mf Z}^d.
\eeqn
Moreover, the corresponding map between the sets of strata
\beqn
(\sigma_{d'd})_*: {\mf Z}^d \to {\mf Z}^{d'}
\eeqn
is bijective.
\end{prop}

\begin{proof}
See Subsection \ref{subsection_34}.
\end{proof}

\subsubsection{Change of groups}

The following proposition is necessary for showing the FOP transversality condition is an open condition on orbifolds. It is also the hardest one to prove. We first recall and introduce some notations. Let $(G, V, W)$ be an NC triple and $H \subset_V G$ be an essential subgroup. Using the basic decomposition $V = V_H \oplus \check V_H$, a point of $V$ is denoted by $(v_H, \check v_H)$. Given an element $P \in \wh {\rm Poly}{}_G^d (V, W)$ viewed as a function from $V$ to $W$, using the basic decomposition $W = \check W_H \oplus W_H$, we also write 
\beqn
P = \check P_H + P_H.
\eeqn
Then one has defined the map in the proof of Lemma \ref{lemma244} 
\beq\label{group_change_map}
\begin{split}
\mu: M^d(G, V, W) &\ \to  M^d(H, V, W)\\
 (v, P) &\ \mapsto \Big( v, \check P_H(v_H, \cdot) + P_H(v)  \Big).
\end{split}
\eeq
Here the term $P_H(v)$ also contains the constant map component with value in $W_G$. Notice that $\mu$ maps $Z^d(G, V, W)$ into $Z^d(H, V, W)$. We shall not expect that the pullback Whitney stratification by $\mu$ coincides with ${\mf Z}^d(G, V, W)$, because the weaker symmetry required by the subgroup $H$ may not see some deeper strata of ${\mf Z}^d(G, V, W)$. However, they do coincide when restricted to the open subset of points whose isotropy group under the $G$-action is a subgroup of $H$.

\begin{prop}\label{prop314}
Given an NC triple $(G, V, W)$, there exists $d_2 \geq d_1$ such that for all $d \geq d_2$ the following is true. For each $H \subset_V G$, consider the open subset 
\beqn
V_H^+:= \{ v \in V\ |\ G_v \subset H \} \subset V.
\eeqn
Denote
\begin{align*}
&\ M^d(G, V, W)_H^+:= V_H^+ \times \wh {\rm Poly}{}_G^d (V, W),\ &\ M^d(H, V, W)_H^+:= V_H^+ \times \wh{\rm Poly}{}_H^d(V, W)
\end{align*}
\beqn
Z^d(G, V, W)_H^+ := Z^d(G, V, W) \cap M^d(G, V, W)_H^+.
\eeqn
Then $\mu|_{M^d(G, V, W)_H^+}$ is transverse to ${\mf Z}^d(H, V, W)$ and 
\beqn
\Big( \mu^* {\mf Z}^d(H, V, W) \Big)|_{Z^d(G, V, W)_H^+} = {\mf Z}^d(G, V, W)|_{Z^d(G, V, W)_H^+}.
\eeqn
Moreover, the induced maps between the sets of strata 
\beqn
\mu_*: {\mf Z}^d(G, V, W)|_{Z^d(G, V, W)_H^+} \to {\mf Z}^d(H, V, W)
\eeqn
is invertible and makes the following diagram commutes (when $d' \geq d$).
\beq\label{eqn33}
\vcenter{ \xymatrix{  {\mf Z}^d( H, V, W) \ar[rr]^-{(\mu_*)^{-1}} \ar[d]_{(\sigma_{d'd})_*} & & {\mf Z}^d(G, V, W) \ar[d]^{(\sigma_{d'd})_*} \\
         {\mf Z}^{d'}(H, V, W) \ar[rr]_-{(\mu_*)^{-1}} & & {\mf Z}^{d'} (G, V, W)   }}
\eeq
\end{prop}

\begin{proof}
See Subsection \ref{subsection37}.
\end{proof}

\subsubsection{Product}

Another property of the canonical Whitney stratification is that it respects direct product. This property is responsible for the {\bf Product Property} stated in Theorem \ref{thm11} which is important in proving the integral Arnold conjecture in \cite{Bai_Xu_Arnold} and in other potential chain-level applications. 

For $i = 1, 2$, let $(G_i, V_i, W_i)$ be NC triples and denote
\beqn
(G, V, W) = (G_1 \times G_2, V_1\oplus V_2, W_1 \oplus W_2).
\eeqn

\begin{prop}\label{prop_product_Z}
For $d$ sufficiently large, consider the map 
\beqn
\begin{split}
 \xi:  M^d(G_1, V_1, W_1) \times M^d(G_2, V_2, W_2) & \to M^d(G, V, W) \\
    (v_1, P_1, v_2, P_2 ) & \mapsto \left( \left[ \begin{array}{c} v_1 \\ v_2  \end{array}\right], \left[ \begin{array}{cc} P_1 & 0 \\ 0 & P_2 \end{array}\right]  \right).
\end{split}
\eeqn
Then $\xi^* {\mf Z}^d(G, V, W) = {\mf Z}^d(G_1, V_1, W_2) \times {\mf Z}^d(G_2, V_2, W_2)$.
\end{prop}

\begin{proof}
See Subsection \ref{subsection_36}.
\end{proof}

\subsubsection{Stabilization property}

The following proposition is responsible for the {\bf Stabilization Property} of Theorem \ref{thm11}. We prove the following proposition. 

\begin{prop}\label{prop_Whitney_stabilization}
Let $(G, V, W)$ be an NC triple and $R$ be an NC representation. For any $d \geq 1$, consider the map 
\beq\label{polynomial_stabilization}
\begin{split}
\eta: M^d(G, V, W) & \to M^d(G, V\oplus R, W\oplus R) \\
(v, P) & \mapsto \left(  \left[ \begin{array}{c} v \\ 0  \end{array} \right], \left[ \begin{array}{cc} P & 0  \\ 0 & {\rm proj}_{\check R_G} 
\end{array}\right]  \right)
\end{split}
\eeq
where ${\rm proj}_{\check R_G}: R \to \check R_G$ 
is the canonical projections induced from the basic decomposition $R = R_G \oplus {\check R_G}$. Then for $d \geq d_1 $ (given by Proposition \ref{prop313}), $\eta$ is transverse to ${\mf Z}^d(G, V\oplus R, W\oplus R)$ and 
\beqn
\eta^\dagger {\mf Z}^d(G, V\oplus R, W\oplus R) = {\mf Z}^d(G, V, W)
\eeqn
(where $\eta^\dagger$ is the refined pullback, see Definition \ref{defn31}). Moreover, the induced map
\beqn
\eta_\dagger: {\mf Z}^d(G, V, W) \to {\mf Z}^d(G, V\oplus R, W \oplus R)
\eeqn
makes the following diagram commutes.
\beq\label{eqn35}
\vcenter{ \xymatrix{ {\mf Z}^d(G, V, W) \ar[r]^-{\eta_\dagger} \ar[d]_{(\sigma_{d'd})_*} & {\mf Z}^d(G, V\oplus R, W \oplus R) \ar[d]^{(\sigma_{d'd})_*} \\
                    {\mf Z}^{d'}(G, V, W) \ar[r]_-{\eta_\dagger} & {\mf Z}^{d'}(G, V\oplus R, W\oplus R)  }}
\eeq
\end{prop}

\begin{proof}
See Subsection \ref{subsection_proof_stabilization_property}.
\end{proof}

\begin{rem}
Proposition \ref{prop313}---\ref{prop_Whitney_stabilization} are stated for general NC triples $(G, V, W)$. However, it is very easy to reduce the proofs to the case when $V_G =  W_G = \{0\}$; in particular, $V$ and $W$ are complex representations. In the subsections proving these propositions, we always make this simplifying assumption. 
\end{rem}

\subsection{The set of universal strata}\label{subsection33}

Now we use the previous technical results about the Whitney stratifications on the various varieties $Z^d(G, V, W)$ to define a set of universal strata ${\mf Z}_k^{\rm univ}$ whose elements can label the strata of the zero locus of an FOP transverse section (see Theorem \ref{thm11}). Roughly, we define ${\mf Z}_k^{\rm univ}$ by taking the disjoint union of ${\mf Z}^d(G, V, W)$ for all NC triples of virtual dimension $k$ modulo a certain equivalence relation. To make this idea more precise, we need to specify a certain category and define ${\mf Z}_k^{\rm univ}$ as a colimit.

\subsubsection{Category of NC triples}

\begin{defn}
\begin{enumerate}

\item 
A {\it morphism of NC triples} from $(H, V, W)$ to $(G, X, Y)$ is a finite composition of the following two types of morphisms
\begin{enumerate}
    \item A {\it stabilization} consists of a group isomorphism $H \cong G$, an NC representation $R$ of $G$, and a pair of isomorphisms $X \cong V \oplus R$ and $Y \cong W \oplus R$ of NC representations of $G$.

    \item An {\it embedding} consists of a group embedding $H \hookrightarrow G$ and a pair of isomorphisms $V \cong X$ and $W \cong Y$ of NC representations of $H$.
\end{enumerate}

\item Two NC triples $(H, V, W)$ and $(G, X, Y)$ are called {\it stably equivalent} if they admit stabilizations that are isomorphic (isomorphism of groups + induced isomorphisms of representations). A stable equivalence class is called a {\it stable NC isotropy type}. The set of stable NC isotropy types is denoted by $\uds{\mb \Gamma}^{\rm NC}$. For each $k \in {\mb Z}$, let $\uds{\mb \Gamma}^{\rm NC}_k\subset \uds{\mb \Gamma}^{\rm NC}$ be the subset of stable NC isotropy types with virtual dimension ${\rm dim}_{\mb R} (V) - {\rm dim}_{\mb R} (W) = k$.

\item Define a partial order in $\uds{\mb\Gamma}_k^{\rm NC}$ as follows. We denote $\uds\gamma \leq \uds\delta$ if there exist a representative $(H, V, W)$ of $\uds\gamma$, a representative $(G, X, Y)$ of $\uds\delta$, and an embedding from $(H, V, W)$ to $(G, X, Y)$. 
\end{enumerate}
\end{defn}

Notice that there is a natural poset map ${\mb\Gamma}_k^{\rm NC} \to \uds{\mb\Gamma}_k^{\rm NC}$. For each NC pair $({\mc U}, {\mc E})$ of orbifolds and vector bundles whose virtual dimension is $k$ and a stable NC isotropy triple $\uds\gamma \in \uds{\mb\Gamma}_k^{\rm NC}$, define
\beqn
{\mc U}_{\uds\gamma}:= \bigcup_{\gamma\mapsto \uds\gamma} {\mc U}_\gamma
\eeqn
where we take the (disjoint) union of all isotropy types $\gamma$ which are sent to $\uds\gamma$.

\subsubsection{A particular category}

The above technical results on Whitney stratifications allow us to gather all strata for different isotropy types. For our purpose related to Theorem \ref{thm11}, we first consider a particular category, denoted by 
\beqn
\uds{\bf Str}_k^{\rm NC}.
\eeqn
Its objects are triples $({\mf Z}, \rho, n)$ where ${\mf Z}$ is a partially ordered set, $\rho: {\mf Z} \to \uds{\mb \Gamma}^{\rm NC}_k$ is a monotone map, and $n: {\mf Z} \to 2{\mb Z}$ is a strictly monotone map. A morphism from $({\mf Z}_1, \rho_1, n_1)$ to $({\mf Z}_2, \rho_2, n_2)$ is a strictly monotone poset map $\iota_{21}: {\mf Z}_1 \to {\mf Z}_2$ such that $\rho_1 = \rho_2 \circ \iota_{21}$ and $n_1 = n_2 \circ \iota_{21}$. 

\begin{lemma}
Colimits exist in $\uds{\bf Str}_k^{\rm NC}$.
\end{lemma}

\begin{proof}
Let $\uds{\bf Poset}$ be the category of posets with morphisms being poset maps and let $\uds{\bf Poset}^s$ be the category of posets with morphisms being strictly monotone poset maps. In both categories colimits exist. On the other hand, there are forgetful functors
\beqn
\xymatrix{ \uds{\bf Str}_k^{\rm NC} \ar[r]^-{\pi} & \uds{\bf Poset}^s \ar[r]^-{\iota} & \uds{\bf Poset}}
\eeqn
where $\pi$ takes a triple to its underlying poset. Let ${\bf F}: \uds{\bf A} \to \uds{\bf Str}_k^{\rm NC}$ be a diagram (functor) where $\uds{\bf A}$ is any category. Then one obtains a colimit
\beqn
\colim_a \pi( {\bf F}(a)) \in {\rm Ob} \uds{\bf Poset}^s.
\eeqn
As the maps $n: {\mf Z} \to 2{\mb Z}$ is a morphism of $\uds{\bf Poset}^s$, it follows that $n$ descends to a strictly monotone map 
\beqn
n: \colim_a \pi({\bf F}(a)) \to 2{\mb Z}.
\eeqn
Moreover, $\iota(\displaystyle \colim_a \pi({\bf F}(a)))$ is a colimit of the induced diagram $\uds{\bf A} \to \uds{\bf Poset}$, hence there exists a natural map 
\beqn
\rho: \colim_a \pi({\bf F}(a)) \to \uds{\mb\Gamma}^{\rm NC}_k.
\eeqn
It is  straightforward to check that $(\displaystyle \colim_a \pi({\bf F}(a)), n, \rho)$ is a colimit of ${\bf F}(a)$.
\end{proof}

\subsubsection{The universal strata}

The sets of strata of the canonical Whitney stratifications we considered give rise to objects of the category $\uds{\bf Str}_k^{\rm NC}$. Let $(G, V, W)$ be an NC triple of virtual dimension $k$. Define a map
\beqn
\rho_k: {\mf Z}^d(G, V, W) \to \uds{\mb \Gamma}^{\rm NC}_k
\eeqn
as follows: if a stratum $Z^d_\alpha \in {\mf Z}^d(G, V, W)$ is contained in $V_H^* \times \wh{\rm Poly}{}_G^d(V, W)$ for some $V$-essential subgroup $H$, then define $\rho(Z^d_\alpha)$ to be the isotropy type represented by the triple $(H, V, W)$. On the other hand, define a map 
\beqn
n_k: {\mf Z}^d(G, V, W) \to 2{\mb Z}
\eeqn
which assigns to a stratum $Z^d_\alpha$ the number
\beq\label{defn:n-map}
n_k (Z^d_\alpha):= {\rm dim}_{\mb R} (Z^d_\alpha) - {\rm dim}_{\mb R} ( \wh{\rm Poly}{}_G^d(V, W)) - k \in 2 {\mb Z}.
\eeq
Then $({\mf Z}^d(G, V, W), \rho_k, n_k)$ is an object of $\uds{\bf Str}_k^{\rm NC}$. Note that Proposition \ref{prop313} allows us to relate this set for different $d$. Indeed, the bijective map (for $d < d'$) induces a morphism
\beqn
\sigma_{d'd}\in {\rm Hom}_{\uds{\bf Str}_k^{\rm NC}} \left( {\mf Z}^d(G, V, W), {\mf Z}^{d'} (G, V, W) \right).
\eeqn
Then define the colimit
\beqn
{\mf Z}(G, V, W):= \colim_{d} {\mf Z}^d (G, V, W) \in \uds{\bf Str}_k^{\rm NC}.
\eeqn

We would like to make it a functor from the category of NC triples to $\uds{\bf Str}_k^{\rm NC}$. 

\begin{lemma}
Suppose $(H, V, W) \to (G, X, Y)$ is a stabilization morphism. Then there is a canonical map 
\beqn
\eta: {\mf Z}(H, V, W) \to {\mf Z}(G, X, Y).
\eeqn
\end{lemma}

\begin{proof}
We identify $H$ with $G$, $V$ with a subspace of $X$, and $W$ with a subspace of $Y$. By definition, there exist a representation $R$ of $G$ and an equivariant decomposition 
\beqn
X \cong V \oplus R, Y \cong W \oplus R.
\eeqn
By Proposition \ref{prop_Whitney_stabilization}, each such decomposition induces a map 
\beqn
\eta_\dagger: {\mf Z}^d(H, V, W) \to {\mf Z}^d(G, X, Y)
\eeqn
which is compatible with the maps $\sigma_{d'd}$ when $d'>d$. Moreover, as the spaces of above equivariant decompositions are connected, the induced map is independent of the choices of the decompositions. Hence there is a well-defined morphism
\beqn
{\mf Z}(H, V, W) \to {\mf Z}(G, X, Y).\qedhere
\eeqn
\end{proof}

Similarly, from Proposition \ref{prop314} one can define a morphism for embeddings of NC triples. 

\begin{lemma}
Suppose $(H, V, W) \to (G, X, Y)$ be an embedding of NC triples. Then there is a canonical morphism 
\beqn
(\mu_*)^{-1}: {\mf Z}(H, V, W) \to {\mf Z}(G, X, Y).
\eeqn
\end{lemma}

Then the association $(G, V, W) \mapsto {\mf Z}(G, V, W)$ is a functor from the category $\uds{\mb \Gamma}_k^{\rm NC}$ to the category $\uds{\bf Str}_k^{\rm NC}$. Then we define 
\beqn
{\mf Z}_k^{\rm univ}:= \colim_{(G, V, W)} {\mf Z}(G, V, W)
\eeqn
which we call the set of \emph{universal strata} (of virtual dimension $k$).

\subsection{Proof of Proposition \ref{prop313}}\label{subsection_34}

First notice that Proposition \ref{prop313} reduces to the case that $V_G = W_G = 0$. Hence within this subsection, $V$ are $W$ are complex representations of $G$ containing no trivial subrepresentations. As a consequence
\beqn
\wh{\rm Poly}{}_G^d(V, W) = {\rm Poly}{}_G^d(V, W).
\eeqn
Our proof closely resembles that of \cite[Lemma 4.11]{BParker_integer}. 

First we prove an algebraic result also used in \cite{BParker_integer} without providing a proof or reference.

\begin{lemma}\label{lemma313}
Suppose $V$ and $W$ are finite-dimensional complex representations of $G$. Then ${\rm Poly}_G(V, W)$ is a finitely generated module over the ring ${\rm Poly}_G(V, {\mb C})$ of $G$-invariant polynomials on $V$.
\end{lemma}

\begin{proof}
The proof follows from \cite{inv-finite}, we present it here for completeness. By Hilbert's basis theorem, given any finite-dimensional complex $G$-representation $V'$, the ring ${\rm Poly}_G(V', {\mb C})$ is finitely generated. Now let $V' = V \oplus W^{\vee}$ where $W^{\vee}$ is the dual to $W$ endowed with the corresponding $G$-action. The ring ${\rm Poly}_G(V', {\mb C})$ has a bi-grading by keeping track of the degree of the $V$-coordinates and $W^{\vee}$-coordinates respectively. Choose $h_1, \dots, h_r, h_{r+1}, \dots, h_{r+m}, \dots, h_n$ which generate ${\rm Poly}_G(V', {\mb C})$ such that $h_1, \dots, h_r$ have $W^{\vee}$-degree $0$, $h_{r+1}, \dots, h_{r+m}$ have $W^{\vee}$-degree $1$, and $h_{r+m+1}, \dots, h_n$ have higher $W^{\vee}$-degrees.

Note that ${\rm Poly}_G(V, W)$ can be identified with the subset of ${\rm Poly}_G(V', {\mb C})$ consisting of elements with $W^\vee$-degree $1$. Then any element in ${\rm Poly}_G(V, W)$ can be written as a linear combination of products of $h_1, \dots, h_r, h_{r+1}, \dots, h_{r+m}$. Observe that $h_1, \ldots, h_r$ are indeed elements of ${\rm Poly}_G(V, {\mb C})$, it shows ${\rm Poly}_G(V, W)$ is generated by $h_{r+1}, \dots, h_{r+m}$ as a module over ${\rm Poly}_G(V, {\mb C})$.
\end{proof}

Next we construct a left inverse of the map $ \sigma_{d'd}$ of \eqref{eqn31}.

\begin{lemma}(cf. \cite[Lemma 4.11]{BParker_integer})\label{lemma323}
Given an NC triple $(G, V, W)$, there exists $d_1$ such that for all $d' > d \geq d_1$, there is a $G$-invariant holomorphic map
\beqn
\rho_{d'd}: M^{d'}(G, V, W) \to M^d(G, V, W)
\eeqn
satisfying the following properties.
\begin{enumerate}
    \item $\rho_{d'd}$ is a submersion with connected fibers and preserves the action stratification.

    \item $\rho_{d'd} \circ \sigma_{d'd}$ is the identity map of $M^d(G, V, W)$.
    
    \item $\ev \circ \rho_{d'd} = \ev: M^{d'}(G, V, W) \to W$.
\end{enumerate}
\end{lemma}

\begin{proof}[Proof of Proposition \ref{prop313}]
Lemma \ref{lemma323}--(3) implies that $Z^{d'} = \rho_{d'd}^{-1}(Z^d)$. By Lemma \ref{lemma323}--(1) and Proposition \ref{prop_pullback_nice}, $\rho_{d'd}^* {\mf Z}^d = {\mf Z}^{d'}$. On the other hand, because $\rho_{d'd} \circ \sigma_{d'd} = {\rm Id}$ is transverse ${\mf Z}^d$, $\sigma_{d'd}$ is transverse to $\rho_{d'd}^* {\mf Z}^d$ (see Lemma \ref{lemmab4}). Hence
\beqn
\sigma_{d'd}^* {\mf Z}^{d'} = \sigma_{d'd}^* \rho_{d'd}^* {\mf Z}^d = (\rho_{d'd} \circ \sigma_{d'd})^* {\mf Z}^d = {\mf Z}^d.
\eeqn
Moreover, because $\rho_{d'd}$ is surjectve and has connected fibers, for each stratum $Z^d_\alpha\in {\mf Z}^d$, $\rho_{d'd}^{-1}(Z^d_\alpha)$ is a single stratum of ${\mf Z}^{d'}$ and $\sigma_{d'd}$ maps $Z_\alpha^d$ into this stratum. Hence $\sigma$ induces a bijection between the sets of strata.
\end{proof}

\begin{proof}[Proof of Lemma \ref{lemma323}] 
As ${\rm Poly}_G(V, W)$ is finitely generated over ${\rm Poly}_G(V, {\mb C})$ (see Lemma \ref{lemma313}), one can find a sufficiently large $d_0$ such that ${\rm Poly}_G^{d_0} (V, W)$ contains a set of generators $Q_1, \ldots, Q_m$. Then when $d \geq d_0$, let ${\rm Homo}_G^d (V, W) \subset {\rm Poly}_G^d (V, W)$ be the subset of $G$-equivariant homogeneous polynomial maps of degree $d$. Then  
\beqn
{\rm Poly}_G^{d'} (V, W) = {\rm Poly}_G^d (V, W) \oplus \bigoplus_{l=d+1}^{d'} {\rm Homo}_G^l (V, W).
\eeqn
Choose a basis $P_1, \ldots, P_k$ of $\bigoplus_{l=d+1}^{d'} {\rm Homo}_G^l (V, W)$, viewed as a ${\mb C}$-vector space. Then there exist polynomials $h_{ij} \in {\rm Poly}_G(V, {\mb C})$ for $1 \leq i \leq k$, $1 \leq j \leq m$ such that
\beqn
P_i = \sum_{j=1}^m h_{ij} Q_j.
\eeqn
Then for each $P \in {\rm Poly}_G^{d'} (V, W)$, write $P = P' + P''$ with $P' \in {\rm Poly}_G^d(V, W)$ and $P'' \in \bigoplus_{l=d+1}^{d'} {\rm Homo}_G^l (V, W)$. We can write  
\beqn
P'' = \sum_{i=1}^k a_i P_i.
\eeqn
Then define
\beq\label{eqn37}
\rho_{d'd} (v, P) = \left( v,   P' + \sum a_i h_{ij}(v) Q_j \right).
\eeq
Abbreviate $\rho_{d'd}$ by $\rho$ and $\sigma_{d'd}$ by $\sigma$. Then $\rho$ is a submersion, preserves the action stratification, and satisfies $\rho \circ \sigma = {\rm Id}_{M^d}$. Moreover, for all $(v, P) \in M^{d'}$ one has
\beqn
\ev( \rho (v, P)) = P'(v) + \sum a_i h_{ij}(v) Q_j(v) = P'(v) + P''(v) = P(v). 
\eeqn
Lastly, we prove that $\rho$ has connected fibers. For any $(v, P) \in M^{d'}$ with $P = P' + P''$ decomposed as above, consider the path
\beqn
(v, P_t) = (v, P' + P_t'') = \left( v, P' + t \sum a_i h_{ij}(v) Q_j + (1-t) P''\right).
\eeqn
Then $\rho (v, P_t) = \rho(v, P)$ for all $t$. Notice that $(v, P_1) = \sigma(\rho (v, P))$. Therefore, every point in the fiber of $\rho(v, P)$ can be connected with $\sigma(\rho (v, P))$. Hence fibers of $\rho$ are all connected.
\end{proof}

\subsection{Proof of Proposition \ref{prop_product_Z}}\label{subsection_36}

Again, it suffices to consider the case that $V_G = W_G = \{0\}$. 

\begin{lemma}
There exists a surjective submersion with connected fibers  
\beqn
\xi': M^d(G_1 \times G_2, V_1 \oplus V_2, W_1 \oplus W_2 ) \to  M^d(G_1, V_1, W_1) \times M^d(G_2, V_2, W_2)
\eeqn
satisfying the following conditions. 
\begin{enumerate}
    \item $\xi' \circ \xi$ is the identity map.
    
    \item For each $(v, P)$ in the domain of $\xi'$, ${\rm ev}( \xi'(v, P)) = {\rm ev}(v, P)$.
    
\end{enumerate}
\end{lemma}

\begin{proof}
For each $P \in {\rm Poly}_{G_1\times G_2}^d ( V_1 \oplus V_2,  W_1 \oplus W_2)$, denote its $W_1$-component by $P_1$ and its $W_2$-component by $P_2$. Then we can regard $P_1$ as a $G_2$-invariant  map 
\beqn
P_1 \in {\rm Poly}_{G_2}^d ( V_2, {\rm Poly}_{G_1}^d ( V_1, W_1))
\eeqn
and regard $P_2$ as a $G_1$-invariant polynomial map
\beqn
P_2 \in {\rm Poly}_{G_1}^d (V_1, {\rm Poly}_{G_2}^d (V_2, W_2)).
\eeqn
Then for $v = (v_1, v_2) \in V_1 \oplus V_2$, define 
\beqn
\xi' (v, P) =\Big( (v_1, P_1(\cdot, v_2)), (v_2, P_2(v_1, \cdot)) \Big).
\eeqn
Then it is easy to verify that $\xi' \circ \xi$ is the identity map and that 
\beqn
{\rm ev}(\xi' (v, P)) = {\rm ev}(v, P). 
\eeqn
Further, it is obvious that $\xi'$ is a surjective submersion. Moreover, if $\xi'(v, P) = \xi'(u, Q)$, then $v = u$ and $\xi'(v, P) = \xi'(v, t P + (1-t) Q)$ for all $t \in [0,1]$. Hence $\xi'$ has connected fibers.
\end{proof}

\begin{proof}[Proof of Proposition \ref{prop_product_Z}]
Notice that the analytic submersion $\xi'$ preserves the action stratification. Then by Proposition \ref{prop_pullback_nice}, $\xi'|_{U^d}$ pulls back the canonical Whitney stratification on $Z_1^d \times Z_2^d$ to the canonical one on $Z^d $. The former is the product ${\mf Z}_1^d \times {\mf Z}_2^d$ by Proposition \ref{prop_product_nice}. Therefore  
\beqn
(\xi')^* ( {\mf Z}_1^d \times {\mf Z}_2^d ) =  {\mf Z}^d.
\eeqn
Then as $\xi' \circ \xi$ is the identity, $\xi$ is transverse to $(\xi')^* ({\mf Z}_1^d \times {\mf Z}_2^d)$. Hence 
\beqn
\xi^* {\mf Z}^d  = \xi^*  {\mf Z}^d  = \xi^* ( (\xi')^* ({\mf Z}_1^d \times {\mf Z}_2^d)) = (\xi' \circ \xi)^* ({\mf Z}_1^d \times {\mf Z}_2^d) = {\mf Z}_1^d \times {\mf Z}_2^d. \qedhere
\eeqn
\end{proof}

\subsection{Proof of Proposition \ref{prop_Whitney_stabilization}}\label{subsection_proof_stabilization_property}

Again, it suffices to consider the case when $V_G = W_G = R_G = 0$; in particular, ${\rm proj}_{\check R_G}$ is the identity map of $R$. Consider equivariant polynomial maps from $V \oplus R$ to $W \oplus R$. A polynomial map from $V \oplus R$ to $W \oplus R$ consists of a polynomial map from $V \oplus R$ to $W$ and a polynomial map from $V \oplus R$ to $R$. For convenience we write them as
\begin{align*}
&\ P \in {\rm Poly}_G^d(V \oplus R, W),\ &\  {\rm Id}_R + Q \in {\rm Poly}_G^d(V\oplus R, R).
\end{align*}
Denote variables of $V \oplus R$ be $(x, y)$ and the image of the map $\eta$ in \eqref{polynomial_stabilization}  by 
\beqn
\wt {\rm Poly}{}_G^d(V\oplus R, W \oplus R) \subset {\rm Poly}{}_G^d(V \oplus R, W \oplus R),
\eeqn
i.e., polynomial maps $(P_0, {\rm Id}_R)$ where $P_0: V \to W$ does not depends on $y\in R$. We introduce an intermediate space
\beqn
\widecheck{\rm Poly}_G(V \oplus R, W \oplus R)
\eeqn
consisting of pairs of the form 
\beqn
(x, y) \mapsto (P_0 (x), y + Q_0 (x)),\ {\rm where}\ P_0 \in {\rm Poly}_G(V, W),\ Q_0 \in {\rm Poly}_G(V, R).
\eeqn
Then consider
\beqn
\widecheck{Z}^d(G, V\oplus R, W \oplus R) = Z^d(G, V\oplus R, W \oplus R) \cap \Big( (V \oplus R) \times \widecheck{\rm Poly}_G( V\oplus R, W \oplus R) \Big)
\eeqn
By Proposition \ref{thm_nice_Whitney}, there are a canonical Whitney stratification $\widecheck{\mf Z}^d(G, V\oplus R, W \oplus R)$ on $\widecheck Z^d(G, V\oplus R, W \oplus R)$ and a canonical Whitney stratification $\wt{\mf Z}^d(G, V\oplus R, W \oplus R)$ on $\wt{Z}^d(G, V\oplus R, W\oplus R)$. The latter can be identified with ${\mf Z}^d(G, V, W)$. Then we can factorize the map $\eta$ as 
\beqn
\eta = \eta_1 \circ \eta_2
\eeqn
where 
\beqn
\eta_1: (V \oplus R) \times \widecheck {\rm Poly}{}_G^d(V\oplus R, W \oplus R) \to (V \oplus R) \times {\rm Poly}_G^d(V \oplus R, W \oplus R)
\eeqn
and 
\beqn
\eta_2: (V \oplus R) \times \wt{\rm Poly}{}_G^d(V\oplus R, W \oplus R) \to (V \oplus R) \times \widecheck {\rm Poly}{}_G^d(V \oplus R, W \oplus R)
\eeqn
are both the natural inclusions. 

\begin{lemma}\label{lemma325}
There exists $d_3$ (depending on $G$, $V$, $W$, and $R$) such that when $d \geq d_2$, $\eta_2$ is transverse to $\widecheck{\mf Z}^d(G, V\oplus R, W \oplus R)$ and 
\beqn
\eta_2^* \widecheck{\mf Z}^d(G, V\oplus R, W\oplus R) = \wt{\mf Z}^d(G, V \oplus R, W \oplus R).
\eeqn
\end{lemma}

\begin{lemma}\label{lemma326}
There exists $d_4$ (depending on $G$, $V$, $W$, and $R$) such that when $d \geq d_3$, $\eta_1$ is transverse to ${\mf Z}^d(G,V\oplus R, W \oplus R)$ and 
\beqn
\eta_1^* {\mf Z}^d(G, V\oplus R, W\oplus R) \equiv \widecheck{\mf Z}^d(G, V\oplus R, W \oplus R).
\eeqn
\end{lemma}

Assuming these two lemmas, we can prove Proposition \ref{prop_Whitney_stabilization}.

\begin{proof}[Proof of Proposition \ref{prop_Whitney_stabilization}]
Let $d$ be no smaller than $d_1$ of Proposition \ref{prop313}, $d_3$ of Lemma \ref{lemma325}, and $d_4$ of Lemma \ref{lemma326}. Lemma \ref{lemma325} and Lemma \ref{lemma326} imply that $\eta$ is transverse to ${\mf Z}^d = {\mf Z}^d(G, V\oplus R, W \oplus R)$ and 
\beqn
\eta^* {\mf Z}^d = (\eta_1 \circ \eta_2)^* {\mf Z}^d =  \eta_2^* \eta_1^* {\mf Z}^d \equiv \eta_2^* \widecheck {\mf Z}^d \equiv \wt {\mf Z}^d.
\eeqn
If $l\geq d_1 $ and $d \geq l$, then by Proposition \ref{prop313}, the composition $\eta \circ \sigma_{dl}$ is also transverse to ${\mf Z}^d$ and 
\beqn
(\eta \circ \sigma_{dl})^* {\mf Z}^d \equiv \sigma_{dl}^* \wt{\mf Z}^d = \wt{\mf Z}^l.
\eeqn
Therefore the transversality claim about $\eta$ is true for $d \geq d_1(G, V, W)$. As $\wt{\mf Z}^d$ has connected strata, it is equal to the refined pullback $\eta^\dagger {\mf Z}^d$. Moreover, when $d'>d$, the arrows of the diagram \eqref{eqn35} are all induced from refined pullbacks of Whitney stratifications. As refined pullback is functorial, the diagram commutes. 
\end{proof}

The proof strategy of Lemma \ref{lemma325} and Lemma \ref{lemma326} is similar to those of Proposition \ref{prop313} and Proposition \ref{prop_product_Z}. We need to construct a submersion in the opposite direction.

\begin{proof}[Proof of Lemma \ref{lemma325}]
Define 
\beqn
\begin{split}
\zeta_2: (V \oplus R) \times \widecheck{\rm Poly}{}_G^d(V\oplus R, W \oplus R) & \to (V\oplus R) \times \wt{\rm Poly}{}_G^d(V\oplus R, W \oplus R)\\
(x, y, P_0, {\rm Id}_R + Q_0) & \mapsto (x, y + Q_0(x), P_0, {\rm Id}_R).
\end{split}
\eeqn
It is straightforward to check that $\zeta_2$ is a surjective submersion, preserves the action stratification, and $\zeta_2^{-1}(\wt{Z}^d) = \widecheck{Z}^d$. Hence by Proposition \ref{prop_pullback_nice}, one has $\zeta_2^* \wt {\mf Z}^d\equiv \widecheck {\mf Z}^d$. Moreover, it is easy to see that $\zeta_2$ has connected fibers. So $\zeta_2^* \wt {\mf Z}^d = \widecheck {\mf Z}^d$. Moreover, since $\zeta_2 \circ \eta_2$ is the identity, one has 
\beqn
\eta_2^* \widecheck {\mf Z}^d = \eta_2^* (\zeta_2^* \wt{\mf Z}^d) = (\zeta_2 \circ \eta_2)^* \wt {\mf Z}^d = \wt{\mf Z}^d.\qedhere
\eeqn
\end{proof}

Now we prove Lemma \ref{lemma326}. The construction of a left inverse of $\eta_1$ is more involved. First we need the following result.

\begin{prop}\label{GWSchwarz}\cite[Proposition 6.8]{GWSchwarz_1980} Let $V_1$ and $V_2$ be finite-dimensional representations of $G$ and let $U_1 \subset V_1$ be a $G$-invariant open subset. Then the space of $G$-invariant holomorphic maps ${\rm Map}_G^h(U_1, V_2)$ from $U_1$ to $V_2$ is a module over the space of $G$-invariant holomorphic functions ${\rm Map}_G^h(U_1, {\mb C})$ generated by ${\rm Poly}_G(V_1, V_2)$. 
\end{prop}

\begin{proof}[Proof of Lemma \ref{lemma326}]
Define a kind of ``left inverse'' to $\eta_1$. For any $d \geq 0$ and any $Q \in {\rm Poly}_G^d (V \oplus R, R)$, consider the equation 
\beqn
y + Q(x,y) = 0.
\eeqn
By the implicit function theorem, there exists a $G$-invariant open neighborhood
\beqn
O^d \subset V \times {\rm Poly}_G^d(V \oplus R, R)
\eeqn
of $V \times {\rm Poly}_G^d(V, R)$, where the second factor is viewed as elements in ${\rm Poly}_G^d(V \oplus R, R)$ not depending on $R$, such that when $(x, Q) \in O^d$, we can always solve $y$ in terms of $x$. More precisely, there is a $G$-equivariant complex analytic map 
\beqn
{\mc F}: O^d \to R
\eeqn
such that 
\beqn
{\mc F}(x, Q) + Q( x, {\mc F}(x, Q)) \equiv 0\ {\rm and}\ {\mc F}(x, Q_0) = - Q_0(x),\ \forall Q_0 \in {\rm Poly}_G^d(V, R).
\eeqn
Notice that ${\mc F}: O^d \to R$ is no longer a polynomial map. However, by Proposition \ref{GWSchwarz}, we can write it as a combination 
\beqn
{\mc F}(x, Q) = \sum_{j=1}^n f_j(x, Q) P_j(x, Q)
\eeqn
where $f_j: O^d \to {\mb C}$ is a $G$-invariant holomorphic function and $P_j$ is $G$-equivariant polynomial map from $V \oplus {\rm Poly}_G^d(V\oplus R, R)$ to $R$. 

Now set
\beqn
U^d:= \Big\{ (x, y, P, {\rm Id}_R + Q) \in (V\oplus R) \times {\rm Poly}_G^d(V \oplus R, W \oplus R)\ |\ {\rm Id}_R + Q \in O^d \Big\}
\eeqn
which is a $G$-invariant open subset of $(V\oplus R) \times {\rm Poly}_G^d(V \oplus R, W \oplus R)$. Define
\beqn
\begin{split}
\widecheck \zeta_1: U^d & \to (V \oplus R) \times \widecheck {\rm Poly}{}_G (V \oplus R, W \oplus R)\\
(x, y, P, {\rm Id}_R + Q) &\ \mapsto \left( x, y, \widecheck P, {\rm Id}_R + \widecheck Q \right).
\end{split}
\eeqn
where
\begin{align*}
&\ \widecheck P(\cdot) = P  \left( \cdot, \sum_{j=1}^n f_j(x, Q) P_j(\cdot, Q) \right),\ &\ \widecheck Q(\cdot) =  Q \left( \cdot, \sum_{j=1}^n f_j(x, Q) P_j( \cdot, Q) \right).
\end{align*}
Indeed, $\widecheck P$ and $\widecheck Q$ are polynomial maps in variable $x \in V$ because $P$, $Q$ are polynomials in $x, y$ and $P_j$ are polynomials in $x$. As $f_j$ is $G$-invariant and $P_j$ is $G$-equivariant, one can also check $\widecheck P$ and $\widecheck Q$ are equivariant. Moreover, there exists $\widecheck d \geq d$ (depending on $d$, $G$, $V$, $W$) such that $\widecheck P$ and $\widecheck Q$ have degrees at most $\widecheck d$. Hence we obtain an analytic map
\beqn
\widecheck \zeta_1: U^d \to (V\oplus R) \times \widecheck{\rm Poly}{}_G^{\widecheck d}(V\oplus R, W\oplus R).
\eeqn

\noindent {\it Claim.} Let 
\beqn
\widecheck \phi: (V\oplus R) \times \widecheck {\rm Poly}{}_G^d (V\oplus R, W\oplus R) \to  (V\oplus R) \times \widecheck {\rm Poly}{}_G^{\widecheck d}(V\oplus R, W\oplus R) 
\eeqn
be the natural inclusion. Then $\widecheck \zeta_1 \circ \eta_1 = \widecheck \phi$.

\vspace{0.1cm}

\noindent {\it Proof of the claim.} Indeed, if we insert $P = P_0$ and $Q = Q_0$ into the above formula of $\widecheck P$ and $\widecheck Q$, as $P_0$ and $Q_0$ only depend on $x$, it follows that $\widecheck P = P_0$ and $\widecheck Q = Q_0$.\ \hfill {\it End of the proof of the claim.}

\vspace{0.1cm}

On the other hand, similar to the proof of Proposition \ref{prop313}, there exists a $G$-invariant left inverse obtained by expressing high degree monomials in terms of low degree ones
\beqn
\widecheck \psi: (V\oplus R) \times \widecheck {\rm Poly}{}_G^{\widecheck d}(V\oplus R, W\oplus R) \to (V\oplus R) \times \widecheck {\rm Poly}{}_G^d (V\oplus R, W\oplus R)
\eeqn
to $\widecheck\phi$. Consider the composition
\beqn
\zeta_1:= \widecheck \psi \circ \widecheck \zeta_1: U^d \to (V\oplus R) \times \widecheck {\rm Poly}{}_G^d(V\oplus R, W\oplus R) 
\eeqn

Now we check the following items.
\begin{enumerate}
    \item $\zeta_1$ preserves the action  stratification.

    \item By the above claim, $\zeta_1$ is a left inverse to $\eta_1$.
    
    \item $\zeta_1^{-1}(\widecheck Z^d)= Z^d\cap U^d$.\footnote{Notice that $\zeta_1$ does not preserve the evaluation but only preserves the zero set.} Indeed, if $\widecheck P(x) = 0$ and $y + \widecheck Q(x) = 0$, then 
    \beqn
    P(x, {\mc F}(x, Q)) = 0,\ y + Q(x, {\mc F}(x, Q)) = 0,
    \eeqn
    then by the definition of ${\mc F}(x, Q)$, it follows that $y = {\mc F}(x, Q)$. Hence $P(x, y) = 0$ and $y + Q(x, y) = 0$. Hence $\zeta_1^{-1}(\widecheck Z^d) \subset Z^d$. On the other hand, if $(x, y, P, {\rm Id}_R + Q)\in Z^d \cap U^d$, then $P(x, y) = 0$ and $y + Q(x,y) = 0$ and $y = {\mc F}(x, Q)$ is the only solution. Hence $\zeta_1( Z^d \cap U^d) = \widecheck Z^d$. 
    \end{enumerate}
Then by Proposition \ref{prop_pullback_nice}, one has 
\beqn
\zeta_1^* \widecheck {\mf Z}^d \equiv {\mf Z}^d|_{Z^d \cap U^d}.
\eeqn
On the other hand, as $\zeta_1 \circ \eta_1$ is the identity map, $\eta_1$ is transverse to ${\mf Z}^d$. Hence 
\beqn
\eta_1^* {\mf Z}^d =  \eta_1^* ({\mf Z}^d|_{Z^d\cap U^d} ) \equiv \eta_1^* (\zeta_1^* \widecheck {\mf Z}^d) =  (\zeta_1 \circ \eta_1)^* \widecheck {\mf Z}^d =  \widecheck {\mf Z}^d.  \qedhere
\eeqn
\end{proof}

\subsection{Proof of Proposition \ref{prop314}}\label{subsection37}

From the previous proofs of various properties of the canonical Whitney stratification, one can see the pattern is to consider the inclusion map of a space of special polynomial maps to a larger space of more general polynomial maps and then create a left inverse, which is typically an analytic submersion. In those situations, the comparisons of Whitney stratifications appeal to Proposition \ref{prop_pullback_nice}. However, the construction for proving Proposition \ref{prop314} here is more complicated. The map we constructed is no long an analytic submersion nor a left inverse. As a result, one needs to appeal to Lemma \ref{lemmab5} which is more technical than Proposition \ref{prop_pullback_nice}. 

\begin{defn}
An {\it averaging map} is an $H$-equivariant smooth map of the form 
\beqn
\kappa:  M^d(H, V, W)_H^+ \to M^d(G, V, W)_H^+,\ \kappa(v, Q) = (v, \uds\kappa (v, Q)).
\eeqn
such that $\ev(\kappa (v, Q)) = Q(v) \in W$.
\end{defn}

Roughly speaking, an averaging map provides a way of producing a (more symmetric) $G$-equivariant map from a (less symmetric) $H$-equivariant map which preserves the evaluation map.

The proof of Proposition \ref{prop314} is based on the following three lemmas.

\begin{lemma}\label{lemma329}
When $d$ is sufficiently large, there exists an averaging map.
\end{lemma}

\begin{proof}
See Subsubection \ref{subsubsection381}.
\end{proof}

The next two lemmas are to verify part of the assumptions of Lemma \ref{lemmab5}.

\begin{lemma}\label{lemma330}
For any averaging map $\kappa$, the following are true.
\begin{enumerate}
\item $\kappa \circ \mu: M^d(G, V, W)_H^+ \to M^d(G, V, W)_H^+$ is transverse to ${\mf Z}^d(G, V, W)$ 
\beqn
(\kappa \circ \mu)^* {\mf Z}^d(G, V, W) = {\mf Z}^d(G, V, W)|_{Z^d(G, V, W)_H^+}.
\eeqn

\item $\mu \circ \kappa:   M^d(H, V, W)_H^+ \to  M^d(H, V, W)_H^+$ is transverse to $  {\mf Z}^d(H, V, W)$ and 
\beqn
(\mu \circ \kappa)^*  {\mf Z}^d(H, V, W) =  {\mf Z}^d(H, V, W)|_{ Z^d(H, V, W)_H^+}.
\eeqn
\end{enumerate}
\end{lemma}

\begin{proof}
See Section \ref{subsubsection373}.
\end{proof}

\begin{lemma}\label{lemma:mu-kappa}
Suppose there exists an averaging map $\kappa$. Then $\mu|_{M^d(G, V, W)_H^+}$ is transverse to ${\mf Z}^d(H, V, W)$ and $\kappa$ is transverse to ${\mf Z}^d(G, V, W)|_{Z^d(G, V, W)_H^+}$.
\end{lemma}

\begin{proof}
By Lemma \ref{lemma330} (cf. Lemma \ref{lemmab4}), $\mu$ is transverse to ${\mf Z}^d(H, V, W)$ at points on ${\rm Im} (\kappa)$ where $\kappa$ is any averaging map. To show that $\mu$ is transverse to ${\mf Z}^d(H, V, W)$, it suffices to show that any $(v_0, P_0) \in Z^d(G, V, W)_H^+$ lies in the image of some averaging map. Let $\kappa$ be an averaging map. Consider the following vector field on $M^d(G, V, W)_H^+$
\beqn
w_{P_0}(v, P) = \Big( 0, P_0 - \uds \kappa (\mu (v, P_0)) \Big).
\eeqn
Its flow is the family of diffeomorphisms
\beqn
F_{P_0, t}(v, P) = \Big( v, P + t (P_0 - \uds \kappa (\mu (v, P_0))) \Big).
\eeqn
It is easy to see that $F_{P_0, t}$ preserves the evaluation map and is $H$-equivariant. Hence the composition 
\beqn
F_{P_0, t} \circ \kappa: M^d(H, V, W)_H^+ \to M^d(G, V, W)_H^+
\eeqn
is an averaging map. Moreover, 
\beqn
F_{P_0, 1} \circ \kappa (\mu ( v_0, P_0)) = \Big( v_0, \uds\kappa (\mu (v_0, P_0)) + (P_0 - \uds\kappa (\mu (v_0, P_0))) \Big) = (v, P_0).
\eeqn
Hence $(v_0, P_0)$ is on the image of some averaging map. Therefore, $\mu$ is transverse to ${\mf Z}^d(H, V, W)$.

The case for $\kappa$ is similar. Lemma \ref{lemma330} (cf. Lemma \ref{lemmab4}) implies that $\kappa$ is transverse to ${\mf Z}^d(G, V, W)$ at points on ${\rm Im}(\mu)$. We would like to show that $\kappa$ is transverse to ${\mf Z}^d(G, V, W)$ at any given $(v_0, Q_0) \in Z^d(H, V, W)_H^+$. Define the diffeomorphism
\beqn
g_{Q_0} (v, Q) = \Big( v, Q +  Q_0 - \uds\mu ( \kappa (v, Q_0 )) \Big).
\eeqn
Then $g_{Q_0}$ preserves the evaluation map. Then the composition $\kappa \circ g_{Q_0}$ is also an averaging map and hence $\kappa \circ g_{Q_0}$ is transverse at points on ${\rm Im}(\mu)$. In particular, it is transverse at the point $\mu(\kappa (v_0, Q_0))$. But this is equivalent to that $\kappa$ is transverse to ${\mf Z}^d(G, V, W)$ at the point 
\beqn
g_{Q_0}(\mu(\kappa (v_0, Q_0))) =  \Big( v_0, \uds\mu (\kappa (v_0, Q_0)) + Q_0 - \uds\mu (\kappa (v_0, Q_0)) \Big) = (v_0, Q_0).
\eeqn
Therefore $\kappa$ is transverse to ${\mf Z}^d(G, V, W)$.
\end{proof}

\begin{proof}[Proof of Proposition \ref{prop314}]
Notice that on the open subset
\beqn
M^d(G, V, W)_H^+ \subset M^d(G, V, W)
\eeqn
the action stratification coming from the $G$-action coincides with the action stratification coming from the $H$-action, as all the strata on both sides are indexed by $V$-essential subgroups of $H$. Denote by
\beqn
\mc{WS}^\infty( Z^d(G, V, W)_H^+; H) \subset \mc{WS}^\infty( Z^d(G, V, W)_H^+)
\eeqn
the set of smooth Whitney stratifications which respect the action stratification. Similarly, one has the corresponding set $\mc{WS}^\infty( Z^d(H,V, W)_H^+; H)$ of Whitney stratifications. Then all the Whitney stratifications contained in the statement of Proposition \ref{prop314} belong to the corresponding subsets. Then by Lemma \ref{lemma330}, \ref{lemma:mu-kappa} and Lemma \ref{lemmab5}, it follows that 
\beqn
\mu^*  {\mf Z}^d(H, V, W)|_{Z^d(G, V, W)_H^+} \equiv {\mf Z}^d(G, V, W)|_{Z^d(G, V, W)_H^+}.
\eeqn
Moreover, as both ${\mf Z}^d(H, V, W)$ and ${\mf Z}^d(G, V, W)|_{Z^d(G, V, W)_H^+}$ have connected strata, one has 
\beqn
\mu^* {\mf Z}^d(H, V, W)|_{Z^d(G, V, W)_H^+} = {\mf Z}^d(G, V, W)|_{Z^d(G, V, W)_H^+}. 
\eeqn
Then consider the induced map 
\beqn
\mu_*: {\mf Z}^d(G, V, W)|_{Z^d(G, V, W)_H^+} \to  {\mf Z}^d(H, V, W)|_{Z^d(G, V, W)_H^+}.
\eeqn
Using an average map $\kappa$ provided by Lemma \ref{lemma329} and its property provided by Lemma \ref{lemma330}, we see $\kappa_*$ is an inverse of $\mu_*$. Lastly, since all arrows of the diagram \eqref{eqn33} are induced by pullbacks of Whitney stratifications and pullback is functorial, the diagram \eqref{eqn33} commutes. 
\end{proof}

\subsubsection{Proof of Lemma \ref{lemma329}}\label{subsubsection381}

\begin{lemma}
Denote $d_0 = |G|$. There exists an $H$-invariant smooth map 
\beqn
V_H^+ \to {\rm Poly}_H^{d_0} (V, {\mb C})
\eeqn
denoted by $v \mapsto L_v$, such that
\beq\label{eqn38}
L_v(\tau v) = \left\{ \begin{array}{cc} 1,\ &\ {\rm if}\ \tau \in H,\\
                                     0,\ &\ {\rm if}\ \tau \notin H.\end{array}\right.
\eeq
\end{lemma}

\begin{proof}
We first construct a map $v \mapsto L_v'' \in {\rm Poly}^{d_0}(V, {\mb C})$ satisfying \eqref{eqn38} without the $H$-invariance condition. For each $v \in V_H^+$, because $G_v$ is contained in $H$, $\tau_1 v \neq \tau_2 v$ if $\tau_1 \in H$ and $\tau_2 \notin H$. Consider a 1-dimensional subspace $S_v \subset V$ and a complement $S_v^\bot$. Let $p_v: V \to S_v$ be the projection induced from the splitting $V \cong S_v \oplus S_v^\bot$. One can choose $S_v$ and $S_v^\bot$ such that for a small open neighborhood $O(v)\subset V_H^+$ of $v$, 
\beqn
p_v(\tau_1 v') \neq p_v(\tau_2 v')\ {\rm if\ } \tau_1 \in H,\ \tau_2 \notin H,\ v' \in O(v).
\eeqn
Then by the Lagrange interpolation method, there exists a complex polynomial $f_{v, v'}: S_v \to {\mb C}$ (given by the exact formula of Lagrange interpolation) of degree at most $d_0 = |G|$ such that
\beqn
f_{v,v'} ( p_v (\tau v')) = \left\{ \begin{array}{cc} 1, & \tau \in H,\\
0, & \tau \notin H. \end{array}\right.
\eeqn
Then define the polynomial
\beqn
L_{v, v'}:= f_{v, v'} \circ p_v: V \to {\mb C},\ \forall v' \in O(v).
\eeqn
Then $L_{v, v'}(\tau v') = 1$ if $\tau \in H$ and $L_{v, v'}(\tau v') = 0$ otherwise. 

For each $v \in V_H^+$, we can find as above an open neighborhood $O(v)$ and a function $v' \mapsto L_{v, v'} \in {\rm Poly}^{d_0}(V, {\mb C})$. As $V_H^+$ is Lindel\"of, one can cover it by a countable subcover $O(v_i)$, $i = 1, \ldots$. Find a subordinate partition of unity $\rho_i$ and define
\beqn
v\mapsto L_v'':= \sum_{i=1}^\infty \rho_i(v) L_{v_i, v}.
\eeqn
The local finiteness of the partition of unity implies that this is a smooth map from $V_H^+$ to ${\rm Poly}^{d_0}(V, {\mb C})$. Moreover, 
\beqn
L_v''(\tau v')|_{v' = v} = \sum_{i=1}^\infty \rho_i(v)  L_{v_i, v}(\tau v')|_{v' = v}  = \left\{ \begin{array}{cc} 1,\ & \tau \in H,\\ 0, \ & \tau \notin H \end{array}\right.
\eeqn
Then define 
\beqn
L_v':= \frac{1}{|H|} \sum_{\tau \in H} L_v'' \circ \tau \in {\rm Poly}_H^{d_0}(V, {\mb C})
\eeqn
which still satisfies \eqref{eqn38}, and define 
\beqn
L_v:= \frac{1}{|H|} \sum_{\tau \in H} L_{\tau v}'\in {\rm Poly}_H^{d_0}(V, {\mb C})
\eeqn
so the map $v \mapsto L_v$ is $H$-invariant satisfying \eqref{eqn38}.
\end{proof}

Now we define a map
\beqn
\begin{split}
\uds\kappa':  M^d(V, H, W)_H^+ & \to {\rm Poly}_G^{d+d_0}(V, W)\\
(v, Q) & \mapsto \frac{1}{|H|} \sum_{\tau \in G} \tau^{-1} \circ \big( L_v Q \big) \circ \tau .
\end{split}
\eeqn
Define
\beqn
\kappa' (v, Q) = (v, \uds\kappa '(v, Q)) \in M^{d + d_0}(G, V, W)_H^+.
\eeqn
On the other hand, when $d \geq d_1$, Lemma \ref{lemma323} provides a $G$-invariant holomorphic map
\beqn
\rho: V \times {\rm Poly}_G^{d+d_0} (V, W) \to V \times {\rm Poly}_G^d ( V, W)
\eeqn
satisfying the two conditions of Lemma \ref{lemma323}. Then define
\beqn
\kappa: = \rho \circ \kappa': M^d(H, V, W)_H^+ \to M^d(G, V, W)_H^+.
\eeqn
As both $\rho$ and $\kappa'$ preserve the evaluation map. $\kappa$ is an averaging map. This concludes the proof of Lemma \ref{lemma329}.

\subsubsection{More preparations for Lemma \ref{lemma330}}

\begin{lemma}\label{lemma333}
Let $(G, F, E)$ be a family NC triple over $X$ (Definition \ref{defn246}). Let $\phi: M^d(G, F, E) \to M^d(G, F, E)$ be a fiberwise diffeomorphism which preserves the action stratification and which preserves the subbundle $Z^d(G, F, E)$ set-wise. Then $\phi^* {\mf Z}^d(G, F, E) = {\mf Z}^d(G, F, E)$.
\end{lemma}

\begin{proof}
Notice that the fiberwise $G$-action gives a family version of the action stratification ${\mf M}$ on $M^d(G, F, E)$. 
We first observe that ${\mf Z}^d(G, F, E)$ is a minimal Whitney stratification which respects ${\mf M}$ (Definition \ref{defn_nice_Whitney_stratification}). Indeed, as the bundles are locally trivial, the restriction of ${\mf Z}^d(G, F, E)$ to the preimage to any small open subset $U \subset X$ is minimal among all Whitney stratifications which respects ${\mf M}$. As the minimality is defined via the comparison of the dimension filtration, which is a local condition, the minimality of ${\mf Z}^d(G, F, E)$ follows. Then by Proposition \ref{propb14}, it follows that $\phi^* {\mf Z}^d(G, F, E) \equiv {\mf Z}^d(G, F, E)$. Moreover, when $X$ is connected, ${\mf Z}^d(G, F, E)$ has connected strata; as $\kappa$ is a diffeomorphism, the pullback $\phi^* {\mf Z}^d(G, F, E)$ also has connected strata over each connected component of $X$. Hence one has the equality as claimed.
\end{proof}

Recall that for any lift ${\mf p}: F \to \wh{\rm Poly}{}_G^d(F, E)$ 
of an NC bundle map $S \in C_G^{\rm NC}(F, E)$ (see Definition \ref{defn244}), its graph is the subset 
\beqn
{\rm graph}({\mf p})= \{(v, {\mf p}(v))\ |\ v \in F \} \subset M^d(G, F, E).
\eeqn
Parker \cite{BParker_integer} observed the following. Two maps ${\mf p}_1, {\mf p}_2: F \to \wh {\rm Poly}{}_G^d(F, E)$ are lifts of the same bundle map $S\in C_G^{\rm NC}(F, E)$ if and only if
\beqn
{\rm graph}({\mf p}_1 - {\mf p}_2) \subset  Z^d(G, F, E).
\eeqn
We reprove a result of Parker \cite{BParker_integer}.

\begin{lemma}(cf. \cite[Lemma 4.10]{BParker_integer}) \label{lemma334}
Let $(G, F, E)$ be as above. Fix an integer $d > 0$. Let $F' \subset F$ be an open subset of the total space $F$ and let ${\mf p}_1, {\mf p}_2: F' \to \wh{\rm Poly}{}_G^d ( F, E)$ be two smooth bundle maps such that
\beqn
{\rm graph} (  {\mf p}_1 - {\mf p}_2) \subset Z^d(G, F, E).
\eeqn
Then ${\rm graph}({\mf p}_1)$ is  transverse ${\mf Z}^d(G, F, E)$ if and only if ${\rm graph}({\mf p}_2)$ is transverse to ${\mf Z}^d(G, F, E)$. Moreover, in this case, the two Whitney stratifications on $S_{{\mf p}_1}^{-1}(0) = S_{{\mf p}_2}^{-1}(0)$ pulled back by ${\rm graph}( {\mf p}_1)$ and by ${\rm graph}({\mf p}_2)$ coincide. 
\end{lemma}

\begin{proof}
Suppose the graph of ${\mf p}_1$ is transverse to ${\mf Z}^d(G, F, E)$. For any $(v_0, {\mf p}_2 (v_0)) \in Z^d(G, F, E)$, we would like to show that the graph of $ {\mf p}_2$ is transverse to $ Z^d(G, F, E)$ at this point $(v_0, {\mf p}_2(v_0))$. Choose a compactly supported cut-off function 
\beqn
\rho_0: F' \to [0, 1]
\eeqn
which is identically $1$ near $v_0$. Consider the smooth vertical vector field on the total space of the vector bundle $F \oplus \wh {\rm Poly}{}_G^d(F, E)$ defined by
\beq\label{eqn39}
f (v, P) = \big( 0,  \rho_0(v) (  {\mf p}_1 (v) - {\mf p}_2 (v)) \big).
\eeq
Accordingly, the flow of $f$ is the 1-parameter family of fiber-preserving diffeomorphisms
\beqn
\phi_t(v, P) = \big( v, P + t \rho_0(v)(  {\mf p}_1(v) -  {\mf p}_2(v)) \big)
\eeqn
of $M^d(G, F, E)$ which exists for all time $t$. It is also easy to see that $\phi_t$ preserves the action stratification on $M^d(G, F, E)$ and the set $Z^d(G, F, E)$. Hence by Lemma \ref{lemma333}, $\phi_t$ pulls back the canonical Whitney stratification on $Z^d(G, F, E)$ to itself. Moreover, $\phi_1$ maps a neighborhood of $(v_0,   {\mf p}_2 (v_0))$ in ${\rm graph}({\mf p}_2 )$ to a neighborhood of $(v_0, {\mf p}_1 (v_0))$ in ${\rm graph}(  {\mf p}_1)$. Hence the graph of $ {\mf p}_2$ is also transverse to ${\mf Z}^d$ at $v_0$.

To show that the two pullback Whitney stratifications agree, we just need to check locally around a point $v_0 \in {\rm graph}( {\mf p}_1)^{-1}( Z^d) = {\rm graph}( {\mf p}_2)^{-1} (  Z^d) \subset F'$. Let $\phi_1$ be the time-1 map of the flow of the vector field \eqref{eqn39}. Then in a neighborhood of $v_0\in F'$, as maps, one has
\beqn
{\rm graph}(  {\mf p}_1)  = \phi_1 \circ {\rm graph}(  {\mf p}_2). 
\eeqn
On the other hand, as $\phi_1$ preserves the canonical Whitney stratification ${\mf Z}^d$, it follows that near $v_0$, ${\rm graph}( {\mf p}_1)^* {\mf Z}^d = {\rm graph}(  {\mf p}_2)^*\phi_1^*   {\mf Z}^d = {\rm graph}(  {\mf p}_2)^* {\mf Z}^d$.
\end{proof}





\subsubsection{Proof of Lemma \ref{lemma330}}\label{subsubsection373}

Consider a family NC triple $(G, F, E)$ over the base $X = \wh{\rm Poly}{}_G^d(V, W)$ with 
\begin{align*}
&\ F:= X \times V,\ &\ E:= X \times W.
\end{align*}
Consider the subbundle 
\beqn
F_H^+:= X\times V_H^+  \subset F.
\eeqn
Then $\kappa \circ \mu$ can be viewed as a bundle map
\beqn
\begin{split}
 {\mf p}: F_H^+ &\ \to \wh{\rm Poly}{}_G^d (F, E) \cong X \times {\rm Poly}_G^d(V, W)\\
(P, v) & \mapsto \Big(P, \uds\kappa (\mu (v, P)) \Big)
\end{split}
\eeqn
On the other hand, there is another natural bundle map 
\beqn
\begin{split}
 {\mf p}_0: F_H^+ &\ \to \wh{\rm Poly}{}_G^d (F, E),\\
(P, v) & \ \mapsto (P, P).
\end{split}
\eeqn
Then the construction implies that 
\beqn
{\rm graph}(  {\mf p} -  {\mf p}_0 ) \subset  Z^d(G, F, E).
\eeqn
As ${\rm graph}(  {\mf p}_0)$ is transverse to the canonical Whitney stratification on $  Z^d (G, F, E)$, by Lemma \ref{lemma334} the graph of $ {\mf p}$ is also transverse to it. Equivalently, it means $\kappa \circ \mu$ is transverse to ${\mf Z}^d (G,V, W)$. As ${\rm graph}(  {\mf p}_0)$ pulls back the canonical Whitney stratification on $  Z^d(G, F, E)$ to itself, it is the same for ${\rm graph}(  {\mf p}_1)$. Equivalently, this means $(\kappa \circ \mu)^* {\mf Z}^d(G, V, W) = {\mf Z}^d(G, V, W)|_{Z^d(G, V, W)_H^+}$.

For the case of $\mu \circ \kappa$, consider another family NC triple $(H, \tilde F, \tilde E)$ over the base $Y = \wh {\rm Poly}{}_H^d(V, W)$ with
\begin{align*}
&\ \tilde F := Y \times V,\ &\  \tilde E:= Y \times W.
\end{align*}
Then there is a subbundle
\beqn
\tilde F_H^+:= Y \times V_H^+ \subset F.
\eeqn
Then $\mu \circ \kappa$ can be viewed as a bundle map
\beqn
\begin{split}
\tilde {\mf q}: \tilde F_H^+ &\  \to \wh{\rm Poly}{}_H^d( \tilde F, \tilde E)\\
(Q, v) &\ \mapsto \Big( Q, \uds\mu (\kappa (v, Q)) \Big).
\end{split}
\eeqn
There is another bundle map 
\beqn
\tilde {\mf q}_0(Q, v) = (Q, Q).
\eeqn
Then ${\rm graph}( \tilde {\mf q} - \tilde {\mf q}_0) \subset   Z^d(H, \tilde F, \tilde E)$. Then by Lemma \ref{lemma334}, ${\rm graph}( \tilde {\mf q})$ is transverse to ${\mf Z}^d(H, \tilde F, \tilde E)$ and pulls it back to the canonical Whitney stratification on $\tilde F_H^+$. It is equivalent to that $\mu \circ \kappa$ is transverse to $ {\mf Z}^d(H, V, W)$ and pulls it back to ${\mf Z}^d(H, V, W)$ over $Z^d(H, V, W)_H^+$.

\section{FOP Transverse Sections}\label{section4}

In this section we define the space of FOP transverse sections and prove Theorem \ref{thm11}. Roughly speaking, a normally complex section is FOP transverse if it induces a certain map that is transverse to the variety $Z$ equipped with the canonical Whitney stratification specified in the previous section. The properties of the canonical Whitney stratifications allow us to deduce properties of FOP transverse sections. Within this section, let $({\mc U}, {\mc E})$ be an NC pair (Definition \ref{defn_normal_complex}), i.e., a normally complex vector bundle ${\mc E}$ over a normally complex orbifold ${\mc U}$.

\subsection{FOP transverse bundle maps}

We first discuss the FOP transversality on the local model. Let $(G, F, E)$ be a family NC triple over a smooth manifold $X$. Recall that one has the space
\beqn
C_G^{\rm NC} (F, E) \subset C_G^\infty(F, E)
\eeqn
of normally complex bundle maps (Definition \ref{defn246}).  

\begin{defn}
Consider an NC bundle map $S \in C_G^{\rm NC} (F, E)$ and a lift ${\mf p}: F \to \wh{\rm Poly}{}_G^d(F, E)$ (Definition \ref{defn244}). 

\begin{enumerate}

\item The lift ${\mf p}$ is said to be {\it FOP transverse} at $(x, v) \in F$ if $d \geq d_1(G, V, W)$ of Proposition \ref{prop313} where $V$ and $W$ are fibers of $F$ and $E$ respectively, and if the graph ${\rm graph}({\mf p}) \subset M^d(G, F, E)$ is transverse to the canonical Whitney stratification of $Z^d(G, F, E)$ at $(x, v, {\mf p}(x, v))$.

\item $S$ is said to be {\it FOP transverse} at $(x, v) \in F$ if there exists a lift ${\mf p}$ which is FOP transverse at $(x, v)$. 
\end{enumerate}
\end{defn}

We show that the transversality condition is independent of the  lift.

\begin{lemma}
Let ${\mf p}: F \to \wh{\rm Poly}{}_G^d (F, E)$ and ${\mf p}': F \to \wh{\rm Poly}{}_G^{d'}(F, E)$ be two lifts of $S: F \to E$, then ${\mf p}$ is FOP transverse at $(x, v)$ if and only if ${\mf p}'$ is FOP transverse at $(x, v)$.
\end{lemma}

\begin{proof}
If $d \leq d'$, we may identify ${\mf p}$ with a map ${\mf q}: F \to \wh{\rm Poly}{}_G^{d'} (F, E)$. Then by Proposition \ref{prop313}, ${\mf p}$ is transverse to ${\mf Z}^d (G, F, E)$ at $(x, v)$ if and only if ${\mf q}$ is transverse to ${\mf Z}^{d'}(G, F, E)$ at $(x, v)$. Hence we may assume $d' = d$. Moreover, as ${\mf p}$ and ${\mf p}'$ induce the same section $S$, it follows that ${\rm graph}({\mf p} - {\mf p}') \subset  Z^d(G, F, E)$. Then by Lemma \ref{lemma334}, the two transversality conditions are equivalent.
\end{proof}

\begin{lemma}\label{lemma_openness}
The set of points $(x, v) \in F$ at which an NC bundle map $S \in C_G^{\rm NC} (F, E)$ is FOP transverse is open in $F$.
\end{lemma}

\begin{proof}
Choose a lift ${\mf p}: F \to \wh{\rm Poly}{}_G^d(F, E)$ for  a sufficiently large $d$. Consider the intersection between ${\rm graph}({\mf p})$ and $Z^d(G, F, E)$. Trotman \cite{Trotman_1978} proved that the set of transverse points to a Whitney stratified subset is open. Therefore, the set of points on ${\rm graph}({\mf p})$ where the intersection is transverse is open. As ${\rm graph}({\mf p})$ is homeomorphic to $F$, the lemma follows. 
\end{proof}

Another crucial feature is related to the change of group. Recall that if $H \subset_F G$ is an essential subgroup, then there is the basic decomposition
\beqn
F = F_H \oplus \check F_H
\eeqn
and the total space $F$ is canonically identified with the total space of $\pi_{F_H}^* \check F_H \to F_H$. Then consider the canonical map defined in the proof of Lemma \ref{lemma244}
\beq\label{mapmu}
\mu: M^d(G, F, E) \to M^d(H, F, E)
\eeq
which induces a natural map between spaces of NC bundle maps, which, by abuse of notation, was denoted by
\beqn
\mu: C_G^{\rm NC} (F, E) \to C_H^{\rm NC} (\pi_{F_H}^* \check F_H, \pi_{F_H}^* E).
\eeqn

\begin{lemma}\label{lemma44}
Suppose $S$ is FOP transverse at a point $(x, v)$ whose stabilizer is contained in $H$, then $\mu(S)$ is FOP transverse at $(x, v) \in F$. 
\end{lemma}

\begin{proof}
Choose a lift ${\mf p}: F \to \wh{\rm Poly}{}_G^d(F, E)$ of $S$. Using the map \eqref{mapmu} one obtains an $H$-invariant bundle map ${\mf q}: F \to \wh{\rm Poly}{}_H^d(F, E)$ so that 
\beqn
\mu( {\rm graph}({\mf p})) = {\rm graph}({\mf q}).
\eeqn
Since the stabilizer of $(x, v)$ is contained in $H$, using the family version (which is a trivial extension) of Proposition \ref{prop314}, one can see that ${\rm graph}({\mf q})$ is transverse to ${\mf Z}^d(H, F, E)$ at $(x, v, {\mf q}(x, v))$. Therefore, as an $H$-equivariant NC bundle map, $S$ is FOP transverse at $(x, v)$. 

On the other hand, under the natural identification 
\beqn
M^d(H, F, E) \cong M^d(H, \pi_{F_H}^* \check F_H, \pi_{F_H}^* E)
\eeqn
${\mf q}$ becomes a lift of the (same) $H$-equivariant bundle map $S: \pi_{F_H}^* \check F_H \to \pi_{F_H}^* E$ and the Whitney stratification ${\mf Z}^d(H, F, E)$ is identical to ${\mf Z}^d(H, \pi_{F_H}^* \check F_H, \pi_{F_H}^* E)$. Therefore, ${\rm graph}({\mf q})$ is transverse to the latter at the corresponding point and hence $\mu (S)$ is FOP transverse at $(x, v)$.
\end{proof}

Next we examine more carefully the FOP transversality condition aiming at proving the existence and density of transverse elements. Let $G, F, E$ be as before. Let ${\rm proj}_2: M^d(G, F, E) \to \wh{\rm Poly}{}_G^d(F, E)$ be the vector bundle projection.

\begin{lemma}\label{lemma45}
Suppose $F_G = 0$ (so $F$ is a complex vector bundle). Then the graph of ${\mf p}: F \to \wh{\rm Poly}{}_G^d(F, E)$ is transverse to a stratum $Z_\alpha\in {\mf Z}^d (G, F, E)$ at a point $(x, 0) \in X \subset F$ if and only ${\mf p}|_X: X \to \wh {\rm Poly}{}_G^d(F, E)$ is transverse to the smooth map ${\rm proj}_2: Z_\alpha \to \wh{\rm Poly}{}_G^d(F, E)$ at $x$.
\end{lemma}

\begin{proof}
We may assume that $F$ and $E$ are trivial with fibers $V$ and $W$ respectively. Then ${\mf p}$ is equivalent to a $G$-invariant map 
\beqn
{\mf p}: X \times V \to \wh{\rm Poly}{}_G^d(V, W).
\eeqn
Given $(x, 0) \in X \subset F$, the tangent space of ${\rm graph}({\mf p})$ at the corresponding point is 
\beqn
\Big\{ \Big( \eta, \xi, D_x {\mf p}(x, 0)(\eta) + D_v{\mf p}(x, 0) (\xi) \Big)\ |\ \xi \in V,\ \eta \in T_x X \Big\} \subset T_x X \oplus V \oplus \wh{\rm Poly}{}_G^d(V, W).
\eeqn
As $V$ contains no trivial representations and ${\mf p}$ is $G$-invariant, the partial derivative $D_v {\mf p}$ vanishes along the zero section. Hence the above tangent space to ${\rm graph}({\mf p})$ is $T_x  X \oplus V \oplus {\rm Im} (D_x {\mf p})$. As ${\rm proj}_2$ collapses the $V$-direction, the transversality is equivalent to the transversality of ${\rm Im} (D_x {\mf p})$ to  ${\rm proj}_2: Z_\alpha \to \wh{\rm Poly}{}_G^d(V, W)$.
\end{proof}

\begin{lemma}\label{lemma46}
Suppose $F_G = 0$ and ${\mf p}_0: F \to \wh{\rm Poly}{}_G^d(F, E)$ is FOP transverse near a closed subset $Y \subset F$, then there exists ${\mf p}: F \to \wh{\rm Poly}{}_G^d(F, E)$ which is FOP transverse near the base $X$ such that ${\mf p} = {\mf p}_0$ near $Y$. Moreover, if we equip $\wh{\rm Poly}{}_G^d(F, E)$ with a norm and are given $\delta >0$, then one can require that $\| {\mf p} - {\mf p}_0\|\leq \delta$.
\end{lemma}

\begin{proof}
By Lemma \ref{lemma45}, ${\mf p}_0|_X$ is transverse to ${\rm proj}_2: Z_\alpha \to \wh{\rm Poly}{}_G^d(F, E)$ near $Y \cap X$ for each stratum $Z_\alpha \in {\mf Z}^d(G, F, E)$. Then the standard transversality result shows that there exists a map ${\mf p}_1: X \to \wh{\rm Poly}{}_G^d(F, E)$ which is transverse to ${\rm proj}_2: Z_\alpha  \to \wh{\rm Poly}{}_G^d(F, E)$ for all $\alpha$, such that ${\mf p}_1$ agrees with ${\mf p}_0|_X$ near $Y \cap X$ and such that ${\mf p}_1$ is sufficiently close to ${\mf p}_0|_X$. Consider the pullback $\pi_F^* {\mf p}_1: F \to \wh{\rm Poly}{}_G^d(F, E)$; by abuse of notation, still denote it by ${\mf p}_1$. 

Now choose a pair of $G$-invariant neighborhoods $O' \subset O''$ of $Y$ such that the closure of $O'$ is contained in $O''$ and such that ${\mf p}_0$ is FOP transverse in $O''$. Choose a $G$-invariant cut-off function $\rho: F \to [0,1]$ supported in a sufficiently small neighborhood of $X \setminus O' \subset F \setminus Y$ which is also identically $1$ near $X \setminus O''$. Define
\beqn
{\mf p} = (1-\rho) {\mf p}_0 + \rho {\mf p}_1: F \to \wh{\rm Poly}{}_G^d(F, E)
\eeqn
which agrees with ${\mf p}_0$ near $Y\subset F$ and which is sufficiently close to ${\mf p}_0$. Moreover, its restriction to the zero section is transverse to each stratum of $Z^d(G, F, E)$ under ${\rm proj}_2$. Then by Lemma \ref{lemma45}, ${\mf p}$ is FOP transverse at each point of $X$. By Lemma \ref{lemma_openness}, it is FOP transverse near $X$.
\end{proof}

\subsection{FOP transverse sections}

Suppose $({\mc U}, {\mc E})$ is equipped with a straightening $({\mc U}^\#, {\mc E}^\#)$ (Definition \ref{defn_straightening}). Let $\hat C = (G, U, E, \hat\psi)$ be a bundle chart for ${\mc E} \to {\mc U}$. Then for each $U$-essential subgroup $H \subset G$, near $U_H$ a section $S: U \to E$ can be identified with a nonlinear bundle map
\beqn
S_H^\#: N^\epsilon U_H \to E|_{U_H}.
\eeqn

\begin{defn}\label{defn_chartwise_transverse}
We say that a $G$-equivariant section $S: U \to E$ is {\it FOP transverse at $x\in U$} if it is an NC section in a neighborhood of $x$ and if the NC bundle map 
\beqn
S_{G_x}^\#: N^\epsilon U_{G_x} \to E|_{U_{G_x}}
\eeqn
is FOP transverse at $(x, 0)$. $S$ is said to be {\it FOP transverse} if it is FOP transverse at all points. 
\end{defn}

It is easy to see that the FOP transversality condition is invariant under chart embeddings.

\begin{lemma}\label{lemma48}
Let $\hat C_i = (G_i, U_i, E_i, \hat\psi_i)$, $i = 1, 2$ be two overlapping bundle charts of ${\mc E} \to {\mc U}$ such that ${\mc S}$ is pulled back to equivariant sections $S_i: U_i \to E_i$. Suppose $\psi_1(x_1) = \psi_2(x_2)$. Then $S_1$ is FOP transverse at $x_1$ if and only if $S_2$ is FOP transverse at $x_2$.
\end{lemma}

\begin{proof}
It follows from the invariance property of the canonical Whitney stratification provided by Proposition \ref{prop_nice_invariance}.
\end{proof}

\begin{lemma}
Let $S: U \to E$ be a $G$-equivariant NC section. Then the set of points $x \in U$ where $S$ is FOP transverse is open.
\end{lemma}

\begin{proof}
We may assume $G_x = G$. Using the straightening, we may assume $S$ is an NC bundle map $S_G^\#\in C_G^{\rm NC} (F, E)$ with $F = NU_G$ and $x$ is a point in the base $X = U_G$. Let ${\mf p}: F \to \wh{\rm Poly}{}_G^d(F, E)$ be a lift of $S_G^\#$. Then by definition, ${\rm graph}({\mf p}) \subset M^d(G, F, E)$ is transverse to $Z^d(G, F, E)$ at $(x, 0, {\mf p}(x, 0))$. By Trotman's theorem, they are transverse at all nearby points $(y, v)$. If $H = G_v \subset G$, then the map $\mu$ in \eqref{mapmu} changes ${\mf p}$ to an $H$-invariant lift ${\mf q}$ of $S_G^\#$ near $(y, v)$. By Lemma \ref{lemma44}, $S_G^\#$ is FOP transverse at $(y, v)$. Therefore, the set of FOP transverse points is open.
\end{proof}

We can then prove the chartwise existence and $C^0$-density of FOP transverse sections.

\begin{prop}\label{prop410}
Let $\hat C = (G, U, E, \hat\psi)$ be a bundle chart of ${\mc E}$ and $S_0: U \to E$ be an NC section (with respect to the fixed straightening), then for each $\delta>0$, there exists an FOP transverse section $S: U \to E$ such that $\| S - S_0 \|_{C^0} \leq \delta$. Moreover, if $Y\subset U$ is a $G$-invariant closed subset and $S$ is FOP transverse near $Y$, then $S$ can be chosen to agree with $S_0$ near $Y$.
\end{prop}

\begin{proof}
Via an induction on essential subgroups, one can see it suffices to prove the following statement. Let $H \subset G$ be a maximal $U$-essential subgroup. Then there exists an NC section $S: U \to E$ which is FOP transverse near $Y \cup G U_H$ and which agrees with $S_0$ near $Y$. To prove this statement, first choose a $G$-invariant
\beqn
\epsilon: GU_H \to {\mb R}_+.
\eeqn
As $H$ is a maximal $U$-essential subgroup, it is possible to choose $\epsilon$ such that for any $H'$ conjugate to $H$,
\beqn
|N^\epsilon U_H| \cap |N^\epsilon U_{H'}| \neq \emptyset \Longrightarrow H = H'.
\eeqn
Next, consider the induced NC bundle map 
\beqn
S_{0, H}^\#: N^\epsilon U_H \to E|_{U_H}.
\eeqn
Choose a lift
\beqn
{\mf p}_0: N^\epsilon U_H \to \wh{\rm Poly}{}_H^d( NU_H, E|_{U_H}).
\eeqn
Then by Lemma \ref{lemma46}, there exists 
\beqn
{\mf p}_1: N^\epsilon U_H \to \wh{\rm Poly}{}_H^d(NU_H, E|_{U_H})
\eeqn
which is sufficiently close to ${\mf p}_0$,  agrees with ${\mf p}_0$ near $Y \cap |N^\epsilon U_H|$, and is FOP transverse near $U_H$. Let $S_{{\mf p}_1}: N^\epsilon U_H \to E|_{|N^\epsilon U_H|}$ be the corresponding section defined in the tubular neighborhood $|N^\epsilon U_H|$ of $U_H$. Via $G$-action, the lift ${\mf p}_1$ induces sections defined on $|N^\epsilon U_{H'}|$ for all $H'$ conjugate to $H$ such that the totality defines a $G$-equivariant NC section
\beqn
S_1: G|N^\epsilon U_H| \to E|_{G|N^\epsilon U_H|}.
\eeqn
Now choose a $G$-invariant cut-off function $\rho: U \to [0, 1]$ supported near $GU_H$. Consider the section 
\beqn
S = (1-\rho) S_0 + \rho S_1.
\eeqn
Then by Lemma \ref{lemma252}, $S$ is an NC section. Moreover, $S$ is FOP transverse near $GU_H$ and agrees with $S_0$ near $Y$.
\end{proof}

\subsection{Proof of Theorem 
\ref{thm11}}

\subsubsection{FOP transverse sections on orbifolds, classical transversality, and locality}

We first formally define the space of FOP transverse sections.
\begin{defn}\label{defn411} 
Let ${\mc U}$ be an NC orbifold (without boundary) and ${\mc E} \to {\mc U}$ be an NC vector bundle. Recall the definition of straightenings (Definition \ref{defn_straightening}).

\begin{enumerate}

\item Let $({\mc U}^\#, {\mc E}^\#)$ be a straightening of $({\mc U}, {\mc E})$, an NC section ${\mc S} \in \Gamma^{\rm NC}({\mc U}^{\#}, {\mc E}^\#)$ is called \emph{FOP transverse} if for each bundle chart $\hat C = (G, U, E, \hat\psi)$, the pullback section $S: U \to E$ is FOP transverse (Definition \ref{defn_chartwise_transverse}).

\item Define $\Gamma^{\rm FOP}({\mc U}^\#, {\mc E}^\#) \subset \Gamma^{\rm NC} ( {\mc U}^\#, {\mc E}^\#)$ to be the subset of sections which are FOP transverse.

\item By taking the union over all straightenings on $({\mc U}, {\mc E})$, define 
\beqn
\Gamma^{\rm FOP}({\mc U}, {\mc E}) = \bigcup_{({\mc U}^\#, {\mc E}^\#)} \Gamma^{\rm FOP}({\mc U}^\#, {\mc E}^\#)
\eeqn
and call its elements {\it FOP transverse sections} of ${\mc E}$.
\end{enumerate}
\end{defn}

The {\bf (Classical Transversality)} and {\bf (Locality)} stated in Theorem \ref{thm11} are straightforward. If ${\mc U}$ is a manifold, then there is a canonical (trivial) NC structure and all smooth sections are NC sections. The equivalence between FOP transversality and classical transversality can be traced via the definition. On the other hand, if ${\mc S} \in \Gamma^{\rm FOP}({\mc U}, {\mc E})$, then it is FOP transverse with respect to a certain straightening $({\mc U}^\#, {\mc E}^\#)$. The straightening restricts to one over any open subset ${\mc U}'$. As the FOP transversality condition is defined pointwise, it follows that ${\mc S}|_{{\mc U}'} \in \Gamma^{\rm FOP}({\mc U}', {\mc E}|_{{\mc U}'})$.

\subsubsection{Density of FOP transverse sections and Extension Property I}

\begin{prop}\label{prop:extension-Y}
Given a smooth section ${\mc S}_0\in \Gamma ({\mc U}, {\mc E})$ which is FOP transverse near a closed subset ${\mc Y} \subset {\mc U}$, for for any $\delta >0$, there exists an FOP transverse section ${\mc S}: {\mc U} \to {\mc E}$ which agrees with ${\mc S}_0$ near ${\mc Y}$ such that $\| {\mc S}_0 - {\mc S} \|_{C^0} \leq \delta$.
\end{prop}

\begin{proof}
By definition, there exists a straightening near ${\mc Y}$ with respect to which ${\mc S}_0$ is normally complex and FOP transverse. By Corollary \ref{cor_straightening_extension}, there exists a straightening $({\mc U}^\#, {\mc E}^\#)$ of $({\mc U}, {\mc E})$ which agrees with the existing one near ${\mc Y}$. As normally complex sections form a submodule over $C^\infty({\mc U})$, by using a cut-off function, we may assume that ${\mc S}_0$ is normally complex with respect to this straightening by the density property stated in Proposition \ref{prop_complex_dense}. For the rest of the proof, we will use this straightening.

On the other hand, as $|{\mc U}|$ is Lindel\"of and paracompact, one can find a countable and locally finite cover of $|{\mc U}|$ by bundle charts $\hat C_i = (G_i, U_i, E_i, \hat\psi_i)$, $i = 1, \ldots, \infty$. We may also find precompact $G_i$-invariant open subsets $U_i' \subset U_i$ such that $\{ \psi_i(U_i') \}$ still covers ${\mc U}$. We prove the following statement by induction on $k = 1, 2, \ldots$.
\begin{itemize}
\item For each $k$, there exists a section ${\mc S}_k \in \Gamma^{\rm NC}({\mc U}^\#, {\mc E}^\#)$ such that ${\mc S}_k$ is FOP transverse near $\bigcup_{i\leq k} \psi_i( \ov{U_i'})$ and which agrees with ${\mc S}_0$ near ${\mc Y}$. Moreover
\beqn
\| {\mc S}_k - {\mc S}_0 \|_{C^0} \leq \left( 1 - \frac{1}{2^k}\right) \delta.
\eeqn
\end{itemize}

For the $k = 1$ case, let $S_{0,1}: U_1 \to E_1$ be the pullback of ${\mc S}_0$ to the chart $\hat C_1$. Then by Proposition \ref{prop410}, there exists $S_1': U_1 \to E_1$ which is FOP transverse on $U_1$ and which agrees with $S_{0,1}$ near $Y_1:= \psi_1^{-1}({\mc Y})$ such that $\| S_1' - S_{0,1} \|_{C^0} \leq \frac{\delta }{2}$. Choose a $G_1$-invariant cut-off function $\rho_1: U_1 \to [0,1]$ supported near $\ov{U_1'}$ such that $\rho_1|_{U_1'} \equiv 1$. Regarding $\rho_1$ as a function on ${\mc U}$, then $\rho_1 S_1'$ defines a normally complex section of ${\mc E}$. Then define 
\beqn
{\mc S}_1 = (1 - \rho_1) {\mc S}_0 + \rho_1 S_1'
\eeqn
which is FOP transverse near $\psi_1(\ov{U_1'})$ and which agrees with ${\mc S}_0$ near ${\mc Y}$. Moreover, 
\beqn
\| {\mc S}_1 - {\mc S}_0 \|_{C^0} \leq \frac{\delta}{2}.
\eeqn 

Suppose we have proved the above statement for $k-1$.  Denote
\begin{align*}
&\ {\mc Y}_k:= {\mc Y} \cup \bigcup_{i \leq k-1} \psi_i( \ov{U_i'} ) \subset {\mc Y} ,\ &\ Y_k:= \psi_k^{-1}({\mc Y}_k) \subset U_k.
\end{align*}
The latter is a $G_k$-invariant closed subset. Let $S_{k-1,k}: U_k \to E_k$ be the pullback of ${\mc S}_{k-1}$ to $\hat C_k$. Then $S_{k-1, k}$ is FOP transverse near $Y_k$. Then by Proposition \ref{prop410}, there exists an FOP transverse section $S_k': U_k \to E_k$ which agrees with $S_{k-1, k}$ near $Y_k$ such that 
\beqn
\| S_k' - S_{k-1, k} \|_{C^0} \leq \frac{\delta}{2^k}.
\eeqn
Choose a $G_k$-invariant cut-off function $\rho_k: U_k \to [0,1]$ supported near $\ov{U_k'}$ such that $\rho_k = 1$ over $U_k'$. Also regard $\rho_k$ as a function on ${\mc U}$. Then $\rho_k S_k'$ is a normally complex section of ${\mc E}$. Define
\beqn
{\mc S}_k = (1-\rho_k) {\mc S}_{k-1} + \rho_k S_k'.
\eeqn
By construction, the section agrees with ${\mc S}_{k-1}$ near ${\mc Y}_k$. Moreover, by construction, it is FOP transverse near $\psi_i(\ov{U_k'})$. Lastly, 
\beqn
\| {\mc S}_k - {\mc S}_0 \|_{C^0} \leq \| {\mc S}_k - {\mc S}_{k-1} \|_{C^0} + \| {\mc S}_{k-1} - {\mc S}_0 \|_{C^0} \leq \left( 1- \frac{1}{2^k} \right) \delta.
\eeqn
This proves the $k$-th case of the induction.

Finally, as the countable cover is locally finite, on any compact subset, the sequence ${\mc S}_k$ stabilizes for $k$ sufficiently large. Hence in the limit one obtains an FOP transverse section ${\mc S}$ satisfying the required conditions.
\end{proof}

\subsubsection{Extension property II}

We consider the problem of extending FOP transverse sections from a suborbifold ${\mc X} \subset {\mc U}$ with an ordinary normal bundle (Definition \ref{defn_ordinary_normal}). Denote ${\mc E}_{\mc X}:= {\mc E}|_{\mc X}$. By the discussion of Subsection \ref{subsec:ncs}, $({\mc X}, {\mc E}_{\mc X})$ is an NC pair. Let ${\mc S}_{\mc X}: {\mc X} \to {\mc E}_{\mc X}$ be an FOP transverse section. Then by definition, there exists a straightening $({\mc X}^\#, {\mc E}_{\mc X}^\#)$ of $({\mc X}, {\mc E}_{\mc X})$ such that ${\mc S}_{\mc X}$ is an FOP transverse section with respect to this straightening. Then by Proposition \ref{prop241}, there exists a straightening $({\mc U}^\#, {\mc E}^\#)$ which extends $({\mc X}^\#, {\mc E}_{\mc X}^\#)$. Moreover, the involved structures are all straightened along $N{\mc X}$. Let $\epsilon: {\mc X} \to {\mb R}_+$ be a sufficiently small smooth function with an associated disk bundle $N^\epsilon {\mc X}$. Then the normal exponential map $\exp^{N{\mc X}}: N^\epsilon {\mc X} \to {\mc U}$ is an open embedding, covered by a bundle map 
\beqn
{\rm par}^{N{\mc X}}: \pi_{N{\mc X}}^* {\mc E}|_{\mc X} \to {\mc E}
\eeqn
which defines a bundle isomorphism when restricting ${\mc E}$ along the image of $\exp^{N{\mc X}}$. We claim that the pullback 
\beq\label{eqn:pull-back section}
{\rm par}^{N{\mc X}} \circ ( \pi_{N{\mc X}}^* {\mc S}_{\mc X}) \circ (\exp^{N{\mc X}})^{-1}: |N^\epsilon {\mc X}| \to {\mc E}
\eeq
is FOP transverse. It suffices to consider this claim locally. Let $\hat C = (G, U, E, \hat\psi)$ be a bundle chart centered at a point $x \in {\mc X}$. Then by definition, $\psi^{-1}({\mc X})$ is a $G$-invariant submanifold $X \subset U$ with a normal bundle $NX \to X$. The section ${\mc S}_{\mc X}$ is pulled back to an NC section 
\beqn
S_X: X \to E|_X
\eeqn
which induces an NC bundle map (for some $\delta: X_{G_x} \to {\mb R}_+$)
\beqn
S_{X, G_x}^\#: N^\delta X_{G_x} \to E|_{X_{G_x}}.
\eeqn
As $G_x$ acts trivialy on the normal bundle $NX \to X$, the trivial factor in the normal direction in $NX$ does not affects the FOP transversality. Hence the pullback 
\beqn
(\pi^{NX})^* S_{X, G_x}^\#
\eeqn
is still FOP transverse. Therefore, after extending the FOP transverse section \eqref{eqn:pull-back section} using a smooth cut-off function, we can apply Proposition \ref{prop:extension-Y} to get the desired extension.

\subsubsection{Product property}

For $i = 1, 2$, let $({\mc U}_i, {\mc E}_i)$ be an NC pair. Given ${\mc S}_i \in \Gamma^{\rm FOP}({\mc U}_i, {\mc E}_i)$. Denote ${\mc U} = {\mc U}_1 \times {\mc U}_2$, ${\mc E} = {\mc E}_1 \boxplus {\mc E}_2$, and ${\mc S} = {\mc S}_1 \boxplus {\mc S}_2$. We prove that ${\mc S}$ is FOP transverse. Indeed, by definition, there exists a straightening $({\mc U}_i^\#, {\mc E}_i^\#)$ on $({\mc U}_i, {\mc E}_i)$ such that 
\beqn
{\mc S}_i \in \Gamma^{\rm FOP}({\mc U}_i^\#, {\mc E}_i^\#).
\eeqn
The two straightenings induce a straightening on $({\mc U}, {\mc E})$ by taking the product. We check that ${\mc S}$ is FOP transverse at each point $p = (p_1, p_2)$ with respect to the product straightening. If $\hat C_i = (G_i, U_i, E_i, \hat \psi_i)$ is a bundle chart of ${\mc E}_i$ centered at $p_i$, then 
\beqn
\hat C = (G, U, E, \hat \psi):= \hat C_1 \times \hat C_2:= (G_1 \times G_2, U_1 \times U_2, E_1 \boxplus E_2, \hat\psi_1 \boxplus \hat\psi_2)
\eeqn
is a bundle chart of ${\mc E}$. Let 
\beqn
S_{G_i}^\#: N^{\epsilon_i} U_{i, G_i} \to E_i|_{U_{G_i}}
\eeqn
be the NC bundle maps induced from the straightenings. Then for $\epsilon < \epsilon_1, \epsilon_2$, the bundle map induced from the product straightening
\beqn
S_G^\# = (S_{G_1}^\#, S_{G_2}^\#): N^\epsilon U_G \to E|_{U_G}
\eeqn
is well-defined, and it is obviously an NC bundle map. Abbreviate $X_i = U_{i, G_i}$, $F_i = NU_{i, G_i}$, $E_i = E_i|_{U_{G_i}}$. Let ${\mf p}_i: F_i \to \wh{\rm Poly}{}_{G_i}^d( F_i, E_i)$ be a lift of $S_{G_i}^\#$. Then 
\beqn
{\mf p}:= ({\mf p}_1, {\mf p}_2): F_1 \boxplus F_2 \to \wh{\rm Poly}{}_G^d( F_1 \boxplus F_2, E_1 \boxplus E_2)
\eeqn
is a lift of $S_G^\#$. Then by Proposition \ref{prop_product_Z}, the FOP transversality of ${\mf p}_1$ and ${\mf p}_2$ implies the FOP transversality of ${\mf p}$. Hence ${\mc S}$ is FOP transverse. 

\subsubsection{Stratified transversality}

Suppose ${\mc S} \in \Gamma^{\rm FOP}({\mc U}^\#, {\mc E}^\#) \subset \Gamma^{\rm FOP}({\mc U}, {\mc E})$. Let $\gamma$ be represented by an NC triple $(G, V, W)$. Let $p \in {\mc U}_\gamma$. Then it suffices to consider on a chart containing $p$. Let $\hat C = (G, U, E, \hat\psi)$ be a bundle chart centered at $p$. Then $\psi: U_G \to {\mc U}_\gamma$ gives a manifold chart of ${\mc U}_\gamma$. Let $S: U \to E$ be the section corresponding to ${\mc S}$. Using the straightening one obtains the bundle map 
\beqn
S_G^\#: N^\epsilon U_G \to E|_{U_G}.
\eeqn
Using the basic decomposition $E|_{U_G} = E_G \oplus \check E_G$, one obtains the two components $S_G$ and $\check S_G$ of $S_G^\#$. Then by the definition of FOP transversality, the section 
\beqn
S_G: N^\epsilon U_G \to E_G
\eeqn
is transverse at $\psi^{-1}(p) \in U_G$. As the derivative of $S_G$ in the normal direction to $U_G$ vanishes, this is equivalent to that $S_G|_{U_G}$ is transverse at $0$. This implies the stratified transversality.

\subsubsection{Stabilization property}

Start with ${\mc S} \in \Gamma^{\rm FOP}({\mc U}^\#, {\mc E}^\#) \subset \Gamma^{\rm FOP}({\mc U}, {\mc E})$. Let ${\mc F} \to {\mc U}$ be an NC vector bundle. Then by Proposition \ref{prop_straightened_connection}, there exist a connection $\nabla^{\mc F}$ straightened with respect to $g^{T{\mc U}}$, an inner product $h^{\mc F}$ preserved by $\nabla^{\mc F}$, and a nearby NC structure ${\mc I}^{\mc F}$ which is straightened with respect to $g^{T{\mc U}}$ and $\nabla^{\mc F}$ such that $(\nabla^{\mc F}, h^{\mc F}, {\mc I}^{\mc F})$ is normally Hermitian. By Proposition \ref{prop_straightening_stabilization}, the bundle metric $g^{T{\mc F}}$ and the bundle NC structure ${\mc I}^{T{\mc F}}$ induce a straightening ${\mc F}^\#$ of ${\mc F}$. There are also straightenings 
\beqn
(\pi_{\mc F}^* {\mc E})^\#\ {\rm and}\ (\pi_{\mc F}^* {\mc F})^\#.
\eeqn
We claim 
\beqn
\pi_{\mc F}^* {\mc S} \oplus \tau_{\mc F} \in \Gamma^{\rm FOP} \Big( {\mc F}^\#, (\pi_{\mc F}^* {\mc E})^\# \oplus ( \pi_{\mc F}^* {\mc F})^\# \Big).
\eeqn
Still, it suffices to consider this problem locally on a single chart. Let $\hat C = (G, U, E, \hat\psi)$ be a bundle chart of ${\mc E}$. Assume ${\mc F}$ is pulled back to a $G$-equivariant bundle $F \to U$ on this chart, which can be achieved by shrinking $U$ if necessary. Then ${\mc S}$ is pulled back to an FOP transverse NC section $S: U \to E$. Consider the section  
\beqn
\tilde S:= \pi_F^* S \oplus \tau_F: F \to \pi_F^* E \oplus \pi_F^* F.
\eeqn
The bundle straightening construction implies that $\tilde S$ is an NC section. Moreover, if ${\mf p}: N^\epsilon U_G \to \wh{\rm Poly}{}_G^d(NU_G, E|_{U_G})$ is a lift of $S_G^\#$, then there is a natural lift  
\begin{multline*}
\tilde {\mf p}: NF_G \to \wh{\rm Poly}{}_G^d( NF_G, \pi_{F_G}^* E \oplus \pi_{F_G}^* F) \\
= \wh{\rm Poly}{}_G^d( NF_G, \pi_{F_G}^* E \oplus \pi_{F_G}^* F_G \oplus \pi_{F_G}^* \check F_G)  \\
= \pi_{F_G}^* \wh{\rm Poly}{}_G^d( NU_G, E|_{U_G}) \oplus \pi_{F_G}^* \wh{\rm Poly}{}_G^d( \check F_G, \check F_G) \oplus \pi_{F_G}^* F_G
\end{multline*}
of $\tilde S_G^\#$ defined as follows. Use $u$ to denote fiber coordinates of $F$; with respect to the basic decomposition $F = F_G \oplus \check F_G$, write $u = (u_G, \check u_G)$. Then 
\beqn
\tilde {\mf p}(x, v, u_G, \check u_G) = \Big( {\mf p}(x, v), {\rm Id}_{\check F_G}, u_G \Big).
\eeqn
In terms of graphs, it follows that 
\beqn
\eta^{-1}({\rm graph}(\tilde {\mf p}))=  {\rm graph}({\mf p})
\eeqn
where $\eta$ is the (family version of the) map \eqref{polynomial_stabilization}. Then by Proposition \ref{prop_Whitney_stabilization}, ${\mf p}$ is FOP transverse at $x \in U_G$ if and only if $\tilde {\mf p}$ is FOP transverse at $(x, 0) \in F$.

\subsubsection{Stratification of zero locus}\label{subsubsec:stratification}

Suppose the virtual dimension of $({\mc U}, {\mc E})$ is $k$. Let ${\mc S}: {\mc U} \to {\mc E}$ be an FOP transverse section. For each universal stratum $\theta \in {\mf Z}_k^{\rm univ}$, we define the subset
\beqn
{\mc S}^{-1}(0)_\theta\subset {\mc S}^{-1}(0)
\eeqn
by identifying its intersection with each chart. Recall that one has the map $\rho_k: {\mf Z}_k^{\rm univ} \to \uds{\mb \Gamma}_k^{\rm NC}$. For a bundle chart $\hat C = (G, U, E, \hat\psi)$, if $\psi(U) \cap {\mc U}_{\rho_k (\theta)} = \emptyset$, then ${\mc S}^{-1}(0)_\theta \cap \psi (U) = \emptyset$. Otherwise, suppose $\rho_k (\theta)$ is represented by a triple $(H, V, W)$ where $H \subset_U G$, $V$ is isomorphic to fibers of $TU|_{U_H}$ and $W$ is isomorphic to fibers of $E|_{U_H}$. Choose a lift  
\beqn
{\mf p}: N^\epsilon U_H \to \wh{\rm Poly}{}_G^d (NU_H, E|_{U_H})
\eeqn
of the induced bundle map $S_H^\#: N^\epsilon U_H \to E|_{U_H}$. Then one has 
\beqn
(S_H^\#)^{-1}(0) = {\rm proj}_1 \Big(  {\rm graph}({\mf p}) \cap Z^d(H, NU_H, E|_{U_H}) \Big)
\eeqn
where 
\beqn
{\rm proj}_1: M^d(H, NU_H, E|_{U_H}) = NU_H \oplus \wh{\rm Poly}{}_H^d(NU_H, E|_{U_H}) \to NU_H
\eeqn
is the projection onto the first factor.
Notice that $\theta$ is represented by an element in $Z_\theta^d \in {\mf Z}^d(H, V, W)$. Then using the local triviality of the canonical Whitney stratification, we can define
\beqn
(S_H^\#)^{-1}(0)_\theta:= {\rm proj}_1 \Big( {\rm graph}({\mf p}) \cap Z^d_\theta \Big).
\eeqn
Then define ${\mc S}^{-1}(0)_\theta$ to be the union of images of all $(S_H^\#)^{-1}(0)_\theta$ under $\psi$. The FOP transversality condition implies that ${\mc S}^{-1}(0)_\theta$ is a smooth submanifold of ${\mc U}_{\rho_k (\theta)}$ of dimension $k + n_k (\theta)$. Therefore we have constructed a decomposition 
\beqn
{\mc S}^{-1}(0) = \bigsqcup_{\theta \in {\mf Z}_k^{\rm univ}} {\mc S}^{-1}(0)_\theta.
\eeqn
We need to prove that after we remove empty pieces, this partition makes ${\mc S}^{-1}(0)$ a Thom--Mather stratified space (Definition \ref{defn_Thom_Mather}). We first prove that it is a stratification. As locally this partition comes from a Whitney stratification on each orbifold chart, the partition is locally finite and each stratum is locally closed. It remains to verify the axiom of frontier (see Definition \ref{defn_stratification}), namely, 
\beqn
{\mc S}^{-1}(0)_\theta \cap \ov{{\mc S}^{-1}(0)_{\theta'}} \neq \emptyset \Longrightarrow {\mc S}^{-1}(0)_{\theta} \subset \ov{{\mc S}^{-1}(0)_{\theta'}}.
\eeqn
If $x \in {\mc S}^{-1}(0)_\theta \cap \ov{{\mc S}^{-1}(0)_{\theta'}}$, choose a bundle chart $\hat C = (G, U, E, \hat\psi)$ centered at $x$. Let $V$ be the fiber of $TU|_{U_G}$ and $W$ be the fiber of $E|_{U_G}$. Then the stable isotropy type $\rho_k (\theta)\in \uds{\mb \Gamma}_k^{\rm NC}$ is represented by $(G, V, W)$. Then for $d$ sufficiently large, $\theta$ corresponds to a $G$-fixed stratum $Z_\alpha \in {\mf Z}^d(G, V, W)$ and $\theta'$ corresponds to the $G$-orbit of some stratum $Z_{\alpha'} \geq Z_\alpha$. Then if $y \in {\mc S}^{-1}(0)_\theta$, choose a bundle chart centered at $y$ and denoted, by abuse of notation, by $\hat C = (G, U, E, \hat\psi)$. Let $S: U \to E$ be the pullback of ${\mc S}$ to $\hat C$. Then choose a lift of $S_G^\#$
\beqn
{\mf p}: N^\epsilon U_G \to \wh{\rm Poly}{}_G^d(NU_G, E|_{U_G}).
\eeqn
Then 
\beqn
y \in \psi \Big( {\rm proj}_1 \Big( {\rm graph}({\mf p}) \cap Z_\alpha^d(G, NU_G, E|_{U_G}) \Big)\Big).
\eeqn
As $Z_\alpha \leq Z_{\alpha'}$ (which means $Z_\alpha \subset \ov{Z_{\alpha'}}$), by Lemma \ref{lemma_pullback_Whitney} and its proof, one sees
\beqn
y \in \psi \Big( {\rm proj}_1 \Big( \ov{ {\rm graph}({\mf p}) \cap Z_{\alpha'}^d (G, NU_G, E|_{U_G}) } \Big)\Big).
\eeqn
It follows that $y \in \ov{{\mc S}^{-1}(0)_{\theta'}}$. 

Lastly we need to show that the decomposition of ${\mc S}^{-1}(0)$ by the smooth manifolds ${\mc S}^{-1}(0)_\theta$ makes the zero locus a Thom--Mather stratified space. Indeed, locally the preimage of ${\mc S}^{-1}(0)$ in each chart is a Whitney stratified space, hence it is a Whitney stratified set in ${\mc U}$ in the sense of Definition \ref{defnb15}. Then by Theorem \ref{thm_b15}, Proposition \ref{prop_b14}, and Definition \ref{defn_Thom_Mather}, ${\mc S}^{-1}(0)$ is a Thom--Mather stratified space.

\section{Integral Euler Cycles}\label{section5}

In this section, we introduce aspects of the bordism theory of derived orbifolds and derive homological consequences from zero loci of FOP transverse sections. We will not discuss such bordism theories in full generality; instead, we only introduce concepts that are needed to state our results. The interested reader can refer to \cite{pardon2022enough} and \cite{pardon2023orbifold} for a comprehensive treatment.

\subsection{Derived orbifold bordism}

\begin{defn}\label{defn_D_orbifold}
A {\it derived orbifold} (with or without boundary) is a triple ${\mc D} = ({\mc U}, {\mc E}, {\mc S})$ where ${\mc U}$ is an effective orbifold (with or without boundary), ${\mc E} \to {\mc U}$ is an orbifold vector bundle, and ${\mc S}: {\mc U} \to {\mc E}$ is a continuous section. 

\begin{enumerate}

\item ${\mc D}$ is called {\it compact} if ${\mc S}^{-1}(0)$ is compact.

\item An {\it orientation} on ${\mc D}$ is an orientation on the orbifold real vector bundle $T{\mc U} \otimes {\mc E}^*$.

\item A {\it normal complex structure} on ${\mc D}$ consists of an NC structure on ${\mc U}$ and an NC structure on ${\mc E}$. Two NC structures are concordant if the corresponding pairs of NC structures on ${\mc U}$ and ${\mc E}$ are concordant.

\item An {\it oriented and normally complex} (ONC for short) derived orbifold chart is a derived orbifold chart ${\mc D}$ together with an orientation and an NC structure.

\item The {\it virtual dimension} of ${\mc D}$ is the integer
\beqn
{\rm dim}_{\mb R}^{\rm vir} ({\mc D})  = {\rm dim}_{\mb R} ( {\mc U} )- {\rm rank}_{\mb R} ( {\mc E} ).
\eeqn

\end{enumerate}
\end{defn}

We also consider a more special structure called {\it stable complex structure}. 

\begin{defn}\label{defn_D_stable_complex}
Let ${\mc D} = ({\mc U}, {\mc E}, {\mc S})$ be a derived orbifold.
\begin{enumerate}

\item A {\it stable complex strucutre} on ${\mc E}$ is an equivalence class of quadruples
\beqn
(k, {\mc F}_0, {\mc F}_1, \phi)
\eeqn
where $k\geq 0$ is an integer, ${\mc F}_0, {\mc F}_1 \to {\mc U}$ are complex vector bundles, and 
\beqn
\phi: \uds{\mb R}^{\oplus k} \oplus {\mc E} \oplus {\mc F}_0 \cong {\mc F}_1
\eeqn
is an isomorphism of orbifold vector bundles. Here $\uds{\mb R}^{\oplus k}$ denotes the trivial rank $k$ bundle. The equivalence relation $\sim$ is generated by the following two relations: 1) we require that \beqn
(k, {\mc F}_0, {\mc F}_1, \phi) \sim (k+2, {\mc F}_0, \uds{\mb C} \oplus {\mc F}_1, \phi_0 \oplus \phi)
\eeqn
where
$\phi_0: \uds{\mb R}^{\oplus 2} \to \uds{\mb C}$ is the map $(x, y) \mapsto x + {\bf i} y$; 2) we require that 
\beqn
(k, {\mc F}_0, {\mc F}_1, \phi) \sim (k, {\mc F}_0 \oplus {\mc F}, {\mc F}_1 \oplus {\mc F}, \phi \oplus {\rm Id}_{\mc F})
\eeqn
where ${\mc F}\to {\mc U}$ is an arbitrary complex vector bundle.

\item A {\it stable complex structure} on ${\mc D}$ consists of a stable complex structure on $T{\mc U}$ and a stable complex structure on ${\mc E}$.
\end{enumerate}
\end{defn}

There are a few standard and basic operations on derived orbifolds. For a derived orbifold with boundary ${\mc D} = ({\mc U}, {\mc E}, {\mc S})$, one can restrict to the boundary 
\beqn
\partial {\mc D} = (\partial {\mc U}, {\mc E}|_{\partial {\mc U}}, {\mc S}|_{\partial {\mc U}}).
\eeqn
If ${\mc D}$ is ONC, so is $\partial {\mc D}$. For derived orbifolds without boundary, there are the operation ${\mc D} \mapsto - {\mc D}$ by reversing the orientation and the operation of taking disjoint union.

Now we consider the derived bordism group of topological spaces. We refrain from using the language of generalized homology theories. The introduction of the bordism group is merely for the convenience of stating the results.

\begin{defn}
Let $X$ be a topological space. The {\it oriented and normally complex derived orbifold bordism group} resp. {\it stably complex derived orbifold group} of $X$, denoted by ${\rm d}\Omega_*^{\rm ONC}(X)$ resp. ${\rm d}\Omega_*^{\rm U}(X)$, is the set of equivalence classes of quadruples $({\mc U}, {\mc E}, {\mc S}, \kappa)$ where ${\mc D} = ({\mc U}, {\mc E}, {\mc S})$ is a compact ONC resp. stably complex derived orbifold without boundary and $\kappa: {\mc U} \to X$ is a continuous map, while the equivalence relation is generated by the following relations
    \begin{enumerate}
    
\item {\bf (Shrinking)} $({\mc U}, {\mc E}, {\mc S}, \kappa) \sim ({\mc U}', {\mc E}', {\mc S}', \kappa')$ if ${\mc U}' \subset {\mc U}$ is an open neighborhood of ${\mc S}^{-1}(0)$, ${\mc E}' = {\mc E}|_{{\mc U}'}$, ${\mc S}' = {\mc S}|_{{\mc U}'}$, and $\kappa' = \kappa|_{{\mc U}'}$.

\item {\bf (Stabilization)} $({\mc U}, {\mc E}, {\mc S}, \kappa) \sim (\hat{\mc U}, \hat{\mc E}, \hat{\mc S}, \hat\kappa)$ if there is a NC resp. complex vector bundle $\pi_{\mc F}: {\mc F} \to {\mc U}$ such that $\hat{\mc U} = {\mc F}$, $\hat{\mc E} = \pi_{\mc F}^* {\mc E} \oplus \pi_{\mc F}^* {\mc F}$, $\hat{\mc S} = \pi_{\mc F}^* {\mc S} \oplus \tau_{\mc F}$, and $\hat\kappa: {\mc F} \to X$ is any extension of the induced map along the zero section.

\item {\bf (Cobordism)} $({\mc U}_0, {\mc E}_0, {\mc S}_0, \kappa_0) \sim ({\mc U}_1, {\mc E}_1, {\mc S}_1, \kappa_1)$ if there is a compact ONC resp. stably complex derived orbifold $\tilde {\mc D} = (\tilde {\mc U}, \tilde {\mc E}, \tilde {\mc S})$ with $\partial \tilde {\mc D} \cong (-{\mc D}_0 ) \sqcup {\mc D}_1$ and a continuous map $\tilde\kappa: \tilde {\mc U} \to X$ such that $\tilde\kappa|_{\partial \tilde {\mc U}}$ coincides with $\kappa_0\sqcup \kappa_1$.
\end{enumerate}
\end{defn}

The bordism sets ${\rm d}\Omega_*^{\rm ONC}(X)$ and ${\rm d}\Omega_*^{\rm U}(X)$ are graded abelian groups with addition induced by disjoint union and grading defined by the virtual dimension. It is easy to see that ${\rm d}\Omega_*^{\rm ONC}$ is a functor from topological spaces to graded abelian groups. Recall that there are the more classical bordism theories, the oriented bordism $\Omega_*^{\rm SO}$ and stable complex bordism $\Omega_*^{\rm U}$. 

\begin{prop}\label{prop54}
There is a natural transformation (the dotted arrow) making the following diagram commute (the other arrows are all naturally defined).
\beqn
\xymatrix{ \Omega_*^{\rm U} \ar[r] \ar[d] &   {\rm d}\Omega_*^{\rm U}   \ar@{.>}[d] \\
\Omega_*^{\rm SO}  \ar[r] & {\rm d}\Omega_*^{\rm ONC}  }
\eeqn
\end{prop}

\begin{proof}
Since complex vector bundles are canonically oriented, a stable complex structure on a vector bundle naturally induces an orientation. To obtain an NC structure, one needs to make certain choices. Suppose ${\mc D} = ({\mc U}, {\mc E}, {\mc S})$ is a stably complex derived orbifold. By definition, there exist two quadruples $(k, {\mc F}_0, {\mc F}_1, \phi)$ and $(l, {\mc E}_0, {\mc E}_1, \eta)$ where 
\begin{align*}
&\ \phi: \uds{\mb R}^{\oplus k} \oplus T {\mc U} \oplus {\mc F}_0 \cong {\mc F}_1,\ &\ \eta: \uds{\mb R}^{\oplus l} \oplus {\mc E} \oplus {\mc E}_0 \cong {\mc E}_1
\end{align*}
are isomorphisms. By taking stabilizations, which induce equivalent quadruples, we can assume ${\mc F}_0 = {\mc E}_0$. The above two isomorphisms can be pulled back to the total space of ${\mc F}_0$ as
\beqn
\pi_{{\mc F}_0}^* \phi: \uds{\mb R}^{\oplus k} \oplus \pi_{{\mc F}_0}^* T{\mc U} \oplus \pi_{{\mc F}_0}^* {\mc F}_0 \cong \pi_{{\mc F}_0}^* {\mc F}_1
\eeqn
and
\beqn
\pi_{{\mc F}_0}^* \eta: \uds{\mb R}^{\oplus l} \oplus \pi_{{\mc F}_0}^* {\mc E} \oplus \pi_{{\mc F}_0}^* {\mc F}_0 \cong \pi_{{\mc F}_0}^* {\mc E}_1.
\eeqn
Upon choosing a connection on ${\mc F}_0$ inducing a bundle isomorphism
\beqn
\Delta: \pi_{{\mc F}_0}^* T{\mc U} \oplus \pi_{{\mc F}_0}^* {\mc F}_0 \cong T{\mc F}_0,
\eeqn
the pullback isomorphism $\pi_{{\mc F}_0}^* \phi$ becomes an isomorphism
\beqn
\uds{\mb R}^{\oplus k} \oplus T{\mc F}_0 \cong \pi_{{\mc F}_0}^* {\mc F}_1.
\eeqn
As $\uds{\mb R}^{\oplus k}$ and $\uds{\mb R}^{\oplus l}$ do not affect the normal direction, the above two isomorphism induce a normal complex structure on the stabilization of ${\mc D}$ by ${\mc F}_0$. Note that the space of connections on ${\mc F}_0$ is contractible; moreover, for different presentations of the stable complex structures on ${\mc D}$, Definition \ref{defn_D_stable_complex} implies that the induced normal complex structures agree by an extension of the connection. Therefore, although the specific construction requires choices, the resulting bordism class in ${\rm d}\Omega_*^{\rm ONC}$ is canonically defined.
\end{proof}

\subsection{Definition of the Euler cycles}

There is a natural transformation from the oriented bordism to the integral homology
\beqn
\Omega_*^{\rm SO}(\cdot) \to H_*(\cdot; {\mb Z})
\eeqn
defined by pushing forward the fundamental classes of compact oriented manifolds. A virtual fundamental class can be viewed as a lift of this natural transformation to the derived orbifold bordism, possibly after changing the coefficients. 

\begin{thm}\label{thm55}
For $k\in {\mb Z}$ and a universal stratum $\theta \in {\mf Z}_k^{\rm univ}$, there is a natural transformation
\beqn
{\rm FOP}_\theta: {\rm d}\Omega_k^{\rm ONC}(\cdot) \to H_{k + n_k (\theta)}(\cdot; {\mb Z})
\eeqn
such that for the special maximal $\theta_k \in {\mf Z}_k^{\rm univ}$ (corresponding to the isotropy free part of the zero locus), the composition
\beqn
\Omega_{k}^{\rm SO}  \to {\rm d}\Omega_{k}^{\rm ONC} \xrightarrow{{\rm FOP}_{\theta_k}} H_{k}(\cdot; {\mb Z})
\eeqn
agrees with the natural transformation $\Omega_k^{\rm SO}(\cdot) \to H_k(\cdot; {\mb Z})$ which takes a closed oriented manifold to its fundamental class.
\end{thm}

The proof is to construct the Euler cycles via FOP transverse perturbations.

\begin{defn}
Let ${\mc D} = ({\mc U}, {\mc E}, {\mc S}_0)$ be a compact derived orbifold. A {\it perturbation} of ${\mc D}$ is a section ${\mc S}: {\mc U} \to {\mc E}$  satisfying the following condition: there exists a precompact open neighborhood ${\mc U}'$ of ${\mc S}_0^{-1}(0)$ and an inner product on ${\mc E}$ such that 
\beqn
| {\mc S}_0(x) - {\mc S}(x) | < |{\mc S}_0(x)|\ \ \forall x\notin \ov{{\mc U}'}.
\eeqn
In particular, ${\mc S}$ does not vanish outside $\ov{{\mc U}'}$.
\end{defn}

\begin{lemma}\label{lemma57}
Let ${\mc D} = ({\mc U}, {\mc E}, {\mc S}_0)$ be a compact and NC derived orbifold. Then there exists a smooth FOP transverse perturbation ${\mc S}: {\mc U} \to {\mc E}$. Moreover, if ${\mc D}$ has boundary and ${\mc S}_{\partial}$ is an FOP transverse perturbation of $\partial {\mc D}$, then one can take ${\mc S}$ such that ${\mc S}|_{\partial {\mc U}} = {\mc S}_{\partial}$.
\end{lemma}

\begin{proof}
We can choose a precompact open neighborhood ${\mc U}'$ of ${\mc S}_0^{-1}(0)$ and an inner product on ${\mc E}$ such that $|{\mc S}_0(x)| \geq 2 \delta$ for some $\delta > 0$ and all $x \notin {\mc U}'$. Then by the $C^0$-density of the space of FOP transverse sections, there exists an FOP transverse ${\mc S}: {\mc U} \to {\mc E}$ such that $\| {\mc S} - {\mc S}_0 \|_{C^0} \leq \delta$, which defines a perturbation of ${\mc D}$. If ${\mc S}_\partial$ is given, notice that $\partial {\mc U} \subset {\mc U}$ can be viewed as a suborbifold with ordinary normal bundle (the collar). Then by ({\bf Extension Property II}) of Theorem \ref{thm11}, there exists an FOP perturbation extending ${\mc S}_\partial$ to the interior. 
\end{proof}

\begin{prop}\label{prop58}
Let ${\mc D} = ({\mc U}, {\mc E}, {\mc S}_0)$ be a compact NC derived orbifold. 
\begin{enumerate}
    
    \item For each FOP transverse perturbation ${\mc S}$ and each universal stratum $\theta$, the closure $\ov{{\mc S}^{-1}(0)_\theta}$ is a compact oriented Thom--Mather stratified pseudomanifold (Definition \ref{defn_TM_pseudomanifold}).

    \item If $\tilde {\mc S}$ is an FOP transverse perturbation on the product $\tilde {\mc D} = [1, 2]\times {\mc D}$ with boundary restrictions denoted by ${\mc S}_1$ and ${\mc S}_2$ respectively, then the closure $\ov{\tilde{\mc S}^{-1}(0)_\theta}$ is a cobordism of oriented Thom--Mather stratified pseudomanifolds (Definition \ref{defnc8}).
\end{enumerate}
\end{prop}

\begin{proof}

For (1), by Theorem \ref{thm11}, the closure $\ov{{\mc S}^{-1}(0)_\theta}$ is a compact oriented Thom--Mather stratified space; furthermore, based on the discussion in Section \ref{subsubsec:stratification}, the frontier of ${\mc S}^{-1}(0)_\theta$ is the union of smooth manifolds of codimension two or higher as the same assertion holds for the universal zero locus. Hence ${\mc S}^{-1}(0)_\theta$ satisfies conditions for pseudomanifolds. For the same reason, in (2) the closure $\ov{\tilde{\mc S}^{-1}(0)_\theta}$ provides a cobordism. 
\end{proof}

\begin{proof}[Proof of Theorem \ref{thm55}]
Given an element $a \in {\rm d}\Omega_*^{\rm ONC}(X)$ represented by $({\mc U}, {\mc E}, {\mc S}_0, \kappa)$. Consider an FOP transverse perturbation ${\mc S}$. By (2) of Proposition \ref{prop58}, the Thom--Mather stratified pseudomanifold $\ov{{\mc S}^{-1}(0)_\theta}$ carries a fundamental class
\beqn
[\ov{{\mc S}^{-1}(0)_\theta}]\in H_{n_k (\theta)}( \ov{{\mc S}^{-1}(0)_\theta}; {\mb Z}),
\eeqn
see Appendix \ref{appendixc}. Via the inclusion into ${\mc U}$ and the map $\kappa: {\mc U} \to X$, we define
\beqn
\chi_\theta^{\rm FOP}(a):= \kappa_* [\ov{{\mc S}^{-1}(0)_\theta}] \in H_{n_k (\theta)}(X; {\mb Z}).
\eeqn
We only need to prove that this element only depends on the element in ${\rm d}\Omega_*^{\rm ONC}(X)$.

We first show that it is independent of the choice of FOP transverse perturbations. Indeed, if ${\mc S}_1$ and ${\mc S}_2$ are two such perturbations, then (3) of Proposition \ref{prop58} and Lemma \ref{lemmac9} imply that the corresponding cycles are equal. 

We then need to show that the class $\chi_\theta^{\rm FOP}(a)$ is invariant under the equivalence relation defining ${\rm d}\Omega_k^{\rm ONC}(X)$. Indeed, the subsequent Lemma \ref{lemma59}, Lemma \ref{lemma510}, and Lemma \ref{lemma511} show that the generators of the equivalence relation, the pushforward class in $X$ is invariant.

\begin{lemma}\label{lemma59}
Let ${\mc D} = ({\mc U}, {\mc E}, {\mc S}_0)$ be a compact ONC derived orbifold without boundary and $\iota: {\mc U}' \hookrightarrow {\mc U}$ be an open embedding onto a neighborhood of ${\mc S}_0^{-1}(0)$ and ${\mc D}'$ be the corresponding shrinking of ${\mc D}$. Then for any $\theta \in \Theta$,
\beqn
\chi_\theta^{\rm FOP}({\mc D}) = \iota_* (\chi_\theta^{\rm FOP}({\mc D}')).
\eeqn
\end{lemma}

\begin{proof}
One can first choose an FOP transverse perturbation ${\mc S}'$ of ${\mc D}'$. Then by using a cut-off function supported in ${\mc U}'$, one can find an FOP transverse perturbation ${\mc S}$ of ${\mc D}$ such that $\iota^* {\mc S}$ coincides with ${\mc S}'$ on a smaller neighborhood of ${\mc S}_0^{-1}(0)$ such that ${\mc S}^{-1}(0) = \iota(({\mc S}')^{-1}(0))$. Then the lemma follows.
\end{proof}

\begin{lemma}\label{lemma510}
Let ${\mc D} = ({\mc U}, {\mc E}, {\mc S}_0)$ be a compact ONC derived orbifold without boundary. Let $\hat {\mc D}$ be the stabilization of ${\mc D}$ by an NC vector bundle ${\mc F} \to {\mc U}$ and let $\iota: {\mc U} \to {\mc F}$ be the zero section. Then for each universal stratum $\theta \in \Theta$, 
\beq\label{eqn51}
\chi_{ \theta}^{\rm FOP}(\hat{\mc D}) = \iota_* \Big( \chi_{ \theta}^{\rm FOP}({\mc D}) \Big).
\eeq
\end{lemma}

\begin{proof}
Choose an FOP transverse perturbation ${\mc S}$ of ${\mc D}$. By the {\bf (Stabilization Property)} of Theorem \ref{thm11}, the stabilization $\hat{\mc S}$ is FOP transverse; it is easy to check it is also a perturbation of $\hat{\mc D}$. Then
\beqn
\iota( {\mc S}^{-1}(0)_\theta) = \hat{\mc S}^{-1}(0)_\theta. 
\eeqn
This is indeed an isomorphism of Thom--Mather stratified spaces and preserves orientation. Hence \eqref{eqn51} holds.
\end{proof}

\begin{lemma}\label{lemma511}
Let ${\mc D}_i = ({\mc U}_i, {\mc E}_i, {\mc S}_{0,i})$, $i = 1, 2$, be compact ONC derived orbifolds without boundary. Let $\tilde {\mc D}$ be a compact ONC cobordism from ${\mc D}_1$ to ${\mc D}_2$. Let $\iota: {\mc U}_1 \sqcup {\mc U}_2 \to \partial {\mc U}$ be the diffeomorphis contained in the cobordism. Then for each universal stratum $\theta \in {\mf Z}^{\rm univ}$, there exists a class $a_\theta \in H_*( \tilde {\mc U}, \partial \tilde {\mc U}; {\mb Z})$ such that 
\beqn
\iota_* \left( \chi_\theta^{\rm FOP}({\mc D}_2) - \chi_\theta^{\rm FOP}({\mc D}_1) \right) = \partial a_\theta.
\eeqn
\end{lemma}

\begin{proof}
Choose FOP transverse perturbations ${\mc S}_i: {\mc U}_i \to {\mc E}_i$ of ${\mc D}_i$. Then by Lemma \ref{lemma57}, one can find an extension $\tilde {\mc S}: \tilde {\mc U} \to \tilde {\mc E}$ as a perturbation of $\tilde {\mc D}$ whose boundary restriction is ${\mc S}_{0, i}$ respectively. Consider the stratum $\tilde {\mc S}^{-1}(0)_\theta \subset \tilde {\mc U}$ which is an oriented manifold with boundary identified with $(-{\mc S}_1^{-1}(0)_\theta) \sqcup ({\mc S}_2^{-1}(0)_\theta)$. Its closure is a cobordism of Thom--Mather stratified pseudomanifolds. Then there is a relative fundamental class
\beqn
[\tilde {\mc S}^{-1}(0)_\theta] \in H_* \Big( \ov{\tilde {\mc S}^{-1}(0)_\theta}, \partial \ov{\tilde {\mc S}^{-1}(0)_\theta} ; {\mb Z} \Big) \to H_* \Big( \tilde {\mc U}, \partial \tilde {\mc U}; {\mb Z} \Big).
\eeqn
By Lemma \ref{lemmac9}, its boundary is the difference of $[{\mc S}_2^{-1}(0)_\theta] - [{\mc S}_1^{-1}(0)_\theta]$. 
\end{proof}

Hence the map $\chi_\theta^{\rm FOP}$ is well-defined. It follows from the construction that this is a natural transformation from ${\rm d}\Omega_k^{\rm ONC}(\cdot)$ to $H_{k + n_k(\theta)} (\cdot; {\mb Z})$. Because any cobordism class of compact oriented manifolds tautologically gives rise to a  cobordism class of ONC derived orbifold, this trivial section is an FOP transverse section and the homology class indexed by $\theta_k$ is indeed the fundamental class.
\end{proof}

\begin{rem}
The product property of the FOP transversality condition allows us to derive certain multiplicative properties of the FOP natural transformations.
\end{rem}

\section{Integer-valued Gromov--Witten Type Invariants}\label{section_GW}

In this section we construct the integer-valued curve counting invariants of a general compact symplectic 
manifold and prove Theorem \ref{thm16}. Let $(X, \omega)$ be a compact symplectic manifold. Let $J$ be an $\omega$-compatible almost complex structure on $X$. Let $A \in H_2(M; {\mb Z})$ be a homology class. Choose nonnegative integers $g$ and $n$. Then there is the moduli space of genus $g$, $n$-marked stable maps in class $A$
\beqn
\ov{\mc M}_{g, n}(X, J, A).
\eeqn
Denote its virtual dimension by $k$. Let $\ov{\mc M}{}_{g,n}$ be the Deligne--Mumford space when $2g + n \geq 3$ and be a single point otherwise. We would like to define symplectic deformation invariants 
\beq\label{GW_bordism_class}
[\ov{\mc M}_{g, n}(X, J, A)_\theta] \in H_{k + n_k(\theta)}  (X^n \times \ov{\mc M}_{g,n}; {\mb Z})\ \forall \theta \in {\mf Z}_k^{\rm univ}
\eeq
using the FOP perturbation method and the global Kuranishi chart construction of Abouzaid--McLean--Smith \cite{AMS, AMS2} and Hirschi--Swaminathan \cite{Hirschi_Swaminathan_2022}. In view of Proposition \ref{prop54} and Theorem \ref{thm55}, to define these invariants, it suffices to construct a canonical element 
\beqn
[\ov{\mc M}_{g, n}(X, J, A)] \in {\rm d}\Omega_k^{\rm U}(X^n \times \ov{\mc M}{}_{g,n}).
\eeqn

\subsection{Stable complex Kuranishi bordism}

The Kuranishi approach of regularizing moduli spaces was originally introduced by Fukaya--Ono \cite{Fukaya_Ono} (see \cite{FOOO_Kuranishi} for a comprehensive discussion). Typically moduli spaces of pseudoholomorphic curves always admit local Kuranishi charts. The {\it global Kuranishi charts} on the stable map moduli spaces were firstly constructed by Abouzaid--McLean--Smith \cite{AMS} in the genus zero case and then extended by Hirschi--Swaminathan \cite{Hirschi_Swaminathan_2022} and Abouzaid--McLean--Smith \cite{AMS2} to the higher genus case. We first recall the notion of stable complex structure in the equivariant setting.

\begin{defn}\label{defn_G_stable_complex}
Let $G$ be a compact Lie group, $U$ be a smooth $G$-manifold, and $E \to U$ be a $G$-equivariant real vector bundle. A {\it $G$-equivariant stable complex structure} on $E$ is an equivalence class of quadruples
\beqn
(k, F_0, F_1, \phi)
\eeqn
where $k$ is a nonnegative integer, $F_0, F_1 \to V$ are $G$-equivariant complex vector bundles, and 
\beqn
\phi: \uds{\mb R}^{\oplus k} \oplus E \oplus F_0 \cong F_1
\eeqn
is an isomorphism of $G$-equivariant real vector bundles. The equivalence relation is generated by the following two relations: 1) we require
\beqn
(k, F_0, F_1, \phi) \sim (k+2, F_0, {\mb C} \oplus F_1, \phi_0 \oplus \phi)
\eeqn
where $\phi_0: {\mb R}^2 \to {\mb C}$ is $(x, y) \mapsto x + {\bf i} y$; 2) we require
\beqn
(k, F_0, F_1, \phi) \sim (k, F_0 \oplus F, F_1 \oplus F, \phi\oplus {\rm Id}_F)
\eeqn
where $F \to V$ is any $G$-equivariant complex vector bundle.
\end{defn}

\begin{defn}(\cite[Definition 4.1, 4.3]{AMS})
\begin{enumerate}

\item A {\it smooth Kuranishi space} (K-space for short) is a quadruple
\beqn
K = (G, U, E, S)
\eeqn
where $G$ is a compact Lie group, $U$ is a smooth $G$-manifold such that the $G$-action is effective and has only finite isotropy groups, $E \to U$ is a smooth $G$-equivariant vector bundle, and $S: U \to E$ is a continuous $G$-equivariant section. The K-space is called {\it compact} if $S^{-1}(0)$ is compact.

\item The {\it virtual dimension} of a K-space $K = (G, U, E, S)$ is 
\beqn
{\rm dim}_{\mb R}^{\rm vir} ( K ):= {\rm dim}_{\mb R} (U) - {\rm dim}_{\mb R} (G) - {\rm rank}_{\mb R} (E).
\eeqn

\item A {\it stable complex structure} on a smooth K-space $K = (G, U, E, S)$ consists of a $G$-equivariant stable complex structure on $TU/ \uds{\mf g}$ (where $\uds{\mf g} \subset TU$ is the trivial subbundle generated by infinitesimal actions, thanks to the condition that the $G$-action has at no continuous stabilizers) and a $G$-equivariant stable complex structure on $E$. A {\it stably complex K-space} consists of a smooth K-space together with a stable complex structure.
\end{enumerate}
\end{defn}

Similar to the case of derived orbifold bordism, we can define the Kuranishi version of the bordism groups.

\begin{defn}
Let $X$ be a topological space. The {\it stable complex Kuranishi bordism group} of $X$, denoted by ${\rm k}\Omega_*^{\rm U}(X)$, is the set of equivalence classes of quintuples $(G, U, E, S, \kappa)$ where $K = (G, U, E, S)$ is a compact and stably complex K-space and $\kappa: U \to X$ is a $G$-invariant continuous map, and where the equivalence relation is generated by the following relations.
\begin{enumerate}
\item {\bf (Shrinking)} $(G, U, E, S, \kappa) \sim (G, U', E', S', \kappa')$ if $U' \subset U$ is a $G$-invariant open neighborhood of $S^{-1}(0)$, $E' = E|_{U'}$, $S' = S|_{U'}$, and $\kappa' = \kappa|_{U'}$.

\item {\bf (Stablization)} $(G, U, E, S, \kappa) \sim (G, \hat U, \hat E, \hat S, \hat \kappa)$ if there is a $G$-equivariant complex vector bundle $\pi_F: F \to U$ and $\hat U = F$, $\hat E = \pi_F^* E \oplus \pi_F^* F$, $\hat S = \pi_F^* S \oplus \tau_F$, and $\hat\kappa = \pi_F^* \kappa$.

\item {\bf (Group enlargement)} $(G, U, E, S, \kappa) \sim (G \times G', U', E', S', \kappa')$ where $G'$ is another compact Lie group, if there is a $G$-equivariant principal $G'$-bundle $\pi_P: P \to U$ such that $U' = P$, $E' = \pi_P^* E$, $S' = \pi_P^* S$, and $\kappa' = \pi_P^* \kappa$.

\item {\bf (Cobordism)} $(G, U_0, E_0, S_0, \kappa_0) \sim (G, U_1, E_1, S_1, \kappa_1)$ if there is a compact and stably complex K-space with boundary $\tilde K = (G, \tilde U, \tilde E, \tilde S)$ whose boundary restriction is isomorphic to the disjoint union $(-K_0) \sqcup K_1$ such that $\tilde \kappa: \tilde U \to X$ is an extension of $\kappa_0 \sqcup \kappa_1$.
\end{enumerate}
\end{defn}

Notice that there is a natural transformation 
\beqn
\Omega_*^{\rm U} \to {\rm k}\Omega_*^{\rm U}
\eeqn
where we view a compact stably complex manifold as a stably complex K-space with trivial symmetry group and trivial obstruction bundle. 

\begin{lemma}
There is a natural transformation (the dotted arrow) making the following diagram commutes.
\beqn
\xymatrix{   &   {\rm k}\Omega_*^{\rm U} \ar@{.>}[d] \\
\Omega_*^{\rm U} \ar[ru] \ar[r] &  {\rm d}\Omega_*^{\rm U} }
\eeqn
\end{lemma}

\begin{proof}
Given a K-space $K = (G, U, E, S)$, one obtains a derived orbifold ${\mc D} = ({\mc U}, {\mc E}, {\mc S})$ by taking the $G$-quotient. More precisely, we set
\begin{align*}
&\ {\mc U}:= U/G, &\  {\mc E}:= E/G
\end{align*}
and ${\mc S}$ is induced from $S$. Moreover, a $G$-invariant map $\kappa: U \to X$ descends to a continuous map from ${\mc U}$ to $X$. 

Now we consider the correspondence of stable complex structures. If $K = (G, U, E, S)$ has a stable complex structure, then it induces a stable complex structure on ${\mc D}$ as follows. By Definition \ref{defn_G_stable_complex}, there exist $(k, F_0, F_1, \phi)$ and $(l, E_0, E_1, \eta)$ such that 
\begin{align*}
&\ \phi: \uds{\mb R}^{\oplus k} \oplus (TU/ \uds{\mf g}) \oplus F_0 \cong  F_1,\ &\ \eta: \uds{\mb R}^{\oplus l} \oplus E \oplus E_0 \cong E_1.
\end{align*}
Notice that $(TU/\uds{\mf g})/G$ is canonically identified with $T{\mc U}$. Therefore the two $G$-equivariant isomorphisms induce isomorphisms of orbifold vector bundles
\begin{align*}
&\ \phi/G: \uds{\mb R}^{\oplus k}\oplus T{\mc U} \oplus (F_0/G) \cong F_1/G,\ &\ \eta/G: \uds{\mb R}^{\oplus l} \oplus {\mc E} \oplus (E_0/G) \cong E_1/G.
\end{align*}
It is straightforward to check that the induced stable complex structure on ${\mc D}$ does not depend on the choices of the quadruples defining the stable complex structure. Moreover, one is readily to check that shrinkings resp. stabilizations of K-charts induce shrinkings resp. stabilizations of D-charts. On the other hand, group enlargements of K-charts induce isomorphic D-charts. Finally, cobordisms between K-spaces descend to cobordisms of derived orbifolds, and the commutativity of the diagram follows from definition.
\end{proof}

\subsection{Proof of Theorem \ref{thm16}}

Global topological Kuranishi charts on the moduli spaces of stable maps are constructed in \cite{AMS}\cite{Hirschi_Swaminathan_2022}\cite{AMS2}. 

\begin{thm}\label{thm_AMS_2}
Let $X, \omega, J, A, g, n$ be as in Theorem \ref{thm16}. There exists a distinguished (nonempty) collection of quintuples $(G, U, E, S, L)$ (called the {\rm AMS charts}) where $K = (G, U, E, S)$ is a stably complex K-space and $L: S^{-1}(0)/G \to \ov{\mc M}{}_{g,n}(X, J, A)$ is a homeomorphism. Each AMS chart is associated to a choice of an auxiliary datum. Moreover, the following is true.
\begin{enumerate}
    \item For each AMS chart $(G, U, E, S, L)$ there exists a $G$-invariant map 
    \beqn
    \kappa: U \to X^n \times \ov{\mc M}{}_{g,n}
    \eeqn
    which extends the evaluation-stabilization map
    \beqn
    (\ev\times {\rm st} ) \circ L^{-1}: S^{-1}(0)/G \to X^n \times \ov{\mc M}{}_{g,n}.
    \eeqn

    \item The element of ${\rm k}\Omega_*^{\rm U}(X^n \times \ov{\mc M}{}_{g,n})$ represented by $(G, U, E, S, \kappa)$ is independent of choices and the extension $\kappa$. Denote this element by $[\ov{\mc M}{}_{g,n}(X, J, A)]$.

    \item Suppose $\omega$ can be deformed via symplectic forms to another symplectic form $\omega'$ on $X$ and $J'$ is $\omega'$-compatible. Then 
    \beqn
    [\ov{\mc M}{}_{g,n}(X, J, A)] = [\ov{\mc M}{}_{g,n}(X, J', A)] \in {\rm k}\Omega_*^{\rm U}(X^n \times \ov{\mc M}{}_{g,n}).
    \eeqn
\end{enumerate}
\end{thm}

\begin{proof}
The construction of the AMS charts are given in \cite[Section 4]{AMS2} and the statementwhere part (1), (2), and the case of fixing $\omega$ and varying $J$ in (3) are summarized in \cite[Corollary 4.70]{AMS2}. The case of deforming both $\omega$ and $J$ can be carried out in the same way as in \cite[Section 4]{AMS2}.
\end{proof}

\begin{proof}[Proof of Theorem \ref{thm16}]
By applying the composition of natural transformations 
\beqn
{\rm k}\Omega_*^{\rm U} \to {\rm d}\Omega_*^{\rm U} \to {\rm d}\Omega_*^{\rm ONC}
\eeqn
to the element $[\ov{\mc M}{}_{g,n}(X, J, A)]$ given by Theorem \ref{thm_AMS_2}, one obtains a canonical element in ${\rm d}\Omega_k^{\rm ONC}(X^n \times \ov{\mc M}{}_{g,n})$. Then for each universal stratum $\theta \in {\mf Z}_k^{\rm univ}$ where $k$ is the virtual dimension, applying the FOP natural transformation, one obtains the $\theta$-th FOP Euler class
\beqn
[\ov{\mc M}{}_{g,n}(X, J, A)]_\theta^{\rm vir} \in H_{n_k (\theta)}(X^n \times \ov{\mc M}{}_{g,n}; {\mb Z}).
\eeqn
Lastly, when $\ov{\mc M}{}_{g,n}(X, J, A)$ contains no elements with nontrivial automorphism groups, for an AMS chart 
$(G, U, E, S, L)$ the $G$-action can be made free. Hence the only nontrivial FOP class is the one corresponding to the special element $\theta_k$ and it agrees with the ordinary virtual fundamental class.  
\end{proof}

\section{A Simple Proof of Abouzaid--McLean--Smith's Splitting Theorem}\label{sec:proof-AMS}

Let $(X, \omega)$ be a closed symplectic manifold. Denote by ${\rm Ham}(X, \omega)$ the group of Hamiltonian diffeomorphisms of $X$. Given an element of $\pi_1{(\rm Ham}(X, \omega))$ represented by a smooth loop $\phi: S^1 \to {\rm Ham}(X, \omega)$ one can construct a fibration $P_{\phi} \rightarrow S^2$ with fiber symplectomorphic to $X$ and the transition function along the equator being the loop $\phi$. We use the construction of integer-valued Gromov--Witten type invariants to give an alternative proof of \cite[Theorem 1.1]{AMS}:

\begin{thm}\label{thm:AMS}
There is an additive isomorphism 
\beqn
H^*(P_\phi; \mathbb{Z}) \cong H^*(X; \mathbb{Z}) \otimes_{\mathbb{Z}} H^*(S^2; \mathbb{Z}).
\eeqn
\end{thm}

For a discussion on the context of such cohomological splitting result, the reader may consult \cite[Section 1]{BPX}. The proof follows from constructing a splitting map in integral cohomology (cf. \cite[Section 3.4]{AMS}), whose proof is provided in this section.
\begin{prop}\label{prop72}
The map ${\rm res}: H^*(P_\phi; {\mb Z}) \to H^*(X; {\mb Z})$ induced by the inclusion of a fiber $X \hookrightarrow P_\phi$ has a right inverse.
\end{prop}

\subsection{The moduli spaces and the idea of the proof}

A typical idea to prove Proposition \ref{prop72} is to consider moduli spaces of holomorphic sections of $P_\phi \to S^2$ with two marked points, whose fundamental classes can be used to define Seidel representations, see \cite{LMP} and \cite{mcduffseidel}. Here we follow the alternate approach provided in \cite[Section 3]{AMS} where we embed the total space $P_\phi$ into a fibration over a 4-dimensional base to bypass the necessity of introducing Seidel elements in quantum cohomology. Our proof essentially follows the same strategy as \cite{AMS} where we replaces the Morava virtual fundamental class by the integral Euler class defined in this paper.

Consider the product $\mb{CP}^1\times \mb{CP}^1$. For convenience we denote the first factor by $B$ and the second factor by $S^2$. Let $M$ be the blow-up of $B \times S^2 \cong \mb{CP}^1\times\mb{CP}^1$ at $(0, 0)$. Let
\beqn
\pi_M: M \to B
\eeqn
be the composition of the blow-down map and the projection to $B$. For each $t \in B$, write $M_t:= \pi_M^{-1}(t)$. Then $M_t \cong S^2$ except at $t = 0 \in B$ where $M_0$ is a reducible rational curve with two components. We may choose a K\"ahler form $\omega_M$ on $M$ such that for an open neighborhood $W_\infty \subset B$ of $\infty \in B \cong \mb{CP}^1$, 
\beqn
\pi_M^{-1}(W_\infty) \cong W_\infty \times S^2
\eeqn
as K\"ahler manifolds. 

Further, there exists a symplectic fibration $\pi_P: P \to M$ equipped with a symplectic form $\omega_P$ on the total space $P$ as well as a compatible almost complex structure $J_P$ satisfying the following conditions (see \cite[Section 3]{mcduffseidel}).

\begin{enumerate}

\item $\pi_P: P \to M$ is pseudo-holomorphic.

\item $\pi_P^{-1}(\pi_M^{-1}(W_\infty))$ is symplectomorphic to the product $W_\infty \times S^2 \times X$ and $J_P$ is the product of given (almost) complex structures on the three factors.

\item $\pi_P^{-1}(M_0)$ is the singular space $P_\phi \cup_X P_{\phi^{-1}}$ such that the restriction of $\omega_P$ to both $P_\phi$ and $P_{\phi^{-1}}$ are in the deformation equivalence class of the symplectic structures determined by the loops $\phi$ and $\phi^{-1}$ (cf. \cite[Section 8.2]{McDuff_Salamon_2004}). 
Denote the two components of $M_0$ by $M_\phi$ and $M_{\phi^{-1}}$, which are the bases of $P_\phi$ and $P_{\phi^{-1}}$ respectively.
\end{enumerate}
\begin{figure}[h]
    \centering
    \includegraphics[scale=1.3]{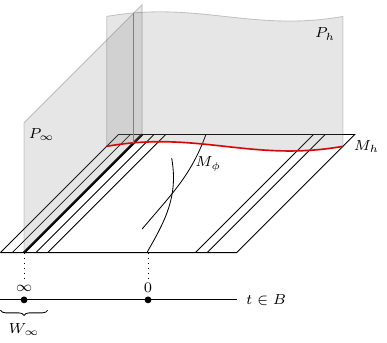}
    \caption{The fibration $P \to M$.}
    \label{figure1}
\end{figure}

Next, choose a holomorphic section $h: B \to M$ whose image $M_h:= h(B) \cong B$ passes through $M_\phi$ but not $M_{\phi^{-1}}$. Over $M_h$, denote
\beqn
P_h:= \pi_P^{-1}(M_h) \subset P.
\eeqn
One can fix a diffeomorphism
\beq\label{eqn71}
P_h \cong B \times X
\eeq
Denote $P_\infty = \pi_P^{-1}(M_\infty)$. See Figure \ref{figure1} for illustration.

We now describe a moduli space of stable maps. We fix a curve class $A \in H_2(P; \mathbb{Z})$ represented by $\{\infty\}\times S^2 \times \{{\rm pt}\}  \subset W_\infty \times S^2 \times X$. We consider the moduli space $\ov{\mc M}{}_{0,2}(P, J_P, A)$ with the evaluation map 
\beqn
{\rm ev} = ({\rm ev}_1, {\rm ev}_2): \ov{\mc M}{}_{0, 2}(P, J_P, A) \to P \times P.
\eeqn
Define
\beqn
\ov{\mc M}{}_h:={\rm ev}^{-1} (P \times P_h)
\eeqn
and its subspace
\beqn
\ov{\mc M}_\infty:= {\rm ev}^{-1} (P_\infty \times P_h).
\eeqn

Notice that there is a special family of curves near $P_\infty:= \pi_P^{-1}(M_\infty)$. As near $\pi_P^{-1}(W_\infty)$, $J_P$ is the product $J_B \times J_{S^2} \times J_X$, there is a subset
\beqn
{\mc M}_{W_\infty} \subset \ov{\mc M}_h
\eeqn
of curves contained in $\pi_P^{-1}(M_t)$ for $t \in W_\infty$ which is constant in the $X$-direction, where the first marked point is mapped into $P_h \cap P_t$ and the second marked point is arbitrary in that line.

\begin{lemma}\label{lemma73}
${\mc M}_{W_\infty}$ is an open neighborhood of $\ov{\mc M}_\infty$ in $\ov{\mc M}_h$ and is transversely cut out.
\end{lemma}

\begin{proof}
First, it is easy to see that $\pi_P: P \to M$ pushes forward the class $A \in H_2(P)$ to the class $\pi_P(A)\in H_2(M)$ represented by $\{t\} \times S^2 \cong M$ for all $ t\in W_\infty$. Then any point of $\ov{\mc M}_h$ close to the  subset $\ov{\mc M}_\infty$ is projected to a holomorphic curve representing $\pi_P(A)$, which must be a sphere contained in $\{t\} \times S^2$ for $t$ close to $\infty$. Such curves are exactly those ones contained in ${\mc M}_{W_\infty}$. The transversality of ${\mc M}_{W_\infty}$ is standard.
\end{proof}

\subsection{Proof using FOP cycles}

We first equip the moduli spaces with derived orbifold charts.

\begin{lemma}
There exists a compact and stably complex derived orbifold ${\mc D} = ({\mc U}, {\mc E}, {\mc S})$ with an isomorphism ${\mc S}^{-1}(0) \cong \overline{\mathcal{M}}_{0,2} (P, J_P, A)$. Moreover, the evaluation map ${\rm ev}$ extends to a smooth submersion
\beqn
\kappa: {\mc U} \to P.
\eeqn
\end{lemma}

\begin{proof}
This is a reformulation of the results from \cite[Section 5.8]{AMS}. The construction of stable complex structure is standard in Gromov--Witten theory.
\end{proof}

\begin{lemma}\label{lemma75}
Define ${\mc U}_h$ by
\beqn
{\mc U}_h:= \kappa^{-1}( P \times P_h) \subset {\mc U}.
\eeqn
Then ${\mc U}_h \subset {\mc U}$ is a suborbifold with ordinary and complex normal bundle. Similarly, define ${\mc U}_\infty$ by
\beqn
{\mc U}_\infty = \kappa^{-1}( P_\infty \times P_h) \subset {\mc U}.
\eeqn
Then ${\mc U}_\infty \subset {\mc U}_h$ is a suborbifold with ordinary and complex normal bundle.
\end{lemma}

\begin{proof}
This follows from \cite[Corollary 5.36]{AMS}. Note that the submersive property of the evaluation maps is used here.
\end{proof}

Now we introduce more notations. From now on all homology and cohomology are assumed to be in ${\mb Z}$-coefficients. If $Y$ is a closed oriented manifold, $\alpha \in H_k(Y)$, and $S \subset Y$ is a codimension $l$ closed oriented submanifold, then denote by
\beqn
\alpha|_S \in H_{k-l}(S)
\eeqn
the class obtained by intersecting with $S$. Geometrically, if $\alpha$ is represented by a pseudocycle in $Y$ which intersects transversely with $S$ (as well as its frontier), then $\alpha|_S$ is represented by the pseudocycle in $S$ obtained by the intersection. 

If $W, Y$ are closed oriented manifolds, $\alpha \in H_*(W \times Y)$, then there is an associated map $\Psi_\alpha: H^*(W) \to H^*(Y)$ defined as the composition
\beqn
\xymatrix{ H^*(W) \ar[r]^{\alpha / \cdot} & H_*(Y) \ar[r]^{\rm PD} & H^*(Y)}
\eeqn
where $\alpha / \cdot$ is the slant product. If $\iota: Y_0 \to Y$ is the inclusion of a closed 
 oriented submanifold and $\alpha_0 = \alpha|_{W \times Y_0}$, then the following diagram commutes.
 \beqn
 \xymatrix{ &   H^*(Y)  \ar[d]^{\iota^*} \\
  H^*(W) \ar[ru]^{\Psi_\alpha} \ar[r]_{\Psi_{\alpha_0}} & H^*(Y_0) }
 \eeqn

Now consider the stably complex derived orbifold 
$${\mc D}_h = ({\mc U}_h, {\mc E}_h, {\mc S}_h) := ({\mc U}_h, {\mc E}|_{{\mc U}_h}, {\mc S}|_{{\mc U}_h})$$
where ${\mc U}_h$ is constructed in Lemma \ref{lemma75}; up to concordance we may regard ${\mc D}_h$ as oriented and normally complex. Let $k$ be the virtual dimension of $\ov{\mc M}_h$. By choosing an FOP transverse perturbation on ${\mc D}_h$, for the stratum $\theta_k \in {\mf Z}_k^{\rm univ}$, there is an integral homology class
\beqn
\chi_{\theta_k}^{\rm FOP}({\mc D}_h) \in H_k ({\mc U}_h; {\mb Z}),
\eeqn
which is pushed forward by the evaluation map to a class
\beqn
\alpha:= {\rm ev}_* (\chi_{\theta_k}^{\rm FOP}({\mc D}_h)) \in H_k ( P_h \times P; {\mb Z}).
\eeqn
Define
\beqn
\alpha_\infty:= \alpha|_{P_h \times P_\infty} \in H_{k-2}(P_h \times P_\infty; {\mb Z})
\eeqn
and 
\beqn
\alpha_\phi:= \alpha|_{P_h \times P_\phi} \in H_{k-2}(P_h \times P_\phi; {\mb Z}).
\eeqn
Then there are the associated maps $\Psi_\alpha$, $\Psi_{\alpha_\infty}$, and $\Psi_{\alpha_\phi}$ which fit into the following commutative diagram
\beq\label{eqn:diagram}
\vcenter{ \xymatrix{
&    & H^*(P_{\infty}) \ar[rd]^{\iota_X^*} &                    \\
H^*(X) \ar[r]^{{\rm proj}_2^*} & H^*(P_h) \ar[r]^{\Psi_\alpha} \ar[ru]^{\Psi_{\alpha_\infty}} \ar[rd]_{\Psi_{\alpha_\phi}} & H^*(P) \ar[r]^{\iota_X^*} \ar[d]_{\iota_\phi^*} \ar[u]^{\iota_\infty^*} & H^*(X) \\
          &    & H^*(P_{\phi}) \ar[ru]_{\iota_X^*}  &    } }
          \eeq
where ${\rm proj}_2: P_h \to X$ is the projection induced from the identification \eqref{eqn71} and the three $\iota_X^*$ are induced from the inclusion of a fiber. 

\begin{lemma}\label{lemma76} 
\begin{enumerate}

\item Near $\ov{\mc M}_\infty$, the original section ${\mc S}_h: {\mc D}_h \to {\mc E}_h$ is transverse in both classical and the FOP sense.

\item $\alpha|_{P_h \times P_\infty} = {\rm ev}_* [\ov{\mc M}_\infty]$.

\item There holds $\iota_X^* \circ \Psi_{\alpha_\infty} \circ {\rm proj}_2^* = {\rm Id}_{H^*(X)}$.
\end{enumerate}
\end{lemma}

\begin{proof}
By Lemma \ref{lemma73}, the open neighborhood of $\ov{\mc M}_\infty$ in ${\mc M}_{W_\infty} \subset \ov{\mc M}_h$ is transversely cut out. Moreover, it contains no curves with nontrivial automorphism groups. Then by the condition {\bf Classical Transversality} of Theorme \ref{thm11}, ${\mc S}_h$ is transverse in the FOP sense near ${\mc M}_{W_\infty}$.

Next, by the condition {\bf Extension Property I} of Theorem \ref{thm11}, after appropriately shrinking the chart ${\mc D}_h$, one can find an FOP transverse perturbation ${\mc S}_h': {\mc U}_h \to {\mc E}_h$ which agrees with ${\mc S}_h$ near $\ov{\mc M}_\infty$. Then we see
\beqn
{\rm ev} (({\mc S}_h')^{-1}(0)_{\theta_k}) \cap (P_h \times P_\infty) = {\rm ev}({\mc S}_h^{-1}(0)) \cap (P_h \times P_\infty) = {\rm ev}(\ov{\mc M}_\infty)
\eeqn
and the intersection is transverse. Homologically this means $\alpha|_{P_h\times P_\infty} = {\rm ev}_*[\ov{\mc M}_\infty]$. 

Lastly, as the moduli space $\ov{\mc M}_\infty$ is identified with $\{\infty\}\times S^2 \times X$, the relation $\iota_X^* \circ \Psi_{\alpha_\infty} \circ {\rm proj}_2^*$ is obvious.
\end{proof}

By the commutativity of the diagram \eqref{eqn:diagram}, Lemma \ref{lemma76} implies that
\beqn
\iota_X^* \circ \Psi_{\alpha_\phi} \circ {\rm proj}_2^* = {\rm Id}_{H^*(X)}.
\eeqn
So Proposition \ref{prop72}, and hence Theorem \ref{thm:AMS} are proved.

\appendix

\section{Technical Results about Straightened Structures}\label{appendixa}

In this appendix we prove Proposition \ref{metric_existence}, Proposition \ref{prop_straightened_connection}, and Proposition \ref{prop_straightening_stabilization}.

\subsection{Proof of Proposition \ref{metric_existence}}\label{subsectiona1}

We first briefly sketch the idea of the proof. We would like to construct the straightened metric and NC structure inductively, first near the deepest strata of the isotropy stratification of the underlying orbifold, then move up until the top stratum. The model for each construction is the bundle metric and bundle NC structure on a total space of a vector bundle with fiberwise group actions. Hence we first prove below that the metric and NC structure in the model case can be straightened. We then carry out the inductive construction using chartwise constructions and partition-of-unity arguments.

\subsubsection{The straightening model}

\begin{defn}\label{defn_straightening_model}
\begin{enumerate}

\item A {\it straightening model} is a triple $(G, X, F)$ where $G$ is a finite group, $X = (X, g^{TX})$ is a Riemannian manifold, $F = (F, h^F, \nabla^F)$ consists of a vector bundle $F \to X$ with a fiberwise linear $G$-action, a $G$-invairant inner product $h^F$, and an $h^F$-preserving connection $\nabla^F$. 

\item Let $F \to X$ be a real vector bundle. A {\it Hermitian triple} is a triple $(\nabla^F, h^F, I^F)$ where $I^F$ is a complex structure on $F$, $h^F$ is a Hermitian inner product on $F$, and $\nabla^F$ is a Hermitian connection on $F$.

\item A {\it complex straightening model} consists of a straightening model $(G, X, F)$ and a $G$-invariant Hermitian triple $(\nabla^F, h^F, I^F)$.
\end{enumerate}

\end{defn}

Associated to a straightening model there is a Riemannian $G$-manifold $(F, g^{TF})$ where $g^{TF}$ is the bundle metric coupled with the Riemannian metric on the base. Moreover, for a complex straightening model, $I^F$ induces a bundle NC structure ${\bm I}^{TF}$. The inductive proof of Proposition \ref{metric_existence} relies on the following technical lemma.

\begin{prop}\label{propa2}
Let $(G, X, F)$ be a straightening model. Then the bundle metric $g^{TF}$ is straightened. In addition, if $(G, X, F)$ is complex, then the bundle NC structure ${\bm I}^{TF}$ is straightened with respect to $g^{TF}$.
\end{prop}

To prove that the bundle metrics are straightened, we need to verify certain elementary facts about such metrics.

\begin{lemma}\label{lemmaa3}
Let $(X, g^{TX})$ be a Riemannian manifold and $\pi_F: F \to X$ be a smooth real vector bundle equipped with an inner product $h^F$ and an $h^F$-preserving connection $\nabla^F$. Let $g^{TF}$ be the bundle metric on $F$. Then 
\begin{enumerate}
    \item fibers of $F$ are totally geodesic and fiberwise geodesics are straight lines.

    \item Let $p^{\pi_F^* F}: TF \to \pi_F^* F$ be the orthogonal projection onto the vertical subbundle $\pi_F^* F \subset TF$. Then the connection on $TF$ induced from the Levi--Civita connection $\nabla^{TF}$ of $g^{TF}$, which in precise terms reads 
    \beqn
    p^{\pi_F^* F} \circ \nabla^{TF} \circ p^{\pi_F^* F},
    \eeqn
    coincides with the pullback connection $\pi_F^* \nabla^F$. 

\end{enumerate}
\end{lemma}

\begin{proof}
Choose local coordinates $(x^1, \ldots, x^m)$ on $X$ and local orthonormal frame $(e_1, \ldots, e_k)$ of $F$. Write  $\nabla^F$ in terms of the connection matrix $\omega_{i,\alpha}^\beta$ as 
\beqn
\nabla_{\partial_i}^F e_\alpha = \sum_{\beta=1}^k \omega_{i,\alpha}^\beta e_\beta,\ 1 \leq i \leq m,\  1 \leq \alpha \leq k.
\eeqn
Let $(y^1, \ldots, y^k)$ be the fiber coordinates dual to the basis $e_1, \ldots, e_k$. Abbreviate
\begin{align*}
&\ \partial_i = \frac{\partial}{\partial x^i},\ &\ \partial_\alpha = \frac{\partial}{\partial y^\alpha}
\end{align*}
Then by the definition of the bundle metric, one has 
\begin{align}\label{eqna1}
&\ \left\langle \partial_i, \partial_\alpha \right\rangle_{g^{TF}} = \sum_{\gamma=1}^k \omega_{i, \gamma}^\alpha y^\gamma,\ &\ \left\langle \partial_\alpha, \partial_\beta  \right\rangle_{g^{TF}} = \delta_{\alpha\beta}.
\end{align}
The two assertions follows from the identities 
\beq\label{eqna2}
\nabla_{\partial_\alpha}^{TF} \partial_\beta  = 0,
\eeq
\beq\label{eqna3}
\left\langle \nabla_{\partial_i}^{TF} \partial_\alpha, \partial_\beta \right\rangle_{g^{TF}} = \omega_{i,\alpha}^\beta.
\eeq
For \eqref{eqna2}, it is easy to see that for $\gamma = 1, \ldots, k$, $\langle \nabla_{\partial_\alpha}^{TF}\partial_\beta, \partial_\gamma \rangle \equiv 0$. Moreover, by the formula of the Levi--Civita connection and \eqref{eqna1},
\beqn
\begin{split}
\langle \nabla_{\partial_\alpha}^{TF} \partial_\beta, \partial_i \rangle = &\ \frac{1}{2} \left( \partial_\alpha \langle \partial_i, \partial_\beta \rangle_{g^{TF}} - \partial_i \langle \partial_\alpha, \partial_\beta \rangle_{g^{TF}} + \partial_\beta \langle \partial_i, \partial_\alpha \rangle \right)\\
= &\ \frac{1}{2} \left( \partial_\alpha \omega_{i, \beta}^\gamma y^\gamma + \partial_\beta \omega_{i, \alpha}^\gamma y^\gamma \right)\\
= &\ \frac{1}{2} ( \omega_{i, \beta}^\alpha + \omega_{i, \alpha}^\beta )
\end{split}
\eeqn
which vanishes as $(e_1, \ldots, e_k)$ is an orthonormal frame and $\nabla^F$ preserves the inner product. This proves \eqref{eqna2}. For \eqref{eqna3}, using the formula for the Levi--Civita connection again, one has
\beqn
\begin{split} \left\langle \nabla_i^{TF} \partial_\alpha, \partial_\beta \right\rangle_{g^{TF}} = &\ \frac{1}{2} \left( \partial_i \left\langle \partial_\alpha, \partial_\beta \right\rangle_{g^{TF}} + \partial_\alpha \left\langle \partial_i, \partial_\beta  \right\rangle_{g^{TF}} - \partial_\beta \left\langle \partial_i, \partial_\alpha \right\rangle_{g^{TF}} \right)\\ = &\ \frac{1}{2} \Big( \partial_\alpha \sum_{\gamma=1}^k \omega_{i, \gamma}^\beta y^\gamma -   \partial_\beta \sum_{\gamma=1}^k \omega_{i, \gamma}^\alpha y^\gamma  \Big)\\
= &\ \frac{1}{2} \big( \omega_{i,\alpha}^\beta - \omega_{i, \beta}^\alpha \big) \\
= &\ \omega_{i,\alpha}^\beta.
\end{split}
\eeqn
This proves \eqref{eqna3}.
\end{proof}

We also need to consider bundles over bundles and bundle metrics obtained in a two-step construction. Let $F$ resp. $\check F$ be equipped with an inner product $h^F$ resp. $h^{\check F}$ and a metric connection $\nabla^F$ resp. $\nabla^{\check F}$. Consider the total space $E: = F \oplus \check F$, which has the bundle metric $g^{TE}$ determined by $(g^{TX}, h^{F} \oplus h^{\check F}, \nabla^F \oplus \nabla^{\check F})$. Notice that there is a natural identification between the total space of $E$ and the total space of the pullback bundle $\pi_F^* \check F \to F$. The second bundle structure gives another bundle metric. We would like to prove that the two bundle metrics agree.

\begin{lemma}\label{lemmaa4}
\begin{enumerate}

\item The bundle isomorphism
\beqn
TE|_F \cong TF \oplus \pi_F^* \check F.
\eeqn
induced from the identification $E \cong \pi_F^* \check F$ and the connection $\pi_F^* \nabla^{\check F}$ is an isometry between the restriction of the bundle metric $g^{TE}$ to $TE|_F$ and the direct sum $g^{TF} \oplus \pi_F^* h^{\check F}$.

\item Let $p^{\pi_F^* \check F}: TE|_F \to \pi_F^* \check F$ be the orthogonal projection. Then
\beqn
p^{\pi_F^* \check F} \circ \nabla^{TE} \circ \pi^{\pi_F^* \check F} = \pi_F^* \nabla^{\check F}.
\eeqn

\item With respect to the natural identification $E \cong \pi_F^* \check F$ of total spaces, the bundle metric $g^{TE}$ coincides with the bundle metric $g^{T \pi_F^* \check F}$ induced from $(g^{TF}, \pi_F^* h^{\check F}, \pi_F^* \nabla^{\check F})$.
\end{enumerate}
\end{lemma}

\begin{proof}
We would like to write down the identification $E \cong \pi_F^* \check F$ of total spaces in local coordinates. Let $(x^i)_{i=1}^m$ be local coordinates on $X$,  $(e_\alpha)_{\alpha=1}^k$ be local orthonormal basis of $F$ with dual fiber coordinates $y^\alpha$, and $(\check e_\gamma)_{\gamma=1}^l$ be local orthonormal basis of $\check F$ with dual fiber coordinates $\check y^\gamma$. Then the identification is 
\beqn
E_x \ni \Big(x, y^\alpha e_\alpha(x) + \check y^\gamma \check e_\gamma(x) \Big) \mapsto \Big( x, y^\alpha e_\alpha(x), \check y^\gamma \pi_F^* e_\gamma(x) \Big) \in \pi_F^* \check F|_{(x, y^\alpha e_\alpha(x))}.
\eeqn
Hence if we regard $\pi_F^* e_\gamma(x)$ as a local basis of $\pi_F^* \check F$, then using the fiber coordinates induced by the given base coordinates and the chosen fiber bases, the map $E \to \pi_F^* \check F$ is just the identity map in $(x^i, y^\alpha, \check y^\gamma)$.

For (1), one only needs to compare along a single fiber of $F$ at some point $x \in X$. One can choose $e_\alpha$ and $\check e_\gamma$ such that the connection matrices at $x$ vanish. Then in the bundle coordinates of $E$, along the fiber $F_x$, one has 
\beqn
g^{TE} = g_{ij} dx^i \otimes dx^j + \sum_{\alpha} dy^\alpha \otimes dy^\alpha + \sum_\gamma d\check y^\gamma \otimes d\check y^\gamma.
\eeqn
Then (1) follows. For (2), we know from (2) of Lemma \ref{lemmaa3} that the restriction of $\nabla^{TE}$ to $\pi_E^* E$ is the pullback $\pi_E^* \nabla^E$. As $\nabla^E$ is the direct sum, it follows that the restriction of $\nabla^{TE}$ to $\pi_E^* \check F$ is the pullback $\pi_E^* \nabla^{\check F}$. Then (2) follows from by further restricting to the submanifold $F \subset E$. Then (3) follows from the explicit form of the bundle metric on $\pi_F^* \check F \to F$.
\end{proof}

\begin{proof}[Proof of Proposition \ref{propa2}] 
First notice that the bundle metric $g^{TF}$ is $G$-invariant. For each $H \subset_F G$, the $H$-fixed point set of $F$ is the total space of the subbundle $F_H \subset F$. Its normal bundle $NF_H$ is then naturally the pullback $\pi_{F_H}^* \check F_H$, which carries the inner product induced $h^{NF_H}$ from $g^{TF}$ and the connection $\nabla^{NF_H}$ induced from $\nabla^{TF}$. Then by definition of straightenedness (see Definition \ref{defn226}), one needs to show that the normal exponential map
\beqn
\exp^{NF_H}: N^\epsilon F_H \to F
\eeqn
is isometric from the bundle metric on $NF_H \cong \pi_{F_H}^* \check F_H$ determined the triple $(g^{TF_H}, h^{NF_H}, \nabla^{NF_H})$. Indeed, by Lemma \ref{lemmaa3} (1), the normal exponential map is just the linear isomorphism $F_H \oplus \check F_H \cong F$ on each fiber of $F$. Hence hence $\exp^{NF_H}$ is just the natural identification $NF_H \cong \pi_{F_H}^* \check F_H \cong F$. Then by Lemma \ref{lemmaa4} (3), this natural identification is an isometry.

Now assume that $(G, X, F)$ is complex. By Definition \ref{defn227}, one needs to show that, for each $H \subset_F G$, there hold (a) the triple $(\nabla^{NF_H}, h^{NF_H}, I^{NF_H})$ is Hermitian and (b) $\exp^{NF_H}$ preserves the NC structure. For (a), we know from the definition of $I^{\check F_H}$, Lemma \ref{lemmaa3} and Lemma \ref{lemmaa4} that 
\beqn
(\nabla^{NF_H}, h^{NF_H}, I^{NF_H}) = (\pi_{F_H}^* \nabla^{\check F_H}, \pi_{F_H}^* h^{\check F_H}, \pi_{F_H}^* I^{\check F_H}).
\eeqn
As the corresponding triple on $\check F_H \to X$ is Hermitian, so is the pullback triple. For (b), notice that all involved normal bundles are tangent to fibers of $F$. As the exponential map is the linear addition in fibers and the NC structures are all pullbacks on fibers, the exponential map preserves all involved complex structures. 
\end{proof}

\subsubsection{Inductive construction}

To prepare the proof of Proposition \ref{metric_existence}, we first work on the case of a single chart. Let $G$ be a finite group and $U$ be a $G$-manifold. 

\begin{lemma}\label{lemmaa5}
Let $g_0^{TU}$ be a $G$-invariant metric on $U$ and ${\bm I}_0^{TU}$ be an NC structure close to ${\bm I}^{TU}$. Let $H \subset_U G$ be a maximal essential subgroup. Let $h_0^{NU_H}$ be the inner product on $NU_H \to U_H$ induced from $g_0^{TU}$ and $\nabla_0^{NU_H}$ be the connection on $NU_H$ induced from the Levi--Civita connection of $g_0^{TU}$. Given a $G$-invariant closed subset $Y$, suppose the following conditions hold.
\begin{enumerate}

\item Near $Y$, $g_0^{TU}$ is straightened and ${\bm I}_0^{TU}$ is straightened with respect to $g_0^{TU}$.

\item Along $NU_H \to U_H$, the triple $(\nabla_0^{NU_H}, h_0^{NU_H}, I_0^{NU_H})$ is Hermitian.
\end{enumerate}
Then there exist a $G$-invariant metric $g_1^{TU}$ and an NC structure ${\bm I}_1^{TU}$ close to ${\bm I}^{TU}$ such that
\begin{enumerate}

\item Near $Y \cup GU_H$, $g_1^{TU}$ is straightened and ${\bm I}_1^{TU}$ is straightened with respect to $g_1^{TU}$.

\item $(g_1^{TU}, {\bm I}_1^{TU})$ coincides with $(g_0^{TU}, {\bm I}_0^{TU})$ near $Y$.
\end{enumerate}
\end{lemma}

\begin{proof} 
We would like to modify the metric near $GU_H$. Indeed, there is a straightening model $(H, U_H, NU_H)$ where $U_H$ is equipped with the metric restricted from $g_0^{TU}$ and $NU_H$ is equipped with the Hermitian triple $(\nabla_0^{NU_H}, h_0^{NU_H}, I_0^{NU_H})$. Hence by Proposition \ref{propa2}, the total space $NU_H$ is equipped with the straightened metric $g^{TNU_H}$ and the straightened bundle NC structure ${\bm I}^{TNU_H}$. Consider the normal exponential map
\beqn
\exp^{NU_H}: N^\epsilon U_H \to U.
\eeqn
As there is no points fixed by groups larger than $H$, we may choose a $G$-invariant function $\epsilon: GU_H \to {\mb R}_+$ small enough such that for any $H' \subset G$ conjugate to $H$, 
\beqn
|N^\epsilon NU_H| \cap |N^\epsilon NU_{H'}| \neq \emptyset \Longrightarrow H = H'.
\eeqn
Then we equip the neighborhood $|N^\epsilon U_H|$ with the pushforward of the the bundle metric $g^{TNU_H}$ and the bundle NC structure ${\bm I}^{TNU_H}$. Using $G$-action one obtains similar pushforward to other conjugate neighborhoods $|N^\epsilon U_{H'}|$. Denote the new metric by $g_1^{TU}$ and the new NC structure by ${\bm I}_1^{TU}$. 

By our assumption, $g_0^{TU}$ is straightened near $Y$ and ${\bm I}_0^{TU}$ is straightened with respect to $g_0^{TU}$ near $Y$; therefore, by construction, near each point $x \in Y \cap GU_H$, $g_0^{TU} = g_1^{TU}$ and ${\bm I}_0^{TU} = {\bm I}_1^{TU}$. Hence there is a $G$-invariant open neighborhood of $Y \cup GU_H$ such that $g_0^{TU}$ and $g_1^{TU}$ glue together to a straightened metric, still denoted by $g_1^{TU}$; on the same neighborhood, ${\bm I}_0^{TU}$ and ${\bm I}_1^{TU}$ glue together to an NC structure straightened with respect to $g_1^{TU}$, which is still denoted by ${\bm I}_1^{TU}$. Notice that since ${\bm I}_1^{TU} = {\bm I}_0^{TU}$ along $Y \cup GU_H$, we may take the neighborhood sufficiently small so that the new one is still close to ${\bm I}^{TU}$. Then by Lemma \ref{lemma_NC_extension}, one can extend ${\bm I}_1^{TU}$ to the whole $U$ without altering its value near $Y \cup GU_H$ and still denoted by ${\bm I}_1^{TU}$, such that it is still close to ${\bm I}^{TU}$. One then use a cut-off function to extend $g_1^{TU}$ to the whole of $U$ without altering its value near $Y \cup GU_H$ and still denoted by $g_1^{TU}$. 
\end{proof}

To extend this argument to an orbifold, one can first get the following easy corollary. Continue the setting of Lemma \ref{lemmaa5}. For each $d \geq -1$, denote
\beqn
U_d:= \bigcup_{{\rm dim}(U_H) \leq d} U_H \subset U.
\eeqn

\begin{cor}\label{cora6}
Suppose $U_{d-1} = \emptyset$ (which implies that $U_d \subset U$ is a disjoint union of fixed point sets and a $G$-invariant submanifold). Let $NU_d \to U_d$ be the normal bundle which carries an induced inner product $h_0^{NU_d}$, an induced connection $\nabla_0^{NU_d}$, and a complex structure $I_0^{NU_d}$. Suppose  
\begin{enumerate}
\item Near $Y$, $g_0^{TU}$ is straightened and ${\bm I}_0^{TU}$ is straightened with respect to $g_0^{TU}$.

\item Along $NU_d \to U_d$, the triple $(\nabla_0^{NU_d}, h_0^{NU_d}, I_0^{NU_d})$ is Hermitian.
\end{enumerate}
Then there exist a $G$-invariant metric $g_1^{TU}$ and an NC structure ${\bm I}_1^{TU}$ close to ${\bm I}^{TU}$ such that 
\begin{enumerate}

\item Near $Y \cup U_d$, $g_1^{TU}$ is straightened and ${\bm I}_1^{TU}$ is straightened with respect to $g_1^{TU}$.

\item $(g_1^{TU}, {\bm I}_1^{TU})$ coincides with $(g_0^{TU}, {\bm I}_0^{TU})$ near $Y$.
\end{enumerate}
\end{cor}

\begin{proof}
Apply the argument of Lemma \ref{lemmaa5} to all fixed point sets $U_K$ with ${\rm dim}(U_K) = d$ (which are disjoint). 
\end{proof}

\begin{proof}[Proof of Proposition \ref{metric_existence}]
For $d \geq -1$, let ${\mc U}_d \subset {\mc U}$ be the union of strata ${\mc U}_\gamma$ of the isotropy stratification whose dimensions are at most $d$. Then ${\mc U}_d$ is closed. We would like to prove the following statement by induction on $d$:

\vspace{0.1cm}

\noindent {\it Claim I.} For each $d$, there exist a Riemannian metric $g_1^{T{\mc U}}$ and NC structure ${\mc I}_1^{T{\mc U}}$ which is close to ${\mc I}^{T{\mc U}}$ such that the pair $(g_1^{T{\mc U}}, {\mc I}_1^{T{\mc U}})$ is a straightened NC structure near ${\mc Y} \cup {\mc U}_d$ and such that $(g_1^{T{\mc U}}, {\mc I}_1^{T{\mc U}}) = (g_0^{T{\mc U}}, {\mc I}_0^{T{\mc U}})$ near ${\mc Y}$.

\vspace{0.1cm}

Notice that the proposition is equivalent to Claim I for $d = {\rm dim}({\mc U})$ case. 

\vspace{0.1cm}

\noindent {\it Claim II.} Suppose Claim I holds for $d-1$ as well for $d$ under the special restriction ${\mc U}_{d-1} = \emptyset$, then Claim I holds for $d$ without this restriction.

\vspace{0.1cm}

\noindent {\it Proof of Claim II.} Let ${\mc U}_{d-1}^+ \subset {\mc U}$ be a closed neighborhood of ${\mc U}_{d-1}$ such that $(g_0^{T{\mc U}}, {\mc I}_0^{T{\mc U}})$ is a straightened NC structure near ${\mc Y} \cup {\mc U}_{d-1}^+$. Then define $\dot {\mc U}:= {\mc U} \setminus {\mc U}_{d-1}$ and $\dot {\mc Y}:= ( {\mc Y} \cup {\mc U}_{d-1}^+ ) \setminus {\mc U}_{d-1}$. We have $\dot {\mc U}_{d-1} = \emptyset$. Then there exists a pair $(\dot g_1^{T{\mc U}}, \dot {\mc I}_1^{T{\mc U}})$ defined on $\dot {\mc U}$ which is a straightened NC structure near $\dot {\mc Y} \cup \dot {\mc U}_d$ and agrees with $(g_0^{T{\mc U}}, {\mc I}_0^{T{\mc U}})$ such that $\dot {\mc I}_1^{T{\mc U}}$ is close to ${\mc I}^{T{\mc U}}$. Together with $(g_0^{T{\mc U}}, {\mc I}_0^{T{\mc U}})$ near ${\mc U}_{d-1}$, one obtains a pair $(g_1^{T{\mc U}}, {\mc I}_1^{T{\mc U}})$ defined on ${\mc U}$ which satisfies the requirement. \qed

\vspace{0.1cm}

Now we heads towards proving Claim I for $d$ under the special restriction ${\mc U}_{d-1} = \emptyset$. We first modify the normal inner products and connections along all $d$-dimensional strata. For each chart $C = (G, U, \psi)$, denote $U_d:= \psi^{-1}({\mc U}_d)$. As ${\mc U}_{d-1} = \emptyset$, it follows that $U_d$ is a closed $G$-invariant submanifold of $U$. Moreover, if $\iota: (G', U', \psi') \to (G, U, \psi)$ is a chart embedding, then $\iota: U_d' \to U_d$ is an open embedding inducing an equivariant bundle map $NU_d' \to NU_d$. Let $NU_d \to U_d$ be the normal bundle, which carries a complex structure $I_0^{NU_d}$, an inner product $h_0^{NU_d}$, and an $h_0^{NU_d}$-preserving connection $\nabla_0^{NU_d}$ such that near $\psi^{-1}({\mc Y}) \cap U_d$, $h_0^{NU_d}$ and $\nabla_0^{NU_d}$ are Hermitian. 

\vspace{0.1cm}

\noindent {\it Claim III.} For all charts $C = (G, U, \psi)$, there exist a $G$-invariant Hermitian inner product $h_1^{NU_d}$ on $NU_d$ such that near $\psi^{-1}({\mc Y}) \cap U_d$, $h_1^{NU_d} = h_0^{NU_d}$, and such that the collection $h_1^{NU_d}$ is invariant under chart embeddings. 

\vspace{0.1cm}

\noindent {\it Proof of Claim III.} We first construct $h_1^{NU_d}$ on each chart. Indeed, take an arbitrary Hermitian inner product $h_1^{NU_d}$ on $NU_d$ which coincides with $h_0^{NU_d}$ near $\psi^{-1}({\mc Y}) \cap U_d$, then take average over $G$. By abuse of notation, let $h_1^{NU_d}$ denote the $G$-invariant average. Then choose a countable cover of ${\mc U}$ by charts $C_i = (G_i, U_i, \psi_i)$ and a subordinate smooth partition of unity $\rho_i: {\mc U} \to [0, 1]$. For each chart $C = (G, U, \psi)$, denote $\rho_i^C:= \rho_i \circ \psi$. If ${\rm supp} \rho_i^C \neq \emptyset$, then $\rho_i^C|_{U_d} h_1^{NU_{i, d}}$ is a Hermitian 2-tensor on the bundle $NU_d$. Then
\beqn
\sum_{i=1}^\infty \rho_i^C|_{U_d} h_1^{NU_{i, d}}
\eeqn
is a Hermitian inner product on $NU_d$ which respects chart embedding.  \qed

\vspace{0.1cm}

\noindent {\it Claim IV.} For each chart $C = (G, U, \psi)$, there exists a $G$-invariant connection $\nabla^{NU_d}$ on $NU_d \to U_d$ which is Hermitian with respect to $I_0^{NU_d}$ and $h_1^{NU_d}$. 

\vspace{0.1cm}

\noindent {\it Proof of Claim IV.} It is similar to the previous case. We first construct the Hermitian connection on a chart $C_i$ and then use a partition of unity to glue together. \qed

\vspace{0.1cm}

Next we construct the new metric and NC structure. Choose a sufficiently small smooth function $\epsilon: {\mc U}_d \to {\mb R}_+$. By abuse of notation, let $\epsilon: U_d \to {\mb R}_+$ the pullback to any chart $C = (G, U, \psi)$. Then the exponential map $\exp^{NU_d}: N^\epsilon U_d \to U$ is an open embedding. Replace $g_0^{TU}|_{|N^\epsilon U_d|}$ by the bundle metric on $N^\epsilon U_d$ determined by $(g_0^{TU_d}, h_1^{NU_d}, \nabla_1^{NU_d})$ and replace the NC structure ${\bm I}_0^{TU}|_{|N^\epsilon U_d|}$ by the bundle NC structure. Denote $|N^\epsilon {\mc U}_d| \subset {\mc U}$ the union of $\psi(|N^\epsilon U_d|)$. Then we obtain a Riemannian metric $g_1^{T{\mc U}}$ and an NC structure ${\mc I}_1^{T{\mc U}}$ defined in $|N^\epsilon {\mc U}_d|$. If $\epsilon$ is sufficiently small, one can guarantee that ${\mc I}_1^{T{\mc U}}$ is close to ${\mc I}^{T{\mc U}}$. One can check that since $(g_0^{T{\mc U}}, {\mc I}_0^{T{\mc U}})$ is straightened near ${\mc Y}$, the new pair $(g_1^{T{\mc U}}, {\mc I}_1^{T{\mc U}})$ coincides with $(g_0^{T{\mc U}}, {\mc I}_0^{T{\mc U}})$ near ${\mc Y} \cap |N^\epsilon {\mc U}_d|$. Then the old and new pairs determines a straightened pair near the closed set ${\mc Y}\cup {\mc U}_d$. By Lemma \ref{lemma_NC_extension}, one can extend ${\mc I}_1^{T{\mc U}}$ to the whose ${\mc U}$ without altering its current value near ${\mc Y} \cup {\mc U}_d$. Then using a cut-off function to extend $g_1^{T{\mc U}}$ to the whole ${\mc U}$ without altering its current value near ${\mc Y}\cup {\mc U}_d$. This proves Claim I for $d$ under the assumption ${\mc U}_{d-1} = \emptyset$. By Claim II, this finishes the proof of the proposition.
\end{proof}

\subsection{Proof of Proposition \ref{prop_straightened_connection}}\label{subsectiona2}

The strategy of proving Proposition \ref{prop_straightened_connection} is similar to that of Proposition \ref{metric_existence}, i.e., inductive construction of straightened structures near strata of the isotropy stratification, from the deepest up. 

\subsubsection{Bundle straightening model}

\begin{defn}\label{defna7} \hfill
\begin{enumerate}

\item A {\it bundle straightening model} is a quadruple $(G, X, F, E)$ where $(G, X, F)$ is a straightening model (Definition \ref{defn_straightening_model}) and $E = (E, h^E, \nabla^E, {\bm I}^E)$ consists of a $G$-equivariant vector bundle $E \to X$ (where $X$ has the trivial $G$-action), $h^E$ is an invariant inner product, $\nabla^E$ is an $h^E$-preserving $G$-invariant connection, and $I^{\check E_G}$ is a $G$-invariant complex structure on $\check E_G \subset E$ such that the triple $(\nabla^{\check E_G}, h^{\check E_G}, I^{\check E_G})$ is Hermitian.

\item A {\it complex bundle straightening model} is a straightening model $(G, X, F, E)$ with $E$ equipped with a $G$-invariant complex structure $I^E$, such that the triple $(\nabla^E, h^E, I^E)$ is Hermitian and such that ${\bm I}^E$ is induced from $I^E$.
\end{enumerate}
\end{defn}

Associated to a bundle straightening model $(G, X, F, E)$ there are the following objects. On the pullback bundle $\pi_F^* E \to F$ there is the pullback connection $\pi_F^* \nabla^E$, the pullback inner product $\pi_F^* h^E$, and the pullback NC structure ${\bm I}^{\pi_F^* E}$. 

\begin{prop}\label{propa8} \hfill
\begin{enumerate}

\item The pullback connection $\pi_F^* \nabla^E$ on $\pi_F^* E \to F$ is straightened with respect to the bundle metric $g^{TF}$. 

\item The pullback inner product $\pi_F^* h^E$ on $\pi_F^* E$ is straightened with respect to $g^{TF}$ and $\pi_F^* \nabla^E$.

\item The pullback NC structure ${\bm I}^{\pi_F^* E}$ on $\pi_F^* E$ is straightened with respect to $g^{TF}$ and $\pi_F^* \nabla^E$.
\end{enumerate}
\end{prop}

\begin{proof}
Choose $H \subset_F G$, whose fixed point set in the total space of $F$ is the total space of $F_H$. Along each geodesic normal to $F_H$ contained in the fiber $F_x\subset F$, the parallel transport of $\pi_F^* E$ with respect to $\pi_F^* \nabla^E$ is simply the identity map of $E_x$. As all involved structures on $\pi_F^* E$ are pullbacks, they are naturally straightened.
\end{proof}

There is a way to obtain a bundle straightening model by looking at restrictions at fixed point sets. Let $(U, g^{TU})$ be a straightened Riemannian $G$-manifold, $E \to U$ be a $G$-equivariant vector bundle equipped with $G$-invariant inner product $h^E$, connection $\nabla^E$, and NC structure ${\bm I}^E$. For $H \subset_U G$, if $(h^{\check E_H}, \nabla^{\check E_H}, I^{\check E_H})$ is a Hermitian triple, then the quadruple 
\beqn
(H, U_H, NU_H, E|_{U_H})
\eeqn
with the naturally induced structures defines a bundle straightening model.

\subsubsection{The inductive construction}
 
Analogous to Lemma \ref{lemmaa5} and Corollary \ref{cora6}, one needs the following lemma for the inductive construction. Recall that $U_d \subset U$ is the union of fixed point loci with dimension at most $d$. Assume $U_{d-1} = \emptyset$. Then $U_d$ is a $G$-invariant closed submanifold. One can decompose
\beqn
E|_{U_d} = E_d \oplus \check E_d
\eeqn
although $E_d$ and $\check E_d$ may have varying dimensions on different components of $U_d$. Then connections, inner products, and NC structures induce corresponding structures on $\check E_d$.

\begin{lemma}\label{lemma_connection_induction}
Let $U$ be a $G$-manifold equipped with a straightened $G$-invariant Riemannian metric $g^{TU}$. Suppose $d \geq 0$ and $U_{d-1} = \emptyset$. Let $Y \subset U$ be a $G$-invariant closed subset.

Let $E \to U$ be a $G$-equivariant vector bundle equipped with an NC structure ${\bm I}^E$. Let $\nabla_0^E$ be a $G$-invariant connection on $E$, $h_0^E$ be a $G$-invariant inner product on $E$, and ${\bm I}_0^E$ is an NC structure on $E$ which is close to ${\bm I}^E$. Suppose
\begin{enumerate}
\item Near $Y$, $\nabla_0^E$ is straightened with respect to $g^{TU}$, $h_0^E$ is preserved by $\nabla_0^E$, and ${\bm I}_0^E$ is straightened with respect to $g^{TU}$ and $\nabla_0^E$.

\item Along $\check E_d \to U_d$, the triple $(\nabla_0^{\check E_d}, h_0^{\check E_d}, I_0^{\check E_d})$ is Hermitian.
\end{enumerate}
Then there exists a $G$-invariant connection $\nabla_1^E$, a $G$-invariant inner product $h_1^E$, and an NC structure ${\bm I}_1^E$ which is close to ${\bm I}^E$ such that 

\begin{enumerate}
\item Near $Y \cup U_d$, $\nabla_1^E$ is straightened with respect to $g^{TU}$, $h_1^E$ is preserved by $\nabla_1^E$, and ${\bm I}_1^E$ is straightened with respect to $g^{TU}$ and $\nabla_1^E$.

\item $(\nabla_1^E, h_1^E, {\bm I}_1^E)$ coincides with $(\nabla_0^E, h_0^E, {\bm I}_0^E)$ near $Y$.
\end{enumerate}
In addition, if $E$ is a complex vector bundle, $h_0^E$ and $\nabla_0^E$ are Hermitian, and ${\bm I}^E$ is induced from the complex structure, and ${\bm I}_0^E = {\bm I}^E$, then we can make ${\bm I}_1^E = {\bm I}_0^E = {\bm I}^E$, $h_1^E = h_0^E$, and $\nabla_1^E$ be Hermitian.
\end{lemma}

\begin{proof}
It suffices to consider the case that $U_d = GU_H$ for an essential subgroup $H$. Consider the quadruple $(H, U_H, NU_H, E|_{U_H})$ where $(H, U_H, NU_H)$ is the induced straightening model with a straightened metric $g^{TNU_H}$ on the total space $NU_H$ and $E|_{U_H}$ is equipped with the induced connection $\nabla_0^{E|_{U_H}}$, the induced inner product $h_0^{E|_{U_H}}$, and a complex structure $I^{\check E_H}$ on $\check E_H \subset E|_{U_H}$. By our assumption, the quadruple $(H, U_H, NU_H, E|_{U_H})$ is a bundle straightening model (Definition \ref{defna7}). Hence by Proposition \ref{propa8}, the pullback bundle $\pi_{NU_H}^* (E|_{U_H})$, equipped with the pullback triple $(\pi_{NU_H}^* \nabla^{E|_{U_H}}, \pi_{NU_H}^* h^{E|_{U_H}}, \pi_{NU_H}^* {\bm I}^{E|_{U_H}})$, is straightened with respect to the bundle metric $g^{TNU_H}$. 

Now let $\epsilon: GU_H \to {\mb R}_+$ be a small enough $G$-invariant function. Using the bundle isomorphism $\pi_{NU_H}^* (E|_{U_H})|_{N^\epsilon U_H} \to E|_{|N^\epsilon U_H|}$ determined by the connection, one obtains a straightened triple $(\nabla_1^E, h_1^E, {\bm I}_1^E)$ on $E$ near $U_d = GU_H$. By the same argument as Lemma \ref{lemmaa5}, this triple agrees near $Y \cap U_d$ with the old one $(\nabla_0^E, h_0^E, {\bm I}_0^E)$, together with which they define a triple which is straightened near $Y \cup U_d$ and which agrees with the old triple near $Y$.

Lastly, suppose $E$ is complex with induced NC structure ${\bm I}^E$, $h_0^E$ and $\nabla_0^E$ are Hermitian, and ${\bm I}_0^E = {\bm I}^E$. Then the bundle isomorphism
\beqn
\pi_{NU_H}^* (E|_{U_H})|_{N^\epsilon U_H} \to E|_{|N^\epsilon U_H|}
\eeqn
is unitary, hence already preserves the NC structures and the inner product. Hence in the procedure of constructing the new triple, one does not need to modify $h_0^E$ and ${\bm I}_0^E$. Moreover, the pullback connection $\nabla_1^E$ is Hermitian.
\end{proof}



\begin{proof}[Proof of Proposition \ref{prop_straightened_connection}] The inductive argument is very similar to the proof of Proposition \ref{metric_existence}. Let ${\mc U}_d \subset {\mc U}$ be the union of strata of the isotropy stratification whose dimensions are at most $d$. Using the same argument as in the proof of Proposition \ref{metric_existence}, it suffices to prove the following claim.

\vspace{0.1cm}

\noindent {\it Claim I.} Suppose ${\mc U}_{d-1} = \emptyset$. Then there exists a triple $(\nabla_1^{\mc E}, h_1^{\mc E}, {\bm I}_1^{\mc E})$ which satisfies the requirement of Proposition \ref{prop_straightened_connection} near ${\mc Y} \cup {\mc U}_d$ and which agrees with the existing triple near ${\mc Y}$.

\vspace{0.1cm}

We first need to modify the inner product and connection in the normal direction of ${\mc U}_d$ to make them Hermitian. 

\vspace{0.1cm}

\noindent {\it Claim II.} There exists a $G$-invariant Hermitian inner product $h_1^{\mc E}$ on ${\mc E}$ which agrees with $h_0^{\mc E}$ near ${\mc Y}$ such that for each chart $\hat C = (G, U, E, \hat \psi)$, the restriction of $h_1^{\mc E}$ to $\check E_d \to U_d$ is Hermitian with respect to $I_0^{\check E_d}$ 

\vspace{0.1cm}

\noindent {\it Proof of Claim II.} We can first modify on each chart. Indeed, let $h_0^{\check E_d}$ be the restriction of $h_0^E$ to $\check E_d$. Then it is Hermitian near $\psi^{-1}({\mc Y}) \cap U_d$. One can find $h_1^{\check E_d}$ which is Hermitian and which agrees with $h_0^{\check E_d}$ near $\psi^{-1}({\mc Y}) \cap U_d$. Then one extend $h_1^{\check E_d}$ to an inner product $h_1^E$ on $E$ which agrees with $h_0^E$ near $\psi^{-1}({\mc Y})$. Then using a partition of unity argument, one can find the desired ones. \qed

\vspace{0.1cm}

\noindent {\it Claim III.} For each chart $\hat C = (G, U, E, \hat\psi)$, there exists a $G$-invariant connection $\nabla^{\check E_d}$ on $\check E_d \to U_d$ which is Hermitian with respect to $I_0^{\check E_d}$ and $h_1^{\check E_d}$. 

\vspace{0.1cm}

\noindent {\it Proof of Claim III.} It is similar to the proof of the above claim. Indeed, on each chart $\hat C = (G, U, E, \hat\psi)$, one first modifies $\nabla_0^{\check E_d}$ to obtain a Hermitian connection $\nabla_1^{\check E_d}$ which agrees with $\nabla_0^{\check E_d}$ near $\psi^{-1}({\mc Y}) \cap U_d$, and then find an $h_1^E$-preserving connection $\nabla_1^E$ extending $\nabla_1^{\check E_d}$ which agrees with $\nabla_0^E$ near $\psi^{-1}({\mc Y})$. Then using a partition of unity argument to find the desired collection of connections. \qed

\vspace{0.1cm}

Then we modify the triple $(\nabla_1^E, h_1^E, I_0^E)$ in a neighborhood of $U_d$ using the bundle isomorphism
\beqn
(\pi_{NU_d}^* (E|_{U_d}))|_{N^\epsilon U_d} \to E|_{|N^\epsilon U_d|}.
\eeqn
induced from the connection $\nabla_1^E$. By abuse of notation, let the new triple be $(\nabla_1^E, h_1^E, {\bm I}_1^E)$, which together with the old triple $(\nabla_0^E, h_0^E, {\bm I}_0^E)$, defines a straightened triple near ${\mc Y} \cup {\mc U}_d$. 

Lastly, suppose ${\mc E}$ is a complex vector bundle, $h_0^{\mc E}$ and $\nabla_0^{\mc E}$ are Hermitian, ${\mc I}_0^{\mc E} = {\mc I}^{\mc E}$ is induced from the complex structure, then the constructions covered by Claim II and Claim III are uncessary. Moreover, the chartwise bundle isomorphism $\Phi^E$ near $U_d$ preserves the NC structure. Therefore one does not need to modify the NC structure and the inner product on ${\mc E}$ during the construction.
\end{proof}

\subsection{Proof of Proposition \ref{prop_straightening_stabilization}}\label{subsectiona3}

One considers the case on a single chart first.

\begin{lemma}\label{lemmaa10}
Let $(U, g^{TU})$ be a straightened Riemannian $G$-manifold. Let $F \to U$ be a $G$-equivariant vector bundle equipped with an inner product $h^F$ and an $h^F$-preserving connection $\nabla^F$. Suppose $\nabla^F$ is straightened with respect to $g^{TU}$ and $h^F$ is straightened with respect to $g^{TU}$ and $\nabla^F$. Then the bundle metric $g^{TF}$ determined by $(g^{TU}, \nabla^F, h^F)$ on the total space of $F$ is straightened.
\end{lemma}

\begin{proof}
Let $H \subset_F G$ be an essential subgroup whose fixed point set in $F$ is the total space of $F_H \to U_H$. As $g^{TU}$ is straightened, we may identify a neighborhood of $U_H$ isometrically with a disk bundle $N^\epsilon U_H$. As $\nabla^F$ and $h^F$ are straightened, we may identify $F|_{|N^\epsilon U_H|}$ with the pullback of $F|_{U_H}$ with the pullback connection and inner products. Then near $U_H$, the total space $F$ can be identified naturally with other total spaces
\beqn
F|_{|N^\epsilon U_H|} \cong NU_H \oplus F_H \oplus \check F_H \cong  \pi_{F_H}^* (NU_H \oplus \check F_H) \cong         NF_H
\eeqn
By Lemma \ref{lemmaa4} (3), bundle metrics obtained from different perspectives are identical. Hence $g^{TF} = g^{TNF_H}$, meaning $g^{TF}$ is straightened along $NF_H$. 
\end{proof}

Next we treat the straightenedness of bundle NC structures. 

\begin{lemma}\label{lemmaa11}
Let $(U, g^{TU})$ be a straightened Riemannian $G$-manifold and ${\bm I}^{TU}$ is a straightened NC structure. Let $F \to U$ be a $G$-equivariant vector bundle equipped with an NC structure ${\bm I}^F$, a Hermitian metric $h^F$, a Hermitian connection $\nabla^F$ such that $\nabla^F$ is straightened with respect to $g^{TU}$ and ${\bm I}^F$ is straightened with respect to $g^{TU}$ and $\nabla^F$.
Then the bundle NC structure ${\bm I}^{TF}$ induced from ${\bm I}^{TU}$, ${\bm I}^F$, and $\nabla^F$ is straightened with respect to the bundle metric $g^{TF}$.
\end{lemma}

\begin{proof}
Using notations from the previous proof, we know near $U_H$, $F$ is naturally identified with the total space of $NU_H \oplus F_H \oplus \check F_H \to U_H$ and the normal exponential map along $NF_H$ is the same as linear additions in fibers of the triple direct sum. As the complex structure $I^{NF_H}$ is the same as the pullback of the direct sum $I^{NU_H} \oplus I^{\check F_H}$ to $F_H$, the normal exponential map preserves the NC structure. 
\end{proof}

Now we turn to straightened structures on pullback vector bundles. 

\begin{lemma}\label{lemmaa12}
Let $(U, g^{TU})$ and $(F, \nabla^F, h^F)$ be as in Lemma \ref{lemmaa10}. Let $E \to U$ be a $G$-equivariant vector bundle equipped with a connection $\nabla^E$ straightend with respect to $g^{TU}$, an inner product $h^E$ preserved by $\nabla^E$, and an NC structure ${\bm I}^E$ straightened with respect to $g^{TU}$ and $\nabla^E$ such that $(\nabla^E, h^E, {\bm I}^E)$ is normally Hermitian. Then 
\begin{enumerate}

\item $\pi_F^* \nabla^E$ is straightened with respect to the bundle metric $g^{TF}$.

\item $\pi_F^* h^E$ is preserved by $\pi_F^* \nabla^E$ and $\pi_F^* {\bm I}^E$ is straightened with respect to $g^{TF}$ and $\pi_F^* \nabla^E$.

\item $(\pi_F^* \nabla^E, \pi_F^* h^E, \pi_F^* {\bm I}^E)$ is normally Hermitian.
\end{enumerate}
\end{lemma}

\begin{proof}
For any essential subgroup $H \subset_F G$, one has the identification of the normal bundle 
\beqn
NF_H \cong \pi_{F_H}^* NU_H \oplus \pi_{F_H}^* \check F_H
\eeqn
which has the equipped bundle metric $g^{TNF_H}$ which agrees with $g^{TF}$ in a tubular neighborhood. As the normal geodesics are straight line segments in fibers of $NF_H$, it is easy to see that $\pi_F^* \nabla^E$ is just the pullback from the restriction to $F_H$. The straightenedness of other structures folllows similarly. 
\end{proof}

\begin{proof}[Proof of Proposition \ref{prop_straightening_stabilization}] The straightenedness of $g^{T{\mc F}}$, ${\mc I}^{T{\mc F}}$, $\pi_{\mc F}^* \nabla^{\mc E}$, $\pi_{\mc F}^* h^{\mc E}$, $\pi_{\mc F}^* {\mc I}^{\mc E}$ can be verified chartwise. These verifications follow from Lemma \ref{lemmaa10}, \ref{lemmaa11}, and \ref{lemmaa12}.
\end{proof}

\section{Whitney Stratifications}\label{appendixb}

In this appendix we prove the existence of a canonical Whitney stratification of complex analytic sets relative to a given stratification (Proposition \ref{thm_nice_Whitney}), properties of such Whitney stratifications (Proposition \ref{propb13} and Proposition \ref{prop_pullback_nice}) as well as Proposition \ref{prop_product_nice}.

\subsection{Most basic results}

Recall the definition of Whitney stratification (Definition \ref{defn31}). A first basic fact is that the partial order relation among strata of a Whitney stratification is compatible with dimensions.

\begin{lemma}\cite[(1.1)]{Topological_stability}\label{lemmab1}
Let $Z_\alpha, Z_\beta$ be two strata of a Whitney stratification on $Z \subset M$. If  $Z_\alpha < Z_\beta$, then ${\rm dim}(Z_\alpha) < {\rm dim}(Z_\beta)$. 
\end{lemma}



Next we consider the pullback operation on Whitney stratifications. 

\begin{lemma}\label{lemma_pullback_Whitney}
Let $Z \subset M$ be equipped with a Whitney stratification ${\mf Z}$ and $f: N \to M$ be a smooth map transverse to ${\mf Z}$, then the partition
\beqn
f^*{\mf Z}:= \Big\{ f^{-1}(Z_\alpha) \neq \emptyset\ |\ Z_\alpha \in {\mf Z} \Big\}
\eeqn
is a Whitney stratification on $f^{-1}(Z) \subset N$.
\end{lemma}

\begin{proof}
The only non-obvious condition for $f^*{\mf Z}$ being a Whitney stratification is the axiom of frontier. Namely, we need to prove
\beqn
f^{-1}(Z_\alpha) \cap \ov{f^{-1}(Z_\beta)} \neq \emptyset \Longrightarrow f^{-1}(Z_\alpha) \subset \ov{f^{-1}(Z_\beta)}.
\eeqn
Suppose $f^{-1}(Z_\alpha) \cap \ov{f^{-1}(Z_\beta)} \neq \emptyset$, then there exist a point $x \in f^{-1}(Z_\alpha)$ and a sequence of points $x_i \in f^{-1}(Z_\beta)$ such that $x_i \to x$. Then $Z_\beta \ni f(x_i) \to f(x) \in Z_\alpha$. Therefore, $f(x) \in Z_\alpha \cap \ov{Z_\beta} \neq \emptyset$. Hence $Z_\alpha < Z_\beta$. On the other hand, $f$ is transverse to ${\mf Z}$ is equivalent to that ${\rm graph}(f) \subset N \times M$ is transverse to each $N \times Z_\alpha$. Moreover, $\{ N \times Z_\alpha\ |\ Z_\alpha \in {\mf Z}\}$ is a Whintey stratification of $N \times Z \subset N \times M$.  Then by \cite[Corollary 10.4]{Mather_1970},
\beqn
{\rm graph}(f) \cap (N \times Z_\alpha) \subset \ov{ {\rm graph} (f) \cap (N \times Z_\beta) }
 \Longrightarrow f^{-1}(Z_\alpha) \subset \ov{f^{-1}(Z_\beta)}. \qedhere
 \eeqn
\end{proof}

It is often helpful to refine a Whitney stratification by taking connected components of all strata. The following lemma guarantees that one still obtains a Whitney stratification.

\begin{lemma}\cite[Proposition 8.7]{Mather_1973}\cite[Corollary 10.5]{Mather_1970}\label{lemmab3}
Let $M$ be a smooth manifold and ${\mf Z}$ be a locally finite partition of a closed subset $Z \subset M$ into locally closed and connected submanifolds such that each  disjoint pair $(Z_\alpha, Z_\beta)$ of ${\mf Z}$ is Whitney regular, then ${\mf Z}$ is a Whitney stratification of $Z$.
\end{lemma}

The following lemma is also frequently used in this paper.



\begin{lemma}\label{lemmab4}
Let $X, M, N$ be smooth manifolds and $h: X \to M$ and $f: M \to N$ be smooth maps. Let $Z \subseteq N$ be a closed subset equipped with a Whitney stratification ${\mf Z}$. Then $f \circ h$ is transverse to ${\mf Z}$ if and only if $f$ is transverse to ${\mf Z}$ at points in ${\rm Im}(h)$. Moreover, suppose $f$ is transverse to ${\mf Z}$ hence pulls back a Whitney stratification $f^* {\mf Z}$ on $f^{-1}(Z)$. Then $h$ is transverse to $f^*{\mf Z}$ if and only if $f \circ h$ is transverse to ${\mf Z}$. 
\end{lemma}

\begin{proof}
The argument is the same as if $Z$ is a smooth closed submanifold and follows from basic linear algebra.
\end{proof}

\subsubsection{Invariance properties of minimal Whitney stratifications}

The purpose of this part is to show certain functorial properties of minimal Whitney stratifications. We do not only consider absolute minimal Whitney stratifications, but also those which are minimal in a restricted class. This is because the Whitney stratification we need in our applications must respect, in a certain sense, the group action on the ambient space (see Definition \ref{defn_action_stratification}). 

\begin{rem}
As the following lemma is rather technical, we would like to try to provide some motivations. Indeed, Lemma \ref{lemmab5} serves as a generalization of Proposition \ref{prop_pullback_nice} which requires more setup in the complex analytic category. Both results deal with the following situation. Let $f: M \to N$ be an embedding of manifolds. Suppose $Z \subset N$ is equipped with a minimal Whitney stratification ${\mf Z}$ and $f$ is transverse to ${\mf Z}$, then $f^* {\mf Z}$ is a Whitney stratification on $Y:= f^{-1}(Z)$ (see Figure \ref{fig:Whitney}). However, we would prefer to have that $f^* {\mf Z}$ is (equivalent to) a  minimal Whitney stratification ${\mf Y}$ (which exists for other reasons), which may not be true in general. This is the typical scenario we have in Section \ref{section3}. In order to show the minimality of $f^* {\mf Z}$, we typically construct a submersive left inverse $g: N \to M$. In the complex analytic category, when the left inverse $g$ is a submersion, one can directly compare $g^* {\mf Y}$ with ${\mf Z}$ and then derive $f^* {\mf Z} \equiv {\mf Y}$. However, there is one case, i.e. in the proof of Proposition \ref{prop314}, where the construction only provides a smooth map $g$ which is not a left inverse either. In that senario, we need to appeal to Lemma \ref{lemmab5} to compare the involved Whitney stratifications. Lastly, we would like to mention that the idea comes from Parker \cite{BParker_integer}.

\begin{figure}[h]
    \centering
    \includegraphics[scale = 1.2]{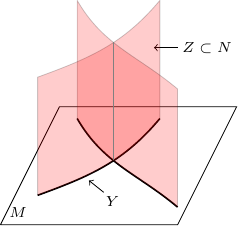}
    \caption{Comparing Whitney stratifications}
    \label{fig:Whitney}
\end{figure}
\end{rem}

Let $M$ be a smooth manifold, $Z \subset M$ be a subset, and $\mc{WS}^\infty(Z)$ be the set of all smooth Whitney stratifications on $Z$. Recall that one can compare different Whitney stratifications using the dimension filtration (Definition \ref{defn_minimal_Whitney}).

\begin{lemma}\label{lemmab5}
Let $M$, $N$ be smooth manifolds and $Y \subset M$, $Z \subset N$ be subsets. Let ${\mc N}(Y) \subset \mc{WS}^\infty(Y)$ resp. ${\mc N}(Z) \subset \mc{WS}^\infty(Z)$ be a subset of smooth Whitney stratifications which is closed under the equivalence relation $\equiv$. Suppose ${\mf Y}\in \mc{N}(Y)$ resp. ${\mf Z} \in \mc{N}(Z)$ is a minimal one. 

Let $f: M \to N$ and $g: N \to M$ be smooth maps satisfying the following.
\begin{enumerate}

\item $f^{-1}(Z) = Y$ and $g^{-1}(Y) = Z$.

\item $f$ resp. $g$ is transverse to ${\mf Z}$ resp. ${\mf Y}$.

\item $f^* {\mf Z} \in {\mc N}(Y)$ and $g^* {\mf Y} \in {\mc N}(Z)$.

\item $g \circ f$ is transverse to ${\mf Y}$ and $(g\circ f)^* {\mf Y} \equiv {\mf Y}$.

\item $f \circ g$ is transverse to ${\mf Z}$ and $(f\circ g)^* {\mf Z} \equiv {\mf Z}$.
\end{enumerate}
Then $f^* {\mf Z} \equiv {\mf Y}$ and $g^* {\mf Y} \equiv {\mf Z}$. 
\end{lemma}

\begin{proof}
Denote $\check{\mf Y}:= f^* {\mf Z}$ and $\check {\mf Z}:= g^* {\mf Y}$. By the minimality of ${\mf Y}$ resp. ${\mf Z}$, one has 
\begin{align*}
&\ {\mf Y} \leq \check{\mf Y},\ &\ {\mf Z} \leq \check{\mf Z}.
\end{align*}
Denote $m = {\rm dim}(M)$, $n = {\rm dim}(N)$.  We prove inductively that for all $l \geq 0$
\begin{align*}
&\ {\mf Y}_{m-l} = \check{\mf Y}_{m-l},\ &\ {\mf Z}_{n-l} = \check {\mf Z}_{n-l}.
\end{align*}
The $l=0$ case is automatically true. Suppose the above is true for all $l<k$. We first show that ${\mf Y}_{m-k} = \check {\mf Y}_{m-k}$. Suppose on the contrary that ${\mf Y}_{m-k} \neq \check {\mf Y}_{m-k}$. Then as ${\mf Y} \leq \check{\mf Y}$, one must have ${\mf Y}_{m-k} \subsetneq \check{\mf Y}_{m-k}$. Hence there exists a point $y \in {\mf Y}_{m-k+1} = \check {\mf Y}_{m-k+1}$ and $y \in \check {\mf Y}_{m-k} \setminus {\mf Y}_{m-k}$. 
Let ${\mf Y}_y$ and $\check {\mf Y}_y$ be the strata through $y$ associated to ${\mf Y}$ and $\check {\mf Y}$ respectively. Then
\begin{align}\label{eqnb1}
&\ {\rm dim}({\mf Y}_y) = m- k+1,\ &\ {\rm dim}(\check {\mf Y}_y ) \leq m-k.
\end{align}
Condition (4) implies that $f$ is transverse to $\check {\mf Z}$ and $f^* \check {\mf Z} \equiv {\mf Y}$. Then \eqref{eqnb1} implies
\beqn
{\rm dim} (\check {\mf Z}_{f(y)} ) = n- m+ {\rm dim} ( {\mf Y}_y ) = n- k + 1.
\eeqn
where $\check {\mf Z}_{f(y)}$ is the stratum containing $f(y)$ associated to the Whitney stratification $\check{\mf Z}$. On the other hand, by \eqref{eqnb1} and the transversality of $f$ to ${\mf Z}$, one has
\beqn
{\rm dim} ({\mf Z}_{f(y)} ) = n - m + {\rm dim}( \check {\mf Y}_y ) \leq n-k.
\eeqn
Therefore, 
\beqn
f(y) \in {\mf Z}_{n-k} \setminus \check {\mf Z}_{n-k} \Longrightarrow {\mf Z}_{n-k} \nsubseteq \check {\mf Z}_{n-k}.
\eeqn
As ${\mf Z}_{n-l} = \check {\mf Z}_{n-l}$ for all $l < k$, this contradicts the fact that ${\mf Z} \leq \check{\mf Z}$. Therefore, ${\mf Y}_{m-k} = \check {\mf Y}_{m-k}$. Similarly, one also obtains that ${\mf Z}_{n-k} = \check {\mf Z}_{n-k}$. Therefore, the inductive step is complete, and we conclude that ${\mf Y} \equiv \check {\mf Y}$ and ${\mf Z} \equiv \check {\mf Z}$. 
\end{proof}

It is useful to consider a special case.

\begin{cor}\label{corb6}
Let $M$ be a smooth manifold, $Y \subseteq M$ be a subset and ${\mc N}(Y) \subset \mc{WS}^\infty(Y)$ be a subset of smooth Whitney stratifications on $Y$ which is closed under equivalences. Suppose ${\mf Y} \in {\mc N}(Y)$ is a minimal element. Suppose $f: M \to M$ is a diffeomorphism such that $f(Y) = Y$, $f^*{\mf Y}\in {\mc N}(Y)$, and $(f^{-1})^* {\mf Y} \in {\mc N}(Y)$. Then $f^* {\mf Y} \equiv {\mf Y}$. 
\end{cor}

\begin{proof}
This follows from Lemma \ref{lemmab5} where $M = N$, $Y = Z$, and $g = f^{-1}$. 
\end{proof}

\subsection{Existence results in the complex analytic case}

The theory of canonical Whitney stratifications can be discussed in many different categories such as semianalytic, semialgebraic, real analytic, complex analytic, or complex algebraic (see \cite[Section 4]{Mather_1973}). In this paper we consider either the complex algebraic category or the complex analytic category. 

\begin{lemma}\label{lemmab7}
A subset $U \subset M$ is strongly analytic (Definition \ref{defn_strongly_analytic}) if and only if there exist closed analytic subsets $V$ and $W$ such that $U = V \setminus W$.
\end{lemma}

\begin{proof}
If $U$ is strongly analytic, then one can take $V = \ov{U}$ and $W = \ov{U}\setminus U$. Conversely, if $U = V\setminus W$, then by \cite[IV.10.Theorem 5]{Lojasiewicz_book}), $\ov{U}$ is analytic. Moreover, as $\ov{U} \subset V$, one has
\beqn
\ov{U} \setminus (\ov{U} \cap W)  = \ov{U} \setminus W \subset V \setminus W = U\Longrightarrow \ov{U}\setminus U \subset \ov{U} \cap W.
\eeqn
As $U \cap W = \emptyset$, it follows that $\ov{U}\setminus U = \ov{U} \cap W$ which is a closed analytic set. Then by definition, $U$ is strongly analytic.
\end{proof}

\begin{defn}
Let $Z \subseteq M$ be a strongly analytic subset. A point $p \in Z$ is called a {\it smooth regular point} resp. {\it analytic regular point} if there is an open neighborhood $N \subseteq M$ of $p$ such that $N \cap Z$ is a smooth resp. complex analytic submanifold.
\end{defn}
By the theorem of Bloom \cite{Bloom_1969}, smooth regular points are also analytic regular points. Hence we do not distinguish these two names. A point $p \in Z$ is {\it singular} if it is not a regular point.
Using the above definition, a regular point $p$ has a well-defined local dimension ${\rm dim}_p ( S)$.\footnote{In this appendix ``dimension'' always means ``real dimension.''} When $p$ is singular define ${\rm dim}_p(Z) = -\infty$. For each dimension $k$, denote 
\beqn
Z_{k, {\rm reg}} \subseteq Z = \{ p\in Z\ {\rm regular}\ |\ {\rm dim}_p(Z)  = k \}.
\eeqn
We also denote
\beqn
{\rm dim}(Z) \leq k\ \Longleftrightarrow {\rm dim}_p (Z) \leq k,\ \forall p \in Z.
\eeqn

\begin{lemma}\label{lemmab9}
Suppose $Z \subset M$ is a closed analytic subset. 
\begin{enumerate}

\item The set of singular points $Z_{\rm sing}$ is a closed analytic subset. 

\item For each $k$, $Z_{k, {\rm reg}}$ is a strongly analytic submanifold.
\end{enumerate}
\end{lemma}

\begin{proof}
For $Z_{\rm sing}$, it follows from \cite[Theorem IV.2.1]{Lojasiewicz_book}. On the other hand, notice that $Z_{k, {\rm reg}}$ is open and closed in $Z \setminus Z_{\rm sing}$, hence is the union of connected components. Indeed, the openness follows from the definition; for closedness, for $p \in \ov{Z}_{k, {\rm reg}}$, if $\dim_p Z \neq k$, then any point sufficiently closed to $p$ must also have local dimension not equal to $k$, so $p \in Z_{\rm sing}$. By \cite[Theorem IV.2.3]{Lojasiewicz_book}, the family of connected components of $Z \setminus Z_{\rm sing}$ is locally finite and each component's closure is analytic. Hence $\ov{Z_{k, {\rm reg}}}$ is analytic. As $Z_{k, {\rm reg}} = \ov{Z_{k, {\rm reg}}} \setminus Z_{\rm sing}$, by Lemma \ref{lemmab7}, $Z_{k, {\rm reg}}$ is a strongly analytic submanifold.
\end{proof}

Now we give a proof of the existence results of minimal Whitney stratifications following Mather \cite{Mather_1973}. We restrict our consideration to the complex analytic case. The construction relies on the following fundamental result, whose original version was proved by Whitney \cite{Whitney_1965} in the complex analytic setting. 

\begin{thm}(see \cite[Theorem 4.1, Addendum 4.3]{Mather_1973} \cite[Section 25]{Lojasiewicz_1965})\label{thmb10} Let $X$ and $Y$ be disjoint strongly analytic submanifolds of a complex manifold $M$. Let $S_b(X, Y)$ be the set of points $y \in Y$ where Whitney's condition (b) fails for $(X, Y)$ at $y$. Then $\ov{S_b(X, Y)}$ is an analytic set. Moreover, if $Y \subset \ov{X}$, then $S_b(X, Y)$ is nowhere dense in $Y$. 
\end{thm}



Now we can construct a canonical Whitney stratification on a closed analytic set. Indeed the following theorem is the ``absolute'' version of Theorem \ref{thmb12} below. However we would like to give a separate proof to show the reader the basic idea of the construction.

\begin{thm}\cite[Theorem 4.9, Addendum 4.11, Addendum 4.12]{Mather_1973}\label{thmb11}
Let $Z$ be a 
closed analytic subset of a smooth 
complex manifold $M$. Then there exists a smooth Whitney stratification ${\mf Z}$ of $Z$ such that for any open subset $U \subseteq M$ (in Euclidean topology), the restriction of ${\mf Z}$ to $Z \cap U$ is minimal. Moreover, ${\mf Z}$ is strongly complex analytic, i.e., its strata are all strongly complex analytic submanifolds.
 
\end{thm}

\begin{proof}
Let $n = {\rm dim}_{\mb R}(M)$. We construct inductively a decreasing sequence of subsets
\beqn
Z = {\mf Z}_n \supseteq {\mf Z}_{n-1} \supseteq \cdots \supseteq {\mf Z}_k
\eeqn
such that for each $l > k$ the following conditions are satisfied.
\begin{enumerate}
    \item ${\mf Z}_l$ is a closed analytic subset of $M$ of real dimension at most $l$.
    
    \item ${\mf Z}_l^*:= {\mf Z}_l \setminus {\mf Z}_{l-1}$ is a strongly analytic submanifold of real dimension $l$ (empty if $l$ is odd).
    
    \item For any $m>l$,  $({\mf Z}_m^*, {\mf Z}_l^*)$ satisfies Whitney's condition (b) at all points $x\in {\mf Z}_l^*$.
    
    \end{enumerate}

We start with ${\mf Z}_n = Z$. Suppose we have constructed ${\mf Z}_l$ for all $l \geq k$ satisfying the induction hypothesis. Then define ${\mf Z}_{k-1}$ to be the closure of points $x \in {\mf Z}_k$ satisfying one of the following conditions. 
\begin{enumerate}
    \item $x$ is a singular point of ${\mf Z}_k$ or a regular point with ${\rm dim}_x ({\mf Z}_k ) < k$;\label{item 1}
    
    \item $x$ is a regular point of ${\mf Z}_k$ with local dimension equal to $k$, and there exists $l > k$ such that the pair $({\mf Z}_l^*, {\mf Z}_{k, {\rm reg}})$ does not satisfy Whitney's condition (b) at $x$ (in particular $x \in \ov{{\mf Z}_l^*}$).\label{item 2}
\end{enumerate}
Then by construction ${\mf Z}_{k-1}$ has dimension at most $k-1$. The set of points satisfying (\ref{item 1}) above is a closed analytic subset. 
On the other hand,  
by Lemma \ref{lemmab9} and Theorem \ref{thmb10} above, the closure of $S_b({\mf Z}_l^*, {\mf Z}_{k, {\rm reg}})$ is a closed analytic set and has dimension at most $k-1$. Hence ${\mf Z}_{k-1}$ is a closed analytic subset of $M$ of dimension at most $k-1$ and ${\mf Z}_k^*:= {\mf Z}_k \setminus {\mf Z}_{k-1}$ is a strongly analytic submanifold of dimension $k$ (by Lemma \ref{lemmab7}). 
Other properties required for the induction hypothesis for ${\mf Z}_{k-1}$ are automatically satisfied. This completes the induction step.

Let ${\mf Z}$ be the collection of connected components of ${\mf Z}_k^*$ for all $k$. 
Then the collection ${\mf Z}$ is locally finite. Then ${\mf Z}$ satisfies all conditions for a Whitney stratification except for the axiom of frontier. The fact that ${\mf Z}$ satisfies the axiom of frontier follows from \cite[Proposition 8.7]{Mather_1973} (cf. Lemma \ref{lemmab3}). Therefore, ${\mf Z}$ is a Whitney stratification on $Z$. We show that ${\mf Z}$ is strongly complex analytic. Indeed, for each $k$, ${\mf Z}_k^*$ is the union of irreducible components, each of which is a strongly analytic submanifold. As ${\mf Z}_k^*$ is nonsingular, its irreducible components do not intersect. On the other hand, each irreducible component is connected as it is regular (see \cite[Section 5.4, Theorem]{Chirka_2012}). Hence all connected components of ${\mf Z}_k^*$ for all $k$, namely all strata of ${\mf Z}$, are strongly complex analytic submanifolds.

One needs to show that for all open subset $U \subset M$, ${\mf Z}|_{U \cap Z}$ is minimal among all smooth Whitney stratifications of $Z \cap U$. Indeed, let $\tilde {\mf Z}$ be any smooth Whitney stratification on $Z \cap U$ with dimension filtration $\tilde{\mf Z}_k$. We need to show that either ${\mf Z}_l \cap U = \tilde {\mf Z}_l$ for all $l$ or it is true for $l > k$ and ${\mf Z}_k \cap U \subsetneq \tilde {\mf Z}_k$. We prove this claim inductively. For $n = {\rm dim}_{\mb R}(M)$, one has ${\mf Z}_n \cap U  = \tilde {\mf Z}_n = Z \cap U$. Suppose we have proved that ${\mf Z}_l \cap U = \tilde {\mf Z}_l$ for all $l >k$ and suppose by contradiction that ${\mf Z}_k \cap U \nsubseteq \tilde {\mf Z}_k$. Then there exists $x \in {\mf Z }_k \cap U$ but $x \notin \tilde {\mf Z}_k$. Hence $x \in \tilde {\mf Z}_{k+1}\setminus \tilde {\mf Z}_k$. As $\tilde {\mf Z}_k$ is closed, it follows $\tilde {\mf Z}_{k+1} = {\mf Z}_{k+1} \cap U$ is locally a smooth submanifold near $x$ of dimension $k+1$. Then $x$ is a regular point and hence $x \in {\mf Z}_{k+1, {\rm reg}}$. Then by the construction of ${\mf Z}_k$, there must be some $m>k$ and a sequence $x_i \in {\mf Z}_{k+1, {\rm reg}}$ such that $x_i \to x$ and $({\mf Z}_{m+1}^*, {\mf Z}_{k+1, {\rm reg}})$ does not satisfy Whitney's condition (b) at $x_i$. As ${\mf Z}_{m+1}^* \cap U = \tilde {\mf Z}_{m+1}^*$, this contradicts the assumption that $\tilde {\mf Z}$ is a smooth Whitney stratification of $Z \cap U$. 
\end{proof}

One can easily extend the above results from the case of closed analytic sets to strongly analytic subsets. Indeed, if $Z$ is strongly analytic, one can first construct the canonical Whitney stratification for the closure $\ov{Z}$ whose strata $\ov{Z}_\alpha$ are strongly analytic submanifolds. Then its restriction to the open subset $Z \subset \ov{Z}$ is a minimal Whitney stratification on $Z$. Moreover, the induced strata $\ov{Z}_\alpha \setminus ( \ov{Z} \setminus Z)$ are again strongly analytic submanifolds (see Lemma \ref{lemmab7}).

\subsection{Relative minimal Whitney stratifications}

To prove Proposition \ref{thm_nice_Whitney}, we need to consider Whitney stratifications which are minimal in some relative sense. Suppose $M$ is a smooth complex manifold. 
Let $Z \subseteq M$ be a closed 
analytic subset and $Z' \subseteq Z$ be an 
open subset such that $Z\setminus Z' $ is a closed analytic set. Suppose $Z'$ is equipped with a
strongly analytic Whitney stratification ${\mf Z}'$. Define
\beqn
\mc{WS}^\infty(Z; {\mf Z}') \subset \mc{WS}^\infty(Z)
\eeqn
to be the subset of smooth Whitney stratifications on $Z$ whose restriction to $Z'$ coincides with ${\mf Z}'$; we say that such Whitney stratifications {\it extend} ${\mf Z}'$. Then for each open subset $U \subset M$, the restriction of Whitney stratifications induces a map 
\beqn
\mc{WS}^\infty(Z; {\mf Z}') \to \mc{WS}^\infty( Z \cap U; {\mf Z}'|_{Z' \cap U}).
\eeqn

\begin{thm}\label{thmb12}
There is a minimal element of ${\mf Z} \in \mc{WS}^\infty(Z; {\mf Z}')$ unique up to equivalence such that for any open subset $U \subset M$, the restriction ${\mf Z}|_{Z \cap U}$ is a minimal element of $\mc{WS}^\infty( Z \cap U; {\mf Z}'|_{Z' \cap U})$. In addition, ${\mf Z}$ is 
strongly  analytic.
\end{thm}

\begin{proof}
We follow the same strategy of the proof of Theorem \ref{thmb11}. As ${\mf Z}'$ is 
strongly complex analytic, for each $k$, ${\mf Z}_k' \setminus {\mf Z}_{k-1}'$ is a 
strongly complex analytic submanifold. Suppose ${\rm dim}_{\mb R}(M) = n$. For each $k \leq n$, we would like to define a decreasing sequence of closed analytic subsets
\beqn
Z\setminus Z' = {\mf Z}_n'' \supseteq {\mf Z}_{n-1}'' \cdots \supseteq {\mf Z}_k''
\eeqn
such that for all $l \geq k$ the following conditions are satisfied.
\begin{enumerate}
\item ${\mf Z}_l''$ is a closed 
analytic set of dimension at most $l$.

\item ${\mf Z}_l'' \setminus {\mf Z}_{l-1}''$ is a 
strongly complex submanifold of real dimension $l$. 

\item ${\mf Z}_l'' \setminus {\mf Z}_{l-1}''$ is disjoint from the closure of ${\mf Z}_l'$.

\item For each $x \in {\mf Z}_l'' \setminus {\mf Z}_{l-1}''$ and $m > l$, the pair $({\mf Z}_m'' \setminus {\mf Z}_{m-1}'' , {\mf Z}_l'' \setminus {\mf Z}_{l-1}'')$ and the pair $({\mf Z}_m' \setminus {\mf Z}_{m-1}', {\mf Z}_l''  \setminus {\mf Z}_{l-1}'')$ satisfy Whitney's condition (b) at $x$. 
\end{enumerate}
Now we start with ${\mf Z}_n'':= Z \setminus Z'$ which is indeed a closed analytic subset of $M$. Suppose we have constructed ${\mf Z}_n'', \ldots, {\mf Z}_k''$ satisfying the above conditions. Now consider the closure of points $x \in {\mf Z}_k''$ satisfying one of the following conditions, denoted by ${\mf Z}_{k-1}'' \subset {\mf Z}_k''$.
\begin{enumerate}

\item $x$ is contained in the closure of ${\mf Z}_k'$.

\item $x$ is either a singular point of ${\mf Z}_k''$ or a regular point with ${\rm dim}_x ({\mf Z}_k'') < k$.

\item $x \in {\mf Z}_{k, {\rm reg}}''$ and there exists $m>k$ such that the pair $({\mf Z}_m'' \setminus {\mf Z}_{m-1}'', {\mf Z}_{k, {\rm reg}}'')$ does not satisfy Whitney's condition (b) at $x$.

\item $x \in {\mf Z}_{k, {\rm reg}}''$ and there exists $m>k$ such that the pair $( {\mf Z}_m' \setminus {\mf Z}_{m-1}', {\mf Z}_{k, {\rm reg}}'')$ does not satisfy Whitney's condition (b) at $x$.
\end{enumerate}
Then following the same argument as the proof of Theorem \ref{thmb11}, ${\mf Z}_{k-1}''$ is a closed analytic set of real dimension at most $k-1$ and hence ${\mf Z}_k'' \setminus {\mf Z}_{k-1}''$ is a 
strongly analytic submanifold. This finishes the inductive step. Then similar to Theorem \ref{thmb11}, one can find a Whitney stratification (of $Z\setminus Z'$) by taking all connected components of ${\mf Z}_k'' \setminus {\mf Z}_{k-1}''$. Together with strata of ${\mf Z}'$, one obtains an extension of ${\mf Z}'$, denoted by ${\mf Z}$.

We need to prove that ${\mf Z}$ satisfies the minimality condition. Let $U \subseteq M$ be an open subset. Let $\tilde {\mf Z}$ be another Whitney stratification on $Z \cap U$ which extends ${\mf Z}'|_{Z' \cap U}$. Then the dimension filtrations of ${\mf Z}|_{Z \cap U}$ and $\tilde {\mf Z}$ are 
\begin{align*}
&\ {\mf Z}_k \cap U = ({\mf Z}_k' \cap U) \sqcup ({\mf Z}_k'' \cap U),\ &\ \tilde {\mf Z}_k = ({\mf Z}_k' \cap U) \sqcup \tilde {\mf Z}_k''
\end{align*}
where $\tilde {\mf Z}_k'' \subset (Z \setminus Z') \cap U$ is the dimension filtration of the restriction of $\tilde {\mf Z}$ to $(Z \setminus S')\cap U$. Then we need to show that either ${\mf Z}_k'' \cap U = \tilde {\mf Z}_k''$ for all $l$ or there exists $k$ such that it is true for all $l>k$ and ${\mf Z}_k'' \cap U \subsetneq \tilde {\mf Z}_k''$. We argue inductively. For $n = {\rm dim}(M)$, one has ${\mf Z}_n'' \cap U = \tilde {\mf Z}_n'' = (Z \setminus Z')\cap U$ by definition. Suppose we have proved that ${\mf Z}_l'' \cap U = \tilde {\mf Z}_l''$ for all $l > k$. Assume in contradiction that ${\mf Z}_k''  \cap U \nsubseteq \tilde {\mf Z}_k''$. Then there exists a point $x \in {\mf Z}_k'' \cap U$ which is not contained in $\tilde {\mf Z}_k''$, namely, $x \in \tilde {\mf Z}_{k+1}'' \setminus \tilde {\mf Z}_k''$. As $\tilde {\mf Z}_k''$ is closed, $x$ is a smooth (equivalently analytic) regular point of $\tilde{\mf Z}_{k+1}'' = {\mf Z}_{k+1}'' \cap U$ of local (real) dimension $k+1$. 

Now because $\tilde {\mf Z}$ is a Whitney stratification on $Z \cap U$ and $x \in \tilde {\mf Z}_{k+1}'' \setminus \tilde {\mf Z}_k''$, it follows from Lemma \ref{lemmab1}, $x\notin \ov{{\mf Z}_{k+1}'} \cap U$. Then as $x \in {\mf Z}_k'' \cap U\subset {\mf Z}_k''$, by the construction of ${\mf Z}$, it follows that there is a sequence $x_i \to x$ and $m > k$ such that either $( {\mf Z}_{m+1}'' \setminus {\mf Z}_{m}'', {\mf Z}_{k+1, {\rm reg}}'')$ or $({\mf Z}_{m+1}' \setminus {\mf Z}_m', {\mf Z}_{k+1, {\rm reg}}'')$ does not satisfy Whitney's condition (b) at $x_i$. This contradicts the hypothesis that $\tilde {\mf Z}$ is a Whitney stratification on $Z \cap U$. Therefore ${\mf Z}_k'' \cap U \subseteq \tilde {\mf Z}_k''$, hence inductively we conclude that ${\mf Z}|_{U \cap Z} \leq \tilde {\mf Z}$.
\end{proof}

\subsection{Proof of Proposition \ref{thm_nice_Whitney}}\label{appendix_nice_Whitney}

We first construct the Whitney stratification. Choose an increasing order of $\{M_\alpha\}$ as $M_{\alpha_1}, \ldots, M_{\alpha_m}$ such that $M_{\alpha_i} \leq M_{\alpha_j}$ implies $i \leq j$. For all $k$, denote 
\beqn
M^{(k)}:= M \setminus \bigcup_{i < k} M_{\alpha_i}
\eeqn
which gives a filtration
\beqn
M = M^{(1)} \supset M^{(2)}  \supset \cdots \supset M^{(m)}.
\eeqn
Denote
\beqn
Z^{(k)}:= Z \cap M^{(k)}
\eeqn
which is a stronlgy complex analytic subset of $M$. Then for each $k$, $\ov{Z^{(k)}}$ is analytic.

As the first step of the induction, by Theorem \ref{thmb11}, there exists a canonical Whitney stratification of the closure $\ov{Z^{(m)}}$ whose strata are strongly complex analytic submanifolds of $M$. Restrict to the open subset $Z^{(m)}$ to obtain a Whitney stratification ${\mf Z}^{(m)}$ which is still strongly complex analytic.

We assume, as the induction hypothesis, that for a $k \leq m$, we have constructed strongly complex analytic Whitney stratifications ${\mf Z}^{(k)}$ on $Z^{(k)}$ such that for each $i \geq k$, $Z_{\alpha_i}$ is the union of strata. Now consider the closure $\ov{Z^{(k-1)}}$ which contains the open subset $Z^{(k)}$ such that $\ov{Z^{(k-1)}} \setminus Z^{(k)}$ is closed analytic. Then by Theorem \ref{thmb12}, there exists a minimal extension of ${\mf Z}^{(k)}$ to $\ov{Z^{(k-1)}}$. Restrict the extension to $Z^{(k-1)}$ (which contains $Z^{(k)}$), one obtains a stronlgy complex analytic Whitney stratification ${\mf Z}^{(k-1)}$. Then the induction can continue until we find a nice strongly complex analytic Whitney stratification on $Z$. 


Lastly, one can combine induction with the proof of the minimality part of Theorem \ref{thmb12} to show that the constructed Whitney stratification is minimal among all smooth Whitney stratifications of $Z$ which respect ${\mf M}$. \qed

We have the following invariance property of the canonical Whitney stratification which respects the given stratification.

\begin{prop}\label{propb13}
Let $M$ be a smooth manifold equipped with a partition ${\mf M}$. Let $Z \subset M$ be a closed subset and let ${\mf Z} \in \mc{WS}^\infty(Z, {\mf M})$ be a minimal element. Let $f$ be a diffeomorphism of $M$ which preserves ${\mf M}$ such that $f(Z) = Z$. Then $f^* {\mf Z} \equiv {\mf Z}$. If ${\mf Z}$ has connected strata, then $f^* {\mf Z} = {\mf Z}$.
\end{prop}

\begin{proof}
The assumption implies that $f^* {\mf Z} \in \mc{WS}^\infty(Z; {\mf M})$. Then by Corollary \ref{corb6}, one has $f^* {\mf Z} \equiv {\mf Z}$. Suppose ${\mf Z}$ has connected strata. As $f$ is a diffeomorphism, for each stratum $Z_\alpha \in {\mf Z}$ which is connected, $f^{-1}(Z_\alpha)$, a stratum of $f^* {\mf Z}$, is connected. Hence $f^* {\mf Z} = {\mf Z}$. 
\end{proof}

\subsection{Proof of Proposition \ref{prop_product_nice}}\label{appendix_product_nice}

In this subsection we prove Proposition \ref{prop_product_nice}. First we gather some basic facts about product Whitney stratifications. Let $M$ and $N$ be smooth manifolds, and $S \subset M$ and $T \subset N$ be subsets equipped with a Whitney stratification ${\mf S}$ resp. ${\mf T}$. Denote 
\begin{align*}
&\ R = S \times T,\ \tilde {\mf R} = {\mf S} \times {\mf T}.
\end{align*}
Then the dimension filtration of $\tilde {\mf R}$ is 
\beqn
\tilde {\mf R}_m = \bigcup_{k + l = m} {\mf S}_k \times {\mf T}_l
\eeqn
while 
\beqn
\tilde {\mf R}_m^* = \bigsqcup_{k + l = m} {\mf S}_k^* \times {\mf T}_l^*
\eeqn
which is the disjoint union.

We first consider the absolute case of Proposition \ref{prop_product_nice}. 

\begin{prop}\label{propb14}
Let $M$ resp. $N$ be a  complex manifold and $S \subset M$ resp. $T \subset N$ be a closed analytic set equipped with the canonical Whitney stratification ${\mf S}$ resp. ${\mf T}$ (provided by Proposition \ref{thm_nice_Whitney}). Then the canonical Whitney stratification of $R:= S \times T \subset M \times M$ is equal to ${\mf S} \times {\mf T}$.
\end{prop}

\begin{proof}
Let ${\mf R}$ be the canonical Whitney stratification on $R$. Then ${\mf R} \leq \tilde {\mf R}$. As both stratifications have connected strata, we only need to prove that they are equivalent, i.e., having the same dimension filtration.  We prove by contradiction. Suppose ${\mf R}$ is not equivalent to $\tilde {\mf R}$. Then by definition, there exists an $m$ such that 
\beqn
{\mf R}_n = \tilde {\mf R}_n \ \forall n \geq m\ {\rm and}\ {\mf R}_{m-1} \subsetneq \tilde {\mf R}_{m-1}.
\eeqn
Therefore, there exists a point 
\beqn
r = (s, t) \in {\mf R}_m = \tilde {\mf R}_m = \bigcup_{k + l = m} {\mf S}_k \times {\mf T}_l,
\eeqn
and $r \in \tilde {\mf R}_{m-1}$, $r \notin {\mf R}_{m-1}$. Then $r = (s, t) \in {\mf R}_m^* \subset {\mf R}_{m, {\rm reg}} = \tilde {\mf R}_{m, {\rm reg}}$. Hence there exist  $k, l$ with $k+l = m$ such that $s \in {\mf R}_{k, {\rm reg}}$ and $t \in {\mf R}_{l, {\rm reg}}$. On the other hand, the condition that $r = (s, t) \in \tilde{\mf R}_{m-1}$ implies that either $s\in {\mf S}_{k-1}$ or $t \in {\mf T}_{l-1}$. Indeed, if $s \in {\mf S}_k^*$ and $t \in {\mf T}_l^*$, then $r \in {\mf S}_k^* \times {\mf T}_l^* \subset \tilde {\mf R}_m^*$, contradicting that $r \in \tilde {\mf R}_{m-1}$. 

Without loss of generality, we assume $s \in {\mf S}_{k-1}$. Then there exist a sequence $s_\mu' \in {\mf S}_{k'}^*$, a sequence $s_\mu \in {\mf S}_{k, {\rm reg}}$, both converging to $s$, such that $T_{s_\mu'} {\mf S}_{k'}^*$ converges to a subspace $H_s \subset T_s M$, the secant line $\ov{ s_\mu s_\mu'}$ converges to a line $L_s \subset T_s M$, and $L_s \nsubseteq H_s$. Now depending on the position of $t$ we discuss in two cases. 
\begin{enumerate}
    \item $t \notin {\mf T}_{l-1}$. Then consider the sequence of points $r_\mu = (s_\mu, t) \in {\mf S}_k \times {\mf T}_l \subset \tilde {\mf R}_{k+l}$ which converges to $r = (s, t) \in \tilde {\mf R}_{k+l, {\rm reg}}$. Then for $\mu$ sufficiently large, $r_\mu \in \tilde {\mf R}_{k+l, {\rm reg}}$. Consider the sequence of points $r_\mu' = (s_\mu', t)\in {\mf S}_{k'}^* \times {\mf T}_l^* \subset \tilde {\mf R}_{k' + l}^*$ which also converges to $r$. Then the sequence of tangent planes $T_{r_\mu'} {\mf S}_{k' + l}^* = T_{s_\mu'} {\mf S}_{k'}^* \oplus T_t {\mf S}_l^*$ converges to $H_s \oplus T_t {\mf S}_l^*$; the sequence of secant lines $\ov{ r_\mu r_\mu'}$ converges to the line $L_s \oplus \{0\}$ which is not contained in $H_s \oplus T_t {\mf T}_l^*$. Hence   condition (b) for $(\tilde {\mf R}_{k' + l}^*, \tilde {\mf R}_{k+l, {\rm reg}})$ fails at $r$. 

    \item $t \in {\mf T}_{l-1}$. Then it follows that for some $l' > l$,  condition (b) for $({\mf T}_{l'}^*, {\mf T}_{l, {\rm reg}})$ fails at $t$. Hence there exist a sequence $t_\mu' \in {\mf T}_{l'}^*$ and a sequence $t_\mu \in {\mf T}_{l, {\rm reg}}$ both converging to $t$ such that $T_{t_\mu'} {\mf T}_{l'}^*$ converges to a subspace $H_t \subset T_t N$, the secant lines $\ov{ t_\mu t_\mu'}$ converge to a line $L_t \subset T_t N$, and $L_t \nsubseteq H_t$. Then consider the sequence $r_\mu = (s_\mu, t_\mu)$ which converges to $r \in \tilde {\mf R}_{k+l, {\rm reg}}$. Then for $\mu$ sufficiently large, $r_\mu \in \tilde {\mf R}_{k+l, {\rm reg}}$. Consider the sequence $r_\mu' = (s_\mu', t_\mu') \in {\mf S}_{k'}^* \times {\mf T}_{l'}^*\subset \tilde {\mf R}_{k' + l'}^*$ which also converges to $r$. The tangent planes $T_{r_\mu'} \tilde {\mf R}_{k' + l'}^*$ converge to $H_s \oplus H_t$ while the secant lines $\ov{r_\mu r_\mu'}$ converge to a line $L_r \subset L_s \oplus L_t$. Then $L_r \nsubseteq H_s \oplus H_t$. Hence condition (b) for $(\tilde {\mf R}_{k' + l'}^*, \tilde {\mf R}_{k+l, {\rm reg}})$ fails at $r$. \qedhere
\end{enumerate}
In both cases, by the induction hypothesis, it follows that $r \in {\mf R}_{k+l-1} = {\mf R}_{m-1}$, contradicting the assumption that $r \in {\mf R}_m^*$. 
\end{proof}

\begin{proof}[Proof of Proposition \ref{prop_product_nice}]
Abbreviate $S_\alpha = S \cap M_\alpha$, $T_\beta = T \cap N_\beta$ with corresponding Whitney stratifications ${\mf S}_\alpha$ and ${\mf T}_\beta$. Abbreviate $R = S \times T$ and $R_{\alpha\beta} = S_\alpha \times T_\beta$. Let the canonical Whintey stratification on $R$ resp. $R_{\alpha\beta}$ be ${\mf R}$ resp. ${\mf R}_{\alpha\beta}$ and the product Whitney stratification be $\tilde {\mf R}$ resp. $\tilde {\mf R}_{\alpha\beta}$. 

We prove this proposition by induction on $(\alpha, \beta)$. First, Proposition \ref{propb14} implies that for any top stratum $R_{\alpha\beta} \subset M \times N$, ${\mf R}_{\alpha\beta} = \tilde {\mf R}_{\alpha\beta}$. Now fix a pair $(\alpha, \beta)$, suppose we have shown that for any $R_{\gamma\delta} > R_{\alpha\beta}$, ${\mf R}_{\gamma\delta} = \tilde {\mf R}_{\gamma\delta}$. We would like to show that ${\mf R}_{\alpha\beta} = \tilde {\mf R}_{\alpha\beta}$. Suppose our claim is false. Then there exists $m \geq 0$ such that 
\beqn
{\mf R}_{\alpha\beta, n} = \tilde {\mf R}_{\alpha\beta, n}\ (\forall n \geq m)\ \ \  {\rm and}\ \ \  {\mf R}_{\alpha\beta, m-1}  \subsetneq \tilde {\mf R}_{\alpha\beta, m-1}.
\eeqn
Then one can choose a point $r = (s, t)$ such that
\begin{align*}
&\ r \in {\mf R}_{\alpha\beta, m}^*,\ &\ r \in \tilde {\mf R}_{\alpha\beta, m-1}.
\end{align*}
Then
\beq\label{yyy}
r = (s, t) \in \tilde {\mf R}_{\alpha\beta, m-1} \subset \tilde {\mf R}_{\alpha\beta, m} = \bigcup_{p + q = m} {\mf S}_{\alpha, p} \times {\mf T}_{\beta, q} = {\mf R}_{\alpha\beta, m}. 
\eeq
Then $r \in {\mf R}_{\alpha\beta, m}^* \subset {\mf R}_{\alpha\beta, m, {\rm reg}} = \tilde {\mf R}_{\alpha\beta, m, {\rm reg}}$. Therefore for some pair $(p, q)$ with $p+q = m$, $(s, t)$ is a regular point of ${\mf S}_{\alpha, p}\times {\mf T}_{\beta, q}$ of local dimension $p + q$. Hence
\begin{align*}
&\ s \in {\mf S}_{\alpha, p, {\rm reg}},\ &\ t \in {\mf T}_{\beta, q, {\rm reg}}.
\end{align*}
Similar to the proof of the absolute case, we know that either $s \in {\mf S}_{\alpha, p-1}$ or $t \in {\mf T}_{\beta, q-1}$. Without loss of generality, assume that $s \in {\mf S}_{\alpha, p-1}$.

\vspace{0.1cm}

\noindent {\it Claim.} For all $r> \alpha$, $p' \leq p$, $\delta > \beta$, $q' \leq q$, 
\begin{align*}
&\ s \notin \partial {\mf S}_{\gamma, p'}^*,\ &\ t \notin \partial {\mf T}_{\beta, q'}^*.
\end{align*}

\noindent {\it Proof of the claim.} Suppose for some $\gamma > \alpha$ and $p' \leq p$, $s \in \partial {\mf S}_{\gamma, p'}^*$. Let $q'\leq q$ be such that $t \in {\mf T}_{\beta, q'}^*$. Then we see that 
\beqn
r = (s, t) \in \partial {\mf S}_{\gamma, p'}^* \times {\mf T}_{\beta, q'}^* \subset \partial \left( {\mf S}_{\gamma, p'}^* \times {\mf T}_{\beta, q'}^* \right) \subset \partial \tilde {\mf R}_{\gamma\beta, p' + q'}^* = \partial {\mf R}_{\gamma\beta, p' + q'}^*.
\eeqn
Here the last equality follows from the induction hypothesis. Notice that $p' + q' \leq m$ and $R_{\gamma\beta} > R_{\alpha\beta}$. This contradicts the fact that $r \notin {\mf R}_{\alpha\beta, m-1}$ and the way to construct the canonical nice Whitney stratification ${\mf R}$. The claimed condition for $t$ can be proved in the same way. \hfill {\it End of the proof of the claim.}

The rest of the argument is similar to the last part of the proof of Proposition \ref{propb14}, except that we need to use the above claim. The assumption $s \in {\mf S}_{\alpha, p-1}$ and the above claim imply that there exists $k' > p$ such that condition (b) for $({\mf S}_{\alpha, k'}^*, {\mf S}_{\alpha, k, {\rm reg}})$ fails at $s$. Then there exist a sequence $s_\mu' \in {\mf S}_{\alpha, k'}^*$, a sequence $s_\mu \in {\mf S}_{\alpha, k, {\rm reg}}$, both converging to $s$, such that $T_{s_\mu'} {\mf S}_{\alpha, k'}^*$ converges to a subspace $H_s \subset T_s M_\alpha$, the secant line $\ov{ s_\mu s_\mu'}$ converges to a line $L_s \subset T_s M_\alpha$, and $L_s \nsubseteq H_s$. Now depending on the position of $t$ we discuss in two cases. 
\begin{enumerate}
    \item Assume $t \notin {\mf T}_{\beta, l-1}$. Then consider the sequence of points $r_\mu = (s_\mu, t) \in {\mf S}_{\alpha,k} \times {\mf T}_{\beta,l} \subset \tilde {\mf R}_{\alpha\beta, k+l}$ which converges to $r = (s, t) \in \tilde {\mf R}_{\alpha\beta, k+l, {\rm reg}}$. Then for $\mu$ sufficiently large, $r_\mu \in \tilde {\mf R}_{\alpha\beta, k+l, {\rm reg}}$. Consider the sequence of points $r_\mu' = (s_\mu', t)\in {\mf S}_{\alpha, k'}^* \times {\mf T}_{\beta,l}^* \subset \tilde {\mf R}_{\alpha\beta, k' + l}^*$ which also converges to $r$. Then the sequence of tangent planes $T_{r_\mu'} {\mf S}_{\alpha, k' + l}^* = T_{s_\mu'} {\mf S}_{\alpha, k'}^* \oplus T_t {\mf S}_l^*$ converges to $H_s \oplus T_t {\mf S}_l^*$; the sequence of secant lines $\ov{ r_\mu r_\mu'}$ converges to the line $L_s \oplus \{0\}$ which is not contained in $H_s \oplus T_t {\mf T}_l^*$. Hence   condition (b) for $(\tilde {\mf R}_{\alpha\beta, k' + l}^*, \tilde {\mf R}_{\alpha\beta, k+l, {\rm reg}})$ fails at $r$. 

    \item $t \in {\mf T}_{\beta, l-1}$. Then the above claim implies that for some $l' > l$,  condition (b) for $({\mf T}_{\beta, l'}^*, {\mf T}_{\beta, l, {\rm reg}})$ fails at $t$. Hence there exist a sequence $t_\mu' \in {\mf T}_{\beta, l'}^*$ and a sequence $t_\mu \in {\mf T}_{\beta, l, {\rm reg}}$ both converging to $t$ such that $T_{t_\mu'} {\mf T}_{\beta, l'}^*$ converges to a subspace $H_t \subset T_t N_\beta$, the secant lines $\ov{ t_\mu t_\mu'}$ converge to a line $L_t \subset T_t N_\beta$, and $L_t \nsubseteq H_t$. Then consider the sequence $r_\mu = (s_\mu, t_\mu)$ which converges to $r \in \tilde {\mf R}_{\alpha\beta, k+l, {\rm reg}}$. Then for $\mu$ sufficiently large, $r_\mu \in \tilde {\mf R}_{\alpha\beta, k+l, {\rm reg}}$. Consider the sequence $r_\mu' = (s_\mu', t_\mu') \in {\mf S}_{\alpha, k'}^* \times {\mf T}_{\beta, l'}^*\subset \tilde {\mf R}_{\alpha\beta, k' + l'}^*$ which also converges to $r$. The tangent planes $T_{r_\mu'} \tilde {\mf R}_{\alpha\beta, k' + l'}^*$ converge to $H_s \oplus H_t$ while the secant lines $\ov{r_\mu r_\mu'}$ converge to a line $L_r \subset L_s \oplus L_t$. Then $L_r \nsubseteq H_s \oplus H_t$. Hence condition (b) for $(\tilde {\mf R}_{\alpha\beta, k' + l'}^*, \tilde {\mf R}_{\alpha\beta, k+l, {\rm reg}})$ fails at $r$. 
\end{enumerate}
In both cases, by the induction hypothesis it follows that $r \in {\mf R}_{\alpha\beta, m-1}$ which contradicts the assumption that $r \in {\mf R}_{\alpha\beta, m}^*$.
\end{proof}

\subsection{Pullback under submersions}

This subsection serves as a technical preparation for proving several statements in Section \ref{section3}. 

First notice that Whitney's condition (b) is invariant under diffeomorphisms. In fact it is preserved under pullbacks by submersions. 

\begin{lemma}
Let $X, Y$ be disjoint smooth submanifolds of $M$. Let $\pi: \tilde M \to M$ be a submersion. Denote $\tilde X = \pi^{-1}(X)$, $\tilde Y = \pi^{-1}(Y)$. Then $(X, Y)$ satisfies Whitney's condition (b) at $y \in \pi(\tilde Y) \subset Y$ if and only if $(\tilde X, \tilde Y)$ satisfies Whitney's condition (b) at $\tilde y \in \tilde Y$ for all $\tilde y \in \pi^{-1}(y)$.
\end{lemma}

\begin{proof}
As $\pi$ is a submersion, it is transverse to both $X$ and $Y$. Then the Whitney condition is preserved under transversal pullbacks (see for example \cite[(1.3)]{Topological_stability}). 
Conversely, suppose $(\tilde X, \tilde Y)$ satifies condition (b) at all points $\tilde y \in \pi^{-1}(y)$. Suppose $x_i \in X$ converges to $y$, $y_i \in Y$ converges to $y$, $\ov{x_i y_i}$ converges to a line $l \subset T_y M$, and $T_{x_i} X$ converges to $\tau \subset T_y M$. Choose a preimage $\tilde y \in \pi^{-1}(y)$. Then there exist nearby points $\tilde x_i \in \pi^{-1}(x_i) \subset \tilde X$ converging to $\tilde y$ and $\tilde y_i \in \pi^{-1}(y_i) \subset \tilde Y$ converging to $\tilde y$. By choosing a subsequence, we may assume $\ov{\tilde x_i \tilde y_i}$ converges to $\tilde l \subset T_{\tilde y} \tilde M$ and $T_{\tilde x_i} \tilde X$ converges to $\tilde \tau \subset T_{\tilde y} \tilde M$. We may also move $\tilde y_i$ in the same fiber of $\pi$ to guarantee that $\tilde l \nsubseteq {\rm ker} (d\pi)$. Then $d\pi(\tilde l) = l$ and 
\beqn
d\pi( \tilde \tau) = d\pi\left( \lim_{i \to \infty} T_{\tilde x_i} \tilde U \right) = \lim_{i \to \infty} d\pi (T_{\tilde x_i} \tilde U) = \lim_{i \to \infty} T_{x_i} U = \tau.
\eeqn
Then $\tilde l \subset \tilde \tau$ which implies $l \subset \tau$. Hence $(X, Y)$ satisfies condition (b) at $y$. 
\end{proof}

One has the following straightforward corollary.

\begin{cor}\label{cora16}
When $\pi$ is a surjective submersion, $S_b( \tilde X, \tilde Y ) = \pi^{-1} (S_b(X, Y))$.
\end{cor}

There are some other easy results about surjective submersions in complex analytic category.

\begin{lemma}\label{lemmab17}
Let $M$ be a complex manifold with ${\rm dim}_{\mb R}(M) = m$ and $U \subset M$. Let $\pi: \tilde M \to M$ be a surjective analytic submersion and ${\rm dim}_{\mb R} (\tilde M) = \tilde m$. Then
\begin{enumerate}
    \item $\pi^{-1} (\ov{U}) = \ov{\pi^{-1}(U)}$.
    
    \item If $U$ is a closed analytic subset, then 
    \beqn
    \pi^{-1}(U)_{\rm sing} = \pi^{-1}(U_{\rm sing})
    \eeqn
    and for each $l$, 
    \beqn
    \pi^{-1}(U)_{\tilde m - l, {\rm reg}} = \pi^{-1}( U_{m-l, {\rm reg}}).
    \eeqn

    \item If $U \subset M$ is strongly analytic, then $\pi^{-1}(U)\subset \tilde M$ is also strongly analytic.
\end{enumerate}
\end{lemma}
\begin{proof}
    For item (1), the inclusion $\ov{\pi^{-1}(U)} \subset \pi^{-1} (\ov{U})$ follows from definition; for the reversed inclusion, it can be proven using the surjective property of $\pi$ and choosing a product neighborhood of limiting points in $\ov{U}$. For item (2), note that $\tilde{p} \in \pi^{-1}(U)$ is a smooth point if and only if $p = \pi(\tilde{p}) \in U$ is a smooth point, so it follows by taking the complement. For item (3), by assumption, both $\ov{U}$ and $\ov{U} \setminus U$ are closed and analytic; it follows from (1) that both $\ov{\pi^{-1}(U)} = \pi^{-1}(\ov{U})$ and $\ov{\pi^{-1}(U)} \setminus \pi^{-1}(U)$ are closed and analytic subsets, so (3) is proved.
\end{proof}

\subsubsection{Pulling back minimal Whitney stratifications}

Let $M$ be a complex manifold and $\pi: \tilde M \to M$ be a surjective analytic submersion. Suppose $Z \subset M$ is a strongly analytic subset and ${\mf Z}$ is a strongly analytic Whitney stratification of $Z$, then by Lemma \ref{lemmab17}, the pullback Whitney stratification $\pi^* {\mf Z}$ is also strongly analytic.

\begin{cor}\label{corb18}
Let $M$ be a complex manifold, $Z \subset M$ be a closed analytic subset, $Z' \subset Z$ be an open subset such that $Z \setminus Z'$ is a closed analytic subset. Suppose $\pi: \tilde M \to M$ is a complex analytic submersion with $Z \subset \pi(\tilde M)$. Let ${\mf Z}'$ be a strongly analytic Whitney stratification of $Z'$ and let ${\mf Z}$ be a minimal extension to $Z$. Then $\pi^* {\mf Z}$ is a minimal extension to $\pi^{-1}(Z)$ of the pullback $\pi^* {\mf Z}'$ on $\pi^{-1}(Z')$.
\end{cor}

\begin{proof}
The proof is based on tracing the steps in the proof of Theorem \ref{thmb12}. Suppose 
${\rm dim}_{\mb R}(M) = m$ and ${\rm dim}_{\mb R} ( \tilde M ) = \tilde m$. Following the recipe of the proof of Theorem \ref{thmb12}, one can construct a minimal extension $\tilde {\mf Z}$ of $\pi^* {\mf Z}'$ by constructing a decreasing sequence
\beqn
\pi^{-1}(Z \setminus Z') = \tilde {\mf Z}_{\tilde m}'' \supseteq \tilde {\mf Z}_{\tilde m-1}'' \cdots \supseteq \tilde {\mf Z}_k''.
\eeqn
We claim that for $l = 0, 1, \ldots$, 
\beqn
\tilde {\mf Z}_{\tilde m - l}'' \equiv \pi^{-1}( {\mf Z}_{m-l}'').
\eeqn
This would imply that $\tilde {\mf Z} \equiv \pi^* {\mf Z}$. Indeed, the set $\tilde {\mf Z}_k''$ is defined as the closure of points in $\tilde {\mf Z}_{k-1}''$ which satisfy one of the conditions listed in the proof of Theorem \ref{thmb12}. Corollary \ref{cora16} and Lemma \ref{lemmab17} imply that these conditions are preserved under the pullback by the submersion $\pi$. Therefore, $\tilde {\mf Z}_{\tilde m-l}'' \equiv \pi^{-1}({\mf Z}_{m-l}'')$. 
\end{proof}

There is another technical situation we need to consider. 

\begin{prop}\label{prop_pullback_nice}
Let $M, {\mf M}, Z$ be as in Proposition \ref{thm_nice_Whitney}. Let ${\mf Z}$ be a minimal Whitney stratification of $Z$ which respects ${\mf M}$. Let $\tilde M$ be a complex manifold and $\pi: \tilde M \to M$ be a surjective analytic submersion. Notice that there is the pullback stratification  
\beqn
\tilde {\mf M} = \{ \tilde M_\alpha \} = \{ \pi^{-1} (M_\alpha)\}
\eeqn
of $\tilde M$ whose strata are strongly analytic submanifolds. Then the pullback Whitney stratification $\pi^* {\mf Z}$ is a minimal one on $\pi^{-1}(Z)$ which respests $\tilde {\mf M}$. 
\end{prop}

\begin{proof}
We go through the construction in the proof of Proposition \ref{thm_nice_Whitney} for both $Z$ and $\tilde Z := \pi^{-1}(Z)$ in parallel. Order the strata $M_\alpha$ increasingly as 
\beqn
M_{\alpha_1}, \ldots, M_{\alpha_m}
\eeqn
such that $M_{\alpha_i} \leq M_{\alpha_j}$ implies $i \leq j$. Define 
\begin{align*}
&\ M^{(k)}:= M \setminus \bigcup_{i < k} M_{\alpha_i},\ &\ \tilde M^{(k)}:= \tilde M \setminus \bigcup_{i < k} \tilde M_{\alpha_i} = \pi^{-1}( M^{(k)}).
\end{align*}
Set 
\begin{align*}
&\ Z^{(k)}:= Z \cap M^{(k)},\ &\ \tilde Z^{(k)}:= \tilde Z \cap \tilde M^{(k)} = \pi^{-1}(Z^{(k)}).
\end{align*}
Then following the proof of Proposition \ref{thm_nice_Whitney} for both $Z$ and $\tilde Z$, we have the strongly analytic Whitney stratifications ${\mf Z}$ and $\tilde {\mf Z}$ respecting ${\mf M}$ and $\tilde{\mf M}$ respectively. For each $k$, let ${\mf Z}^{(k)}$ resp. $\tilde {\mf Z}^{(k)}$ be the restriction of ${\mf Z}$ resp. $\tilde {\mf Z}$ on the open subset $Z^{(k)}$ resp. $\tilde Z^{(k)}$. We prove by a reversed induction on $k$ that for each $k$, one has 
\beq\label{eqnb3}
\tilde {\mf Z}^{(k)} = \pi^* {\mf Z}^{(k)}.
\eeq
When $k = m$, recall that ${\mf Z}^{(m)}$ is the restriction of the canonical Whitney stratification on $\ov{Z^{(m)}}$ to $Z^{(m)}$. Then the absolute case (i.e. the $Z' = \emptyset$ case) of Corollary \ref{corb18} shows that the canonical Whitney stratification on $\ov{\tilde Z^{(m)}} = \pi^{-1}( \ov{Z^{(m)}})$ is equivalent to the pullback of the canonical Whitney stratification on $\ov{Z^{(m)}}$. Then the restriction to $\tilde Z^{(m)}$ is equivalent to $\pi^* {\mf Z}^{(m)}$. The relative case of Corollary \ref{corb18} allows us to prove \eqref{eqnb3} inductively. 
\end{proof}

\section{Thom--Mather Stratified Spaces}\label{appendixc}

In this appendix we review the notion of Thom--Mather stratified spaces (or called abstract stratified spaces), especially their fundamental classes. It was proved by Mather \cite{Mather_1970} that Whitney stratified sets in smooth manifolds are Thom--Mather stratified spaces. The zero locus of an FOP transverse section is a stratified space inside an effective orbifold which satisfies similar properties as a Whitney stratified set in a manifold. In order to construct a fundamental class, we extend Mather's result to this kind of stratified subsets in orbifolds.

\subsection{Thom--Mather stratified spaces}

Thom--Mather stratified spaces are spaces stratified by smooth manifolds with certain structures regularizing the normal behavior near each stratum. Since the definitions in standard references (\cite{Thom_1964, Thom_1969}\cite{Mather_1970}) require the control data as part of the structure, which is not canonically given in relevant constructions, we modify the definition in the following way.

\begin{defn}\label{defn_Thom_Mather}
A {\it Thom--Mather stratified space} is a stratified space $(X, {\mf X})$ satisfying the following conditions. 
\begin{enumerate}

\item $X$ is a locally compact, Hausdorff, and second countable space.

\item ${\mf X}$ is a stratification such that each stratum $X_\alpha \in {\mf X}$ is a smooth manifold.

\item There exists a set of {\it control data}, which is a collection 
\beqn
{\mf J} = \Big\{(N_\alpha, \pi_\alpha, \rho_\alpha)\ |\ X_\alpha \in {\mf X} \Big\}
\eeqn
where $N_\alpha \subset X$ is an open neighborhood of $X_\alpha$, $\pi_\alpha: N_\alpha \to X_\alpha$ is a continuous retraction, and $\rho_\alpha: N_\alpha \to [0, +\infty)$ is a continuous function. Moreover, these triples satisfy the following conditions.
\begin{enumerate}
\item $X_\alpha = \{ v \in N_\alpha \ |\ \rho_\alpha (v) = 0\}$.

\item For each pair of strata $X_\alpha, X_\beta \in {\mf X}$, define $N_{\beta\alpha} = N_\alpha \cap X_\beta$ (which is an open subset of the manifold $X_\beta$), $\pi_{\beta\alpha} = \pi_\alpha |_{N_{\beta\alpha}}$, and $\rho_{\beta\alpha} = \rho_\alpha |_{N_{\beta\alpha}}$. We require that the map 
\beqn
(\pi_{\beta\alpha}, \rho_{\beta\alpha}): N_{\beta\alpha} \to X_\alpha \times (0, +\infty)
\eeqn
is a smooth submersion. 

\item For any three strata $X_\alpha, X_\beta, X_\gamma$ one has 
\begin{align*}
&\ \pi_{\beta\alpha}\circ \pi_{\gamma\beta} = \pi_{\gamma\alpha},\ &\ \rho_{\beta\alpha} \circ \pi_{\gamma\beta} = \rho_{\gamma\alpha}
\end{align*}
whenever both sides of the equations are defined.
\end{enumerate}
\end{enumerate}
\end{defn}

The {\it dimension} of a Thom--Mather stratified space $X$ is the maximum of dimensions of its strata. The {\it dimension filtration} of $X$ is the filtration
\beqn
X_n \supseteq X_{n-1} \supseteq \cdots
\eeqn
where $X_k$ is the union of strata of dimension at most $k$. Notice that the union of any subset of strata is naturally a Thom--Mather stratified space. On the other hand, any open subset of $X$ is also a Thom--Mather stratified space.

\begin{defn}
An {\it isomorphism} of Thom--Mather stratified spaces from $X$ to $Y$ is a homeomorphism $f: X \to Y$ whose restriction to each stratum $X_\alpha \subset X$ is a diffeomorphism onto a stratum of $Y$. An {\it open embedding} of Thom--Mather stratified space is an isomorphism of Thom--Mather stratified space onto an open subset of the codomain.
\end{defn}

Goresky \cite{Goresky_1978} proved that Thom--Mather stratified spaces admits triangulations. 

\begin{thm} \cite[5. Proposition]{Goresky_1978} \label{thm_triangulation} Let $(X, {\mf X})$ be a Thom--Mather stratified space. Then there exists a polyhedral complex $K$ and a homeomorphism $f: |K| \to X$ such that for each stratum $X_\alpha$, $f^{-1}(\ov{X_\alpha})$ is a subcomplex and $f: f^{-1}(X_\alpha) \to X_\alpha$ is a smooth triangulation.\footnote{Goresky's construction {\it a priori} depends on the control data.}
\end{thm}

We would like to show that certain Thom--Mather stratified spaces admit fundamental classes. Let $X$ be a compact $n$-dimensional Thom--Mather stratified space. 

\begin{lemma}\label{lemmac4}
If $k > l$, then $H_k(X_l; {\mb Z}) = 0$.
\end{lemma}

\begin{proof}
Each $X_l$ itself is a Thom--Mather stratified space. Then the lemma follows from Goresky's theorem. 
\end{proof}

\begin{lemma}\cite[0.5.1 Proposition]{Goresky_thesis} Given a control data on $(X, {\mf X})$, for $\epsilon>0$ sufficiently small, let $N_\alpha(\epsilon):= \rho_\alpha^{-1}([0, \epsilon)) \subset N_\alpha \subset X$. Define
\beqn
O_{n-1}(\epsilon):= \bigcup_{{\rm dim} (X_\alpha) \leq n-1} N_\alpha(\epsilon).
\eeqn
Then for $\epsilon$ sufficiently small, the inclusion map $X_{n-1} \to O_{n-1}(\epsilon)$ is a homotopy equivalence.
\end{lemma}

Now suppose each top stratum of an $n$-dimensional Thom--Mather stratified space $X$ is oriented. Abbreviate $O_{n-1}(\epsilon)$ by $O_{n-1}$. Then there exists a fundamental class 
\beqn
[X] \in H_n(X, O_{n-1})
\eeqn
(in integer coefficients). The homotopy equivalence $X_{n-1} \simeq O_{n-1}$ implies
\beqn
H_n(X, X_{n-1}) \cong H_n(X, O_{n-1}).
\eeqn
Hence we regard $[X]$ as a class in $ H_n(X, X_{n-1})$. We call it the {\it fundamental class} of $X$ (with respect to the orientations on the top strata). If we are given a triangulation of $X$, then the fundamental class is represented by the singular cycle which is the sum of all the oriented top-dimensional cells of the triangulation.

\begin{defn}\label{defn_TM_pseudomanifold}
An $n$-dimensional {\it Thom--Mather stratified pseudomanifold} is an $n$-dimensional Thom--Mather stratified space $X$ such that $X_{n-1} = X_{n-2}$.
\end{defn}

If $X$ is a compact oriented Thom--Mather stratified pseudomanifold, the exact sequence
\beqn
H_n(X) \to H_n(X, X_{n-1}) \to H_{n-1}(X_{n-1}) = H_{n-1}(X_{n-2}) = \{0\},
\eeqn
where the last equality follows from Lemma \ref{lemmac4}, implies that $[X]$ lives in $H_n(X; {\mb Z})$. We record this discussion as follows.

\begin{prop}
    Any $n$-dimensional Thom--Mather stratified pseudomanifold $X$ with oriented top-dimensional strata has a well-defined fundamental class 
    \beqn
    [X] \in H_n(X; {\mb Z}).
    \eeqn
\end{prop}

\begin{defn}\label{defnc8}
A {\it cobordism} of compact oriented $n$-dimensional Thom--Mather stratified pseudomanifolds from $Y$ to $Z$ consists of a compact oriented $n+1$-dimensional Thom--Mather stratified space $X$, a closed union of strata $\partial X \subset X$, and an isomorphism $Y \sqcup Z \cong \partial X$ satisfying the following conditions.
\begin{enumerate}

\item $X \setminus \partial X$ is an $n+1$-dimensional Thom--Mather stratified pseudomanifold.

\item There exists an open embedding $\partial X \times [0, \epsilon) \to X$ whose restriction to $\partial X \times \{0\}$ is the identity map of $\partial X$. This implies that top strata of $\partial X$ have naturally induced orientations.

\item The isomorphism $Y \sqcup Z \cong \partial X$ is orientation preserving if $Y \sqcup Z$ is oriented as $(-Y) \sqcup Z$.
\end{enumerate}
\end{defn}

In the situation of the above definition, $X_n = X_{n-1} \cup \partial X$. Consider the fundamental class $[X] \in H_{n+1}(X, X_n)$ and the exact sequence 
\beqn
H_{n+1} (X_n, \partial X) \to H_{n+1} ( X, \partial X) \to H_{n+1} (X, X_n ) \to H_n (X_n, \partial X).
\eeqn
By the collar structure, excision, and the fact that $X_n \setminus \partial X = X_{n-1} \setminus \partial X$, it follows that $H_{n+1} (X_n, \partial X) = H_n ( X_n, \partial X) = 0$. Hence $[X]$ lives in $H_{n+1} (X, \partial X)$.

\begin{lemma}\label{lemmac9}
The image of $[X]$ under the boundary map $H_{n+1} (X, \partial X) \to H_n (\partial X)$ is $[\partial X] = [Z] - [Y]$.
\end{lemma}

\begin{proof}
By Theorem \ref{thm_triangulation}, one can choose a triangulation of $X$. We see that the codimension one faces of the top-dimensional cells, if not cancelled in pairs, are exactly the top-dimensional faces of the boundary $\partial X$.
\end{proof}

\subsection{Tubular neighborhoods in orbifolds}

\subsubsection{Tubular neighborhoods in manifolds}

We recall several useful notions and results about tubular neighborhoods in smooth manifolds towards constructing control data on Whitney stratified subsets of orbifolds. The main reference is \cite[Section 6]{Mather_1970}. Let $(M, g)$ be a smooth manifold and $S \subset M$ be a closed submanifold. Recall the notion of normal bundle $NS \to S$ and disk bundle $N^\epsilon S \subset NS$ of radius $\sqrt{\epsilon}$ (defined with respect to the inner product on $NS$ induced from $g$) associated to a smooth function $\epsilon: S \to {\mb R}_+$ (cf. \eqref{disk_bundle}). We need to consider identifications of disk bundles in $NS$ and tubular neighborhoods beyond those obtained by normal exponential maps.

\begin{defn}
Let $(M, g)$ and $S$ be as above.
\begin{enumerate}
\item A {\it tubular neighborhood} of $S$ is a pair $T_S = (\epsilon, \varphi)$ where $\epsilon: S \to {\mb R}_+$ is a smooth function and $\varphi: N^\epsilon S \to M$ is an open embedding satisfying the following condition. For any $x \in S$, $\varphi(x) = x$ and composed linear map 
\beqn
\xymatrix{  N_x S \ar[r]^{d\varphi_x} & T_x M \ar[r] & N_x S}
\eeqn
is the identity map of $N_x S$. In particular, $d\varphi|_S$ gives a splitting of the canonical exact sequence
\beqn
\xymatrix{
0 \ar[r] & TS \ar[r] & TM|_S \ar[r] & NS \ar[r] & 0
}.
\eeqn

\item Associated to each tubular neighborhood $T_S$, denote 
\beqn
|T_S| = \varphi(N^\epsilon S),
\eeqn
over which there is a {\it tubular projection} 
\beqn
\pi_S:= \pi_{NS} \circ \varphi^{-1}: |T_S| \to S
\eeqn
and the {\it tubular function}
\beqn
\rho_S: |T_S| \to [0, +\infty),\ \rho_S(y) = \| \varphi^{-1}(y) \|^2.
\eeqn

\item The {\it restriction} of a tubular neighborhood $T_S$ of $S$ to an open subset $U \subset S$ is the tubular neighborhood $T_S|_U:= (\epsilon|_U, \varphi|_{N^\epsilon S|_U})$ of $U$.

\item Two tubular neighborhoods $T_S = (\epsilon, \varphi)$ and $T_S' = (\epsilon', \varphi')$ of $S$ are {\it equivalent} if 
\beqn
(\pi_S, \rho_S) = (\pi_S', \rho_S')
\eeqn
as germs of maps defined near $S$. Equivalently, $T_S$ is equivalent to $T_S'$ if there exist a smooth function $\epsilon'': S \to {\mb R}_+$ satisfying $\epsilon'' \leq \min (\epsilon, \epsilon')$ and a (nonlinear) bundle map 
\beqn
\tau: N^{\epsilon''} S \to N^{\epsilon'' } S
\eeqn
which is fiberwise diffeomorphism such that $\| \tau(y) \| = \| y\|$, such that the linearization along the zero section 
\beqn
d\tau|_{S}: NS \to NS
\eeqn
is the identity map of $NS$, and such that 
\beqn
\varphi|_{N^{\epsilon''} S} = \varphi' \circ \tau |_{N^{\epsilon''} S}.
\eeqn
If $T_S$ and $T_S'$ are equivalent, we denote $T_S \sim T_S'$.
\end{enumerate}
\end{defn}

\begin{defn}\label{compatible_neighborhood}
Let $M, P$ be smooth manifolds and $f: M \to P$ be a smooth map. Let $S \subset M$ be a submanifold. We say a tubular neighborhood $T_S$ is {\it compatible} with $f$ if $f \circ \pi_S = f$ in $|T_S|$.
\end{defn}

Below is a lemma useful in the induction argument. As its argument will be used in the orbifold setting, we provide an alternate self-contained proof.

\begin{lemma}\cite[Proposition 6.2]{Mather_1970} \label{lemma_tubular_extension}
In the situation of Definition \ref{compatible_neighborhood}, suppose restriction $f|_S: S \to P$ is a submersion. Let $O \subset S$ be an open subset and $D \subset S$ is a closed subset contained in $O$. Suppose $T_0 = (\epsilon_0, \varphi_0)$ is a tubular neighborhood of $O$ in $M$ which is compatible with $f$, then there exists a tubular neighborhood $T_S = (\epsilon, \varphi)$ of $S$ which is compatible to $f$ such that $T_S$ and $T_0$ are equivalent near $D$. In particular, the old and new tubular projections and tubular functions agree near $D$.
\end{lemma}

\begin{proof}
For each $p \in P$, denote $M_p = f^{-1}(p)$ and $S_p = f^{-1}(p) \cap S$. By the condition that $f|_S: S \to P$ is a submersion, $S_p$ is a smooth submanifold and $M_p$ is a smooth submanifold near $S_p$. Moroever, the normal bundle of $S_p$ inside $M_p$ is isomorphic to $NS|_{S_p}$.

We first extend the tubular neighborhood in the linear level. Consider the exact sequence of vector bundles with a splitting $\tau: NS \to TM|_S$
\beqn
\xymatrix{ 0 \ar[r] & TS \ar[r] & TM|_S \ar[r] & NS  \ar@/_1pc/[l]_{\tau} \ar[r] & 0 }.
\eeqn
A splitting $\tau$ is said to be compatible with $f$ if for each $x \in S$, $\tau (N_x S) \subset T_x M_{f(x)}$. For the existing tubular neighborhood $T_0 = (\epsilon_0, \varphi_0)$, the linearization 
\beqn
\tau_0:= d\varphi_0|_O: NS|_O \to TM|_O
\eeqn
gives a splitting of the exact sequence over $O$ which is compatible with $f$. Then one can find a splitting $\tau_1: NS \to TM|_S$ compatible with $f$ such that $\tau_1$ agrees with $\tau_0$ over an open neighborhood $O' \subset O$ of $D$.

We use the extended splitting $\tau_1$ to define a tubular neighborhood. For each $p\in P$, let $\exp_{M_p}$ be the exponential map inside the submanifold $M_p$ with respect to the restriction of $g^{TM}$ to $M_p$. For $\epsilon: S \to {\mb R}_+$ sufficiently small, define
\beqn
\varphi_1: N^\epsilon S \to M,\ \varphi_1 (x, v) = \exp_{M_{f(x)}} \big( \tau_1 (x, v) \big) \in M_{f(x)} \subset M.
\eeqn
It is easy to check that this is a tubular neighborhood which is compatible with $f$. However it does not necessarily match with the existing one.

We interpolate $\varphi_0$ and $\varphi_1$ over the region $O' \setminus D$. Choose a cut-off function $\nu_{O', D}: S \to [0, 1]$ supported in $O'$ such that $\nu_{O', D} \equiv 1$ near $D$. For each $x \in O' \setminus D$, compare the two fibers 
\beqn
\varphi_1 (N_x^\epsilon S)\ {\rm and}\ \varphi_0 (N_x^\epsilon S).
\eeqn
As $\tau_1$ and $\tau_0$ has the same linearization along $O'$, we see that for all $(x, v) \in NS|_{O'}$, one has 
\beqn
d_{M_{f(x)}} \Big( \varphi_0 (x, v), \varphi_1 (x, v) \Big) = O(\|v\|^2).
\eeqn
Then if $\epsilon$ is sufficiently small, for all $(x, v) \in N^\epsilon O'$, there exists a unique shortest geodesic
\beqn
\gamma_{(x, v)}: [0, 1] \to M_{f(x)}
\eeqn
such that 
\begin{align*}
&\ \gamma_{(x, v)}(0) = \varphi_0 (x, v),\ &\ \gamma_{(x, v)}(1) = \varphi_1 (x, v).
\end{align*}
Then define $\varphi: N^\epsilon S \to M$ by 
\beqn
\varphi(x, v) = \left\{ \begin{array}{rl}  \varphi_0 (x, v),\ &\ x \in D,\\
     \gamma_{(x, v)} \big( 1- \nu_{O', D}(x) \big),\ &\ x \in O' \setminus D,\\
     \varphi_1 (x, v),\ &\ x \notin O'.
     \end{array}\right.
\eeqn
(See Figure \ref{figure_tubular_extension} for an illustration of this interpolation.)

\begin{figure}[h]
    \centering
    \includegraphics[width=0.5\linewidth]{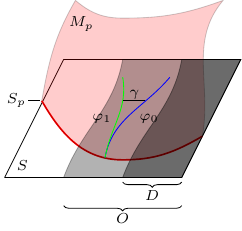}
    \caption{The fibers of $\varphi_1$ and $\varphi_0$ over the same point are contained in the same preimage $M_p = f^{-1}(p)$. Their tangents at the zero section are the same so their distance is small compared to the distance from the zero section. One interpolates between them along the shortest geodesics $\gamma\subset M_p$ in a region near $D$.}
    \label{figure_tubular_extension}
\end{figure}
One can easily check that $\varphi$ is a tubular neighborhood, compatible with $f$, and agrees with $\varphi_0$ near $D$. Hence $T_S = (\epsilon, \varphi)$ is a tubular neighborhood of $S$ satisfying the desired conditions.
\end{proof}

\subsubsection{Normal bundles of the isotropy stratification}

We would like to consider tubular neighborhoods of smooth manifolds inside an orbifold (without boundary). We first discuss the normal bundle to strata of the isotropy stratification (Lemma \ref{lemma_isotropy_stratification}). Here we do not consider the vector bundle ${\mc E} \to {\mc U}$ and the isotropy stratification of ${\mc U}$ means that for the rank zero vector bundle.

Let $\gamma$ be an isotropy type represented by $(G_0, V_0)$ and ${\mc U}_\gamma \subset {\mc U}$ be the stratum of the isotropy stratification. In the normal direction, locally, a neighborhood of ${\mc U}_\gamma$ is a fiberwise quotient of a vector bundle by the group $G_0$. It is tempting to define a normal bundle to ${\mc U}_\gamma$ as a fiberwise quotient of a global vector bundle over ${\mc U}_\gamma$ such that for each orbifold chart $C = (G, U, \psi)$ containing a point of ${\mc U}_\gamma$ the vector bundle is the normal bundle over $U_G \subset U$. However, there is no canonical transition functions. Indeed, the correct viewpoint is that ${\mc U}_\gamma$ is a non-effective orbifold and the local pieces $NU_G$ define an orbifold vector bundle. However, we want to refrain from discussing non-effective orbifolds as we always do in this paper.

Let ${\mc A}_\gamma$ be the set of charts of ${\mc U}$ of the form $(G, U, \psi)$ such that $G \cong G_0$ and the fiber of $TU|_{U_G}$ is isomorphic to $V_0$ as representations. Then each such chart provides a manifold chart of ${\mc U}_\gamma$. Then define a space
\beqn
{\mc N}_\gamma:= \left( \bigsqcup_{C = (G, U, \psi) \in {\mc A}_\gamma} NU_G \right)/\sim
\eeqn
where the equivalence relation is defined as follows: a point $(x_1, v_1) \in NU_{1, G_1}$ (where $x_1 \in U_{1, G_1}$ and $v_1 \in NU_{1, G_1}|_{x_1}$) coming from a chart $C_1 = (G_1, U_1, \psi_1)$ and a point $(x_2, v_2)  \in NU_{2, G_2}$ coming from a chart $C_2 = (G_2, U_2, \psi_2)$ are equivalent if there exists a chart $C = (G, U, \psi) \in {\mc A}_\gamma$, $(x, v) \in NU_G$, and chart embeddings
\beqn
\iota_i: U \to U_i,\ i = 1, 2
\eeqn
such that $\iota_i(x) = x_i$ and $d\iota_i(x)(v) = v_i$. We claim that ${\mc N}_\gamma$ is naturally an effective orbifold. Indeed, for each $C = (G, U, \psi)\in {\mc A}_\gamma$, the natural map
\beqn
d\psi: NU_G \to {\mc N}_\gamma
\eeqn
is $G$-invariant and a homeomorphism onto an open subset. Hence $(G, NU_G, d\psi)$ is an orbifold chart of ${\mc N}_\gamma$. Moreover, the compatibility of charts in ${\mc A}_\gamma$ guarantees the compatibility of the collection of orbifold charts of ${\mc N}_\gamma$. Therefore, ${\mc N}_\gamma$ is an effective orbifold. Notice that there is a canonical projection map ${\mc N}_\gamma \to {\mc U}_\gamma$ and a canonical zero section ${\mc U}_\gamma \to {\mc N}_\gamma$. We call ${\mc N}_\gamma$ the {\it normal bundle} to ${\mc U}_\gamma$. Notice that ${\mc N}_\gamma$ has its isotropy stratification whose strata are indexed by isotropy types $\delta \geq\gamma$. 

Moreover, suppose ${\mc U}$ is equipped with a straightened Riemannian metric $g^{T{\mc U}}$, then there is a smooth function on ${\mc N}_\gamma$ whose pullback to each chart $(G, NU_G, d\psi)$ is the function $(x, v) \mapsto \| v\|^2$ where $\| v\|$ is calculated by the Riemannian metric.

\subsubsection{Tubular neighborhoods in orbifolds}

Let ${\mc U}$ be an effective orbifold equipped with a straightened Riemannian metric $g^{T{\mc U}}$. Within this appendix, a {\it submanifold} of ${\mc U}$ means a submanifold of a stratum of the isotropy stratification of ${\mc U}$. On the other hand, let $P$ be a smooth manifold. Recall that a continuous map $f: {\mc U} \to P$ is called smooth if its pullback to each orbifold chart is a smooth map to $P$; it is a submersion if its pullback to each orbifold chart is a submersion  (see Section \ref{sec:orbifolds}). As an example, the projection ${\mc N}_\gamma \to {\mc U}_\gamma$ is a submersion.

\begin{defn}
Let $Z \subset {\mc U}$ be a submanifold contained in ${\mc U}_\gamma$. Let $ \pi_{N_\gamma S}: N_\gamma Z \to Z$ be the normal bundle of $Z$ inside ${\mc U}_\gamma$. 
The {\it normal bundle} of a submanifold $Z$ in ${\mc U}$ is the orbifold 
\beqn
NZ:= \pi_{N_\gamma Z}^* {\mc N}_\gamma|_Z
\eeqn
with a natural projection map $\pi_{NZ}: NZ \to Z$. A point of $NZ$ is denoted by $(x, v, w)$ where $x \in Z$, $v \in N_\gamma Z|_x$, and $w \in {\mc N}_\gamma|_x$. Define
\beqn
\rho_{NZ}(x, v, w):= \| v\|^2 + \| w\|^2.
\eeqn
For each smooth function $\epsilon: Z \to {\mb R}_+$, denote
\beqn
N^\epsilon Z = \Big\{ (x, v, w) \in NZ \ |\ \rho_{NZ}(x, v, w) < \epsilon (x) \Big\}.
\eeqn
\end{defn}

\begin{defn}
A {\it tubular neighborhood} of $Z$ inside ${\mc U}$ is a pair $T_Z = (\epsilon, \varphi)$ where $\epsilon: Z \to {\mb R}_+$ is a smooth function and 
\beqn
\varphi: N^\epsilon Z \to {\mc U}
\eeqn
is an orbifold open embedding satisfying the following conditions.
\begin{enumerate}
\item $\varphi|_Z$ is the inclusion of $Z$ into ${\mc U}$.

\item For each $x \in Z$, if we lift $\varphi$ to a chart $C = (G, U, \psi)$ containing $x$, then $\varphi$ is a $G$-equivariant tubular neighborhood of $\psi^{-1}(Z)$ inside $U$.

\item Associated to each tubular neighborhood $T_Z$, define
\beqn
|T_Z|:= \varphi( N^\epsilon Z)
\eeqn
over which there is a tubular projection $\pi_Z: |T_Z| \to Z$ and a tubular function $\rho_Z: |T_Z| \to [0, +\infty)$. Both are smooth maps on orbifolds.

\item Two tubular neighborhoods $T_Z = (\epsilon, \varphi)$ and $T_Z' = (\epsilon', \varphi')$ are equivalent if there exist a smooth function $\epsilon'': Z \to {\mb R}_+$ such that $\epsilon'' \leq \min ( \epsilon, \epsilon')$, an orbifold isomorphism
\beqn
\tau: N^{\epsilon''}Z \to N^{\epsilon''} Z
\eeqn
satisfying 1) $\pi_{NZ} \circ \tau = \pi_{NZ}$, 2) $\rho_{NZ} \circ \tau = \rho_{NZ}$, 3), for each chart $C = (G, U, \psi) \in {\mc A}_\gamma$, $\tau$ is lifted to a $G$-invariant map $\tau: NU_G \to NU_G$ such that the linearization along the zero section is the identity map of $NU_G$, and 4)
\beqn
\varphi|_{N^{\epsilon''} Z} = \varphi' \circ \tau|_{N^{\epsilon''} Z}.
\eeqn

\item A tubular neighborhood $T_Z$ is {\it compatible} with a smooth map $f:{\mc U} \to P$ if $f\circ \pi_Z = f$ on $|T_Z|$. 

\end{enumerate}
\end{defn}

Now we can generalize the tubular neighborhood theorem to the orbifold setting.

\begin{prop}\label{prop_orbifold_tubular_extension}
Let $Z \subset {\mc U}$ be a submanifold. Let $f: {\mc U} \to P$ be a smooth map such that $f|_Z: Z \to P$ is a submersion. Let $O \subset Z$ be an open subset and $D \subset Z$ be a closed subset contained in $O$. Suppose $T_0 = (\epsilon_0, \varphi_0)$ is a tubular neighborhood of $O$ which is compatible with $f$. Then there exists a tubular neighborhood $T_Z = (\rho, \varphi)$ of $Z$ compatible with $f$ such that $T_Z$ is equivalent to $T_0$ near $D$.
\end{prop}

\begin{proof}
The argument is analogous to the manifold case. We construct the tubular neighborhood on the chart level and prove that the chartwise constructions are compatible with chart embeddings. 

We introduce similar notations first. For each $x \in Z$, choose an orbifold chart $C = (G, U, \psi)$ centered at $x$. Denote $\tilde Z = \psi^{-1}(Z) \subset U_G$, $\tilde O = \psi^{-1}(O)$ and $\tilde D = \psi^{-1}(D)$. Then by definition $\tilde f:= f\circ \psi: U \to P$ is smooth. For each $p$, denote $U_p = \tilde f^{-1}(p)$ and $\tilde Z_p = U_p \cap \tilde Z$. Any smooth function $\epsilon: Z \to {\mb R}_+$ is lifted to a smooth function $\epsilon: \tilde Z \to {\mb R}_+$. Then the existing tubular neighborhood $T_0$ is lifted to 
\beqn
\tilde \varphi_0: N^\epsilon \tilde O \to U
\eeqn
after appropriately shrinking the function $\epsilon$.

Consider the $G$-equivariant exact sequence
\beqn
\xymatrix{   0 \ar[r] & T\tilde Z \ar[r] & TU|_{\tilde Z} \ar[r] & N\tilde Z \ar[r] & 0 }.
\eeqn
A $G$-equivariant splitting $\tau: N\tilde Z \to TU|_{\tilde Z}$ is said to be compatible with $f$ if for each $x \in \tilde Z$, $\tau ( N_x \tilde Z) \subset TU_{f(\psi(x))}$. Then the linearization of $\tilde \varphi_0$ defines a splitting over $\tilde O$ which is compatible with $f$, denoted by $\tau_0$. Then we can choose a $G$-equivariant splitting 
\beqn
\tau: N\tilde Z \to TU|_{\tilde Z}
\eeqn
which coincides with $\tau_0$ over an open neighborhood $\tilde O' \subset \tilde O$ of $\tilde D$. Then define the map 
\beqn
\tilde \varphi_1: N^\epsilon \tilde Z \to U,\ \tilde \varphi_1 (x, v) = \exp_{U_{\tilde f(x)}} \big( \tau (x, v) \big) \in U_{\tilde f(x)}.
\eeqn

Next, similar to the proof of Lemma \ref{lemma_tubular_extension}, we interpolate $\tilde \varphi_1$ and $\tilde \varphi_0$ within $\tilde O' \setminus \tilde D$. Indeed, if $x \in \tilde O'$, then
\beqn
d_{U_{\tilde f(x)}} \left( \tilde \varphi_1(x, v), \tilde \varphi_0(x, v) \right) = O( \| v \|^2 ).
\eeqn
Hence there is a shortest geodesic $\gamma_{(x, v)}: [0, 1] \to U_{\tilde f(x)}$ connecting $\tilde \varphi_0(x, v)$ and $\tilde \varphi_1 (x, v)$. Choose a smooth cut-off function $\nu_{\tilde O', \tilde D}: \tilde Z \to {\mb R}_+$ supported inside $\tilde O'$ such that $\nu_{\tilde O', \tilde D} \equiv 1$ near $\tilde D$. Then define
\beqn
U_{\tilde f(x)}\ni \tilde \varphi(x, v) = \left\{ \begin{array}{rl} \tilde \varphi_0(x, v),\ &\ x \in \tilde D,\\
     \gamma_{(x, v)}\big( 1-  \nu_{\tilde O', \tilde D} (x)  \big) 
  ,\ &\ x \in \tilde O \setminus \tilde D,\\
     \tilde \varphi_1(x, v),\ &\ x \notin \tilde O. \end{array}\right.
\eeqn

As the submanifolds $U_{\tilde f(x)}$ are $G$-invariant and the Riemannian metric is $G$-invariant, the chartwise tubular neighborhoods are $G$-equivariant. Also easy to check that these tubular neighborhoods are compatible with chart embedding. Hence this is a tubular neighborhood of $Z$. From the construction we can see it agrees with $T_0$ near $D$.
\end{proof}

\subsection{Whitney stratified sets in orbifolds are Thom--Mather sets}

\begin{defn}\label{defnb15}
Let ${\mc U}$ be an effective orbifold (without boundary). A {\it Whitney stratified subset} of ${\mc U}$ is a subset $Z \subset {\mc U}$ equipped with a stratification ${\mf Z}$ satisfying the following conditions.
\begin{enumerate}
    \item Each stratum $Z_\alpha \in {\mf Z}$ is a submanifold contained in a stratum ${\mc U}_{\gamma(\alpha)}$ of the isotropy stratification.

    \item For each orbifold chart $C = (G, U, \psi)$, the above conditions imply that for each stratum $Z_\alpha \in {\mf Z}$, $\psi^{-1}(Z_\alpha) \subset U$ is a (possibly empty) smooth submanifold. Then for each pair $Z_\alpha < Z_\beta$, $(\psi^{-1}(Z_\beta), \psi^{-1}(Z_\alpha))$ is Whitney regular. 
\end{enumerate}
\end{defn}

\begin{defn}\label{defn_control_data}
Let $(Z, {\mf Z})$ be a Whitney stratified subset of ${\mc U}$. A set of {\it ambient control data} for ${\mf Z}$ is a collection ${\mf I} = \{ T_\alpha\ |\ Z_\alpha \in {\mf S}\}$ of  tubular neighborhoods $T_\alpha$ of $Z_\alpha$ such that whenever $Z_\alpha < Z_\beta$, one has 
\begin{align*}
&\ \pi_\alpha \circ \pi_\beta = \pi_\alpha,\ &\ \rho_\alpha \circ \pi_\beta = \rho_\alpha
\end{align*}
as germs of maps defined near $Z_\alpha$.
\end{defn}

The following is an obvious consequence of Definition \ref{defn_Thom_Mather} and Definition \ref{defn_control_data}.

\begin{prop}\label{prop_b14}
Suppose one has a set of ambient control data ${\mf I} = \{ (T_\alpha)\ |\ Z_\alpha \in {\mf Z} \}$ for $(Z, {\mf Z})$. Let ${\mf J}$ be the data
\beqn
{\mf J} = \Big\{ (\pi_\alpha|_Z, \rho_\alpha|_Z)\ |\ Z_\alpha \in {\mf Z} \Big\}.
\eeqn
Then the triple $(Z, {\mf Z}, {\mf J})$ is a Thom--Mather stratified space.
\end{prop}

Our theorem is 
\begin{thm}\label{thm_b15}
Let ${\mc U}$ be an orbifold (without boundary) and $(Z, {\mf Z})$ is a closed Whitney stratified subset. Then there exists a set of ambient control data for ${\mf Z}$.
\end{thm}

The proof resembles that of \cite[Proposition 7.1]{Mather_1970}. Since the proof (as well as Mather's proof) is complicated, we would like to demonstrate the main idea for the simple case when there are only three strata $Z_1 < Z_2 < Z_3$. We can first construct a tubular neighborhood $T_1 = (\epsilon_1, \varphi_1)$ of $Z_1$, which comes with the tubular projection 
\beqn
\pi_1: |T_1| \to Z_1
\eeqn
and tubular function 
\beqn
\rho_1: |T_1| \to [0, +\infty).
\eeqn

\begin{lemma}
There exists an open neighborhood $O_1 \subset |T_1 |$ of $Z_1$ such that 
\beqn
(\pi_1, \rho_1): Z_2 \cap O_1 \to Z_1 \times (0, +\infty)
\eeqn
is a submersion.
\end{lemma}

\begin{proof}
It suffices to prove the assertion locally. For each $p \in Z_1$, choose an orbifold chart $C = (G, U, \psi)$ containing $p$ so that $\psi^{-1}(Z_1) \subset U_G$ is a submanifold. The two maps $\pi_1$ and $\rho_1$ are lifted to $G$-invariant maps near $\psi^{-1}(Z_1)$. Denote them by $\tilde \pi_1$ and $\tilde \rho_1$. On the other hand, $\psi^{-1}(Z_2) \subset U$ is a $G$-invariant smooth submanifold. By definition, $(\psi^{-1}(Z_2), \psi^{-1}(Z_1))$ satisfies Whitney's condition (b) at each point of $\psi^{-1}(Z_1)$. This implies (see \cite[Lemma 7.3]{Mather_1973}) that there exists a $G$-invariant neighborhood $\tilde O_1$ of $\psi^{-1}(Z_1)$ such that the projection 
\beqn
(\tilde \pi_1, \tilde \rho_1): \psi^{-1}(Z_2) \cap \tilde O_1 \to \psi^{-1}(Z_1) \times (0, +\infty)
\eeqn
is a submersion. This is equivalent to say that $(\pi_1, \rho_1)$ is a submersion in $Z_2 \cap \psi(\tilde O_1)$. Then $O_1$ can be constructed by taking union of these local pieces. 
\end{proof}

Then by the relative case of Proposition \ref{prop_orbifold_tubular_extension}, if we shrink $|T_1|$ appropriately (by taking a function $\epsilon_1' \leq \epsilon_1$), then there exists a tubular neighborhood $T_2$ of $Z_2 \cap |T_1|$ which is compatible with $(\pi_1, \rho_1)$. Then use the absolute case of Proposition \ref{prop_orbifold_tubular_extension} again, we can extend it to a tubular neighborhood $T_2$ of $Z_2$ which is compatible with $T_1$ near $Z_1$.

Now we consider the third stratum $Z_3$. Notice that we can freely shrink $|T_1|$ and $|T_2|$ without altering the compatibility between $T_1$ and $T_2$. First consider the open subset $Z_3 \cap |T_2|\subset Z_3$. In the same way as before, after shrinking $|T_2|$ appropriately, there exists a tubular neighborhood $T_3$ of $Z_3 \cap |T_2|$ inside $|T_2|$ which is compatible with $T_2$. We need to show that this is also compatible with $T_1$. Indeed, 
\beqn
\pi_1 \circ \pi_3 = (\pi_1 \circ \pi_2) \circ \pi_3 = \pi_1 \circ (\pi_2 \circ \pi_3) = \pi_1 \circ \pi_2 = \pi_1,
\eeqn
where the equality holds near $|T_1| \cap Z_2$. Similarly $\rho_1 \circ \pi_3 = \rho_1$ over the same region. Therefore, the tubular neighborhood $T_3$ is compatible with $(\pi_1, \rho_1)$. 

Then by Proposition \ref{prop_orbifold_tubular_extension}, after appropriately shrinking $|T_1|$ and $|T_2|$, one can find a tubular neighborhood $T_3$ of $Z_3$ within $|T_1| \cup |T_2|$. Using this proposition again, we can construct a global tubular neighborhood $T_3$ which is compatible with both $T_1$ and $T_2$. This finishes the proof in this sample case.

\begin{proof}[Proof of Theorem \ref{thm_b15}]
We follow the proof of Mather \cite[Section 7]{Mather_1970}. Suppose $Z$ has dimension $n$ and let $Z_k \subset Z$ be the union of strata of dimension at most $k$ and ${\mf Z}_k$ be the restricted stratification. We construct inductively for $k$. The induction hypothesis is that there exists a set of ambient control data ${\mf I}_k$ for $(Z_k, {\mf Z}_k)$ (which is itself a Whitney stratified set in ${\mc U}$). Now we construct a set of ambient control data ${\mf I}_{k+1}$ for $(Z_{k+1}, {\mf Z}_{k+1})$ which extends ${\mf I}_k$. Notice that it is always the freedom to shrink the existing tubular neighborhood without breaking the compatibility conditions.

We construct ${\mf I}_{k+1}$. As strata of the same dimension do not intersect, one can construct a tubular neighborhood for each individual strata of dimension $k+1$ separately. Without loss of generality, assume that there is only one stratum $Z_\beta$ of dimension $k+1$. The construction of the  tubular neighborhood $T_\beta$ uses a reversed induction on $l = 0, 1, \ldots, k+1$. For each $l$, denote
\beqn
O_l = \bigcup_{l \leq {\rm dim}( Z_\alpha) \leq k} |T_\alpha|
\eeqn
which is an open subset of ${\mc U}$. We construct inductively a tubular neighborhood $T_{\beta, l}$ of $Z_\beta \cap O_l$ inside $O_l$ which is compatible with all existing tubular neighborhoods. The $l = k+1$ case has nothing to prove. Suppose we have constructed a tubular neighborhood for $Z_\beta \cap T_{l+1}$ satisfying the requirement. Then for each $Z_\alpha$ of dimension $l$, we can check that 
\beqn
\pi_\alpha \circ \pi_\beta = \pi_\alpha,\ \rho_\alpha \circ \pi_\beta = \rho_\alpha
\eeqn
near $Z_\alpha$. Therefore, $T_{\beta, l+1}$ is compatible with $(\pi_\alpha, \rho_\alpha)$. Then one can extend $T_{\beta, l+1}$ to a tubular neighborhood $T_{\beta, l}$ of $Z_\beta \cap O_l$ inside $O_l$ using Proposition \ref{prop_orbifold_tubular_extension}. Then finally the tubular neighborhood $T_{\beta, 0}$ is a tubular neighborhood of $Z_\beta$ inside ${\mc U}$ which is compatible with all existing tubular neighborhoods. The induction on $k$ can be carried on and stops after finitely many stages.
\end{proof}

\subsubsection{The case for orbifolds with boundary}

We extend Definition \ref{defnb15} and Theorem \ref{thm_b15} to the case when the orbifold has boundary. Let ${\mc U}$ be an orbifold with boundary and $Z \subset {\mc U}$ be a subset. Denote
\begin{align*}
    &\ {\rm Int}Z:= Z \cap {\rm Int} {\mc U},\ &\ \partial Z:= Z \cap \partial {\mc U}.
\end{align*}
If ${\mf Z}$ is a stratification on $Z$ with strata $Z_\alpha$, then denote
\begin{align*}
    &\ {\rm Int} Z_\alpha:= Z_\alpha \cap {\rm Int} {\mc U},\ \partial Z_\alpha:= Z_\alpha \cap \partial {\mc U}.
\end{align*}

\begin{defn}
A {\it collared Whitney stratified subset} of ${\mc U}$ is a subset $Z \subset {\mc U}$ equipped with a stratification ${\mf Z}$ satisfying the following conditions.
\begin{enumerate}

    \item The interior ${\rm Int} Z$ with the induced stratification is a Whitney stratified subset of ${\rm Int}{\mc U}$. Let ${\mc U}_{\gamma(\alpha)}$ be the stratum of ${\mc U}$ whose interior contains ${\rm Int} Z_\alpha$.

    \item $\partial Z_\alpha$ is a submanifold of $\partial {\mc U}$ contained in $\partial {\mc U}_{\gamma(\alpha)}$.

    \item There is an open embedding $\iota: \partial {\mc U} \times [0, \epsilon) \to {\mc U}$ (a collar neighborhood of the boundary) sending $\partial {\mc U} \times  \{0\}$ identically onto $\partial {\mc U}$ such that the restriction of $\iota$ onto $\partial Z_\alpha \times [0, \epsilon)$ is an open embedding into $Z_\alpha$.
\end{enumerate}
\end{defn}

If $Z \subset {\mc U}$ is a collared Whitney stratified subset, then one can refine the stratification as
\beqn
Z = \bigsqcup_\alpha {\rm Int} Z_\alpha \sqcup \bigsqcup_{\alpha} \partial Z_\alpha.
\eeqn

\begin{cor}
Let $Z \subset {\mc U}$ be a collared Whitney stratified subset. Then $Z$ with the refined stratification specified above admits a set of ambient control data.
\end{cor}

\begin{proof}
The proof can be carried out in the same way as that of Theorem \ref{thm_b15}. The collared structure allows one to choose product type ambient control data near boundary strata. 
\end{proof}

\bibliography{mathref}

\bibliographystyle{amsalpha}

\end{document}